%% file: main.tex
\documentclass[a4paper,12pt,leqno,openbib]{book}
\pdfoutput=1

\selectfont
\usepackage[english]{babel}
\usepackage{bbold}
\usepackage{stmaryrd}
\usepackage{verbatim}
\usepackage[utf8x]{inputenc}
\usepackage{amsmath}
\numberwithin{equation}{subsection}
\usepackage{amssymb}
\usepackage{amsthm}
\usepackage{graphicx}
\usepackage{mathtools}
\usepackage[T1]{fontenc}
\usepackage{amsthm}
\usepackage{setspace}
\usepackage{amsfonts}
\usepackage{eucal}
\usepackage{mathrsfs}
\usepackage{pictexwd,dcpic}
\usepackage{tikz}
\usepackage{tqft}
\usepackage{pigpen}
\usepackage{xcolor}
\usepackage{pdfpages}
\usepackage{geometry}

\usepackage[colorinlistoftodos]{todonotes} 
             
\usepackage{fancyhdr} 

\usepackage{hyperref}
\hypersetup{
    colorlinks=true,
    linkcolor=blue,
    filecolor=magenta,      
    urlcolor=cyan,
}

\newcommand*\circled[1]{\tikz[baseline=(char.base)]{
            \node[shape=rectangle,draw,inner sep=2.5pt] (char) {#1};}}

\newlength\FHoffset
\fancyheadoffset{\FHoffset}

\newlength\FHleft
\newlength\FHright

\newbox\FHline
\setbox\FHline=\hbox{\hsize=\paperwidth%
  \hspace*{\FHleft}%
  \rule{\dimexpr\headwidth-\FHleft-\FHright\relax}{\headrulewidth}\hspace*{\FHright}%
}

\renewcommand{\headrulewidth}{0.4pt}

\fancypagestyle{plain}{%
\fancyhf{} 
\fancyhead[LO,RE]{\circled{\thepage}}
\renewcommand{\headrulewidth}{0pt}
}

\newtheorem*{theorem-non}{Theorem}
\newtheorem*{theorem-nonf}{Théorème}
\newtheorem*{corollary-non}{Corollary}
\newtheorem*{corollary-nonf}{Corollaire}
\newtheorem{lmm}[subsubsection]{Lemma}
\newtheorem{thm}[subsubsection]{Theorem}
\newtheorem{prop}[subsubsection]{Proposition}
\newtheorem{cor}[subsubsection]{Corollary}
\newtheorem{caution}[subsubsection]{Caution}
\theoremstyle{definition}
\newtheorem{defn}[subsubsection]{Definition}
\newtheorem{context}[subsubsection]{Context}
\newtheorem{ex}[subsubsection]{Example}
\newtheorem{rmk}[subsubsection]{Remark}

\newtheorem{notation}[subsubsection]{Notation}
\newtheorem{construction}[subsubsection]{Construction}
\title{thesis title}
\author{Massimo Pippi}
\date{ date dissertation }

\newcommand{\drl}[0]{\llparenthesis}
\newcommand{\drr}[0]{\rrparenthesis}
\newcommand{\dsl}[0]{\llbracket}
\newcommand{\dsr}[0]{\rrbracket}
\newcommand{\ital}[1]{\textit{#1}} 

\newcommand{\x}[0]{\times}

\newcommand{\hm}[0]{\textup{Hom}}
\newcommand{\heart}{\ensuremath\heartsuit}
\newcommand{\oo}{\infty}
\newcommand{\Z}[0]{\mathbb{Z}}
\newcommand{\N}[0]{\mathbb{N}}
\newcommand{\C}[0]{\mathbb{C}}
\newcommand{\Q}[0]{\mathbb{Q}}
\newcommand{\aff}[2]{\mathbb{A}^{#1}_{#2}}
\newcommand{\proj}[2]{\mathbb{P}^{#1}_{#2}}
\newcommand{\sch}[0]{\textup{\textbf{Sch}}_{S}}
\newcommand{\Cc}[0]{\mathcal{C}}

\newcommand{\sm}[0]{\textup{\textbf{Sm}}_S}
\newcommand{\ag}[1]{\mathbb{A}^1_{#1}/\mathbb{G}_{m,#1}}
\newcommand{\gm}[1]{\mathbb{G}_{m,#1}}
\newcommand{\bgm}[1]{\textup{B}\mathbb{G}_{m,#1}}
\newcommand{\Ec}[0]{\mathcal{E}}
\newcommand{\Lc}[0]{\mathcal{L}}
\newcommand{\Lb}[0]{\mathbb{L}}
\newcommand{\Mc}[0]{\mathcal{M}}
\newcommand{\Oc}[0]{\mathcal{O}}
\newcommand{\Rc}[0]{\mathcal{R}}
\newcommand{\Fc}[0]{\mathcal{F}}

\newcommand{\Map}[0]{\textup{Map}}
\newcommand{\rhom}[0]{\mathbb{R}\textup{Hom}}
\newcommand{\bu}[0]{\textup{B}\mathbb{U}}
\newcommand{\Ql}[1]{\mathbb{Q}_{\ell#1}}
\newcommand{\Rl}[0]{\mathcal{R}^{\ell}}
\newcommand{\Rlv}[0]{\mathcal{R}^{\ell,\vee}}
\newcommand{\Hl}[0]{\mathbb{H}_{\mathbb{Q}_{\ell}}}
\newcommand{\rbul}[1]{\mathcal{R}^{\ell}_{#1}(\textup{B}\mathbb{U}_{#1})}
\newcommand{\Hr}[0]{\mathbb{H}^{\rbul{}}_{\mathbb{Q}_{\ell}}}
\newcommand{\rperf}[0]{\mathcal{R}_{\textup{\textbf{Perf}}}}
\newcommand{\hk}[0]{\textup{HK}}
\newcommand{\Mv}[0]{\mathcal{M}^{\vee}}
\newcommand{\HQ}[0]{\textup{H}\mathbb{Q}}
\newcommand{\MB}[1]{\textup{M}\mathbb{B}_{#1}}

\newcommand{\dbarrow}[4]{
\begindc{\commdiag}[30]
\obj(0,1)[]{$#1$}
\obj(35,1)[]{$#2$}
\mor(7,0)(28,0){$#3$}[\atright,\solidarrow]
\mor(28,2)(7,2){$#4$}[\atright,\solidarrow]
\enddc
}

\newcommand{\Sp}[0]{\textup{\textbf{Sp}}}
\newcommand{\ch}[0]{\textup{Ch}_A}

\newcommand{\shvl}[0]{\textup{\textbf{Shv}}_{\mathbb{Q}_{\ell}}}
\newcommand{\shv}[0]{\textup{\textbf{Shv}}}
\newcommand{\md}[0]{\textup{\textbf{Mod}}}
\newcommand{\fclgm}[0]{\textup{LG}}
\newcommand{\lgm}[1]{\textup{LG}_{S}(#1)}
\newcommand{\tlgm}[1]{\textup{LG}_{(#1,\Lc_{#1})}}
\newcommand{\tlgmproj}[1]{\textup{LG}_{(\mathbb{P}^{n-1}_{#1},\mathcal{O}(1))}}

\newcommand{\dgcat}{\textup{dgCat}_S}
\newcommand{\dgcatlf}{\textup{dgCat}^{\textup{lf}}_S}
\newcommand{\dgcatdk}{\textup{\textbf{dgCat}}_S}
\newcommand{\dgcatm}{\textup{\textbf{dgCat}}^{\textup{idm}}_S}
\newcommand{\dgcatlfo}{\textup{dgCat}^{\textup{lf},\otimes}_S}
\newcommand{\dgcatdko}{\textup{\textbf{dgCat}}^{\otimes}_S}
\newcommand{\dgcatmo}{\textup{\textbf{dgCat}}^{\textup{idm},\otimes}_S}
\newcommand{\dgcatft}{\textup{\textbf{dgCat}}_S^{\textup{ft}}}

\newcommand{\cohb}[1]{\textup{\textbf{Coh}}^{b}(#1)}
\newcommand{\perf}[1]{\textup{\textbf{Perf}}(#1)}
\newcommand{\qcoh}[1]{\textup{\textbf{QCoh}}(#1)}
\newcommand{\cohm}[1]{\textup{\textbf{Coh}}^{-}(#1)}
\newcommand{\cohbp}[2]{\textup{\textbf{Coh}}^b(#1)_{\textup{\textbf{Perf}}(#2)}}
\newcommand{\cohmp}[2]{\textup{\textbf{Coh}}^{-}(#1)_{\textup{\textbf{Perf}}(#2)}}
\newcommand{\cohs}[1]{\textup{Coh}^{s}(#1)}
\newcommand{\sing}[1]{\textup{\textbf{Sing}}(#1)}
\newcommand{\pdgcat}{\textup{Pairs-dgCat}^{\textup{lf}}_S}
\newcommand{\pdgcato}{\textup{Pairs-dgCat}^{\textup{lf},\otimes}_S}
\newcommand{\sring}{\textup{\textbf{sRing}}}
\newcommand{\dscaff}{\textup{\textbf{dSch}}^{\textup{aff}}}
\newcommand{\dsc}{\textup{\textbf{dSch}}}

\newcommand{\mf}[2]{\textup{\textbf{MF}}(#1,#2)}
\newcommand{\schet}[1]{\textup{\textbf{Sch}}_{\textup{et},#1}}
\newcommand{\uH}[1]{\textup{\textbf{H}}_{#1}}
\newcommand{\sh}[0]{\textup{\textbf{SH}}_S}
\newcommand{\shm}[0]{\textup{\textbf{SH}}_S^{\otimes}}
\newcommand{\shnc}[0]{\textup{\textbf{SH}}^{\textup{nc}}_S}

\newcommand{\shncm}[0]{\textup{\textbf{SH}}^{\textup{nc}, \otimes}_S}
\newcommand{\shncmop}[0]{\textup{\textbf{SH}}^{\textup{nc,op}, \otimes}_S}
\newcommand{\nc}[1]{\mathcal{N}c_{#1}}
\begin{document}

\includepdf{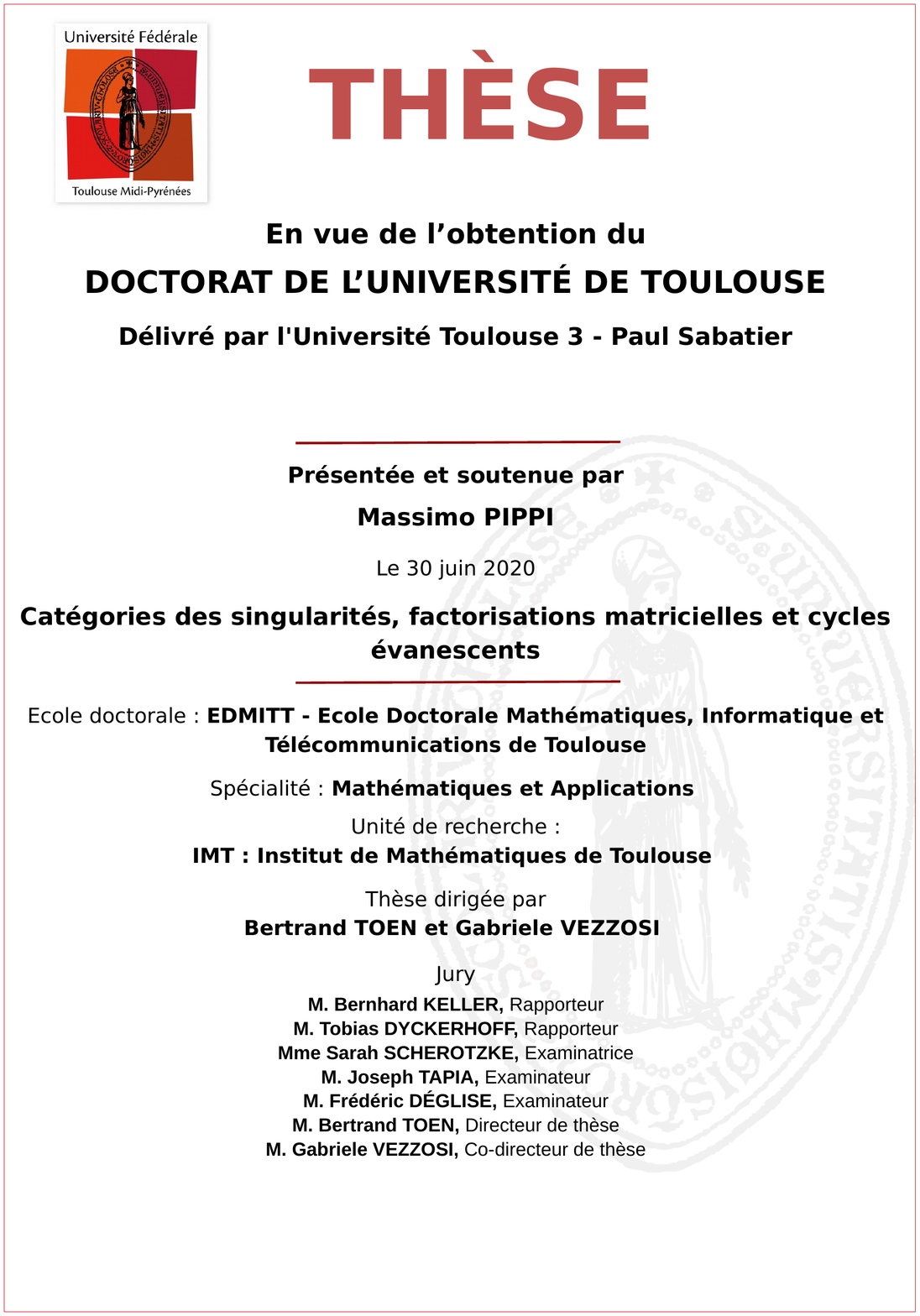}
\leavevmode\thispagestyle{empty}\newpage

\thispagestyle{empty}
\input{chapters/resume}

\newpage

\input{chapters/abstract}

\leavevmode\thispagestyle{empty}\newpage

\chapter*{Introduction}
\input{chapters/introduction_french}

\chapter*{Introduction}
\input{chapters/introduction}

\tableofcontents

\chapter{\textsc{Preliminaries and notation}}
\input{chapters/preliminaries}

\chapter{\textsc{dg-categories of relative singularities and matrix factorizations}}
\input{chapters/chapter02}

\chapter{\textsc{The \texorpdfstring{$\ell$}{l}-adic realization of dg-categories of singulairities of a twisted LG model}}
\input{chapters/chapter03}

\chapter{\textsc{Iterated vanishing cycles}}
\input{chapters/chapter04}

\end{document}

%% file: chapters/resume.tex
\begin{center}
    \large{\textbf{Catégories des singularités, factorisations matricielles et cycles évanescents}}
\end{center}
\begin{center}
    \textbf{Resumé}
\end{center}

Le but de cette thèse est d'étudier les dg-catégories de singularités
$\sing{X,s}$, associées à des couples $(X,s)$, où $X$ est un schéma et $s$ est une section
d'un fibré vectoriel sur $X$. La dg-catégorie $\sing{X,s}$ est définie comme
le noyau du dg-foncteur de $\sing{X_0}$ vers $\sing{X}$ induit par l'image
directe le long de l'inclusion du lieu de zéros (dérivé) $X_0$ de $s$ dans
$X$. 

Dans une première partie, nous supposons que le fibré
vectoriel est trivial de rang $n$. On démontre alors un théorème de structure pour
$\sing{X,s}$ dans le cas où $X=Spec(B)$ est affine. Cet énoncé affirme
que tout objet de $\sing{X,s}$ est représenté par un complexe de
$B$-modules concentré dans $n+1$ degrés. Lorsque
$n=1$, cet énoncé généralise l'équivalence d'Orlov , qui identifie
$\sing{X,s}$ avec la dg-catégorie des factorisations matricielles
$\mf{X}{s}$, au cas où $s \in \Oc_X(X)$ n'est pas nécessairement plat. 

Dans une seconde
partie, nous étudions la cohomologie $\ell$-adique de $\sing{X,s}$ (définie par A.~Blanc -
M.~Robalo - B.~Toën and G.~Vezzosi), où $s$ est une section globale d'un
fibré en droites. Pour cela, on introduit le faisceau $\ell$-adique des
cycles évanescents invariantes par monodromie. En utilisant un
théorème de D.~Orlov généralisé par J.~Burke et M.~Walker, on calcule
la réalisation $\ell$-adique de $\sing{Spec(B),(f_1,..,f_n)}$ pour
$(f_1,..,f_n) \in B^n$. 

Dans le dernier chapitre, nous introduisons les
faisceaux $\ell$-adiques des cycles évanescents itérés pour un schéma sur
un anneau de valuation discrète de rang $2$. On relie ces faisceaux
$\ell$-adiques à la réalisation $\ell$-adique des dg-catégories de singularités
des fibres prises sur certains sous-schémas fermés de la base.

\vspace{1cm}
\textbf{Mots-Clés:} géométrie algébrique dérivée, géométrie non-commutative, cycles évanescents, dg-catégories des singularités, factorisations matricielles, réalisations motivique et $\ell$-adique des dg-catégories

%% file: chapters/abstract.tex
\begin{center}
\large{\textbf{Categories of singularities, matrix factorizations and vanishing cycles}}
\end{center}
\begin{center}
    \textbf{Abstract}
\end{center}

The aim of this thesis is to study the dg categories of singularities $\sing{X,s}$ of pairs $(X,s)$, where $X$ is a scheme and $s$ is a global section of some vector bundle over $X$. $\sing{X,s}$ is defined as the kernel of the dg functor from $\sing{X_0}$ to $\sing{X}$ induced by the pushforward along the inclusion of the (derived) zero locus $X_0$ of $s$ in $X$.

In the first part, we restrict ourselves to the case where the vector bundle is trivial. We prove a structure theorem for $\sing{X,s}$ when $X=Spec(B)$ is affine. Roughly, it tells us that every object is $\sing{X,s}$ is represented by a complex of $B$-modules concentrated in $n+1$ consecutive degrees (if $s\in B^n$). By specializing to the case $n=1$, we generalize Orlov's theorem, which identifies $\sing{X,s}$ with the dg category of matrix factorizations $\mf{X}{s}$, to the case where $s \in \Oc_X(X)$ is not flat.

In the second part, we study the $\ell$-adic cohomology of $\sing{X,s}$ (as defined by A.~Blanc - M.~Robalo - B.~Toën and G.~Vezzosi) when $s$ is a global section of a line bundle. In order to do so, we introduce the $\ell$-adic sheaf of monodromy-invariant vanishing cycles. Using a theorem of D.~Orlov generalized by J.~Burke and M.~Walker, we compute the $\ell$-adic realization of $\sing{Spec(B),(f_1,..,f_n)}$ for $(f_1,..,f_n) \in B^n$.

In the last chapter, we introduce the $\ell$-adic sheaves of iterated vanishing cycles of a scheme over a discrete valuation ring of rank $2$. We relate one of these $\ell$-adic sheaves to the $\ell$-adic realization of the dg category of singularities of the fiber over a closed subscheme of the base.

\vspace{1cm}
\textbf{Keywords:} derived algebraic geometry, non-commutative geometry, vanishing cycles, dg categories of singularities, matrix factorizations, motivic and $\ell$-adic realizations of dg categories

%% file: chapters/introduction_french.tex
\pagenumbering{roman}
La géométrie non-commutative, dans la version envisagée par M.~Kontsevich, est l'étude de certaines dg-catégories qui ont des bonnes propriétés. Ces bonnes propriétés sont celles que l'on retrouve dans les dg-catégories attachées aux schémas. En fait, pour tout-schéma $X$, on peut considérer la dg-catégorie des complexes parfaits $\perf{X}$. Cette construction peut être promue en $\oo$-foncteur 
\begin{equation*}
 \perf{\bullet}:    \sch^{\textup{op}} \rightarrow \dgcatm
\end{equation*}
Cet $\oo$-foncteur est le pont qui relie le monde commutatif avec celui non-commutatif. $\perf{X}$ est parfois appelée l'ombre non-commutative de $X$\footnote{par exemple, voir \cite{pa11}}. Plus précisément, il faudrait l'appeler \emph{une} ombre non-commutative de $X$. En fait, il y a des autres dg-catégories qu'on peut associer à un schéma, par exemple $\qcoh{X}$ et $\cohb{X}$. Ils forment une chaîne d'inclusions de dg-catégories
\begin{equation*}
    \perf{X}\subseteq \cohb{X} \subseteq \qcoh{X}
\end{equation*}
Bien que les constructions $X\mapsto \perf{X}$ et $X\mapsto \qcoh{X}$ soient différentes a priori, elles donnent les mêmes informations. On peut récupérer la dg-catégorie des faisceaux quasi-cohérents sur $X$ a partir de celle des complexes parfaits par la formule 
\begin{equation*}
    \qcoh{X}=\textup{Ind}(\perf{X})
\end{equation*}
Dans la direction opposée, on peut retrouver $\perf{X}$ en considérant les objets compacts dans $\qcoh{X}$. Par conséquent, la différence entre ces deux dg-catégories ne contient pas d'information de nature géométrique. Autrement dit, la différence entre $\perf{X}$ et $\qcoh{X}$ est purement de nature catégorielle.

En revanche, les deux ombres non-commutatives $\perf{X}$ et $\cohb{X}$ sont différentes d'une manière plus intéressante: une conséquence d'un théorème de
M.~Auslander-D.A.~Buchsbaum (\cite[Theorem 4.1]{ab56}) et de J.P.~Serre (\cite[Théorème 3]{se55}) affirme que ces deux dg-catégories sont équivalentes si et seulement si $X$ est un schéma régulier. Ainsi, on peut considérer une nouvelle ombre non-commutative de $X$, connue comme la dg-catégorie des singularités de $X$, qui est définie comme le dg-quotient 
\begin{equation*}
    \sing{X}:=\cohb{X}/\perf{X}
\end{equation*} 
Elle mesure (dans un certain sens) à quel point $X$ est singulier. Ça explique aussi le choix de la terminologie. 

Plus précisément, on va considérer la dg-catégorie des singularités attachée aux couples $(X,s)$, où $X$ est un schéma et $s$ est une sectionne globale d'un fibré vectoriel. La dg-catégorie de $(X,s)$, que nous noterons $\sing{X,s}$, est définie comme le noyau du dg-foncteur 
\begin{equation*}
    \mathfrak{i}_* : \sing{X_0}\rightarrow \sing{X}
\end{equation*}
induit par l'image directe le long du plongement $\mathfrak{i}:X_0:=V(s)\rightarrow X$ du lieu (dérivé) des zéros de $s$ dans $X$. On peut y penser comme à un invariant non-commutatif qui mesure comment les singularités de $X_0$ sont pires que celles de $X$. Le cas le plus simple est celui ou $X$ est un schéma régulier et donc $\sing{X,s}$ coïncide avec la dg-catégorie des singularités de $X_0$ (comme $\sing{X}\simeq 0$). Nous utiliserons ce fait élémentaire très souvent dans cette thèse, même sans le mentionner explicitement.

On va étudier deux aspects de ces dg-catégories
\begin{itemize}
    \item dans la première partie on essaye de mieux comprendre comment ces dg-catégories sont faites. Par exemple, nous démontrerons une généralisation du théorème d'Orlov qui nous dit que, si $s$ est une sectionne globale de $\Oc_X$, alors $\sing{X,s}$ est équivalente à $\mf{X}{s}$, une deuxième ombre non-commutative de $(X,s)$ connue comme la dg-catégorie des factorisations matricielles de $(X,s)$. L'auteur aime bien penser à cette équivalence en termes d'une analogie avec l'isomorphisme de groupes $\pi_1(S^1)\simeq \Z$. Si on adopte le point de vue d'un topologue, on peut interpréter $\pi(S^1)$ comme un objet qui correspond à une définition conceptuelle et $\Z$ comme un joli modèle. De façon équivalente, nous pensons à $\sing{X,s}$ comme l'objet qui sort d'une définition conceptuelle et à $\mf{X}{s}$ comme un joli modèle
    \item dans la deuxième partie on va étudier la cohomologie $\ell$-adique de $\sing{X,s}$, où $s$ est une sectionne globale d'un fibré en droites. Pour cela, on aura besoin d'introduire une généralisation des cycles évanescents invariants par inertie, où le rôle du disque est joué par l'espace total du fibré en droites. En outre, on va étudier aussi la cohomologie $\ell$-adique d'une certaine dg-catégorie dans le cas où l'on se donne un schéma sur un anneau de valuation discrète de rang $2$. Ces deux approches sont des tentatives de généraliser le théorème principal dans \cite{brtv}, et elles correspondent aux deux généralisations en dimension supérieure de la définition de DVR
        \begin{equation*}
        \begindc{\commdiag}[18]
        \obj(0,20)[1]{$\begin{minipage}[c]{3.5cm}
        anneau local noethérien régulier de dimension $1$
        \end{minipage}=\begin{minipage}[c]{3.5cm}
        anneau de valuation discrète de rang $1$ 
        \end{minipage}$}
        \obj(-80,-20)[2]{$\begin{minipage}[c]{3.5cm}
        anneau local noethérien régulier de dimension $n$
        \end{minipage}$}
        \obj(80,-20)[3]{$\begin{minipage}[c]{3.5cm}
        anneau de valuation discrète de rang $n$ 
        \end{minipage}$}
        \mor(-10,10)(-70,-10){$$}
        \mor(10,10)(70,-10){$$}
        \enddc
    \end{equation*}
\end{itemize}
\section*{\textsc{Dans cette thèse}}
Cette thèse est divisée en $4$ chapitres. Au début de chacun (exception faite pour le premier), il y a une section où les théorèmes principaux sont énoncés. Ces sections et la présente se chevauchent presque entièrement. L'auteur s'excuse pour cette redondance, mais il la considère inoffensive, alors qu'une réminiscence latine lui suggère que "melius abundare est quam deficere".
\subsection*{\textsc{Chapitre \texorpdfstring{$1$}{1}}}
Dans ce premier chapitre on va introduire les outils mathématiques dont on aura besoin dans la suite. Aucune section dans ce chapitre n'a la prétention d'être exhaustive. Dans tous les cas, on a cherché à souligner les idées principales des théories qui sont mentionnées. On profitera de ce chapitre pour fixer la notation et la terminologie que nous utiliserons. Aucun résultat, ni aucune construction, de ce chapitre n'est due à l'auteur.
\subsection*{\textsc{Chapitre \texorpdfstring{$2$}{2}}}
On introduit la catégorie monoïdale symétrique (ordinaire) des modèles de Landau-Ginzburg de dimension $n$ sur un schéma de base $S$. Les objets de cette catégorie sont des couples $(X,\underline{f})$ où $X$ est un schéma plat sur $S$ et $\underline{f}$ est une sectionne globale de $\Oc_X^{n}$. Nous étudions la dg-catégorie des singularités $\sing{X,\underline{f}}$. Plus précisément, à l'aide de calculs explicites sur le morphisme counité
\begin{equation*}
    \mathfrak{i}^*\mathfrak{i}_*\Ec\rightarrow \Ec \hspace{1cm} \Ec \in \cohbp{X_0}{X}
\end{equation*}
où $\mathfrak{i}:X\times^h_{\aff{n}{S}}S\rightarrow X$ et $\cohbp{X_0}{X}=\{\Ec \in \cohb{X_0}: \mathfrak{i}_*\Ec \in \perf{X}\}$, on peut montrer le théorème suivant: soit $K(B,\underline{f})$ l'algèbre de Koszul, i.e.la dg-algèbre $B[\varepsilon_1,\dots,\varepsilon_n]$ où les $\varepsilon_i$ sont des générateurs libres en degré $-1$ tels que 
\begin{equation*}
    \varepsilon_i^2=0 \hspace{0.5cm}
    d(\varepsilon_i)=f_i \hspace{0.5cm}
    \varepsilon_i \varepsilon_j+\varepsilon_j 
\varepsilon_i=0
\end{equation*}
\begin{theorem-nonf}{\textup{\ref{structure theorem}}}
Soit $(X=Spec(B),\underline{f})$ un modèle LG affine $n$-dimensionnel sur $S$. Alors, chaque objet dans  $\sing{B,\underline{f}}$ est un rétracte d'un objet représenté par un $K(B,\underline{f})$-dg module concentré en $n+1$ degrés. 
\end{theorem-nonf}
Note que $n+1$ est l'amplitude de l'algèbre de Koszul $K(B,\underline{f})$.
Successivement, on se spécialise au cas $n=1$ et on démontre le résultat suivant 
\begin{theorem-nonf}{\textup{\ref{structure theorem n=1}}}
Soit $(B,f)$ un modèle de Landau-Ginzburg (de dimension $1$) sur $S$. Alors, pour chaque $K(B,f)$-dg module $E$ dont le $B$-dg module sous-jacent est parfait, il y a une équivalence naturelle dans $\sing{B,f}$
\begin{equation*}
    E \simeq \begindc{\commdiag}[18]
    \obj(-40,0)[1]{$\Bigl ( \bigoplus_{i\in \Z}E_{2i}$}
    \obj(40,0)[2]{$\bigoplus_{i\in \Z}E_{2i+1}\Bigr )$}
    \mor(-30,-2)(30,-2){$d+h$}[\atright,\solidarrow]
    \mor(30,2)(-30,2){$d+h$}[\atright,\solidarrow]
    \enddc
\end{equation*}
où $d$ est le différentiel, $h$ est le morphisme induit par le générateur en degré $-1$ de l'algèbre de Koszul $K(B,f)$ et le complexe à droite est placé dans le degrés $0$ et $1$.
\end{theorem-nonf}
Ce théorème nous permet de généraliser l'équivalence d'Orlov (\cite{brtv}, \cite{buch}, \cite{bw12}, \cite{bw15}, \cite{efpo15},  \cite{orl04}, \cite{orl12}) à n'importe quel modèle de Landau-Ginzburg de dimension $1$ (on ne demande pas que le potentiel $f$ soit plat). Soit $\mf{B}{f}$ la catégorie des factorisations matricielles de $(B,f)$ (\cite{eis80}), i.e. la dg-catégorie $2$-périodique dont les objets sont des quadruplés $(E_0,E_1,\phi_0,\phi_1)$, où $E_0, E_1$ sont des $B$-modules projectifs de type fini et $\phi_0:E_0\rightarrow E_1, \phi_1:E_1\rightarrow E_0$ sont des morphismes $B$-linéaires tels que $\phi_0\circ \phi_1 = f\cdot id_{E_1}$ et $\phi_1\circ \phi_0=f\cdot id_{E_0}$.
\begin{corollary-nonf}{\textup{\ref{Orlov equivalence}}}
L'$\oo$-transformation naturelle
\begin{equation*}
    Orl^{-1,\otimes}:\sing{\bullet,\bullet}^{\otimes}\rightarrow \mf{\bullet}{\bullet}^{\otimes}
\end{equation*}
construite dans \cite[\S 2.4]{brtv}, qui est définie pour un modèle LG affine comme
\begin{equation*}
    Orl^{-1}:\sing{B,f}\rightarrow \mf{B}{f}
\end{equation*}
\begin{equation*}
    E\mapsto \begindc{\commdiag}[18]
    \obj(-40,0)[1]{$\Bigl ( \bigoplus_{i\in \Z}E_{2i}$}
    \obj(40,0)[2]{$\bigoplus_{i\in \Z}E_{2i+1}\Bigr )$}
    \mor(-30,-2)(30,-2){$d+h$}[\atright,\solidarrow]
    \mor(30,2)(-30,2){$d+h$}[\atright,\solidarrow]
    \enddc
\end{equation*}
est une $\oo$-équivalence naturelle.
\end{corollary-nonf}
On note qu'une factorisation matricielle définit un $K(B,f)$-dg module.
Cela définit (approximativement) la transformation naturelle inverse 
\begin{equation*}
    Orl: \mf{\bullet}{\bullet} \rightarrow \sing{\bullet,\bullet}
\end{equation*}
Le Théorème \ref{structure theorem n=1} nous dit que $Orl \circ Orl^{-1}$ est équivalent à l'identité. À la fin de ce chapitre nous suggérons qu'il pourrait être intéressant de considérer la localisation de $\sing{X,\underline{f}}$ par rapport aux opérateurs d'Eisenbud, qui ont été introduits par la première fois dans \cite{eis80}. Nous pensons que l'on peut aussi les décrire comme suit. Soit $X_k$ le lieu des zéros dérivés de $(f_1,\dots,f_{k-1},f_{k+1},\dots,f_n)$ et soit $\mathfrak{i}_{0k}:X_0 \rightarrow X_k$. Alors les opérateurs d'Eisenbud devraient correspondre aux morphismes naturels 
\begin{equation*}
    \chi_k:\Ec\rightarrow cofib(\mathfrak{i}_{0k}^*\mathfrak{i}_{0k*}\Ec\rightarrow \Ec)\simeq \Ec[2]
\end{equation*}
D'après des théorèmes de D.~Orlov (\cite{orl06}) et de J.~Burke-M.~Walker (\cite{bw15}), il y a une équivalence des dg-catégories (au moins quand $X$ est affine et $\underline{f}$ est une suite régulière, mais nous croyons que cela est vrai en général)
\begin{equation*}
    \sing{X,\underline{f}}\simeq \sing{\proj{n-1}{B},W_{\underline{f}}}
\end{equation*}
où $W_{\underline{f}}=f_1\cdot T_1+\dots +f_n\cdot T_n$. De plus, dans \cite{bw15}, les auteurs démontrent que les opérateurs d'Eisenbud  $\chi_k$ correspondent dans le terme de droite à la multiplication par les éléments $T_k$. On s'attend donc d'avoir l'équivalence suivante
\begin{equation*}
    \sing{X,\underline{f}}[\chi_1^{-1},\dots,\chi_n^{-1}]\simeq \mf{\gm{X}^{n-1}}{W_{\underline{f}|\gm{X}^{n-1}}}
\end{equation*}
Notons que ces deux dg-catégories sont $2$-périodiques, dotées des $n-1$ autoéquivalences: les $\chi_k\circ \chi_{1}^{-1}$ à gauche et les $T_k\circ T_1^{-1}$ à droite ($k\geq 2$). En d'autres termes, l'équivalence ci-dessus devrait être vraie dans la catégorie des $\perf{\Oc_S[u_1,\dots,u_n,u_1^{-1},\dots,u_n^{-1}]}$-modules, où les $u_i$'s sont des générateurs libres en degré $2$.

Ces faits sont encore investigués par l'auteur aujourd'hui.
\subsection*{\textsc{Chapitre \texorpdfstring{$3$}{3}}}
Le lien entre les dg-catégories des factorisations matricielles et la théorie des cycles évanescents est bien connu grâce à A.~Prygel (\cite{pr11}), T.~Dyckerhoff (\cite{dy11}), A.~Efimov (\cite{efi18}) et bien d'autres. Dans leur récent article \cite{brtv}, les auteurs démontrent une formule qui redonne la cohomologie $\ell$-adique de la dg-catégorie de la fibre spéciale d'un schéma régulier, propre et plat sur un trait hénselien strict excellent (voir la Section \ref{the l-adic realization of dg categories}) avec la partie invariante par inertie de la cohomologie $\ell$-adique évanescente $2$-périodique. Plus précisément, soit $S$ un trait hénselien strict excellent et soit $p:X\rightarrow S$ un schéma régulier, propre et plat sur $S$. Soit $\Phi_p(\Ql{,X}(\beta))$ le faisceau $\ell$-adique des cycles évanescents de $\Ql{,X}=\bigoplus_{i\in \Z}\Ql{,X}(i)[2i]$ et soit $(-)^{\textup{h}I}$ l'$\oo$-foncteur des points fixes à homotopie prés, où $I$ est le groupe d'inertie.
\begin{theorem-nonf}{\textup{\cite[Théorème 4.39]{brtv}}}
Il y a une équivalence de  $i_{\sigma}^*\Rlv_S(\sing{S,0})\simeq \Ql{,\sigma}^{\textup{h}I}\otimes_{\Ql{,\sigma}}\Ql{,\sigma}(\beta)$-modules
\begin{equation*}
    i_{\sigma}^*\Rlv_S(\sing{X_{\sigma}})\simeq p_{\sigma*}\Phi_p(\Ql{,X}(\beta))^{\textup{h}I}[-1]
\end{equation*}
\end{theorem-nonf}
Ce théorème, au-delà d'être intéressant, est crucial pour la stratégie de B.~Toën et G.~Vezzosi pour prouver la conjecture du conducteur de Bloch formulée dans \cite{bl85}. Pour plus de détails sur la stratégie de B.~Toën et G.~Vezzosi pour résoudre cette conjecture, voir \cite{tv17}, \cite{tv19a} and \cite{tv19b}.

Ce chapitre est né comme une tentative de généraliser cette formule au cas où $S$ est un anneau régulier, local et noethérien de dimension $n\geq 1$. Notons que le choix des générateurs $\pi_1,\dots,\pi_n$ de l'idéal maximal de $S$ définit une suite régulière de $X$. Donc, d'après le théorème de D.~Orlov et J.~Burke-M.~Walker que nous avons déjà mentionné, on a
\begin{equation*}
    \sing{X,\underline{\pi}\circ p}\simeq \sing{\proj{n-1}{X},W_{p\circ\underline{\pi}}}
\end{equation*}
L'objet à droite à l'avantage d'être la dg-catégorie des singularités d'un sous-schéma de $\proj{n-1}{X}$ de codimension $1$, mais $W_{p\circ \underline{\pi}}$ est une section globale de $\Oc(1)$, et donc il n'est pas clair à priori quelle est la bonne généralisation de $\Phi_p(\Ql{,X}(\beta))$. 

On définit la catégorie monoïdale symétrique des modèles LG tordus sur $(S,\Lc_S)$, où $\Lc_S$ est un fibré en droites sur $S$,  comme la catégorie dont les objets sont les couples $(X,s_X)$ où $X$ est un schéma (plat) sur $S$ et $s_X$ est une section globale de l'image réciproque de $\Lc_S$ à $X$. On considère le diagramme 
\begin{equation*}
    \begindc{\commdiag}[18]
    \obj(-60,10)[1]{$X_0$}
    \obj(0,10)[2]{$X$}
    \obj(60,10)[3]{$X_{\mathcal{U}}$}
    \obj(-60,-10)[4]{$X$}
    \obj(0,-10)[5]{$V_X=Tot(\Lc_X)$}
    \obj(60,-10)[6]{$\mathcal{U}_X=V_X-X$}
    \mor{1}{2}{$i$}
    \mor{3}{2}{$j$}
    \mor{1}{4}{$s_0$}
    \mor{2}{5}{$s_X$}
    \mor{3}{6}{$s_{\mathcal{U}}$}
    \mor{4}{5}{$i_0$}
    \mor{6}{5}{$j_0$}
    \enddc
\end{equation*}
Nous définissons le \emph{faisceau} $\ell$-\emph{adique des cycles évanescents invariants par monodromie} comme
\begin{equation*}
    \Phi_{(X,s_X)}^{\textup{mi}}(\Ql{}(\beta)):=cofib(i^*s_X^*j_{0*}\Ql{,\mathcal{U}_X}(\beta)\rightarrow i^*j_*s_{\mathcal{U}}^*\Ql{,\mathcal{U}_X}(\beta))
\end{equation*}
Si le fibré en droites est trivial et la section est induite par l'uniformisant d'un trait hénselien strict excellent, cette définition retrouve la partie invariante par inertie du faisceau $\ell$-adique des cycles évanescents.

Alors on démontre
\begin{theorem-nonf}{\textup{\ref{main theorem}}}
Soit $(X,s_X)$ un modèle LG tordu sur $(S,\Lc_S)$. On suppose que $X$ est un schéma régulier. Il y a une équivalence de $i^*s_X^*j_{0*}\Ql{,\mathcal{U}_X}(\beta)\simeq s_0^*\Rlv_X(\sing{X,0})$-modules dans $\shvl(X_0)$
\begin{equation*}
    i^*\Rlv_X(\sing{X,s_X})\simeq \Phi^{\textup{mi}}_{(X,s_X)}(\Ql{}(\beta))[-1]
\end{equation*}
\end{theorem-nonf}
On déduit le corollaire suivant du théorème précédent et du théorème de D.~Orlov et J.~Burke-M.~Walker, qui nous donne une réponse aux problème initial: généraliser la formule de \cite{brtv} au cas où la base est un anneaux régulier, local et noethérien de dimension $n\geq 1$.
\begin{corollary-nonf}{\textup{\ref{answer intial question}}}
Si $g:X\rightarrow S$ est un schéma affine régulier plat sur un anneau local noethérien de dimension $n$ $S$, il y a une équivalence de $\Rlv_X(\sing{X,0})$-modules dans $\shvl(X)$
\begin{equation*}
    \Rlv_X(\sing{X,\underline{\pi}\circ g})\simeq p_*i_*\Phi^{\textup{mi}}_{(\proj{n-1}{X},W_{\underline{\pi}\circ g})}[-1]
\end{equation*}
où $p:\proj{n-1}{X}\rightarrow X$, $i:V(W_{\underline{\pi}\circ g})\rightarrow \proj{n-1}{X}$ et $\underline{\pi}:S \rightarrow \aff{n}{S}$ est le morphism induit par une suite de générateurs du point fermé de $S$.
\end{corollary-nonf}
On termine avec deux sections supplémentaires. Dans la première, on explique qu'il devrait être possible de définir un formalisme de cycles évanescents sur $\ag{S}$. Le diagramme de base est, dans ce cas
\begin{equation*}
    \bgm{S}\hookrightarrow \ag{S} \hookleftarrow S
\end{equation*}
Pour construire une action qui joue le rôle de l'action du groupe d'inertie, on considère les morphismes de champs
\begin{equation*}
    \Theta_n:\ag{S} \rightarrow \ag{S}
\end{equation*}
\begin{equation*}
    (X,\Lc,s)\mapsto (X,\Lc^{\otimes n},s^{\otimes n})
\end{equation*}
Le faisceau $\Phi^{\textup{mi}}_{(X,s_X)}(\Ql{}(\beta))$ devrait être retrouvé si l'on applique l'analogue du foncteur des points fixes dans ce contexte. Dans la deuxième section, on donne quelques remarques sur le fait que l'on doit supposer $X$ régulier.

L'auteur a l'intention d'étudier ces aspects dans un avenir proche. 
\subsection*{\textsc{Chapitre \texorpdfstring{$4$}{4}}}
Dans ce chapitre final on essaye aussi de généraliser la formule de \cite{brtv}. Plutôt que se donner un schéma régulier sur un anneau noetherien de dimension $n$ comme dans le chapitre $3$, on considère la situation où le schéma vit au dessus de $R=\C \dsl x \dsr + y\cdot \C \drl x \drr \dsl y \dsr \subseteq \C \drl x \drr \drl y \drr$. Il s'agit de l'anneau des entiers de la valuation de rang $2$
\begin{equation*}
    v: \C \drl x \drr \drl y \drr^{\x}\rightarrow \Z \x \Z
\end{equation*}
\begin{equation*}
\sum_{i=n}^{\oo} \bigl ( \sum_{j=m_i}^{\oo} a_{ij}x^j\bigr )y^i \mapsto (n,m_n)
\end{equation*}
où on a mis l'ordre lexicographique sur $\Z\times \Z$. Le diagramme de base est
\begin{equation*}
    \begindc{\commdiag}[40]
\obj(0,15)[1]{$\sigma=Spec(\C)$}
\obj(30,15)[2]{$S_1=Spec(\C \dsl x \dsr)$}
\obj(60,15)[3]{$\varepsilon=Spec(\C \drl x \drr)$}
\obj(90,15)[4]{$\bar{\varepsilon}=Spec(\overline{\C \drl x \drr})$}
\obj(30,0)[5]{$S=Spec(R)$}
\obj(60,0)[6]{$U_0=Spec(\C \drl x \drr \dsl y \dsr)$}
\obj(90,0)[7]{$\bar{U}_0=Spec(\overline{\C \drl x \drr} \dsl y \dsr)$}
\obj(30,-15)[8]{$\eta=Spec(\C \drl x \drr \drl y \drr)$}
\obj(60,-15)[9]{$\eta^{mu}=Spec(\overline{\C \drl x \drr} \drl y \drr)$}
\obj(90,-15)[10]{$\bar{\eta}=Spec(\overline{\C \drl x \drr \drl y \drr})$}
\mor{1}{2}{$i_{01}$}[\atleft,\injectionarrow]
\mor{1}{5}{$i_0$}[\atright,\injectionarrow]
\mor{2}{5}{$i_1$}[\atleft,\injectionarrow]
\mor{3}{6}{$k$}[\atleft,\injectionarrow]
\mor{4}{7}{$\bar{k}$}[\atleft,\injectionarrow]
\mor{3}{2}{$l$}[\atright,\injectionarrow]
\mor{6}{5}{$j_{0}$}[\atright,\injectionarrow]
\mor{8}{5}{$j_{1}$}[\atright,\injectionarrow]
\mor{8}{6}{$j_{10}$}[\atright,\injectionarrow]
\mor{9}{7}{$j_{10}^{mu}$}[\atright,\injectionarrow]
\mor{4}{3}{$$}
\mor{7}{6}{$$}
\mor{9}{8}{$$}
\mor{10}{9}{$$}
\mor{10}{7}{$\bar{j_{10}}$}[\atright,\solidarrow]
\cmor((90,17)(90,21)(88,23)(60,23)(32,23)(30,21)(30,17)) \pdown(60,25){$\bar{l}$}
\cmor((88,1)(88,3)(86,5)(60,5)(34,5)(32,3)(32,1)) \pdown(75,7){$\bar{j}_0$}
\cmor((90,-17)(45,-22)(15,-21)(15,-5)(20,0)) \pright(60,-24){$\bar{j}_1$}
\enddc
\end{equation*}
Notons que la ligne en haut et la colonne à droite sont les diagrammes correspondants à la situation classique d'un trait hénselien strict. Comme le point fermé de $\overline{U}_0$ coïncide avec le point ouvert géométrique de $S_1$, on peut itérer la construction des cycles proches. Plus précisément, si $X$ est un $S$-schéma et $\Fc \in \shvl(X)$, on définit
\begin{itemize}
    \item $\Psi^{(1)}(\Fc):=i_{01}^*\bar{l}_*\Fc_{|X_{\bar{\varepsilon}}} \in \shvl(X_{\sigma})^{Gal(\bar{\eta}/\eta)}$ (cycles proches d'ordre $1$)
    \item $\Psi^{(2)}(\Fc):=i_{01}^*\bar{l}_*\bar{k}^*\bar{j}_{10*}\Fc_{|X_{\bar{\eta}}} \in \shvl(X_{\sigma})^{Gal(\bar{\eta}/\eta)}$ (cycles proches d'ordre $2$)
\end{itemize}
où $\shvl(X_{\sigma})^{Gal(\bar{\eta}/\eta)}$ est l'$\oo$-catégorie des faisceaux $\ell$-adiques sur $X_{\sigma}$ avec une action continue de $Gal(\bar{\eta}/\eta)$. Les morphismes que nous avons utilisés dans la définition sont ceux qu'on obtient par pullback le long du morphism $X\rightarrow S$.
Notons qu'on dispose d'un triangle $Gal(\bar{\eta}/\eta)$-équivariant
\begin{equation*}
    \begindc{\commdiag}[18]
    \obj(-20,10)[1]{$\Fc_{|X_{\sigma}}$}
    \obj(20,10)[2]{$\Psi^{(1)}(\Fc)$}
    \obj(0,-10)[3]{$\Psi^{(2)}(\Fc)$}
    \mor{1}{2}{$$}
    \mor{1}{3}{$$}
    \mor{2}{3}{$$}
    \enddc
\end{equation*}

Définissons 
\begin{itemize}
    \item $\Phi^{(1)}(\Fc):=cofib(\Fc_{|X_{\sigma}}\rightarrow \Psi^{(1)}(\Fc)) \in \shvl(X_{\sigma})^{Gal(\bar{\eta}/\eta)}$ (cycles évanescents d'ordre $1$)
    \item $\Phi^{(2)}(\Fc):=cofib(\Fc_{|X_{\sigma}}\rightarrow \Psi^{(2)}(\Fc)) \in \shvl(X_{\sigma})^{Gal(\bar{\eta}/\eta)}$ (cycles évanescents d'ordre $2$)
\end{itemize}
On relie les cycles évanescents invariants par rapport à $Gal(\bar{\eta}/\eta)$ de $\Ql{,X}(\beta)$ avec la cohomologie $\ell$-adique de $\sing{Y}$, où $Y=X\times_SSpec(R/y \cdot R)$, en supposant que $X$ soit régulier, propre, plat et de présentation finie sur $S$, que la réalisation $\ell$-adique existe aussi dans le cadre non-noetherian et que les résultats dans \cite[\S B.4]{pr11} valent dans notre contexte (pour plus de détails à ce propos, voir le paragraphe \ref{caution}).
\begin{theorem-nonf}\ref{theorem chpt 4} 
Il y a un suite fibre-cofibre dans $\shvl(\C)$
\begin{equation*}
    i_{0}^*\Rlv_{S}(\sing{Y})\rightarrow p_{\C*}\Phi^{(2)}(\Ql{}(\beta))[-1]^{\textup{h}Gal(\bar{\eta}/\eta)}\rightarrow cofib(\xi)
\end{equation*}
où $\xi$ est un morphism explicite. 
\end{theorem-nonf}

%% file: chapters/introduction.tex
Non-commutative geometry, as envisioned by M.~Kontsevich, is the study of dg categories that have some good properties. These good properties are those that we find if we look at the dg category associated to a scheme. Indeed, for any scheme $X$, we can consider the dg category of perfect complexes $\perf{X}$. This assignment can be enhanced to an $\oo$-functor
\begin{equation*}
 \perf{\bullet}:    \sch^{\textup{op}} \rightarrow \dgcatm
\end{equation*}
This $\oo$-functor is the bridge that connects the commutative world to the non-commutative one. Sometimes, $\perf{X}$ is therefore called the non-commutative shadow of $X$\footnote{see \cite{pa11} for instance}. It should be actually called \emph{a} non-commutative shadow of $X$. Indeed, there are other dg categories that we can associate to the same scheme, for example $\qcoh{X}$, $\cohb{X}$. They form a chain of inclusions of dg categories
\begin{equation*}
    \perf{X}\subseteq \cohb{X} \subseteq \qcoh{X}
\end{equation*}
Even if the assignments $X\mapsto \perf{X}$ and $X\mapsto \qcoh{X}$ are different a priori, they give us the same amount of information. We can recover the dg category of quasi-coherent sheaves on $X$ from that of perfect complexes by the formula
\begin{equation*}
    \qcoh{X}=\textup{Ind}(\perf{X})
\end{equation*}
In the opposite way, we can recover $\perf{X}$ from $\qcoh{X}$ by taking the compact objects in the latter dg category. Whence, the difference between these two dg categories does not contain any geometrical information. We could say that the difference between these two dg categories is purely of categorical nature.

However, the two non-commutative shadows $\perf{X}$ and $\cohb{X}$ differ in a more interesting way: it follows from a theorem due to M.~Auslander-D.A.~Buchsbaum (\cite[Theorem 4.1]{ab56}) and J.P.~Serre (\cite[Théorème 3]{se55}) that these two dg categories are equivalent if and only if $X$ is a regular scheme. Therefore, we can consider another non-commutative shadow of $X$, known as the dg category of singularities of $X$, that is defined as the dg quotient 
\begin{equation*}
    \sing{X}:=\cohb{X}/\perf{X}
\end{equation*} 
and measures (in some sense) how singular the space $X$ is. This also explains the choice of its name.

We will actually consider the dg category of singularities of pairs $(X,s)$ where $X$ is a scheme and $s$ is a global section of some vector bundle. The dg category of $(X,s)$, that we denote $\sing{X,s}$, is defined as the kernel of the dg functor
\begin{equation*}
    \mathfrak{i}_* : \sing{X_0}\rightarrow \sing{X}
\end{equation*}
induced by the pushforward along the embedding $\mathfrak{i}:X_0:=V(s)\rightarrow X$ of the (derived) zero locus of $s$ inside $X$. Roughly, we could think of it as a non-commutative invariant that measures how worse the singularities of $X_0$ are compared to those of $X$. The simplest case is when $X$ is a regular scheme, in which case $\sing{X,s}$ is just the dg categories of singularities of $X_0$ (as $\sing{X}\simeq 0$). We will use this basic fact very often in this thesis, even without explicitly mentioning it. 

We will investigate two different aspects of these dg categories
\begin{itemize}
    \item in the first part we will try to better understand how they look like. For example, we will show a generalization of Orlov's theorem which tells us that, whenever $s$ is a global section of $\Oc_X$, then $\sing{X,s}$ is equivalent to $\mf{X}{s}$, a second non-commutative shadow of the pair $(X,s)$ known as the dg category of matrix factorizations. The author likes to think about this equivalence in terms of an analogy with the isomorphism $\pi_1(S^1)\simeq \Z$. If we adopt a topologist's point of view, we may think of $\pi_1(S^1)$ as an object corresponding to a conceptual definition, while we interpret $\Z$ as a nice model of the former. Similarly, we think of $\sing{X,s}$ as the conceptual definition and of $\mf{X,s}$ as a nice model
    \item in the second part we will study the $\ell$-adic cohomology of $\sing{X,s}$ for a global section of a line bundle $\Lc$. In order to do so, we will need to introduce a generalization of inertia-invariant vanishing cycles, where the role of the disk is played by the total space of the line bundle. We will also investigate the $\ell$-adic cohomology of a certain dg category of singularities in the setting of a scheme over a discrete valuation ring of rank $2$. Both these approaches arise as attempts to generalize the main theorem in \cite{brtv}, and correspond to the following two higher dimensional generalizations of the definition of a DVR
    \begin{equation*}
        \begindc{\commdiag}[18]
        \obj(0,20)[1]{$\begin{minipage}[c]{3cm}
        noetherian local regular ring of dimension $1$
        \end{minipage}=\begin{minipage}[c]{3cm}
        rank $1$ discrete valuation ring 
        \end{minipage}$}
        \obj(-80,-20)[2]{$\begin{minipage}[c]{3cm}
        noetherian local regular ring of dimension $n$
        \end{minipage}$}
        \obj(80,-20)[3]{$\begin{minipage}[c]{3cm}
        rank $n$ discrete valuation ring 
        \end{minipage}$}
        \mor(-10,10)(-70,-10){$$}
        \mor(10,10)(70,-10){$$}
        \enddc
    \end{equation*}
\end{itemize}
\section*{\textsc{What is in this thesis}}
This thesis is divided into $4$ chapters. At the beginning of each chapter (except chapter $1$) there is a section where the main results are stated. These sections and the present one overlap almost completely. The author apologises for this redundancy, but he considers it harmless, while a latin reminescence of his suggests him that "melius abundare est quam deficere".
\subsection*{\textsc{Chapter \texorpdfstring{$1$}{1}}}
In this first chapter we introduce the mathematical tools that we are going to need. None of the sections in this chapter claims exhaustiveness. However, we have tried to highlight some of the main ideas lying behind the sophisticated theories that are mentioned. In any case, the existence of this chapter is justified as it allows us to fix some notation and terminology. No new results or constructions due to the author appear in this chapter.
\subsection*{\textsc{Chapter 2}}
We introduce the symmetric monoidal (ordinary) category of $n$-dimensional Landau-Ginzburg models over a base scheme $S$. The objects of this category are pairs $(X,\underline{f})$ where $X$ is a flat scheme over $S$ and $\underline{f}$ is a global section of $\Oc_X^{n}$. We study the dg category of relative singularities $\sing{X,\underline{f}}$. More precisely, by explicit computations on the counit morphism
\begin{equation*}
    \mathfrak{i}^*\mathfrak{i}_*\Ec\rightarrow \Ec \hspace{1cm} \Ec \in \cohbp{X_0}{X}
\end{equation*}
where $\mathfrak{i}:X\times^h_{\aff{n}{S}}S\rightarrow X$ and $\cohbp{X_0}{X}=\{\Ec \in \cohb{X_0}: \mathfrak{i}_*\Ec \in \perf{X}\}$, we are able to show the following theorem: recall that the Koszul algebra $K(B,\underline{f})$ is the dg algebra $B[\varepsilon_1,\dots,\varepsilon_n]$ where $\varepsilon_i$ are free generators in degree $-1$ and satisfy 
\begin{equation*}
    \varepsilon_i^2=0 \hspace{0.5cm}
    d(\varepsilon_i)=f_i \hspace{0.5cm}
    \varepsilon_i \varepsilon_j+\varepsilon_j 
\varepsilon_i=0
\end{equation*}
\begin{theorem-non}{\textup{\ref{structure theorem}}}
Let $(X=Spec(B),\underline{f})$ be an affine $n$-dimensional LG model over $S$. Then, every object in $\sing{B,\underline{f}}$ is a retract of an object represented by a $K(B,\underline{f})$-dg module concentrated in $n+1$ degrees. 
\end{theorem-non}
Notice that $n+1$ is exactly the amplitude of the Koszul algebra $K(B,\underline{f})$.

We then specialize to the case $n=1$ and prove
\begin{theorem-non}{\textup{\ref{structure theorem n=1}}}
Let $(B,f)$ be a ($1$-dimensional) LG model over $S$. Then, for each $K(B,f)$-dg module $E$ whose underling $B$-dg module is perfect, there exists a natural equivalence in $\sing{B,f}$
\begin{equation*}
    E \simeq \begindc{\commdiag}[18]
    \obj(-40,0)[1]{$\Bigl ( \bigoplus_{i\in \Z}E_{2i}$}
    \obj(40,0)[2]{$\bigoplus_{i\in \Z}E_{2i+1}\Bigr )$}
    \mor(-30,-2)(30,-2){$d+h$}[\atright,\solidarrow]
    \mor(30,2)(-30,2){$d+h$}[\atright,\solidarrow]
    \enddc
\end{equation*}
where $d$ is the differential, $h$ is the morphism induced by the generator in degree $-1$ of the Koszul algebra $K(B,f)$ and the complex on the r.h.s. lies in degrees $0$ and $1$.
\end{theorem-non}
This theorem allows us to generalize Orlov's equivalence (\cite{brtv}, \cite{buch}, \cite{bw12}, \cite{bw15}, \cite{efpo15},  \cite{orl04}, \cite{orl12}) to all $1$-dimensional LG models (that is, the potential $f$ does not need to be flat). Let $\mf{B}{f}$ be the category of matrix factorizations of $(B,f)$ (\cite{eis80}), i.e. the $2$-periodic dg category whose objects are quadruplets $(E_0,E_1,\phi_0,\phi_1)$, where $E_0, E_1$ are projective $B$-modules of finite type and $\phi_0:E_0\rightarrow E_1, \phi_1:E_1\rightarrow E_0$ are $B$-linear morphisms such that $\phi_0\circ \phi_1 = f\cdot id_{E_1}$ and $\phi_1\circ \phi_0=f\cdot id_{E_0}$.
\begin{corollary-non}{\textup{\ref{Orlov equivalence}}}
The lax-monoidal $\oo$-natural transformation
\begin{equation*}
    Orl^{-1,\otimes}:\sing{\bullet,\bullet}^{\otimes}\rightarrow \mf{\bullet}{\bullet}^{\otimes}
\end{equation*}
constructed in \cite[\S 2.4]{brtv}, which is defined for an affine LG model as
\begin{equation*}
    Orl^{-1}:\sing{B,f}\rightarrow \mf{B}{f}
\end{equation*}
\begin{equation*}
    E\mapsto \begindc{\commdiag}[18]
    \obj(-40,0)[1]{$\Bigl ( \bigoplus_{i\in \Z}E_{2i}$}
    \obj(40,0)[2]{$\bigoplus_{i\in \Z}E_{2i+1}\Bigr )$}
    \mor(-30,-2)(30,-2){$d+h$}[\atright,\solidarrow]
    \mor(30,2)(-30,2){$d+h$}[\atright,\solidarrow]
    \enddc
\end{equation*}
is an $\oo$-natural equivalence.
\end{corollary-non}
Notice that a matrix factorization can be used to define a $K(B,f)$-dg module by placing $E_0$ in degree $0$, $E_1$ in degree $-1$ and by declaring $\phi_1$ to be the differential and $\phi_0$ to be the morphism induced by the generator in degree $-1$ of the Koszul algebra. This (roughly) defines the inverse natural equivalence
\begin{equation*}
    Orl: \mf{\bullet}{\bullet} \rightarrow \sing{\bullet,\bullet}
\end{equation*}
Then, the content of Theorem \ref{structure theorem n=1} is that $Orl \circ Orl^{-1}$ is equivalent to the identity.

At the end of this chapter we suggest that it might be interesting to consider the localization of $\sing{X,\underline{f}}$ with respect to the Eisenbud operators. These were first introduced in \cite{eis80}, but we believe that they can also be described as follows. Let $X_k$ denote the derived zero locus of $(f_1,\dots,f_{k-1},f_{k+1},\dots,f_n)$ and label $\mathfrak{i}_{0k}:X_0 \rightarrow X_k$. The Esinbud operators should correspond to the natural morphisms
\begin{equation*}
    \chi_k:\Ec\rightarrow cofib(\mathfrak{i}_{0k}^*\mathfrak{i}_{0k*}\Ec\rightarrow \Ec)\simeq \Ec[2]
\end{equation*}
After theorems of D.~Orlov (\cite{orl06}) and J.~Burke-M.~Walker (\cite{bw15}), there is an equivalence of dg cateories (at least when $X$ is affine and $\underline{f}$ is a regular sequence, even if we believe that this holds true in general)
\begin{equation*}
    \sing{X,\underline{f}}\simeq \sing{\proj{n-1}{B},W_{\underline{f}}}
\end{equation*}
where $W_{\underline{f}}=f_1\cdot T_1+\dots +f_n\cdot T_n$. Moreover, in \cite{bw15}, the authors prove that Eisenbud operators $\chi_k$ correspond to multiplication by the elements $T_k$ on the r.h.s. We then expect that the following equivalence holds
\begin{equation*}
    \sing{X,\underline{f}}[\chi_1^{-1},\dots,\chi_n^{-1}]\simeq \mf{\gm{X}^{n-1}}{W_{\underline{f}|\gm{X}^{n-1}}}
\end{equation*}
Notice that both these dg categories are two-periodic, endowed with $n-1$ autoequivalences: the $\chi_k\circ \chi_{1}^{-1}$ on the l.h.s. and $T_k\circ T_1^{-1}$ on the r.h.s. ($k\geq 2$). In other words, the equivalence above should hold in the category of $\perf{\Oc_S[u_1,\dots,u_n,u_1^{-1},\dots,u_n^{-1}]}$-modules, where the $u_i$'s are free generators in degree $2$.

These facts are currently under investigation.

\subsection*{\textsc{Chapter 3}}
The link between the dg category of matrix factorizations and the theory of vanishing cycles is known after the work of A.~Preygel (\cite{pr11}), T.~Dyckerhoff (\cite{dy11}), A.~Efimov (\cite{efi18}) and many others. In their recent paper \cite{brtv}, the authors provide a formula that relates the $\ell$-adic cohomology of the dg category of singularities of the special fiber of a regular, proper and flat scheme over an excellent strictly henselian trait (see Section \ref{the l-adic realization of dg categories}) with the inertia-invariant, $2$-periodic $\ell$-adic vanishing cohomology. More precisely, let $S$ be a strictly henselian excellent trait and let $p:X\rightarrow S$ a proper, flat and regular $S$-scheme. Let $X_{\sigma}$ denote the special fiber. Let $\Phi_p(\Ql{,X}(\beta))$ denote the $\ell$-adic sheaf of vanishing cycles of $\Ql{,X}(\beta)=\bigoplus_{i\in \Z}\Ql{,X}(i)[2i]$ and let $(-)^{\textup{h}I}$ denote the homotopy fixed point $\oo$-functor, where $I$ is the inertia group.
\begin{theorem-non}{\textup{\cite[Theorem 4.39]{brtv}}}
There is an equivalence of $i_{\sigma}^*\Rlv_S(\sing{S,0})\simeq \Ql{,\sigma}^{\textup{h}I}\otimes_{\Ql{,\sigma}}\Ql{,\sigma}(\beta)$-modules
\begin{equation*}
    i_{\sigma}^*\Rlv_S(\sing{X_{\sigma}})\simeq p_{\sigma*}\Phi_p(\Ql{,X}(\beta))^{\textup{h}I}[-1]
\end{equation*}
\end{theorem-non}
This theorem, in addition to having its own interest, is of major importance for B.~Toën and G.~Vezzosi's strategy to prove the Bloch's conductor conjecture, formulated in \cite{bl85}. For more details about the B.~Toën and G.~Vezzosi's strategy to solve this conjecture, we refer to \cite{tv17}, \cite{tv19a}, \cite{tv19b}.

This chapter arises as an attempt to generalize this formula to the situation in which $S$ is a noetherian regular local ring of dimension $n\geq 1$. Notice that the choice of generators $\pi_1,\dots,\pi_n$ of the maximal ideal of $S$ defines a regular sequence on $X$. Therefore, by the theorem due to D.~Orlov and J.~Burke-M.~Walker already mentioned above, we have
\begin{equation*}
    \sing{X,\underline{\pi}\circ p}\simeq \sing{\proj{n-1}{X},W_{p\circ\underline{\pi}}}
\end{equation*}
While the object on the r.h.s. has the advantage of being the dg category of singularities of a closed subscheme of $\proj{n-1}{X}$ of codimension $1$, $W_{p\circ\underline{\pi}}$ is a global section of $\Oc(1)$ and therefore it is not clear a priori which is the correct generalization of $\Phi_p(\Ql{,X}(\beta))$.

We define the symmetric monoidal category of twisted LG models over $(S,\Lc_S)$, $\Lc_S$ being a line bundle over $S$, as the category whose objects are pairs $(X,s_X)$, where $X$ is a (flat) $S$-scheme and $s_X$ is a global section of the pullback of $\Lc_S$ to $X$. Consider the diagram
\begin{equation*}
    \begindc{\commdiag}[18]
    \obj(-60,10)[1]{$X_0$}
    \obj(0,10)[2]{$X$}
    \obj(60,10)[3]{$X_{\mathcal{U}}$}
    \obj(-60,-10)[4]{$X$}
    \obj(0,-10)[5]{$V_X=Tot(\Lc_X)$}
    \obj(60,-10)[6]{$\mathcal{U}_X=V_X-X$}
    \mor{1}{2}{$i$}
    \mor{3}{2}{$j$}
    \mor{1}{4}{$s_0$}
    \mor{2}{5}{$s_X$}
    \mor{3}{6}{$s_{\mathcal{U}}$}
    \mor{4}{5}{$i_0$}
    \mor{6}{5}{$j_0$}
    \enddc
\end{equation*}
We define the $\ell$-\emph{adic sheaf of monodromy-invariant vanishing cycles} of $(X,s_X)$ as
\begin{equation*}
    \Phi_{(X,s_X)}^{\textup{mi}}(\Ql{}(\beta)):=cofib(i^*s_X^*j_{0*}\Ql{,\mathcal{U}_X}(\beta)\rightarrow i^*j_*s_{\mathcal{U}}^*\Ql{,\mathcal{U}_X}(\beta))
\end{equation*}
When the line bundle is trivial and the section is induced by the uniformiser of a strictly henselian excellent trait, this definition recovers the $\ell$-adic sheaf of inertia invariant vanishing cycles.

We are then able to prove
\begin{theorem-non}{\textup{\ref{main theorem}}}
Let $(X,s_X)$ be a twisted LG model over $(S,\Lc_S)$ and assume that $X$ is a regular scheme. There is an equivalence of $i^*s_X^*j_{0*}\Ql{,\mathcal{U}_X}(\beta)\simeq s_0^*\Rlv_X(\sing{X,0})$-modules in $\shvl(X_0)$
\begin{equation*}
    i^*\Rlv_X(\sing{X,s_X})\simeq \Phi^{\textup{mi}}_{(X,s_X)}(\Ql{}(\beta))[-1]
\end{equation*}
\end{theorem-non}
The following corollary can be deduced by the previous theorem and by the theorem due to D.~Orlov and J.Burke-M.~Walker. It provides us an answer to the initial problem: to generalise the formula in \cite{brtv} to the case where the base is a noetherian local regular ring of dimension $n\geq 1$.
\begin{corollary-non}{\textup{\ref{answer intial question}}}
If $g:X\rightarrow S$ is a regular, affine flat scheme over a noetherian regular local ring of dimension $n$ $S$, there is an equivalence of $\Rlv_X(\sing{X,0})$-modules in $\shvl(X)$
\begin{equation*}
    \Rlv_X(\sing{X,\underline{\pi}\circ g})\simeq p_*i_*\Phi^{\textup{mi}}_{(\proj{n-1}{X},W_{\underline{\pi}\circ g})}[-1]
\end{equation*}
where $p:\proj{n-1}{X}\rightarrow X$, $i:V(W_{\underline{f}})\rightarrow \proj{n-1}{X}$ and $\underline{\pi}:S\rightarrow \aff{n}{S}$ is the morphism induced by a collection of generators of the closed point of $S$.
\end{corollary-non}
We conclude this chapter with two extra sections. In the first one we explain that it should be possible to define a formalism of vanishing cycles over $\ag{S}$. The base diagram then is
\begin{equation*}
    \bgm{S}\hookrightarrow \ag{S} \hookleftarrow S
\end{equation*}
To construct an action, which replaces of the action of the inertia group in the classical picture, one then considers the morphisms of stacks
\begin{equation*}
    \Theta_n:\ag{S} \rightarrow \ag{S}
\end{equation*}
\begin{equation*}
    (X,\Lc,s)\mapsto (X,\Lc^{\otimes n},s^{\otimes n})
\end{equation*}
Then $\Phi^{\textup{mi}}_{(X,s_X)}(\Ql{}(\beta))$ should be recovered if we apply the analogous of taking inertia fixed points.
In the second section, we comment the regularity hypothesis that we, as well as the authors in \cite{brtv}, impose on $X$.

The author intends to investigate more accurately these aspects in the near future.
\subsection*{\textsc{Chapter 4}}
In this final chapter we also try to generalize the formula proved in \cite{brtv}. Instead of considering schemes over a noetherian regular local ring of dimension $n$ as in Chapter 3, we consider the situation where the scheme lies over $R=\C \dsl x \dsr + y\cdot \C \drl x \drr \dsl y \dsr \subseteq \C \drl x \drr \drl y \drr$. This is the ring of integers of the rank $2$-valuation
\begin{equation*}
    v: \C \drl x \drr \drl y \drr^{\x}\rightarrow \Z \x \Z
\end{equation*}
\begin{equation*}
\sum_{i=n}^{\oo} \bigl ( \sum_{j=m_i}^{\oo} a_{ij}x^j\bigr )y^i \mapsto (n,m_n)
\end{equation*}
where we have endowed $\Z \times \Z$ with the lexicographic order. The base diagram is
\begin{equation*}
        \begindc{\commdiag}[40]
\obj(0,15)[1]{$\sigma=Spec(\C)$}
\obj(30,15)[2]{$S_1=Spec(\C \dsl x \dsr)$}
\obj(60,15)[3]{$\varepsilon=Spec(\C \drl x \drr)$}
\obj(90,15)[4]{$\bar{\varepsilon}=Spec(\overline{\C \drl x \drr})$}
\obj(30,0)[5]{$S=Spec(R)$}
\obj(60,0)[6]{$U_0=Spec(\C \drl x \drr \dsl y \dsr)$}
\obj(90,0)[7]{$\bar{U}_0=Spec(\overline{\C \drl x \drr} \dsl y \dsr)$}
\obj(30,-15)[8]{$\eta=Spec(\C \drl x \drr \drl y \drr)$}
\obj(60,-15)[9]{$\eta^{mu}=Spec(\overline{\C \drl x \drr} \drl y \drr)$}
\obj(90,-15)[10]{$\bar{\eta}=Spec(\overline{\C \drl x \drr \drl y \drr})$}
\mor{1}{2}{$i_{01}$}[\atleft,\injectionarrow]
\mor{1}{5}{$i_0$}[\atright,\injectionarrow]
\mor{2}{5}{$i_1$}[\atleft,\injectionarrow]
\mor{3}{6}{$k$}[\atleft,\injectionarrow]
\mor{4}{7}{$\bar{k}$}[\atleft,\injectionarrow]
\mor{3}{2}{$l$}[\atright,\injectionarrow]
\mor{6}{5}{$j_{0}$}[\atright,\injectionarrow]
\mor{8}{5}{$j_{1}$}[\atright,\injectionarrow]
\mor{8}{6}{$j_{10}$}[\atright,\injectionarrow]
\mor{9}{7}{$j_{10}^{mu}$}[\atright,\injectionarrow]
\mor{4}{3}{$$}
\mor{7}{6}{$$}
\mor{9}{8}{$$}
\mor{10}{9}{$$}
\mor{10}{7}{$\bar{j_{10}}$}[\atright,\solidarrow]
\cmor((90,17)(90,21)(88,23)(60,23)(32,23)(30,21)(30,17)) \pdown(60,25){$\bar{l}$}
\cmor((88,1)(88,3)(86,5)(60,5)(34,5)(32,3)(32,1)) \pdown(75,7){$\bar{j}_0$}
\cmor((90,-17)(45,-22)(15,-21)(15,-5)(20,0)) \pright(60,-24){$\bar{j}_1$}
\enddc
\end{equation*}
Notice that the top row and the column on the right are the diagrams that correspond to the classical setup of a strictly henselian trait. Since the closed point of $\bar{U}_0$ coincides with the geometric generic point of $S_1$, we can iterate the nearby cycles construction. More precisely, if $X$ is an $S$-scheme and $\Fc \in \shvl(X)$, we can define
\begin{itemize}
    \item $\Psi^{(1)}(\Fc):=i_{01}^*\bar{l}_*\Fc_{|X_{\bar{\varepsilon}}} \in \shvl(X_{\sigma})^{Gal(\bar{\eta}/\eta)}$ (first order nearby cycles)
    \item $\Psi^{(2)}(\Fc):=i_{01}^*\bar{l}_*\bar{k}^*\bar{j}_{10*}\Fc_{|X_{\bar{\eta}}} \in \shvl(X_{\sigma})^{Gal(\bar{\eta}/\eta)}$ (second order nearby cycles)
\end{itemize}
where $\shvl(X_{\sigma})^{Gal(\bar{\eta}/\eta)}$ is the $\oo$-category of $\ell$-adic sheaves on $X_{\sigma}$ endowed with a continuous action of $Gal(\bar{\eta}/\eta)$. The morphisms we used in the previous definition are those that arise by pulling back the diagram above along the structure morphism $X\rightarrow S$.
Notice that we have a $Gal(\bar{\eta}/\eta)$-equivariant triangle
\begin{equation*}
    \begindc{\commdiag}[18]
    \obj(-20,10)[1]{$\Fc_{|X_{\sigma}}$}
    \obj(20,10)[2]{$\Psi^{(1)}(\Fc)$}
    \obj(0,-10)[3]{$\Psi^{(2)}(\Fc)$}
    \mor{1}{2}{$$}
    \mor{1}{3}{$$}
    \mor{2}{3}{$$}
    \enddc
\end{equation*}

Define
\begin{itemize}
    \item $\Phi^{(1)}(\Fc):=cofib(\Fc_{|X_{\sigma}}\rightarrow \Psi^{(1)}(\Fc)) \in \shvl(X_{\sigma})^{Gal(\bar{\eta}/\eta)}$ (first order vanishing cycles)
    \item $\Phi^{(2)}(\Fc):=cofib(\Fc_{|X_{\sigma}}\rightarrow \Psi^{(2)}(\Fc)) \in \shvl(X_{\sigma})^{Gal(\bar{\eta}/\eta)}$ (second order vanishing cycles)
\end{itemize}
We then link $Gal(\bar{\eta}/\eta)$-invariant vanishing cycles of $\Ql{,X}(\beta)$ with the $\ell$-adic realization of $\sing{Y}$, where $Y=X\times_SSpec(R/y\cdot R)$, under the assumption that $X$ is a regular, proper, finitely presented flat $S$-scheme and, that the $\ell$-adic realization $\oo$-functor exists in the non-noetherian case too and that the results in \cite[\S B.4]{pr11} extend to our context (see Caution \ref{caution} for more details).
\begin{theorem-non}{\textup{\ref{theorem chpt 4}}}
There is a fiber-cofiber sequence in $\shvl(\C)$
\begin{equation*}
    i_{0}^*\Rlv_{S}(\sing{Y})\rightarrow p_{\C*}\Phi^{(2)}(\Ql{}(\beta))[-1]^{\textup{h}Gal(\bar{\eta}/\eta)}\rightarrow cofib(\xi)
\end{equation*}
where $\xi$ is an explicit morphism.
\end{theorem-non}

\newpage
\section*{Acknowledgements}
My deepest gratitude goes to my supervisors, Bertrand Toën and Gabriele Vezzosi. They taught me almost everything I know and they showed their enthusiasm for every little goal I achieved, which was determinant as this gave me the confidence I needed in the most stressful moments of the past years. Their generosity in sharing their ideas and insights is the backbone upon which all this work relies on. 

It is a pleasure to thank Bernhard Keller and Tobias Dyckerhoff, who honored me by accepting to referee my thesis and to be part of the jury together with Sarah Scherotzke, Joseph Tapia and Frédéric Déglise, to whom I also express my gratitude for having agreed to be part of the committee and for all the beautiful and stimulating mathematical conversations I had the pleasure to have with them over the years.

Many people contributed to this work, one way or another. I will not be able to mention them all, and I ask who will not appear to forgive me.

During the last few years I have benefited many times from the precious advice of Mathieu Anel, Ben Antieau, Dario Beraldo, Marcello Bernardara, Denis-Charles Cisinski, Luca Dall’Ava, Benjamin Hennion, Moritz Kerz, Marc Levine, Étienne Mann, Valerio Melani, Alberto Merici, Joan Milles, Tasos Moulinous, Tony Pantev, Riccardo Pengo, Mauro Porta, Arthur Renaudineau, Marco Robalo, Ed Segal, Nicolò Sibilla and Michel Vaquié. I thank all of them for their interest in my work, for the things they explained to me along the way and for their friendship.

This adventure wouldn’t have been the same if not shared with my “brothers and sisters in thesis”, Jorge Antonio, Elena Dimitriadis Bermejo, Ludovic Monier and Leyth Akrout Tabib.

Merci beaucoup à, Tamara Azaiez, Martine Labruyère, Agnès Requis et Céline Rocca pour leur aide et pour leur compréhension, qu’elles ne m’ont jamais fait manquer.

Thanks to Florian Bertuol, Laurence Fourrier, Valentin Huguin, Van Tu Le, Jules Martel, Dominique Mattei, Jasmin Raissy and Sonny Willets, who contributed to make me feel at home at the IMT.

The people of the Mathematical Sciences Research Institut provided me with the perfect environment during my stay in January and February 2019, and I warmly thank everyone of them for this.

Thank you to my lifetime friends Leonardo, Francesco, Giulio, Davide, Lorenzo, Pietro, Victoria and Francesco for always being there and for all the laughs we had together through the years.

It is hard to tell how much I owe to Alice. Her support and her love have been, and still are, crucial to me.

Infine, un ringraziamento speciale va alla mia famiglia e in particolare a mia mamma Elena, a mio babbo Andrea, a mio fratello Alessandro, a mia nonna Giuliana e a mio nonno Giuliano: senza di voi tutto questo non sarebbe stato possibile.

%% file: chapters/preliminaries.tex
\pagenumbering{arabic}
In this preliminary chapter we will briefly recall some of the mathematical tools that we will need later. We will also fix the notation that we will use in the other chapters of this thesis. 
\section{\textsc{Some notation}}
\begin{itemize}
    \item we will ignore set-theoretical issues, at least in the first Chapter. However, the reader should not worry about this: these were considered in the original sources, where the approach of Grothendieck's universes is adopted
    \item $\sm$ denotes the category of smooth schemes of finite type over $S$
    \item $\sch$ denotes the category of schemes of finite type over $S$
    \item $\oo$-category will always mean $(\oo,1)$-category
    \item $\mathcal{S}$ denotes the $\oo$-category of spaces
\end{itemize}
\section{\textsc{Reminders on the theories of  \texorpdfstring{$\oo$}{oo}-categories and of dg categories}}
We will freely use the language of $\oo$-categories which has been developed in \cite{htt} and \cite{ha}, i.e. their incarnation as quasi-categories (in Joyal's terminology). Of major importance for us are dg categories, that should be thought as $\oo$
-enhancements of triangulated categories. For a more detailed account on these objects we refer the reader to \cite{htt}, \cite{ha}, \cite{gro10}, \cite{ke06}, \cite{to08} and/or \cite{ro14}.
\subsection{\textsc{\texorpdfstring{$\oo$}{oo}-categories}}
Recall that an $\oo$-category is a simplicial set that satisfies the inner horns condition:
\begin{defn}
Let $\Cc$ be a simplicial set. It is said to be an $\oo$-category if for every $n\geq2$ and for every $0<i<n$ there exists a dotted arrow making the diagram below commute
\begin{equation*}
    \begindc{\commdiag}[18]
    \obj(-10,10)[1]{$\Lambda_i^n$}
    \obj(10,10)[2]{$\Cc$}
    \obj(-10,-10)[3]{$\Delta^n$}
    \mor{1}{2}{$$}
    \mor{1}{3}{$$}
    \mor{3}{2}{$$}[\atright,\dashArrow]
    \enddc
\end{equation*}
\end{defn}

One of the main merits of this definition is that it puts under the same roof ordinary categories and topological spaces. Indeed, one can see that those $\oo$-categories for which there exists an unique dotted arrow satisfying the lifting property are those which arise as nerves of ordinary categories. On the other hand, if the lifting property is satisfied for all $0\leq i \leq n$, then $\Cc$ is a Kan complex, which is just a simplicial terminology for topological space. 

One should think of $n$-simplices of an $\oo$-category as the $n$-morphisms of a higher category. In particular, $0$-simplices correspond to the objects and $1$-simplices to the morphisms. It is worth mentioning that $\oo$-categories model those higher categories for which all $n$-morphism are invertible whenever $n\geq 2$. 
Every $\oo$-category $\Cc$ can be used to define an ordinary category $\textup{h}\Cc$, whose objects coincide with those of $\Cc$ and morphisms are homotopy classes of $1$-simplices. It should be thought as a (brutal) truncation of $\Cc$.

An $\oo$-functor is just a morphism of simplicial sets in this model.

Lots of efforts have been made to make sense, in the theory of $\oo$-categories, of the usual constructions one has at his disposal in the world of ordinary categories. In particular, one disposes of the Yoneda embedding, of the adjoint functor theorem, of a solid theory of symmetric monoidal $\oo$-categories etc. To give an account of this machinery goes beyond the scope of this work. Many excellent accounts on the subject already exist and therefore it doesn't seem necessary to spend more words on this matter here. We just point out some of the sources we have in mind: \cite{htt}, \cite{ha}, \cite{ro14}, \cite{gro10}, \cite{jo08}.

\subsection{\textsc{dg categories}}
\begin{rmk}
For more details on the theory of dg categories, we invite the reader to consult \cite{ke06}, \cite{to11} and/or \cite{ro14}.
\end{rmk}
Let $A$ be a commutative ring.
\begin{defn}
An $A$-linear dg category $\Cc$ is a category such that for every pair of objects $x,y$, $\hm_{\Cc}(x,y)$ is a cochain complex of $A$-modules. Moreover, the composition morphisms are required to be morphisms of cochain complexes.

A dg functor $F:\Cc \rightarrow \Cc'$ between two dg categories is a functor such that, for every $x,y\in \Cc$, $\hm_{\Cc}(x,y)\rightarrow \hm_{\Cc'}(F(x),F(y))$ is a morphism of complexes. 
\end{defn}

\begin{ex}
The basic example of dg category is $\ch$, the category of cochain complexes over $A$.
\end{ex}

To any dg category $\Cc$ we can associate one ordinary category $\textup{H}^0(\Cc)$: its objects are those of $\Cc$ and the sets of morphisms are defined as $\hm_{\textup{H}^0(\Cc)}(x,y)=\textup{H}^0(\hm_{\Cc}(x,y))$. This construction is the analogue of the one mentioned in the paragraph above in the context of $\oo$-categories.
$\textup{H}^0(\Cc)$ is called the homotopy category of $\Cc$.

As usual, one should think of a cochain complexes as a vessel for carrying information about its cohomology\footnote{this is a quote: \cite[p.27]{gl18}}. Following this philosophy, the correct notion of equivalence in the context of dg categories is the following one:
\begin{defn}
A Dwyer-Kan (DK for short) equivalence of dg categories is a dg functor $F:\Cc\rightarrow \Cc'$ such that:
\begin{itemize}
    \item for every $x,y \in \Cc$, $\hm_{\Cc}(x,y)\rightarrow \hm_{\Cc'}(F(x),F(y))$ is a quasi-isomorphism
    \item the functor induced on the homotopy categories $\textup{H}^0(F):\textup{H}^0(\Cc)\rightarrow \textup{H}^0(\Cc')$ is essentially surjective
\end{itemize}
\end{defn}
Let $\dgcat$ be the category of small $A$-linear dg categories together with $A$-linear dg functors. This category can be endowed with a  combinatorial model category structure, where weak equivalences are DK-equivalences (see \cite{tab05}). The underlying $\oo$-category of this model category coincides with the $\oo$-localization of $\dgcat$ with respect to the class of DK-equivalences. We will denote this $\oo$-category by $\dgcatdk$. 

Recall the notion of modules over a dg category $\Cc$:
\begin{defn}
A $\Cc$-module is a dg functor $F:\Cc \rightarrow \ch$.
\end{defn}
Another very important class of dg functors is that of Morita-equivalences: 
\begin{defn}
We say that a dg functor $\Cc\rightarrow \Cc'$ is a Morita equivalence if the induced dg functor at the level of dg categories of perfect dg modules is a DK-equivalence.
\end{defn} 
Every DK-equivalences is a Morita equivalences. We can therefore endow $\dgcat$ with a second combinatorial model category structure by using the theory of Bausfield localizations. In this case weak-equivalences are Morita equivalences. Similarly to the previous case, the underlying $\oo$-category of this model category coincides with the $\oo$-localization of $\dgcat$ with respect to Morita equivalences. We will label this $\oo$-category by $\dgcatm$. 

Let  $\hat{\Cc}_c$ denote the dg category of perfect $\Cc^{\textup{op}}$-dg modules. Then $\dgcatm$ is equivalent to the full subcategory of $\dgcatdk$ spanned by dg categories $\Cc$ for which the Yoneda embedding $\Cc \hookrightarrow \hat{\Cc}_c$ is a DK-equivalence.

To summerize, we have the following pair of composable $\oo$-localizations
\begin{equation}
\dgcat \rightarrow \dgcatdk \rightarrow \dgcatm
\end{equation}
$\dgcatdk \rightarrow \dgcatm$ is a left adjoint to the inclusion $\dgcatm\hookrightarrow \dgcatdk$ under the identification mentioned above. At the level of objects, it is defined by the assignment $T\mapsto \hat{T}_c$.

It is possible to enhance $\dgcatdk$ and $\dgcatm$ with symmetric monoidal structures. Furthermore, if we restrict to the full subcategory $\dgcatlf \subseteq \dgcat$ of locally-flat (small) dg categories, we get two composable symmetric monoidal $\oo$-functors
\begin{equation}
    \dgcatlfo \rightarrow \dgcatdko \rightarrow \dgcatmo
\end{equation}
More details on the Morita theory of dg categories can be found in \cite{to07}.

One of the most recurrent operations that occur in this work is that of forming quotients of dg categories:
given a dg category $\Cc$ together with a full sub dg category $\Cc'$, both of them in $\dgcatm$, we can consider  the pushout $\Cc \amalg_{\Cc'}0$ in $\dgcatm$. Here $0$ denotes the final object in $\dgcatm$, i.e. the dg category with only one object whose endomorphisms are given by the zero hom-complex. We denote this pushout by $\Cc/\Cc'$ and refer to it as the dg quotient of $\Cc' \hookrightarrow \Cc$. Equivalently, The dg category $\Cc/\Cc'$ can also be obtained as the image in $\dgcatm$ of the pushout $\Cc \amalg_{\Cc'}0$ formed in $\dgcatdk$. Its homotopy category coincides with (the idempotent completion of) the Verdier quotient of $\textup{H}^0(\Cc)$ by the full subcategory $\textup{H}^0(\Cc')$ (see \cite{dri}). 

Compact objects in $\dgcatm$ are dg categories of finite-type over $A$, as defined in \cite{tv07}. In particular,
\begin{equation}
    \textup{Ind}(\dgcatft)\simeq \dgcatm
\end{equation}

Moreover, $\dgcatm$ is equivalent to the $\oo$-category of small, idempotent complete, $A$-linear stable $\oo$-categories (\cite{co13}).

Here there are some dg categories that we will use: let $X$ be a derived scheme (stack).
\begin{itemize}
    \item $\qcoh{X}$ will denote the dg category of quasi-coherent $\Oc_X$-modules. It can be defined as follows. If $X=Spec(B)$, then $\qcoh{X}=\md_B$, the dg category of dg modules over the dg algebra associated to $B$ via the Dold-Kan equivalence. In the general case, $\qcoh{X}=\varprojlim_{Spec(B)\rightarrow X}\md_B$ (the functoriality being that induced by base change)
    \item $\perf{X}$ will denote the full sub dg category of $\qcoh{X}$ spanned by perfect complexes. If $X=Spec(B)$, then $\perf{B}$ is the smallest subcategory of $\md_B$ which contains $B$ and that is stable under the formation of finite colimits and retracts. More generally, an object $\Ec \in \qcoh{X}$ is perfect if, for any $g:Spec(B)\rightarrow X$, the pullback $g^*\Ec \in \perf{B}$. Perfect complexes coincide with dualizable objects in $\qcoh{X}$ and, under some mild additional assumptions (that will always be verified in our examples), with compact objects (see \cite{bzfn})
\end{itemize}
Moreover, if we assume to work in the noetherian setting, we can consider
\begin{itemize}
    \item $\cohb{X}$ will denote the full sub dg category of $\qcoh{X}$ spanned by those cohomologically bounded complexes $\Ec$ (i.e. $\textup{H}^i(\Ec)\neq 0$ only for a finite number of indexes) such that $\textup{H}^*(\Ec)$ is a coherent $\textup{H}^0(\Oc_X)$-module
    \item $\cohm{X}$ will denote the full dg category of $\qcoh{X}$ spanned by those cohomologically bounded above complexes $\Ec$ (i.e. $\textup{H}^i(\Ec)=0$ for $i>>0$) such that $\textup{H}^*(\Ec)$ is a coherent $\textup{H}^0(\Oc_X)$-modules. These are also known as pseudo-coherent complexes
    \item Let $p : X\rightarrow Y$ be a proper morphism locally almost of finite type. By \cite[Cpt.4 Lemma 5.1.4]{gr17}, we have an induced $\oo$-functor $p_*:\cohb{X}\rightarrow \cohb{Y}$. We denote $\cohbp{X}{Y}$ the full subcategory of $\cohb{X}$ spanned by those objects $\Ec$ such that $p_*\Ec\in \perf{Y}$
\end{itemize}
\section{\textsc{Derived algebraic geometry}}
Derived algebraic geometry is a broad generalization of algebraic geometry whose building blocks are simplicial commutative rings, rather then commutative rings. It is usually better behaved in the situations that are typically defined bad in the classical context, e.g. non-transversal intersections. The main idea is to develop algebraic geometry in an homotopical context: instead of saying that two elements are equal we rather say that they are homotopic, the homotopy being part of the data. In this thesis derived schemes appear exclusively as (homotopy) fiber products of ordinary schemes. It is necessary to allow certain schemes to be derived, as if we restrict ourselves to work with discrete (i.e. classical) schemes, certain statements are false as stated here (e.g. Corollary \ref{Orlov equivalence}) and some important characters do not appear (e.g. the algebra which acts on certain dg categories of (relative) singularities).

The ideas and motivations that led to derived algebraic geometry go back to J.P.~Serre (Serre's intersection formula, \cite{se65}), P.~Deligne (algebraic algebraic geometry in a symmetric monoidal category, \cite{de90}), L.~Illusie, A.~Grothendieck, M.~André, D.~Quillen (the cotangent complex, \cite{an74}, \cite{ill71}, \cite{qu70}, \cite{sga6}) $\dots$ however, the theory nowadays relies on solid roots thanks to the work of J.~Lurie (\cite{sag}) and B.~Toën - G.~Vezzosi (\cite{tv05}, \cite{tv08}).
We shall briefly recall the main definition and features of the theory following \cite{to14}.

Let $\sring$ be the $\oo$-category of simplicial rings. It is equivalent to the $\oo$-localization of the category of simplicial rings with respect to the class of morphism that are weak homotopy equivalences on the underlying simplicial sets.

For any topological space $X$ we can consider the $\oo$-category $\sring(X)$ of $\oo$-functors $Op(X)^{\textup{op}}\rightarrow \sring$ satisfying the descent condition\footnote{$Op(X)$ is the category associated to the partially ordered set of open subsets of $X$.}. Similarly to the classical case, for any continuous map $\phi:X\rightarrow Y$ there is an adjunction (of $\oo$-functors)
\begin{equation*}
    \phi^{-1}:\sring(Y)\rightleftarrows \sring(X): \phi_* 
\end{equation*}
One can define an $\oo$-category $\textup{\textbf{dRingSp}}$ of derived ringed spaces, whose objects are pairs $(X,\Oc_X)$ consisting of a topological space $X$ and a sheaf of simplicial rings $\Oc_X$ on it. For two such objects $(X,\Oc_X)$, $(Y,\Oc_Y)$, the mapping space is given by the formula
\begin{equation}
    \Map_{\textup{\textbf{dRingSp}}}\bigl ((X,\Oc_X),(Y,\Oc_Y) \bigr )\simeq \amalg_{\phi:X\rightarrow Y}\Map_{\sring(Y)}(\Oc_Y,\phi_*\Oc_X)
\end{equation}
To any derived ringed space $(X,\Oc_X)$ we can associate a ringed space $(X,\pi_0(\Oc_X))$ by considering the sheaf of connected components of $\Oc_X$. It is called the truncation of $(X,\Oc_X)$. 

This allows us to restrict our attention to derived locally ringed spaces, i.e. to the subcategory $\textup{\textbf{dRingSp}}^{\textup{loc}}$ of $\textup{\textbf{dRingSp}}$: its objects are the objects of $\textup{\textbf{dRingSp}}$ whose truncations are locally ringed spaces and its maps are those inducing local morphisms on the truncations.

\begin{defn}
The $\oo$-category $\dsc$ of derived schemes is the full subcategory of $\textup{\textbf{dRingSp}}^{\textup{loc}}$ spanned by pairs $(X,\Oc_X)$ such that 
\begin{itemize}
    \item the truncation $(X,\pi_0(\Oc_X))$ is a scheme
    \item each $\pi_i(\Oc_X)$ is a quasi-coherent $\pi_0(\Oc_X)$-module
\end{itemize}
\end{defn}
By definition, each derived scheme $X$ has an underlying scheme $t_0(X)$ (its truncation). Indeed, the assignment $X\mapsto t_0(X)$ is part of an adjunction 
\begin{equation}
  \iota: \textup{\textbf{Sch}} \rightleftarrows \dsc :  t_0 
\end{equation}
A derived scheme is affine if its underlying scheme is so. There is an equivalence of $\oo$-categories
\begin{equation}
    \dscaff\simeq \sring^{\textup{op}}
\end{equation}
where $\dscaff$ is the full subcategory of $\dsc$ spanned by affine derived schemes. For any simplicial commutative ring $A$, we denote $Spec(A)$ the associated derived scheme.

The $\oo$-category $\dsc$ has all finite limits. In particular, it has fiber products. For example, if we consider a diagram
\begin{equation*}
    Spec(B)\rightarrow Spec(A) \leftarrow Spec(C)
\end{equation*}
of affine derived schemes, the fiber product is equivalent to $Spec(B\otimes^{\mathbb{L}}_AC)$, the spectrum of the derived tensor product. 

Consider two (underived) $S$-schemes $X,Y$ ($S$ being underived itself). Then the fiber product computed in $\dsc$ (denoted $X\times^h_SY$) might differ from the one computed in $\textup{\textbf{Sch}}$ (denoted $X\times_S Y$). However, they are related by the formula
\begin{equation}
    t_0(X\times_S^hY)\simeq X\times_SY
\end{equation}

\section{\textsc{\texorpdfstring{$\ell$}{l}-adic sheaves}}\label{l-adic sheaves}
We shall briefly introduce the $\oo$-category of $\ell$-adic sheaves, following \cite{gl18}.

Fix a prime number $\ell$, which is invertible in each residue field of our base scheme $S$.

For $X$ a scheme of finite type over $S$, let $\schet{X}$ be the category of $X$-schemes whose structure morphism is etale. Notice that the morphisms in this category are all etale.
\begin{defn}
An etale covering of $V\in \schet{X}$ is a family of morphisms $\{U_{\alpha}\rightarrow V\}$ in $\schet{X}$ such that every point in $V$ has a preimage in some $U_{\alpha}$, i.e. it is a family of étale morphism which is jointly surjective.
\end{defn}
\begin{rmk}
Also notice that, since all the schemes involved are quasi-compact, for every etale covering of $V$ there exists a finite number of morphisms $U_{\alpha_i}\rightarrow V$ such that $\amalg U_{\alpha_i}\rightarrow V$ is surjective.
\end{rmk}

Let $A$ be a commutative ring and let $\md_A$ be the $\oo$-category of $A$-modules. Recall that this is the simplicial set defined as follows
\begin{itemize}
\item $\md_{A,n}$ is the set of pairs $\bigl ( \{M(i)\}_{i=0,\dots,n}, \{f_{I}\}_{I\subseteq \{0,\dots,n\}} \bigr )$, where $M(0),\dots,M(n)$ are $A$-projective cochain complexes (i.e. the differential increases the degree) and the $f_I$ are cochain maps, compatible under composition up to coherent homotopy: for any subset $I=\{i_{-}<i_m< \dots <i_1<i_+ \}\subseteq \{0,\dots,n\}$, we are given a collection of maps $f_I: M(i_-)_k\rightarrow M(i_+)_{k+m}$ such that
\begin{equation}
d(f_I(x))= (-1)^mf_I(dx)+ \sum_{j=1}^m (-1)^j \bigl (f_{I-\{i_j\}}(x)-f_{\{i_j,\dots, i_+\}}\circ f_{\{i_-,\dots,i_j\}}(x) \bigr )
\end{equation}
\item For $\alpha : [m] \rightarrow [n]$ a map of ordered sets, the map $\alpha^* :\md_{A,n}\rightarrow \md_{A,m}$ is defined by
\begin{equation}
\alpha^*(\{M(i)\}_{i=0,\dots,n}, \{f_{I}\}_{I\subseteq \{0,\dots,n\}})=(\{M(\alpha(i))\}_{i=1,\dots,n},\{g_J\}_{J\subseteq\{0,\dots,m\}}) 
\end{equation}
where
\[
g_J(x)=
\begin{cases}
f_{\alpha(J)} \hspace{0.5cm} \text{ if $\alpha_{|J}$ is injective} \\
x \hspace{0.5cm} \hspace{0.5cm} \text{ if } J=\{j,j'\} \text{ and } \alpha(j)=\alpha(j') \\
0 \hspace{0.5cm} \text{ otherwise}
\end{cases}
\]
\end{itemize} 
\begin{rmk}
The simplicial set $\md_A$ is the dg nerve of the dg category of $A$-projective cochain complexes (see \cite[Construction 1.3.1.6]{ha}), and therefore is an $\oo$-category (\cite[Proposition 1.3.1.10]{ha}).
\end{rmk}
\begin{defn}
A $\md_A$-valued presheaf on $X$ is an $\oo$-functor
\[
\mathcal{F} : (\schet{X})^{op}\rightarrow \md_A
\]
\end{defn}

Therefore, if $\mathcal{F}$ is a $\md_A$-valued presheaf on $X$, we are given an $A$-projective cochain complex $\mathcal{F}(U)$ for each étale $X$-scheme $U$. In particular, for each $n \in \Z$, we can consider the $\md_A^{\heart}$-valued presheaf $U\mapsto \textup{H}^{-n}(\mathcal{F}(U))$. Denote by $\pi_{n}(\mathcal{F})$ the associated sheaf.

\begin{defn}
We say that a $\md_A$-valued presheaf on $X$ is locally acyclic if, for any $n \in \Z$, the sheaf $\pi_n(\mathcal{F})$ vanishes.
\end{defn}
\begin{defn}
A $\md_A$-valued (hypercomplete) sheaf $\mathcal{F}$ on $X$ is a presheaf such that, if $\mathcal{G}$ is locally acyclic, any morphism $\mathcal{G}\rightarrow \mathcal{F}$ is nullhomotopic. We will label $\shv(X,A)$ the full subcategory of $\textup{Fun}((\schet{X})^{op}, \md_A)$ spanned by such objects.
\end{defn}
\begin{rmk}
The inclusion functor $\shv(X,A)\subseteq \textup{Fun}\bigl ((\schet{X})^{op},\md_A \bigr )$ admits a left adjoint functor, called the sheafification functor.
\end{rmk}

\begin{ex}
\begin{itemize}
\item If $X=Spec(\bar{k})$ is the spectrum of a separably closed field, then $\shv(X,A)$ coincides with $\md_A$
\item If $X=Spec(k)$, then $\shv(X,A)$ is, roughly, the $\oo$-category of $A$-projective cochain complexes endowed with a continuous action of the absolute Galois group of $k$. We shall investigate this more accurately later on.
\end{itemize}
\end{ex}

The $\oo$-category $\shv(X,A)$ is stable and admits a t-structure. For an integer $n$, let $\shv(X,A)^{\leq n}$ be the full subcategory of $\shv(X,A)$ spanned by those objects $\mathcal{F}$ such that, for each $m>n$, $\pi_{-m}(\mathcal{F})\simeq 0$. Similarly, define $\shv(X,A)^{\geq n}$ to be the full subcategory of $\shv(X,A)$ spanned by those objects $\mathcal{G}$ such that, for each $m< n$, $\pi_{-m}(\mathcal{G})\simeq 0$. Then $(\shv(X,A)^{\geq 0},\shv(X,A)^{\leq 0})$ determines a t-structure on $\shv(X,A)$. Moreover, the heart $\shv(X,A)^{\heart}$ with respect to this t-structures is equivalent to the (more classical) abelian category of étale sheaves on $X$ (see \cite[Remark 2.1.2.1]{sag}). Moreover, the abelian category $\shv(X,A)^{\heart}$ is Grothendieck. Therefore, we can consider its unbounded derived category $\mathcal{D}\bigl (\shv(X,A)^{\heart} \bigr )$ (see \cite[\S 1.3.5]{ha}). Then, \cite[Theorem 2.1.2.2]{sag} tells us that we have a t-exact equivalence
\begin{equation}
\mathcal{D}\bigl ( \shv(X,A)^{\heart} \bigr )\simeq \shv(X,A)
\end{equation}
For any morphism of $S$-schemes $f: X\rightarrow Y$ it is possible to define a functor
\begin{equation}
\schet{Y} \rightarrow \schet{X}
\end{equation}
\[
(U\rightarrow Y)\mapsto (U\times_Y X\rightarrow X)
\]
Therefore, given a $\md_A$-valued presheaf $\mathcal{F}$ on $Y$, we can define a $\md_A$-valued presheaf $f_*\mathcal{F}$ by the composition of $\mathcal{F}$ with the functor above: informally, $f_*\mathcal{F}$ is defined by the assgnment
\begin{equation}
U\mapsto f_*\mathcal{F}(U)= \mathcal{F}(U\times_Y X)
\end{equation}
This defines an $\oo$-functor
\begin{equation}
f_*:\textup{Fun}\bigl ( (\schet{X})^{op},\md_A \bigr)\rightarrow \textup{Fun}\bigl ( (\schet{Y})^{op},\md_A \bigr)
\end{equation}
which induces a functor on the $\oo$-categories of sheaves
\begin{equation}
f_* : \shv(X,A)\rightarrow \shv(Y,A)
\end{equation}
called the pushforward $\oo$-functor, which admits a left adjoint, the pullback $\oo$-functor 
\begin{equation}\label{pullback oo}
f^* : \shv(Y,A)\rightarrow \shv(X,A)
\end{equation}
Assume that $Spec(k)$ is the spectrum of a field. Let $X$ be a $k$-scheme (of finite type). Recall from \cite[Exposé XIII]{sga7ii} that if $\mathcal{F}$ is a $\md_A^{\heart}$-valued étale sheaf (in the classical sense), an action of $Gal(\bar{k}/k)$ on $\mathcal{F}$ is compatible with the natural action of $Gal(\bar{k}/k)$ on $\bar{X}=X\times_{Spec(k)}Spec(\bar{k})$ if for any $g \in Gal(\bar{k}/k)$, we are given an isomorphism 
\begin{equation}
\sigma(g): g_*\mathcal{F}\rightarrow \mathcal{F}
\end{equation}
and these isomorphisms are compatible with the group multiplication: 
\begin{equation}
\sigma(gh)=\sigma(g)\sigma(h) \hspace{0.5cm} \text{ for all } g,h \in Gal(\bar{k}/k)
\end{equation}
Consider the ordinary pullback functor
\begin{equation}
\shv(X,A)^{\heart} \rightarrow \shv(\bar{X},A)^{\heart}
\end{equation}
\[
\mathcal{F}\mapsto \bar{\mathcal{F}}
\]
\begin{rmk}
The $\oo$-functor (\ref{pullback oo}) is an $\oo$-enhancement of the derived functor of the classical pullback considered above.
\end{rmk}
It is a well known fact that the classical pullback functor induces an equivalence 
\begin{equation}
\shv(X,A)^{\heart}\simeq \shv(\bar{X},A)^{Gal(\bar{k}/k),\heart}
\end{equation}
where the category on the right is the category of $\md_A^{\heart}$-valued sheaves on $\bar{X}$, endowed with a continuous action of $Gal(\bar{k}/k)$, compatible with the action of $Gal(\bar{k}/k)$ on $\bar{X}$ (see \cite[Rappel 1.1.3]{sga7ii}). In particular, this category is Grothendieck abelian, as $\shv(X,A)^{\heart}$ is so. Therefore, we can consider its derived category (in its $\oo$-categorical incarnation)  $\mathcal{D}\bigl ( \shv(\bar{X},A)^{Gal(\bar{k}/k),\heart} \bigr )$, which is equivalent to $\shv(X,A)$.
\begin{notation}
We will write $\shv(\bar{X},A)^{Gal(\bar{k}/k)}$ instead of $\mathcal{D}\bigl ( \shv(\bar{X},A)^{Gal(\bar{k}/k),\heart} \bigr )$.
\end{notation} 
We will recall, following \cite[\S 1.3.5]{ha} how this is defined.

Consider the category $Ch(\shv(\bar{X},A)^{Gal(\bar{k}/k),\heart})$ of cochain complexes constructed out of $\shv(\bar{X},A)^{Gal(\bar{k}/k),\heart}$. This category admits a left proper combinatorial model structure (\cite[Proposition 1.3.5.3]{ha}), where
\begin{itemize}
\item (W) weak equivalences are quasi-isomorphisms
\item (C) a morphism is a cofibration if it is a levelwise monomorphism
\item (F) a morphism is a fibration if it has the right lifting property with respect to morphisms which are both weak equivalences and cofibrations
\end{itemize}
Notice that, as the zero object in this category is the sheaf which associates to each étale $\bar{X}$-scheme the zero complex, every object is cofibrant.
Moreover, $Ch(\shv(\bar{X},A)^{Gal(\bar{k}/k),\heart})$ is naturally a dg category. Let $Ch(\shv(\bar{X},A)^{Gal(\bar{k}/k),\heart})^{\circ}$ be the full subcategory of $Ch(\shv(\bar{X},A)^{Gal(\bar{k}/k),\heart})$ spanned by fibrant-cofibrant objects. then $\shv(\bar{X},A)^{Gal(\bar{k}/k)}$ is defined as the dg nerve (see \cite[Construction 1.3.1.6]{ha}) of $Ch(\shv(\bar{X},A)^{Gal(\bar{k}/k),\heart})^{\circ}$. That is:
\begin{itemize}
\item For $n\geq 0$, $\shv(\bar{X},A)^{Gal(\bar{k}/k)}_n$ is the set of pairs $\bigl ( \{X(i)\}_{i=0,\dots,n},\{f_I\}_{I\subseteq [n]}\bigr )$, where each $X(i)$ is an object in $Ch(\shv(\bar{X},A)^{Gal(\bar{k}/k),\heart})^{\circ}$ and, for each $I=\{i_-<i_m<\dots<i_1<i_+\}\subseteq [n]$, $f_I\in \Map(X(i_-,X_+))_m$ and the following equalities hold
\[
df_I= d\circ f_I - (-1)^m f_I\circ d =\sum_{j=1}^m (-1)^j \bigl (f_{I-\{i_j\}}(x)-f_{\{i_j,\dots, i_+\}}\circ f_{\{i_-,\dots,i_j\}}(x) \bigr )
\]
\item For $\alpha : [m] \rightarrow [n]$ a map of ordered sets, the map 
\begin{equation*}
    \alpha^* \shv(\bar{X},A)^{Gal(\bar{k}/k)}_n\rightarrow \shv(\bar{X},A)^{Gal(\bar{k}/k)}_m
\end{equation*} is defined by
\begin{equation}
\alpha^*(\{X(i)\}_{i=0,\dots,n}, \{f_{I}\}_{I\subseteq \{0,\dots,n\}})=(\{X(\alpha(i))\}_{i=1,\dots,n},\{g_J\}_{J\subseteq\{0,\dots,m\}}) 
\end{equation}
where
\[
g_J(x)=
\begin{cases}
f_{\alpha(J)} \hspace{0.5cm} \text{ if $\alpha_{|J}$ is injective} \\
x \hspace{0.5cm} \hspace{0.5cm} \text{ if } J=\{j,j'\} \text{ and } \alpha(j)=\alpha(j') \\
0 \hspace{0.5cm} \text{ otherwise}
\end{cases}
\]
\end{itemize}
Therefore, given a $k$-scheme $X$, we can regard each $\md_A$-valued sheaf $\mathcal{F}$ on $X$ as a $\md_A$-valued sheaf on $\bar{X}$ with an action of $Gal(\bar{k}/k)$.  

Let us stick back to the case of $S$-schemes, where $S$ is an affine noetherian regular local scheme.
The $\oo$-category $\shv(X,A)$ is compactly-generated, its compact objects being constructible sheaves (see \cite[Proposition 3.38]{brtv}). In what follows, we will denote $\shv^c(X,A)$ the full subcategory of $\shv(X,A)$ spanned by compact objects.
 
For our purposes, it will be also important the functorial behaviour of $\shv(X,A)$ with respect to the ring $A$. For a ring homomorphism $A\rightarrow B$ we get a forgetful $\oo$functor
\begin{equation}
\shv(X,B)\rightarrow \shv(X,A)
\end{equation}
which admits a left adjoint, the extension of scalars $\oo$-functor
\begin{equation}
\shv(X,A)\rightarrow \shv(X,B)
\end{equation}
roughly given by the assignment 
\[
\mathcal{F}\mapsto \mathcal{F}\otimes_A^{\mathbb{L}}B
\]
We will now focus on the case $A=\Z/\ell^d\Z$ and on the ring homomorphisms $\Z/\ell^d\Z \rightarrow \Z/\ell^{d-1}\Z $. We get a sequence of $\oo$-functors
\begin{equation}
\shv(X,\Z/\ell^{d}\Z)\rightarrow \shv(X,\Z/\ell^{d-1}\Z)
\end{equation}
and it follows from \cite[Proposition 2.2.8.4]{gl18} that the image of a constructible sheaf is again constructible, yielding 
\begin{equation}
\shv^c(X,\Z/\ell^{d}\Z)\rightarrow \shv^c(X,\Z/\ell^{d-1}\Z)
\end{equation}
We can then consider the limit of the diagram of $\oo$-categories
\begin{equation}
\shv^c_{\ell}(X):=\varprojlim \shv^c(X,\Z/\ell^{d}\Z)
\end{equation}
This $\oo$-category can be identified with the full subcategory of $\shv(X,\Z)$ spanned by $\ell$-complete constructible sheaves, namely by those objects $\mathcal{F} \in \shv(X,\Z)$ such that
\begin{enumerate}
\item $\mathcal{F}\simeq \varprojlim \mathcal{F}/\ell^d\mathcal{F}$;
\item for any $d\geq 1$, $\mathcal{F}/\ell^d\mathcal{F}$ is constructible. 
\end{enumerate}
We will refer to this $\oo$-category with the $\oo$-category of constructible $\ell$-adic sheaves.

The pushforward for a morphism $f : X\rightarrow Y$ of $k$-schemes of finite type induces an $\oo$-functors at the level of constructible $\ell$-adic
sheaves
\begin{equation}
f_* : \shv^c_{\ell}(X)\rightarrow \shv^c_{\ell}(Y)
\end{equation}
which admits a left adjoint
\begin{equation}
f^*:\shv^c_{\ell}(Y)\rightarrow \shv^c_{\ell}(X)
\end{equation}
which, at the level of objects, takes a constructible $\ell$-adic sheaf to the $\ell$-completion of its pullback.

We next consider the ind-completion of such categories:
\begin{equation}
\shv_{\ell}(X):= \textup{Ind}\bigl ( \shv^c_{\ell}(X) \bigr )
\end{equation}
the $\oo$-category of $\ell$-adic sheaves. It is then a formal fact that we have a couple of adjoint functors, also called the pushforward and the pullback, defined at the level of $\ell$-adic sheaves.

Finally, we consider the localization of $\shv_{\ell}(X)$ with respect to the class of morphisms $\{ \mathcal{F}\rightarrow \mathcal{F}[\ell^{-1}] \}$, obtaining the $\oo$-category of $\mathbb{Q}_{\ell}$-adic sheaves $\shvl(X)$.

Notice that, in the case we consider schemes of finite type over a field, by applying the same construction to the $\oo$-categories $\shv^c(\bar{X},\Z/\ell^{d}\Z)^{Gal(\bar{k}/k)}$, we obtain an equivalent $\oo$-category $\shvl(\bar{X})^{Gal(\bar{k}/k)}$.
\section{\textsc{Stable homotopy categories}}
In this section we will briefly recall the constructions and main properties of $\sh$ and of $\shnc$, the stable $\oo$-category of schemes and the stable $\oo$-category of non-commutative spaces (a.k.a. dg categories).
\subsection{\textsc{The stable homotopy category of schemes}}
The stable homotopy category of schemes was first introduced by F.~Morel and V.~Voevodsky in their celebrated paper \cite{mv99}. The main idea is to develop an homotopy theory for schemes, where the role of the unit interval - which is not available in the world of schemes - is played by the affine line. It was first developed using the language of model categories. We will rather use that of $\oo$-categories, following \cite{ro14} and \cite{ro15}. The two procedures are compatible, as shown in (\cite{ro14} and \cite{ro15}).

Let $S=Spec(A)$ denote an affine scheme. Recall that the Nisnevich topology on $\sm$ is the Grothendieck topology associated to the pretopology whose covering families are collections of étale morphisms $\{U_i\rightarrow X\}$ such that, for every point $x\in X$, there exist an index $i$ and a point $u_i\in U_i$ verifying $k(u_i)\simeq k(x)$. Roughly, we are excluding Galois coverings: if $X=Spec(k)$ is the spectrum of a field, a (non trivial) finite Galois extension $Spec(l)\rightarrow Spec(k)$ is an étale covering, but it is not a Nisnevich one. Recall that an elementary Nisnevich square is a square
\begin{equation}
    \begindc{\commdiag}[18]
    \obj(-15,10)[1]{$U\times_X V$}
    \obj(15,10)[2]{$V$}
    \obj(-15,-10)[3]{$U$}
    \obj(15,-10)[4]{$X$}
    \mor{1}{2}{$$}
    \mor{1}{3}{$$}
    \mor{2}{4}{$p$}
    \mor{3}{4}{$j$}
    \enddc
\end{equation}
such that 
\begin{itemize}
    \item it is a pullback square
    \item $p$ is an étale morphism and $j$ is an open embedding
    \item $p$ induces an isomorphism between $(X-U)\times_XV$ and $X-U$ (both endowed with the reduced scheme structure)
\end{itemize}
Notice that every elementary Nisnevich square determines a Nisnevich covering family $\{U\rightarrow X, V\rightarrow X\}$. Indeed, this is all we need: a presheaf on $\sm$ is a Nisnevich sheaf if and only if it sends Nisnevich elementary squares to pushout squares. 
Then one can produce the unstable homotopy category of schemes as follows: one considers the $\oo$-category $\textup{Fun}(\sm^{\textup{op}},\mathcal{S})$ of presheaves (of spaces) on $\sm$. 
Its full subcategory spanned by Nisnevich sheaves, $\textup{Sh}_{Nis}(\sm)$ is an example of an $\oo$-topos (in the sense of \cite{htt}). 
Then one has to consider its hypercompletion $\textup{Sh}_{Nis}(\sm)^{hyp}$, which coincides with the localization of $\textup{Fun}(\sm^{\textup{op}},\mathcal{S})$ spanned by objects that are local with respect to Nisnevich hypercovers. 
If we further localize with respect to the projections $\{\aff{1}{X}\rightarrow X\}$, we obtain the unstable homotopy $\oo$-category of schemes, which we will denote by $\uH{S}$. This presentable $\oo$-category is characterized by the following universal property (see \cite[Theorem 2.30]{ro14}): for every $\oo$-category with small colimits, the $\oo$-functor
\begin{equation}
    \textup{Fun}^{\textup{L}}(\uH{S},\mathcal{D})\rightarrow \textup{Fun}(\sm^{\textup{op}},\mathcal{D})
\end{equation}
induced by the canonical map $\sm \rightarrow \uH{S}$ is fully faithful and its essential image coincides with those $\oo$-functors $\sm^{\textup{op}}\rightarrow \mathcal{D}$ that satisfy Nisnevich descent and $\aff{1}{}$-invariance. Here $\textup{Fun}^{\textup{L}}(\uH{S},\mathcal{D})$ denotes the $\oo$-category of colimit-preserving $\oo$-functors.
It is also important that the canonical $\oo$-functor $\sm \rightarrow \uH{S}$ can be promoted to a symmetric monoidal $\oo$-functor with respect to the Cartesian structures.
One then considers the pointed version of $\uH{S}$, $\uH{S*}$: it comes equipped with a canonical symmetric monoidal structure $\uH{S*}^{\wedge}$ and there is a symmetric monoidal $\oo$-functor $\uH{S}^{\times}\rightarrow \uH{S*}^{\wedge}$.
The final step consists in stabilization: in classical stable homotopy theory one forces stabilization by inverting $S^1$. However, in this context, there exist two circles, the topological circle $S^1:=\Delta^1/\partial \Delta^1$ and the algebraic circle $\gm{S}$. One the stabilizes $\uH{S*}$ by inverting $S^1\wedge \gm{S}\simeq(\proj{1}{S},\oo):=cofib(S\xrightarrow{\oo}\proj{1}{S})\footnote{clearly, this cofiber is taken in $\uH{S*}$}$ - this can be done using the machinery developed in \cite[\S2.1]{ro15}. As a result, we obtain the presentable, symmetric monoidal, stable $\oo$-category $\sh^{\otimes}:=\uH{S*}^{\wedge}[(\proj{1}{S},\oo)^{-1}]$, called the stable homotopy $\oo$-category of schemes. It is moreover characterized by the following universal property (see \cite[Corollary 2.39]{ro15}): there is a symmetric monoidal $\oo$-functor $\Sigma_+^{\oo}:\sm^{\times}\rightarrow \sh^{\otimes}$ and for any presentable, symmetric monoidal pointed $\oo$-category $\mathcal{D}^{\otimes}$, the map
\begin{equation}
    \textup{Fun}^{\otimes,\textup{L}}(\sh^{\otimes},\mathcal{D}^{\otimes})\rightarrow \textup{Fun}^{\otimes}(\sm^{\times},\mathcal{D}^{\otimes})
\end{equation}
induced by $\Sigma_+^{\oo}$ is fully faithful and its image coincides with those symmetric monoidal $\oo$-functors $F:\sm^{\times}\rightarrow \mathcal{D}^{\otimes}$ that satisfy
\begin{itemize}
    \item Nisnevich descent
    \item $\aff{1}{}$-invariance
    \item $cofib(F(S\xrightarrow{\oo}\proj{1}{S}))$ is an invertible object in $\mathcal{D}$
\end{itemize}
\begin{rmk}
It can be shown, using results of Ayoub (\cite{ay08i}, \cite{ay08ii}), Cisinski-Déglise (\cite{cd19}) and the machinery developed by Gaitsgory-Rozenblyum (\cite{gr17}) and Liu-Zheng (\cite{lz15}, \cite{lz17i}, \cite{lz17ii}), that the assignment $S\mapsto \sh$ defines a sheaf of $\oo$-categories enhanced with a Grothendieck $6$-functors formalism. See \cite{ro14} and \cite{ro15} and the appendix in \cite{brtv}.
\end{rmk}
Among the objects of $\sh$, the spectrum of homotopy invariant non-connective K-theory $\bu$\footnote{more commonly, this spectrum is denoted $KGL_S$, see \cite{ci13}. We use the notation $\bu_S$, which comes from topology, following the lead of \cite{brtv}, which is our main source of inspiration.} will play a crucial role in what follows. Recall that it is an object in $\sh$ such that, for every object $X\in \sm$,
\begin{equation}
    \Map_{\sh}(\Sigma^{\oo}_+X,\bu_S)\simeq \hk(X)
\end{equation}
Moreover, $\bu_S$ satisfies the algebraic Bott periodicity:
\begin{equation}
    \bu_S\simeq \bu_S(1)[2]
\end{equation}
\subsection{\textsc{The stable homotopy category of non-commutative spaces}}
\begin{rmk}
This section is an exposition of the ideas and results of \cite{ro14} and \cite{ro15}.
\end{rmk}
It is possible to mimic the procedure described in the previous paragraph to construct a presentable, symmetric monoidal, stable $\oo$-category $\shnc$ defined starting with non-commutative spaces rather than schemes. 
The first step is to understand which category should play the role of $\sm$.
It comes up that the right notion of smooth non-commutative spaces is that of dg category of finite type. 
\begin{defn}
The $\oo$-category of smooth non-commutative spaces is $\nc{S}:=\dgcat^{\textup{ft,op}}$
\end{defn}
Next, one needs to define an appropriate substitute for the Nisnevich topology in this $\oo$-category.
\begin{notation}
If $X\in \nc{S}$ we will write $T_X$ for the corresponding object in $\dgcat^{\textup{ft}}$.
\end{notation}
One starts with the definition of open embedding:
\begin{defn}
Let $U\rightarrow X$ be a morphism of smooth non-commutative spaces. We say that it is an open embedding if there exists a smooth dg category with a compact generator $K_{X-U}$ and a fully faithful embedding $K\rightarrow T_X$ that fits in an exact triangle of dg categories
\begin{equation}
    K\hookrightarrow T_X \rightarrow T_U
\end{equation}
\end{defn}
\begin{ex}
The previous definition is inspired by the following: let $U\hookrightarrow X$ be an open embedding of smooth $S$-schemes. Then $\perf{U}\rightarrow \perf{X}$ is an open embedding of smooth non-commutative spaces.
\end{ex}
\begin{defn}
A commutative square in $\nc{S}$
\begin{equation}\label{non commutative Nisnevich square}
    \begindc{\commdiag}[18]
    \obj(-15,10)[1]{$W$}
    \obj(15,10)[2]{$V$}
    \obj(-15,-10)[3]{$U$}
    \obj(15,-10)[4]{$X$}
    \mor{1}{2}{$$}
    \mor{1}{3}{$$}
    \mor{2}{4}{$$}
    \mor{3}{4}{$$}
    \enddc
\end{equation}
is Nisnevich if
\begin{itemize}
    \item $U\rightarrow X$ and $W\rightarrow V$ are open embeddings
    \item $T_X\rightarrow T_V$ sends the compact generator of $K_{X-U}$ to the compact generator of $K_{V-W}$ and induces an equivalence $K_{X-U}\simeq K_{V-W}$
    \item the square is a pushout in $\dgcatm$
\end{itemize}
\end{defn}
\begin{rmk}\label{convention emptyset}
By convention, $\emptyset$ is a Nisnevich covering of the zero object in $\dgcat^{\textup{ft}}$.
\end{rmk}
\begin{rmk}
The notion on Nisnevich squares for smooth non-commutative spaces is compatible with the classical one: in
\begin{equation}
    \begindc{\commdiag}[18]
    \obj(-15,10)[1]{$V\times_X U$}
    \obj(15,10)[2]{$V$}
    \obj(-15,-10)[3]{$U$}
    \obj(15,-10)[4]{$X$}
    \mor{1}{2}{$$}
    \mor{1}{3}{$$}
    \mor{2}{4}{$$}
    \mor{3}{4}{$$}
    \enddc
\end{equation}
is an elementary Nisnevich square, then its image via $\perf{\bullet}$ is a Nisnevich square of non-commutative spaces (\cite[Proposition 3.21]{ro15}).
\end{rmk}
\begin{rmk}
It is worth mentioning that Nisnevich squares do not define a Grothendieck topology on $\dgcat^{\textup{ft}}$. Indeed, they are not stable under pullbacks, see \cite[Remark 3.23]{ro15}.
\end{rmk}
One then considers the category of presheaves $\textup{Fun}(\nc{S}^{\textup{op}},\mathcal{S})$ and its localization with respect to the class of morphisms $\{j(U)\amalg_{j(X)}j(V)\rightarrow j(W)\}$ determined by Nisnevich squares of non-commutative spaces as in $(\ref{non commutative Nisnevich square})$, where $j:\nc{S} \rightarrow \textup{Fun}(\nc{S}^{\textup{op}},\mathcal{S})$ is the Yoneda embedding. We further localize with respect to the morphisms $\{X\otimes \perf{S}\rightarrow X\otimes \perf{\aff{1}{S}}\}$\footnote{Since the Yoneda embedding is symmetric monoidal (when we consider the Day convolution product on the $\oo$-category of presheaves), then it suffices to localize with respect to the morphism $j(\perf{S}\rightarrow j(\perf{\aff{1}{S}}))$.} and obtain a presentable symmetric monoidal $\oo$-category $\uH{S}^{\textup{nc},\otimes}$. Pursuing the analogy with the commutative case, we should now force the existence of a zero object in $\uH{S}^{\textup{nc}}$ and then stabilize with respect to the topological and algebraic circles. It comes out that the situation is simpler in the non-commutative case: let 
\begin{equation*}
    \psi^{\otimes}:\nc{S}^{\otimes}\rightarrow \uH{S}^{\textup{nc},\otimes}\rightarrow \shncm:= \uH{S}^{\textup{nc},\otimes}[(S^1)^{-1}]
\end{equation*} 
denote the $\oo$-functor that we obtain if we force the topological circle to be invertible. Then 
\begin{itemize}
    \item $\shnc$ is pointed due to the convention of Remark \ref{convention emptyset}
    \item The non-commutative motive $\psi(\proj{1}{S},\oo)=cofib(\psi(\perf{S} \xrightarrow{\oo^*} \perf{\proj{1}{S}}))$ is invertible in $\shncm$ (\cite[Proposition 3.24]{ro15}): indeed, it is equivalent to $\perf{S}\simeq 1^{\textup{nc}}$ and therefore it is not necessary to force $\gm{S}$ to be invertible
\end{itemize}
\subsection{\textsc{The bridge Between motives and non-commutative motives}}\label{the bridge between motives and non commutative motives}
We have now at our disposal the following picture:
\begin{equation}
    \begindc{\commdiag}[18]
    \obj(-20,20)[1]{$\sm^{\times}$}
    \obj(20,20)[2]{$\nc{S}^{\otimes}$}
    \obj(-20,-20)[3]{$\sh^{\otimes}$}
    \obj(20,-20)[4]{$\shncm$}
    \mor{1}{2}{$\perf{\bullet}$}
    \mor{1}{3}{$\Sigma^{\oo}_+$}[\atright,\solidarrow]
    \mor{2}{4}{$\psi^{\otimes}$}
    \mor{3}{4}{$\rperf$}[\atright,\dashArrow]
    \enddc
\end{equation}
The existence of the dotted map, which we name perfect realization, is granted by the universal property of $\sh^{\otimes}$. Indeed
\begin{itemize}
    \item $\psi \circ \perf{\bullet}$ sends (ordinary) Nisnevich squares to pushout squares in $\shnc$: this is a consequence of the compatibility of non-commutative Nisnevich squares with the classical ones and of the definition of $\shnc$
    \item $\aff{1}{}$-invariance is forced by construction
    \item $\psi(\proj{1}{S},\oo)$ is an invertible object in $\shncm$
\end{itemize}
The fact that $\rperf$ commutes with colimits is also guaranteed by the universal property of $\sh^{\otimes}$. As both $\sh$ and $\shnc$ are presentable $\oo$-categories, the adjoint functor theorem \cite[Corollary 5.5.2.9]{htt} implies the existence of a lax-monoidal right adjoint
\begin{equation}
    \mathcal{M}^{\otimes}:\shncm \rightarrow \sh^{\otimes}
\end{equation}
For our purposes, it will be very important the following result
\begin{thm}{\textup{\cite[Corollary 4.10]{ro15}}}
The image of the monoidal unit $1^{\textup{nc}}_S$ via $\mathcal{M}$ is equivalent to $\bu_S$.
\end{thm}
In particular, the $\oo$-functor $\mathcal{M}$ factors trough the full subcategory of $\sh$ spanned by modules over $\bu_S$:
\begin{equation}
    \mathcal{M}^{\otimes}: \shncm\rightarrow \md_{\bu_S}(\sh)^{\otimes}
\end{equation}
We will now introduce a dual version of this $\oo$-functor.
Consider the endomorphism of $\shnc$ induced by the internal hom:
\begin{equation}
    \rhom_{\shnc}(-,1^{\textup{nc}}_S):\shncmop\rightarrow \shncm
\end{equation}
Then consider
\begin{equation}
    \begindc{\commdiag}[18]
    \obj(-50,10)[1]{$\dgcat^{\textup{ft},\otimes}$}
    \obj(0,10)[2]{$\shncmop$}
    \obj(90,10)[3]{$\shncm$}
    \obj(140,10)[4]{$\md_{\bu_S}(\sh^{\otimes})$}
    \obj(-50,-10)[5]{$\dgcatmo$}
    \mor{1}{2}{$\psi^{\otimes}$}
    \mor{2}{3}{$\rhom_{\shnc}(-,1^{\textup{nc}_S})$}
    \mor{3}{4}{$\mathcal{M}^{\otimes}$}
    \mor{1}{5}{$$}[\atright,\injectionarrow]
    \mor{5}{2}{$\psi^{\otimes}$}[\atright,\solidarrow]
    \cmor((-30,-10)(80,-10)(100,-10)(130,-8)(140,0)) \pup(50,-14){$\mathcal{M}^{\vee,\otimes}$}
    \enddc
\end{equation}

where the vertical map on the left is given by the inclusion of $\dgcat^{\textup{ft},\otimes}$ in its Ind-completion $\dgcatmo$ and the oblique $\psi^{\otimes}$ is induced by the universal property of the Ind-completion.
Let $\Cc \in \dgcatm$. Then $\Mv(T)$ is the sheaf of spectra $X \in \sm \mapsto \hk(\perf{X}\otimes_{S} T)$:
\begin{equation}
    \Mv(T)(X)=\Map_{\sh}(\Sigma^{\oo}_+X,\Mv(T))\simeq \Map_{\shnc}(\rperf(\Sigma^{\oo}_+X),\rhom_{\shnc}(T,1^{\textup{nc}}_S))
\end{equation}
\begin{equation*}
    \simeq \Map_{\shnc}(\rperf(\Sigma^{\oo}_+X)\otimes_{S} T,1^{\textup{nc}}_S)\simeq \hk(\perf{X}\otimes_{S} T)
\end{equation*}
Moreover, $\Mv$ has the following nice properties
\begin{itemize}
    \item it is lax monoidal (it is a composition of lax monoidal $\oo$-functors)
    \item it commutes with filtered colimits (see \cite[Remark 3.4]{brtv})
    \item it sends exact sequences of dg categories to fiber-cofiber sequences in $\md_{\bu_S}(\sh)^{\otimes}$
\end{itemize}
\subsection{\textsc{The \texorpdfstring{$\ell$}{l}-adic realization of dg categories}}\label{the l-adic realization of dg categories}
We will need a way to associate an $\ell$-adic sheaf to a dg category. We follow the construction given in \cite[\S3.6, \S3.7]{brtv}, which relies on results of J.~Ayoub and D.C.~Cisinski-F.~Déglise.

Recall that for any commutative algebra $A$ in the $\oo$-category of spectra $\Sp$ there is a universal $\oo$-functor
\begin{equation}
  -\otimes A :  \sh^{\otimes}\rightarrow \md_{A}(\sh)^{\otimes}
\end{equation}
In particular, if we set $A=\HQ$, the Eilenberg-MacLane spectrum of rational homotopy theory, then we get
\begin{equation}
    -\otimes \HQ:\sh^{\otimes}\rightarrow \md_{\HQ}(\sh)^{\otimes}
\end{equation}
which identifies the $\oo$-category on the right hand side with non-torsion objects of $\sh$.
Similarly, if one puts $\bu_{S,\Q}:=\HQ\otimes \bu_S$, then one gets
\begin{equation}
    -\otimes \HQ:\md_{\bu_S}(\sh)^{\otimes}\rightarrow \md_{\bu_{S,\Q}}(\sh)^{\otimes}
\end{equation}
where the right hand side identifies with non-torsion $\bu_S$-modules. This $\oo$-functor is strongly compatible with the $6$-functors formalism.
Moreover, we state for future reference the following crucial fact:
\begin{thm}{\textup{\cite[Theorem 5.3.10]{ri10}, \cite[\S14.1]{cd19}, \cite[Proposition 3.35]{brtv}}}\label{thm riou}
Let $X$ be a scheme of finite Krull dimension an let $\MB{X}\in \textup{CAlg}(\textup{\textbf{SH}}_X)$ denote the spectrum of Beilinson motivic cohomology. Then the canonical morphism $1_X(1)[2]\rightarrow \bu_X\otimes \HQ=\bu_{X,\Q}$ induces an equivalence of commutative algebra objects
\begin{equation}
    \MB{X}(\beta):=Sym(\MB{X}(1)[2])[\nu^{-1}]\simeq \bu_{X,\Q}=\bu_X\otimes \HQ
\end{equation}
where $\nu$ is the free generator in degree $(1)[2]$.
\end{thm}
By using the theory of h-motives developed by the authors in \cite{cd16} one can define an $\ell$-adic realization $\oo$-functor
\begin{equation}\label{l-adic realization from beilinson motives}
    \Rl: \md_{\MB{}}(\textup{\textbf{SH}})\rightarrow \shvl(-)
\end{equation}
strongly compatible with the $6$-functors formalism, at least for noetherian schemes of finite Krull dimension.
Then, using the equivalence
\begin{equation}
    \md_{\bu}(\md_{\MB{}}(\textup{\textbf{SH}}))\simeq \md_{\bu}(\textup{\textbf{SH}})
\end{equation}
we obtain an $\ell$-adic realization $\oo$-functor
\begin{equation}
    \Rl: \md_{\bu}(\textup{\textbf{SH}})\xrightarrow{-\otimes \HQ} \md_{\bu_{\Q}}(\textup{\textbf{SH}})\rightarrow \md_{\Rl(\bu)}(\shvl(-))\simeq \md_{\Ql{}(\beta)}(\shvl(-))
\end{equation}
strongly compatible with the $6$-functors formalism. The last equivalence above holds as (\ref{l-adic realization from beilinson motives}) is symmetric monoidal and commutes with Tate twists (see \cite[Remark 3.43]{brtv}):
\begin{equation}
    \Rl_X(\bu_X)\underbracket{\simeq}_{(\ref{thm riou})}\Rl_X(Sym(\MB{X}(1)[2])[\nu^{-1}])\simeq Sym(\Rl_X(\MB{X})(1)[2])[\nu^{-1}]
\end{equation}
\begin{equation*}
\simeq \Ql{,X}(\beta)=Sym(\Ql{,X}(1)[2])[\nu^{-1}]
\end{equation*}
In this thesis, we will use the notation
\begin{equation}
    \Rlv_S:=\Rl_S \circ \Mv_S : \dgcatmo \rightarrow \md_{\Ql{,S}(\beta)}(\shvl(S))
\end{equation}

%% file: chapters/chapter02.tex
In this chapter we will investigate the connection between dg categories of singularities and dg categories of matrix factorizations. Since we wish to use a functorial point of view, we will first introduce some variants of categories of Landau-Ginzburg models. We will then define the dg categories of relative singularities and the dg categories of matrix factorizations, following the lead of \cite{brtv}. We will prove a structure theorem for the objects in a dg category of relative singularities. We will next use it in the particular case $n=1$ to show that dg categories of relative singularities of a single potential and matrix factorizations are related by a natural equivalence, improving the results existing in literature.
\section{\textsc{What is in in this chapter}}
A matrix factorizations of a pair $(B,f)$, where $B$ is a ring and $f\in B$ is the datum of two projective finitely-generated modules $(E_0,E_1)$ together with two morphisms $d_0:E_0\rightarrow E_1$, $d_1:E_1\rightarrow E_0$ such that $d_1\circ d_0=f\cdot id_{E_0}$ and $d_0\circ d_1=f\cdot id_{E_1}$. These objects, introduced by D. Eisenbud in \cite{eis80}, can be organized in a $2$-periodic dg category $\mf{B}{f}$ in a natural way. On the other hand, given such a pair $(B,f)$, we can define another dg category $\sing{B,f}$, called the dg category of relative singularities of the pair. The pushforward along the inclusion $\mathfrak{i}:Spec(B)\times^h_{\aff{1}{S}}S\rightarrow Spec(B)$ induces a dg functor
\begin{equation*}
    \mathfrak{i}_*:\sing{Spec(B)\times^h_{\aff{1}{S}}S}\rightarrow \sing{Spec(B)}
\end{equation*}
where $\sing{Z}$ stands for $\cohb{Z}/\perf{Z}$. Then $\sing{B,f}$ is defined as the fiber of this dg functor.
The connection between dg categories of relative singularities and dg categories of matrix factorizations has been first envisioned by R.O. Buchweitz and D. Orlov (see \cite{buch} and \cite{orl04}), who showed that if $B$ is regular ring and $f$ is a regular section, then (the homotopy-categories of) $\mf{B}{f}$ and $\sing{B,f}$ are equivalent. Notice under these hypothesis $Spec(B)\times^h_{\aff{1}{S}}S=Spec(B/f)$ and $\sing{B,f}\simeq \sing{B/f}$\footnote{Indeed, if $X$ is an underived (Noetherian) scheme, $\sing{X}=0$ if and only if $X$ is regular.}. The dg category of relative singularities was first introduced by J. Burke and M. Walker in \cite{bw12} in order to remove the regularity hypothesis on $B$. 

In the recent paper \cite{brtv} the authors show, along the way, that these equivalences are part of a lax-monoidal $\oo$-natural transformation
\begin{equation*}
    Orl^{-1,\otimes} : \sing{\bullet,\bullet} \rightarrow \mf{\bullet}{\bullet} : \lgm{1}^{op,\boxplus}\rightarrow \dgcatmo
\end{equation*}
and suggest that, in order to remove the regularity hypothesis on $f$, one should consider the derived zero locus $Spec(B)\times^h_{\aff{1}{S}}S$ instead of classical one. This remark comes from the observation that if $f$ is regular the two notion coincide and that if $B$ is regular and $f=0$, one can compute that both $\mf{B}{0}$ and $\sing{B,0}\simeq \sing{Spec(B)\times^h_{\aff{1}{S}}S}$ are equivalent to $\perf{B[u,u^{-1}]}$, where $u$ sits in cohomological degree $2$, while the classical zero locus of $f$ coincides with $B$ and thus the associated dg category of singularities is zero.

More generally, one can consider the dg categories of relative singularities of any pair $(B,\underline{f})$, where $\underline{f}\in B^n$ with $n\geq 1$, defined analogously to the case where $n=1$:
\begin{equation}
    \sing{B,\underline{f}}:=fiber\bigl ( \mathfrak{i}_*: \sing{Spec(B)\times^h_{\aff{n}{S}}S}\rightarrow \sing{Spec(B)}\bigr )
\end{equation}
There exists an algorithm which shows that this dg category is built up from $K(B,\underline{f})$-dg modules concentrated in $n+1$ degrees:
\begin{theorem-non}{\textup{(\ref{structure theorem})}}
Let $(Spec(B),\underline{f})$ be a $n$-dimensional affine Landau-Ginzburg model  over $S$. Then every object in the dg category of relative singularities $\sing{B,\underline{f}}$ is a retract of an object represented by a $K(B,\underline{f})$-dg module concentrated in $n+1$ degrees.
\end{theorem-non}
Moreover, when $n=1$, the algorithm mentioned above can be used to show that
\begin{theorem-non}{\textup{(\ref{structure theorem n=1})}}
Let
\begin{equation*}
\begindc{\commdiag}[20]
\obj(-10,0)[]{$(E,d,h)= 0$}
\obj(30,0)[]{$E_m$}
\obj(60,0)[]{$E_{m+1}$}
\obj(90,0)[]{$\dots$}
\obj(120,0)[]{$E_{m'-1}$}
\obj(150,0)[]{$E_{m'}$}
\obj(180,0)[]{$0$}
\mor(3,0)(27,0){$$}
\mor(33,-1)(57,-1){$d_{m}$}[\atright, \solidarrow]
\mor(57,1)(33,1){$h_{m+1}$}[\atright, \solidarrow]
\mor(63,-1)(87,-1){$d_{m+1}$}[\atright, \solidarrow]
\mor(87,1)(63,1){$$}[\atright, \solidarrow]
\mor(93,-1)(117,-1){$$}[\atright, \solidarrow]
\mor(117,1)(93,1){$h_{-1}$}[\atright, \solidarrow]
\mor(123,-1)(147,-1){$d_{-1}$}[\atright, \solidarrow]
\mor(147,1)(123,1){$h_0$}[\atright, \solidarrow]
\mor(153,0)(177,0){$$}
\enddc
\end{equation*}
be an object in $\textup{Coh}^s(B,f)$. Then the following equivalence holds in $\sing{B,f}$:
\begin{equation*}
(E,d,h)\simeq 
\begindc{\commdiag}[25]
\obj(0,0)[]{$\underbracket{\bigoplus_{i\in \mathbb{Z}}E_{2i-1}}_{-1}$}
\obj(40,0)[]{$\underbracket{\bigoplus_{i\in \mathbb{Z}}E_{2i}}_{0}$}
\mor(8,-1)(32,-1){$d+h$}[\atright, \solidarrow]
\mor(32,1)(8,1){$d+h$}[\atright, \solidarrow]
\enddc
\end{equation*}
Moreover, it is natural in $(E,d,h)$.
\end{theorem-non}
It is then possible to deduce the following
\begin{corollary-non}{\textup{(\ref{Orlov equivalence})}}
The lax monoidal $\oo$-natural transformation
\begin{equation*}
Orl^{-1,\otimes} : \sing{\bullet,\bullet} \rightarrow \mf{\bullet}{\bullet} : \lgm{1}^{op,\boxplus}\rightarrow \dgcatmo
\end{equation*}
constructed in \cite[\S2.4]{brtv} defines a lax-monoidal $\oo$-natural equivalence.
\end{corollary-non}
Recall that, for an affine LG model $(Spec(B),f)$, 
\begin{equation*}
    Orl^{-1}_{(B,f)} : \sing{B,f}\xrightarrow{\simeq} \mf{B}{f}
\end{equation*}
is defined as the dg functor
\begin{equation*}
    \begindc{\commdiag}[20]
    \obj(0,30)[1]{$(E,d,h)$}
    \obj(70,30)[2]{$\bigoplus_{i\in \mathbb{Z}}E_{2i-1}$}
    \obj(120,30)[3]{$\bigoplus_{i\in \mathbb{Z}}E_{2i}$} 
    \mor{1}{2}{$$}[\atright,\aplicationarrow]
    \mor(78,29)(112,29){$d+h$}[\atright,\solidarrow]
    \mor(112,31)(78,31){$d+h$}[\atright,\solidarrow]
    \obj(0,0)[4]{$(E',d',h')$}
    \obj(70,0)[5]{$\bigoplus_{i\in \mathbb{Z}}E'_{2i-1}$}
    \obj(120,0)[6]{$\bigoplus_{i\in \mathbb{Z}}E'_{2i}$} 
    \mor{4}{5}{$$}[\atright,\aplicationarrow]
    \mor(78,-1)(112,-1){$d'+h'$}[\atright,\solidarrow]
    \mor(112,1)(78,1){$d'+h'$}[\atright,\solidarrow]
    \mor{1}{4}{$\phi$}
    \mor{2}{5}{$\oplus \phi_{2i-1}$}
    \mor{3}{6}{$\oplus \phi_{2i}$}
    \enddc
\end{equation*}

The corollary above improves all the previous results on the equivalence between the dg categories of singularities and the dg category of matrix factorizations as it removes the regularity assumption on the potential.

In the last section we discuss how one can force $2$-periodicity on $\sing{X,\underline{f}})$ in the case $n\geq 1$

\begin{rmk}
This chapter corresponds to the author's paper \cite{pi19}.
\end{rmk}
\section{\textsc{Preliminaries}}
\subsection{\textsc{Higher dimensional Landau-Ginzburg models}}
\begin{context}
Assume that $A$ is a local, Noetherian regular ring of finite Krull dimension.
\end{context}
Recall that the category of flat $S$-schemes of finite type together with a potential (i.e. a map to $\aff{1}{S}$) is often referred to as the category of Landau-Ginzburg models. Its morphisms are those morphisms of $S$-schemes which are compatible with the potential. Moreover, this category has a natural symmetric monoidal enhancement due to the fact that $\aff{1}{S}$ is a scheme in abelian groups. It is very easy to generalize this category to the case where schemes are provided with multipotentials, i.e. with maps to $\aff{n}{S}$, for any $n\geq 1$.
\begin{defn}
Fix $n \geq 1$. Define the category of $n$-dimensional Landau-Ginzburg models over $S$ ($n$-LG models over $S$ for brevity) to be the full subcategory of $\textup{\textbf{Sch}}_{S/\mathbb{A}^{n}_S}$ spanned by those objects
\[
\begindc{\commdiag}[20]
\obj(0,20)[1]{$X$}
\obj(60,20)[2]{$\aff{n}{S}$}
\obj(30,0)[3]{$S$}
\mor{1}{2}{$\underline{f}=(f_1,\dots,f_n)$}
\mor{1}{3}{$p$}[\atright,\solidarrow]
\mor{2}{3}{$\textup{proj.}$}
\enddc
\]
where $p$ is a flat morphism. Denote this category by $\lgm{n}$ and its objects by $(X,\underline{f})$. 

For convenience, we also introduce the following (full) subcategories of $\lgm{n}$:
\begin{itemize}
\item $\lgm{n}^{\text{fl}}\subseteq \lgm{n}$, the category of flat Landau-Ginzburg models of order $n$ over $S$, spanned by those objects $(X,\underline{f})$ such that $\underline{f}$ is flat and by $(S,\underline{0})$
\item $\lgm{n}^{\text{aff}}\subseteq \lgm{n}$, the category of affine Landau-Ginzburg models of order $n$ over $S$, spanned by those objects $(X,\underline{f})$ such that $X$ is affine
\item $\lgm{n}^{\text{aff,fl}}\subseteq \lgm{n}$, the category of flat, affine Landau-Ginzburg models of order $n$ over $S$, spanned by those objects $(X,\underline{f})$ such that $X$ is affine and $\underline{f}$ is flat and by $(S,\underline{0})$
\end{itemize}
\end{defn}
\begin{construction}
As in \cite{brtv}, we can enhance $\lgm{n}$ (and its variants) with a symmetric monoidal structure. Consider the "sum morphism"\footnote{notice that it is flat}
\begin{equation}
+ : \aff{n}{S}\times_S \aff{n}{S}\rightarrow \aff{n}{S}
\end{equation}
on $\aff{n}{S}$, corresponding to 
\[
A[T_1,\dots,T_n]\rightarrow A[X_1,\dots,X_n]\otimes_A A[Y_1,\dots,Y_n]
\]
\[
T_i\mapsto X_i\otimes 1 + 1 \otimes Y_i \hspace{0.5cm} i=1,\dots,n
\]
Then define 
\begin{equation}
\boxplus : \lgm{n}\times \lgm{n} \rightarrow \lgm{n}
\end{equation}
by the formula
\[
\bigl ( (X,\underline{f}), (Y,\underline{g})\bigr )\mapsto (X\times_S Y, \underline{f}\boxplus \underline{g})
\]
Here, $\underline{f}\boxplus \underline{g}$ is the following composition
\[
X\times_S Y\xrightarrow{\underline{f}\times \underline{g}} \aff{n}{S}\times_S \aff{n}{S}\xrightarrow{+} \aff{n}{S}
\]
Notice that $X\times_S Y$ is still flat over $S$, whence this functor is well defined. It is also easy to remark that $\boxplus$ is associative - i.e. there exist natural isomorphisms $\bigl ( (X,\underline{f})\boxplus (Y,\underline{g}) \bigr) \boxplus (Z,\underline{h})\simeq  (X,\underline{f})\boxplus \bigl ((Y,\underline{g}) \boxplus (Z,\underline{h})\bigr) $ - and that for any object $(X,\underline{f})$, $(S,\underline{0})\boxplus (X,\underline{f})\simeq (X,\underline{f})\simeq (X,\underline{f})\boxplus (S,\underline{0})$. More briefly, $\bigl ( \lgm{n}, \boxplus, (S,\underline{0})\bigr )$ is a symmetric monoidal category. It is not hard to see that this construction works on $\lgm{n}^{\text{fl}}$, $\lgm{n}^{\text{aff}}$ and $\lgm{n}^{\text{aff,fl}}$ too. Indeed,  this is clear for $\lgm{n}^{\text{aff}}$  and if $\underline{f}$ and $\underline{g}$ are flat morphisms, so is $\underline{f}\times \underline{g}$ and therefore $\underline{f}\boxplus \underline{g}$ is a composition of flat morphisms. 
\end{construction}
\begin{notation}
We will denote by $\lgm{n}^{\boxplus}$ (resp. $\lgm{n}^{\text{fl},\boxplus}$, $\lgm{n}^{\text{aff},\boxplus}$ and  $\lgm{n}^{\text{aff,fl}, \boxplus}$) this symmetric monoidal categories.
\end{notation}
\begin{rmk}
Notice that $\lgm{1}^{\boxplus}$ is exactly the symmetric monoidal category
$\textup{LG}_S^{\boxplus}$ defined in \cite[\S 2]{brtv}.
\end{rmk}
\begin{rmk}
Fix $n\geq 1$. Notice that the symmetric group $\mathcal{S}_n$ acts on the category of $n$-LG models over $S$. Indeed, for any $\sigma \in \mathcal{S}$ and for any $(X,\underline{f}) \in \lgm{n}$, we can define
\[
\sigma \cdot (X,\underline{f}) := (X,\sigma \cdot \underline{f})
\]
\end{rmk}

\subsection{\textsc{dg category of singularities of an \texorpdfstring{$n$}{n}-LG model}}\label{dg Category of singularities of an n-LG model over S}
It is a classic theorem due to due to M.~Auslander-D.A.~Buchsbaum (\cite[Theorem 4.1]{ab56}) and J.P.~Serre (\cite[Théorème 3]{se55}) that a Noetherian local ring $R$ is regular if and only if it has finite global dimension. This extremely important fact can be rephrased by saying that every object in $\cohb{R}$ is equivalent to an object in $\perf{R}$. In particular, $R$ is regular if and only if $\cohb{R}/\perf{R}$ is zero. This explains why the quotient above is called category of singularities. 

Before going on with the precise definitions, let us fix some notation.

Let $(X,\underline{f})$ be a $n$-LG model over $S$. Then consider the (derived) zero locus of $\underline{f}$, i.e. the (derived) fiber product

\begin{equation}
\begindc{\commdiag}[20]
\obj(0,15)[1]{$X_0$}
\obj(30,15)[2]{$X$}
\obj(0,-15)[3]{$S$}
\obj(30,-15)[4]{$\aff{n}{S}$}
\mor{1}{2}{$\mathfrak{i}$}
\mor{1}{3}{$$}

\mor{2}{4}{$\underline{f}$}
\mor{3}{4}{$\text{zero}$}
\enddc
\end{equation}

\begin{rmk}
Notice that $X_0\simeq X\times_{\aff{n}{S}} S$ coincides with the classical zero locus whenever $(X,\underline{f})$ belongs to $\lgm{n}^{\text{fl}}$. In general, we always have a closed embedding $t : X\times_{\aff{n}{S}} S=\pi_0(X_0) \rightarrow X_0$.
\end{rmk}
\begin{rmk}
As $S\xrightarrow{\text{zero}} \aff{n}{S}$ is lci and this class of morphism is closed under derived fiber products, we get that $\mathfrak{i}: X_0\rightarrow X$ is a lci morphism of derived schemes.
\end{rmk}

\begin{rmk}
Analogously to \cite[\S2]{brtv}, we have the following inclusions
\[
\perf{X}\subseteq \cohb{X} \subseteq \cohm{X}\subseteq \qcoh{X}
\]
\[
\perf{X_0}\subseteq \cohbp{X_0}{X}\subseteq \cohb{X_0}\subseteq \cohm{X_0}\subseteq \qcoh{X_0}
\]
Indeed, being $X$ and $X_0$ eventually coconnective (see \cite[\S4, \textbf{Definition 1.1.6}]{gr17}), we have the inclusions $\perf{X}\subseteq \cohb{X}$ and $\perf{X_0}\subseteq \cohb{X_0}$. Moreover, as $\mathfrak{i}$ is lci, by \cite{to12}, we have that $\mathfrak{i}_*$ preserves perfect complexes. Thus, the inclusion $\perf{X_0}\subseteq \cohbp{X_0}{X}$ holds.
\end{rmk}
\begin{rmk}
As explained in \cite[Remark 2.14]{brtv}, the dg categories $\perf{X}$, $\perf{X_0}$, $\cohb{X}$, $\cohb{X_0}$ and $\cohbp{X_0}{X}$ are idempotent complete. Indeed, the same argument provided in \textit{loc.cit.} for the case $n=1$ works in general.
\end{rmk}
Notice that all the results in \cite[\S 2.3.1]{brtv} are not specific of the monopotential case and they remain valid in our situation. We will recall these statements for the reader's convenience and refer to \ital{loc. cit.} for the proofs, which remain untouched.
\begin{prop}{\textup{\cite[Proposition 2.17]{brtv}}}
Let $(X,\underline{f})\in \lgm{n}$. Then the inclusion functor induces an equivalence 
\begin{equation}
\cohbp{X_0}{X}\simeq \cohmp{X_0}{X}
\end{equation}
In particular, the following square is cartesian in $\dgcatm$
\begin{equation}
\begindc{\commdiag}[20]
\obj(0,15)[1]{$\cohm{X_0}$}
\obj(60,15)[2]{$\cohm{X}$}
\obj(0,-15)[3]{$\cohbp{X_0}{X}$}
\obj(60,-15)[4]{$\perf{X}$}
\mor{1}{2}{$\mathfrak{i}_*$}
\mor{3}{1}{$$}
\mor{4}{2}{$$}
\mor{3}{4}{$\mathfrak{i}_*$}
\enddc
\end{equation}
\end{prop}
We now give definitions for the relevant dg categories of singularities. The reader should be aware that there are plenty of this objects that one can consider, and we will define some of them later on. The following category, known as category of absolute singularities, first appeared in \cite{orl04}. The following definition is a dg enhancement of the original one. It already appears in \cite{brtv}.
\begin{defn}\label{absolute sing}
Let $Z$ be a derived scheme of finite type over $S$ whose structure sheaf is cohomologically bounded. The dg category of absolute singularities of $Z$ is the dg quotient (in $\dgcatm$)
\begin{equation}
\sing{Z}:= \cohb{Z}/\perf{Z}
\end{equation}
\end{defn}
\begin{rmk}
Notice that the finiteness hypothesis on $Z$ in \textbf{Definition} (\ref{absolute sing}) are absolutely indispensable, as otherwise $\perf{Z}$ may not contained in $\cohb{Z}$.
\end{rmk}
\begin{rmk}
It is well known that, for an underived (Noetherian) scheme $Z$, the dg category $\sing{Z}$ is zero if and only if the scheme is regular. On the other hand, when we allow $Z$ to be a derived scheme, $\sing{Z}$ may be non trivial even if the underlying scheme is regular. For example, consider $Z=Spec(A\otimes^{\mathbb{L}}_{A[T]}A)$.
\end{rmk}

Following \cite{brtv} we next consider the dg category of singularity associated to an $n$-dimesional LG-model.
\begin{defn}
Let $(X,\underline{f})\in \lgm{n}$. The dg category of singularities of $(X,\underline{f})$ is the following fiber in $\dgcatm$
\begin{equation}
\sing{X,\underline{f}}:= Ker\bigl ( \mathfrak{i}_* : \sing{X_0}\rightarrow \sing{X}\bigr )
\end{equation}
\end{defn}
\begin{rmk}
Notice that $\sing{X,\underline{f}}$ is a full sub-dg category of $\sing{X_0}$ (see \cite[Remark 2.24]{brtv}). Moreover, these two dg categories coincide whenever $X$ is a regular $S$-scheme.
\end{rmk}
\begin{prop}{\textup{(See \cite[Proposition 2.25]{brtv})}}
Let $(X,\underline{f})$ be a $n$-dimensional LG model over $S$. Then there is a canonical equivalence
\begin{equation}
\cohbp{X_0}{X}/\perf{X_0}\simeq \sing{X,\underline{f}}
\end{equation}
where the quotient on the left is taken in $\dgcatm$. 
\end{prop}

We now re-propose, for the multi-potential case, the strict model for $\cohbp{X_0}{X}$ which was first introduced in \cite{brtv}.

\begin{construction}\label{strict model derived tensor product}
Let $(Spec(B),\underline{f})\in \lgm{n}^{\textup{aff}}$. Consider the Koszul complex $K(B,\underline{f})$
\begin{equation}
0\rightarrow \bigwedge^n (B\varepsilon_1\oplus \dots \oplus B\varepsilon_n) \rightarrow \dots \rightarrow \bigwedge^2 (B\varepsilon_1\oplus \dots \oplus B\varepsilon_n) \rightarrow (B\varepsilon_1\oplus \dots \oplus B\varepsilon_n) \rightarrow B \rightarrow 0
\end{equation}
concetrated in degrees $[-n,0]$. The differential is given by
\[
\bigwedge^k (B\varepsilon_1\oplus \dots \oplus B\varepsilon_n)\rightarrow \bigwedge^{k-1} (B\varepsilon_1\oplus \dots \oplus B\varepsilon_n)
\]
\[
v_1\wedge \dots \wedge v_k\mapsto \sum_{i=1}^k (-1)^{1+i}\phi(v_i) v_1\wedge \dots \wedge \hat{v_i} \wedge \dots \wedge v_k
\]
where $\phi : B^n \rightarrow B$ is the morphism corresponding to the matrix $[f_1\dots f_n]$. Multiplication is given by concatenation. Notice that $K(B,\underline{f})$ is a cofibrant $B$-module and that we always have a truncation morphism $K(B,\underline{f})\rightarrow B/\underline{f}$, which is a quasi-isomorphism whenever $\underline{f}$ is a regular sequence.

Therefore, we can present $K(B,\underline{f})$ as the cdga $B[\varepsilon_1,\dots,\varepsilon_n]$, where the $\varepsilon_i$'s sit in degree $-1$ and are subject to the following conditions:
\[
d(\varepsilon_i)=f_i \hspace{0.5cm} i=1,\dots,n
\]
\[
\varepsilon_i^2=0 \hspace{0.5cm} i=1,\dots,n
\]
\[
\varepsilon_{i_1}\dots \varepsilon_{i_k}= (-1)^{\sigma} \varepsilon_{i_{\sigma(1)}}\dots \varepsilon_{i_{\sigma(k)}} \hspace{0.5cm} \{i_1,\dots,i_k\}\subseteq \{1,\dots,n\}, \sigma \in \mathcal{S}_k
\]
\end{construction}
\begin{ex}
For instance, when $n=1$, $K(B,f)$ is the cdga $B\varepsilon \xrightarrow{f} B$ concentrated in degrees $[-1,0]$ and, when $n=2$, $K(B,(f_1,f_2))$ is the cdga $B\varepsilon_1 \varepsilon_2\xrightarrow{\begin{bmatrix} -f_2 \\ f_1 \end{bmatrix}} B\varepsilon_1 \oplus B\varepsilon_2 \xrightarrow{\begin{bmatrix} f_1 & f_2 \end{bmatrix}} B$ concentrated in degrees $[-2,0]$.
\end{ex}
\begin{rmk}
Notice that $K(B,\underline{f})$ provides a model for the cdga associated to the simplicial commutative algebra $B\otimes^{\mathbb{L}}_{A[T_1,\dots,T_n]}A$. Indeed, this can be computed explicitly for $n=1$ and the general case follows from the compatibility of the Dold-Kan correspondence with (derived) tensor products.
\end{rmk}
This strict model for the derived zero locus of an affine LG model of order $n$ over $S$ gives us strict models for the relevant categories of modules too. Following \cite{brtv}:
\begin{itemize}
\item There is an equivalence of $A$-linear dg categories between $\qcoh{X_0}$ and the dg category (over $A$) of cofibrant $K(B,\underline{f})$-dg modules, which we will denote $\widehat{K(B,\underline{f})}$. A $K(B,\underline{f})$-dg module is the datum of a cochain complex of $B$-modules $(E,d)$, together with $n$ morphisms $h_1,\dots,h_n : E\rightarrow E[1]$ of degree $-1$ such that
\begin{equation}\label{identities h_i}
\begin{cases}
h_i^2=0 \hspace{0.5cm} i=1,\dots,n \\
[d,h_i]=f_i \hspace{0.5cm} i=1,\dots,n \\
[h_i,h_j]=0 \hspace{0.5cm} i,j=1,\dots,n \\
\end{cases}
\end{equation}
\item Under this equivalence, $\cohb{X_0}\subseteq \qcoh{X_0}$ corresponds to the full sub-dg category of $\widehat{K(B,\underline{f})}$ spanned by those modules of cohomologically bounded amplitude and whose cohomology is coherent over $B/\underline{f}$
\item Under this equivalence, $\perf{X_0}\subseteq \qcoh{X_0}$ corresponds to the full sub-dg category of $\widehat{K(B,\underline{f})}$ spanned by those modules which are homotopically finitely presented
\end{itemize}
\begin{rmk}\label{commutation d and h}
Notice that, for any $K(B,\underline{f})$-dg module, for any $1\leq k \leq n$ and for any $\{ i_1,\dots,i_k\}\subseteq \{1,\dots,n\}$ (where the $i_j$'s are pairwise distinguished),  the following formula holds:
\[
[d,h_{i_1}\dots h_{i_k}]=\sum_{j=1}^k(-1)^j f_{i_j} h_{i_1}\circ \dots \widehat{h_{i_j}} \dots \circ h_{i_k}
\]
\end{rmk}
\begin{rmk}
As in the mono-potential case (see \cite[Remark 2.30]{brtv}), $\mathfrak{i}_* : \qcoh{X_0}\rightarrow \qcoh{X}$   
corresponds, under these equivalences, to the forgetful functor $\widehat{K(B,\underline{f})}\rightarrow \widehat{B}$ ($K(B,\underline{f})$ is a cofibrant $B$-module). 
\end{rmk}
We propose the following straightforward generalization of \cite[Construction 2.31]{brtv} as a strict model for $\textup{Coh}^b(X_0)_{\textup{Perf}(X)}$:
\begin{construction}
Let $\textup{Coh}^s(B,\underline{f})$ be the $A$-linear sub-dg category of $K(B,\underline{f})$ spanned by those modules whose image along the forgetful functor $K(B,\underline{f})-\text{dgmod} \rightarrow B-\text{dgmod}$ is a strictly perfect complex of $B$-modules. In particular, an object $E$ in $\textup{Coh}^s(B,\underline{f})$ is a degree-wise projective cochain complex of $B$-modules together with $n$ morphisms $h_1,\dots,h_n$ of degree $-1$ satisfying the identities (\ref{identities h_i}). As $A$ is a local ring, it is clear that $\textup{Coh}^s(B,\underline{f})$ is a locally flat $A$-linear dg category.
\end{construction}
\begin{lmm}\label{explicit model sing}
Let $(X,\underline{f})=(Spec(B),\underline{f})$ be a $n$-dimensional affine LG model over $S$. Then the cofibrant replacement dg functor induces an equivalence
\begin{equation}
\textup{Coh}^s(B,\underline{f})[q.iso^{-1}]\simeq \textup{Coh}^s(B,\underline{f})/\textup{Coh}^{s,acy}(B,\underline{f})\simeq \cohbp{X_0}{X}
\end{equation}
where $\textup{Coh}^{s,acy}(B,\underline{f})$ is the full sub-dg category of $\textup{Coh}^s(B,\underline{f})$ spanned by acyclic complexes. In particular, this implies that we have equivalences of dg categories
\begin{equation}
\textup{Coh}^{s}(B,\underline{f})/\textup{Perf}^s(B,\underline{f})\simeq \cohbp{X_0}{X}/\perf{X_0}\simeq \sing{X,\underline{f}}
\end{equation}
where $\textup{Perf}^s(B,\underline{f})$ is the full sub-dg category of $\textup{Coh}^s(B,\underline{f})$ spanned by those modules which are perfect over $K(B,\underline{f})$.
\end{lmm}
\begin{proof}
See \cite[Lemma 2.33]{brtv}. The same proof  holds true in our situation too.
\end{proof}

We now exhibit the functorial properties of $\textup{Coh}^s(\bullet,\bullet)$. Let $u:(Spec(C),\underline{g})\rightarrow (Spec(B,\underline{f}))$ be a morphism in $\lgm{n}^{\textup{aff}}$. Define the dg functor
\begin{equation}
u^* : \textup{Coh}^s(B,\underline{f})\rightarrow \textup{Coh}^s(C,\underline{g})
\end{equation}
by the law
\[
E\mapsto E\otimes_B C
\]
Notice that this dg functor is well defined as $E\otimes_B C$ is strictly bounded and degree-wise $C$-projective. It is clear that if two composable morphisms 
\[
(Spec(B),\underline{f})\xrightarrow{u}(Spec(B'),\underline{f}')\xrightarrow{u'}(Spec(B''),\underline{f}'')
\]
are given, $u'^*\circ u^*\simeq (u'\circ u)^*$ are equivalent dg functors $\textup{Coh}^s(B'',\underline{f}'')\rightarrow \textup{Coh}^s(B,\underline{f})$. It is also clear that $id_{(Spec(B),\underline{f})}^*\simeq id_{\textup{Coh}^s(B,\underline{f})}$ and that this law is (weakly) associative and (weakly) unital. In other words,
\begin{equation}\label{Coh^s}
\textup{Coh}^s(\bullet,\bullet) : \lgm{S}^{\textup{aff, op}}\rightarrow \dgcatlf
\end{equation}
has the structure of a pseudo-functor.
We next produce a lax-monoidal structure on this pseudo-functor. We begin by producing a map
\begin{equation}\label{right lax map Coh^s}
\textup{Coh}^s(B,\underline{f})\otimes \textup{Coh}^s(C,\underline{g})\rightarrow \textup{Coh}^s(B\otimes_A C,\underline{f}\boxplus \underline{g})
\end{equation}
Consider the following diagram
\begin{equation}
\begindc{\commdiag}[20]
\obj(0,30)[1]{$Spec(K(B\underline{f}))\times_S^h Spec(C,\underline{g})$}
\obj(80,30)[2]{$Spec(K(B\otimes_A C,\underline{f}\boxplus \underline{g}))$}
\obj(140,30)[3]{$Spec(B\otimes_A C)$}
\obj(0,0)[4]{$S$}
\obj(80,0)[5]{$\aff{n}{S}$}
\obj(140,0)[6]{$\aff{n}{S}\times_S \aff{n}{S}$}
\obj(80,-30)[7]{$S$}
\obj(140,-30)[8]{$\aff{n}{S}$}
\mor{1}{2}{$\phi$}
\mor{2}{3}{$\psi$}
\mor{1}{4}{$$}

\mor{2}{5}{$$}

\mor{3}{6}{$\underline{f}\times \underline{g}$}
\mor{4}{5}{$\text{zero}$}
\mor{5}{6}{$\begin{bmatrix}id & -id\end{bmatrix}$}
\mor{5}{7}{$$}

\mor{6}{8}{$+$}
\mor{7}{8}{$\text{zero}$}
\enddc
\end{equation}

Notice that all the squares in this diagram are (homotopy) cartesian and that all the horizontal maps are lci morphisms of (derived) schemes. Write
\[
K(B,\underline{f})=B[\varepsilon_1,\dots,\varepsilon_n]
\]
\[
K(C,\underline{g})=C[\delta_1,\dots,\delta_n]
\]
\[
K(B\otimes_A C,\underline{f}\boxplus \underline{g})=B\otimes_A C[\gamma_1,\dots,\gamma_n]
\]
where all the $\varepsilon_i$'s, $\delta_i$'s and $\gamma_i$'s sit in degree $-1$ and are subject to the relations (\ref{strict model derived tensor product}). Then $\phi$ corresponds to the morphism of cdga's 
\[
K(B\otimes_A C,\underline{f}\boxplus \underline{g})\rightarrow K(B,\underline{f})\otimes_A K(C,\underline{g})
\]
\[
\gamma_i \mapsto \varepsilon_i\otimes 1 + 1\otimes \delta_i \hspace{0.5cm} i=1,\dots, n
\]
which is the identity in degree zero, while $\psi$ and $\psi \circ \phi$ are just the obvious inclusion of $B\otimes_A C$. 

Then, we define (\ref{right lax map Coh^s}) by 
\[
(\mathcal{F},\mathcal{G})\mapsto \mathcal{F}\boxtimes \mathcal{G}:= \phi_*\bigl (pr_1^*\mathcal{F}\otimes_{K(B,\underline{f})\otimes_A K(C,\underline{g})} pr_2^*\mathcal{G}\bigr )
\]
where $pr_1$ and $pr_2$ are the projections from the (homotopy) fiber product 
\[
Spec(K(B,\underline{f}))\times^h_S Spec(K(C,\underline{g}))
\]
to $Spec(K(B,\underline{f}))$ and $Spec(K(C,\underline{g}))$ respectively. We need to show that $\mathcal{F}\boxtimes \mathcal{G}$ lies in $\textup{Coh}^s(B\otimes_A C,\underline{f}\boxplus \underline{g})$. This is equivalent to the statement that the underlying complex of $pr_1^*\mathcal{F}\otimes_{K(B,\underline{f})\otimes_A K(C,\underline{g})} pr_2^*\mathcal{G}$ is perfect over $B\otimes_A C$. Consider the (homotopy) cartesian square
\begin{equation}
\begindc{\commdiag}[20]
\obj(-60,15)[1]{$Spec(K(B,\underline{f}))\times_S^h Spec(K(C,\underline{g}))$}
\obj(60,15)[2]{$Spec(K(C,\underline{g}))$}
\obj(-60,-15)[3]{$Spec(K(B,\underline{f}))$}
\obj(60,-15)[4]{$S$}
\mor{1}{2}{$pr_2$}
\mor{1}{3}{$pr_1$}
\mor{2}{4}{$q$}
\mor{3}{4}{$p$}
\enddc
\end{equation}
and notice that we have the following chain of equivalences
\begin{equation}
p_*pr_{1*} \bigl ( pr_1^*\mathcal{F}\otimes_{K(B,\underline{f})\otimes_A K(C,\underline{g})} pr_2^*\mathcal{G} \bigr ) \underbracket{\simeq}_{\text{proj. form.}} p_* \bigl ( \mathcal{F}\otimes_{K(B,\underline{f})} pr_{1*}pr_2^*\mathcal{G} \bigr ) 
\end{equation}
\[
\underbracket{\simeq}_{\text{base change}} p_* \bigl ( \mathcal{F}\otimes_{K(B,\underline{f})} p^*q_*\mathcal{G} \bigr ) \underbracket{\simeq}_{\text{proj. form.}} p_*\mathcal{F}\otimes_A q_*\mathcal{G}
\]
As $\mathcal{F} \in \textup{Coh}^s(B,\underline{f})$, $\mathcal{G}\in \textup{Coh}^s(C,\underline{g})$ and perfect complexes are stable under tensor product, we conclude.
 
We next exhibit the lax unit\footnote{$\underline{A}$ denotes the $\otimes$-unit in $\dgcatlfo$, i.e. the dg category with one object $\bullet$ whose complex of endomorphisms $\hm_{\underline{A}}(\bullet,\bullet)$ is just $A$ in degree $0$.}
\begin{equation}\label{lax unit Coh^s}
\underline{A} \rightarrow Coh^s(A,\underline{0})
\end{equation}
This is simply the dg functor defined by
\[
\bullet \mapsto A
\]
where $A$ (concentrated in degree $0$) is seen as a module over $K(A,\underline{0})$ in the obvious way, i.e. the $\varepsilon_i$'s act via zero. 

This defines a (right) lax monoidal structure on (\ref{Coh^s})
\begin{equation}\label{Coh^s monoidal}
\textup{Coh}^{s,\boxtimes}(\bullet,\bullet): \lgm{n}^{\textup{aff,op}, \boxplus}\rightarrow \dgcatlfo
\end{equation}
\begin{rmk}
Notice that the same structure defines a lax monoidal structure on the functor
\begin{equation}
\textup{Perf}^{s,\boxtimes}(\bullet,\bullet): \lgm{n}^{\textup{aff,op}, \boxplus}\rightarrow \dgcatlfo
\end{equation}
\end{rmk}
By the same technical arguments of \cite[Construction 2.34, Construction 2.37]{brtv} we produce a (right) lax monoidal $\oo$-functor
\begin{equation}\label{oo Coh^s monoidal}
\cohbp{\bullet}{\bullet}^{\otimes}: \lgm{n}^{\textup{aff,op},\boxplus}\rightarrow \dgcatmo
\end{equation}
In order to define the lax monoidal $\oo$-functor
\begin{equation}\label{functor sing}
\sing{\bullet,\bullet}^{\otimes} : \lgm{n}^{\text{aff,op},\boxplus}\rightarrow \dgcatmo
\end{equation}
consider the category $\pdgcat$ whose objects are pairs $(T,S)$, where T is an $A$-linear dg category and $S$ a class of morphisms in $T$. Given two objects $(T,S)$ and $(T',S')$, morphisms $(T,S)\rightarrow (T',S')$ are those dg functors $F:T\rightarrow T'$ such that $S$ is sent into $S'$. Composition and identities are defined in the obvious way.
Given a morphism $(T,S)\rightarrow (T',S')$, we say that it is a Dwyer-Kan equivalence if the underlying dg functor is so (i.e. it  is a quasi-equivalence). We denote the class of Dwyer-Kan equivalences in $\pdgcat$ by $W_{DK}$.

Notice that $\pdgcat$ inherits a symmetric monoidal structure from $\dgcatlfo$ by setting $(T,S)\otimes (T',S')=(T\otimes T',S\otimes S')$. We will refer to this symmetric monoidal category by $\pdgcato$. As we are considering locally flat dg categories, it is immediate that this tensor structure is compatible with DK equivalences.
For any $n$-dimensional affine LG model $(Spec(B),\underline{f})$ over $S=Spec(A)$, define $W_{\textup{Perf}^s(B,\underline{f})}$ as the class of morphisms $\bigl ( 0\rightarrow E \bigr )_{E\in \textup{Perf}^s(B,\underline{f})}$ in $\textup{Coh}^s(B,\underline{f})$. Consider the functor
\begin{equation}\label{canonical pair dg cat}
\lgm{n}^{\textup{aff,op}} \rightarrow \pdgcat
\end{equation}
\[
(Spec(B,\underline{f}))\mapsto (\textup{Coh}^s(B,\underline{f}),W_{\textup{Perf}^s(B,\underline{f})})
\]
If $E\in \textup{Perf}^s(B,\underline{f})$ and $F\in 
\textup{Perf}^s(C,\underline{g})$, then $(0\rightarrow E)\otimes (0\rightarrow F)\in W_{\textup{Perf}^s(B,\underline{f})}\otimes W_{\textup{Perf}^s(C,\underline{g})}$ is sent to $0\rightarrow E\boxtimes F$ via (\ref{right lax map Coh^s}), which belongs to $W_{\textup{Perf}^s(B\otimes_A C,\underline{f}\boxplus \underline{g})}$. Then the functor (\ref{canonical pair dg cat}) has a lax monoidal enhancement 
\begin{equation}\label{lax canonical pair dg cat}
\lgm{n}^{\textup{aff,op}, \boxplus} \rightarrow \pdgcato
\end{equation}

By \cite[Construction 2.34 and Construction 2.37]{brtv}, there is a symmetric monoidal $\oo$-functor
\begin{equation}\label{loc dg}
loc_{dg}^{\otimes} : \pdgcato [W_{DK}^{-1}] \rightarrow \dgcatdko
\end{equation}
sending a pair $(T,S)$ to the dg localization $T[S^{-1}]_{dg}$.

We finally define (\ref{functor sing}) as the following composition
\begin{equation}\label{sing affine}
\lgm{n}^{\textup{aff,op}, \boxplus} \xrightarrow{loc. \circ (\ref{lax canonical pair dg cat})}\pdgcato[W_{DK}^{-1}]\xrightarrow{(\ref{loc dg})} \dgcatdko \xrightarrow{\text{loc.}} \dgcatmo
\end{equation}
Notice that $(Spec(B),\underline{f})\in \lgm{n}^{\textup{aff}}$ is sent to $\sing{B,\underline{f}}$ by Lemma \ref{explicit model sing} and by the fact that the quotient $\textup{Coh}^s(B,\underline{f})/\textup{Perf}^s(B,\underline{f})$ is, by definition, the dg localization $\textup{Coh}^s(B,\underline{f})[W_{\textup{Perf}^s(B,\underline{f})}]$ (see \cite[\S8.2]{to07}).

Moreover, we use the same notation
\begin{equation}
    \sing{\bullet,\bullet}^{\otimes} :\lgm{n}^{\textup{op},\boxplus}\rightarrow \dgcatmo
\end{equation}
for the extension to (non-affine) LG models. This can be done by considering the monoidal Kan extension of (\ref{sing affine}) along the inclusion $\lgm{n}^{\textup{aff,op},\boxplus}\subseteq \lgm{n}^{\textup{op},\boxplus}$.
\begin{rmk}
If $n=1$, the lax monoidal structure on the $\oo$-functor $\sing{\bullet,\bullet}^{\otimes}$ identifies with the lax monoidal structure on the $\oo$-functor defined in \cite[Proposition 2.45]{brtv}.
\end{rmk}
\section{\textsc{Orlov's theorem}}
\subsection{\textsc{The structure of} Sing(B,\underline{f})}

In this section we will prove that, in the category of relative singularities $\sing{B,\underline{f}}$ associated to a $n$-dimensional affine Landau-Ginzburg model over $S$, every object is a retract of an object that can be represented by a $K(B,\underline{f})$-dg module concentrated in $n+1$ degrees. We begin with the following observation:
\begin{lmm}\label{cone K(B,f) dg mod}
Let $\phi : (E,d,h)\rightarrow (E',d',h')$ be a cocycle-morphism of $K(B,f)$-dg modules\footnote{Here $d$ (resp. $d'$) stands for the differential and $h^i$ (resp.$h'^i$ ) stands for the action of $\varepsilon_i$, where 

$K(B,\underline{f})=0\rightarrow \underbracket{B \varepsilon_1\dots \varepsilon_n}_{-n}\rightarrow \dots \rightarrow \underbracket{B\varepsilon_1\oplus \dots \oplus B\varepsilon_n}_{-1}\rightarrow \underbracket{B}_{0}$}.  Then the cone of $\phi$ is given by
\begin{equation}
\begindc{\commdiag}[20]
\obj(0,0)[]{$\underbracket{E_{n+1}\oplus E'_{n}}_{n}$}
\obj(80,0)[]{$\underbracket{E_{n+2}\oplus E'_{n+1}}_{n+1}$}
\obj(160,0)[]{$\underbracket{E_{n+3}\oplus E'_{n+2}}_{n+2}$}
\mor(10,2)(65,2){$\begin{bmatrix}
-d_{n+1} & 0 \\
\phi_{n+1} & d'_{n}
\end{bmatrix}$}[\atright, \solidarrow]
\mor(95,2)(145,2){$\begin{bmatrix}
-d_{n+2} & 0 \\
\phi_{n+2} & d'_{n+1}
\end{bmatrix}$}[\atright, \solidarrow]
\mor(145,4)(95,4){$\begin{bmatrix}
-h_{n+3}^{i} & 0 \\
0 & h_{n+2}^{'i}
\end{bmatrix}$}[\atright, \solidarrow]
\mor(65,4)(10,4){$\begin{bmatrix}
-h_{n+2}^{i} & 0 \\
0 & h_{n+1}^{'i}
\end{bmatrix}$}[\atright, \solidarrow]
\enddc
\end{equation}
\end{lmm}
\begin{proof}
Note that the underlying complex of $B$-modules is the classical cone. It only remains to check that all the morphisms involved in the proof of the fact that this complex of $B$ modules is the cone are compatible with the action of $\varepsilon$. This is a tedious but elementary verification.
\end{proof}

Consider an object $(E,d,\{h^i\}_{i\in \{1,\dots,n\}}) \in \textup{Coh}^s(B,\underline{f})$. Then its underlying $B$-dg module $(E,d)$ is strictly perfect. As the (derived) pullback preserves perfect complexes, $(E,d)\otimes_B K(B,\underline{f})$ lies is $\textup{Perf}^s(B,\underline{f})$. This is the $K(B,\underline{f})$-dg module which, in degree $m$ and $m+1$ has the shape
\begin{equation}
\bigoplus_{k=0}^n E_{m+k}\otimes_B\bigwedge^k \bigl (B\varepsilon_1\oplus \dots \oplus B\varepsilon_n \bigr ) \xrightarrow{\partial_m} \bigoplus_{k=0}^n E_{m+k+1}\otimes_B\bigwedge^k \bigl ( B\varepsilon_1\oplus \dots \oplus B\varepsilon_n \bigr ) 
\end{equation}
Moreover, $\partial_m$ is defined as follows: for any $x\in E_{m+k}$
\begin{equation}
\partial_m(x\otimes \varepsilon_{i_1}\wedge \dots \wedge \varepsilon_{i_k})= (-1)^k d_{m+k}(x)\otimes \varepsilon_{i_1}\wedge \dots \wedge \varepsilon_{i_k}+\sum_{j=1}^n \bigl ((-1)^jf_{i_j}x \otimes \varepsilon_{i_1}\wedge \dots \wedge \widehat{\varepsilon_{i_j}} \wedge \dots \wedge \varepsilon_{i_k} \bigr )
\end{equation}
The $-1$ degree morphisms 
\begin{equation}
\eta^j_{m} :\bigoplus_{k=0}^n E_{m+k}\otimes_B\bigwedge^k \bigl (B\varepsilon_1\oplus \dots \oplus B\varepsilon_n \bigr ) \rightarrow \bigoplus_{k=0}^n E_{m+k-1}\otimes_B\bigwedge^k \bigl (B\varepsilon_1\oplus \dots \oplus B\varepsilon_n \bigr ) 
\end{equation}
are defined, for $x\in E_{m+k}$, by
\begin{equation}
\eta_m^j(x\otimes \varepsilon_{i_1}\wedge \dots \wedge \varepsilon_{i_k})=x\otimes \varepsilon_j\wedge \varepsilon_{i_1}\wedge \dots \wedge \varepsilon_{i_k}
\end{equation}
Notice that we have a morphism of $B$-dg modules 
\begin{equation*}
    \phi : (E,d)\otimes_B K(B,\underline{f})\rightarrow (E,d,\{h^i\}_{i\in \{1,\dots,n\}})
\end{equation*} which is defined in degree $m$ by ($x\in E_{m+k}$)
\begin{equation}
\bigoplus_{k=0}^n E_{m+k}\otimes_B \bigwedge^k \bigl (B\varepsilon_1\oplus \dots \oplus B\varepsilon_n \bigr ) \rightarrow E_m
\end{equation}
\[
x\otimes \varepsilon_{i_1}\wedge \dots \wedge \varepsilon_{i_k} \mapsto h_{m-1}^{i_1}\circ \dots \circ h_{m+k}^{i_k}(x)
\]
where with this notation, when $k=0$, we just mean the identity morphism.
\begin{lmm}
$\phi$ is a cocycle morphism of $K(B,\underline{f})$-dg modules.
\end{lmm}
\begin{proof}
It is clear that $\phi$ is a morphism of $K(B,\underline{f})$-dg modules, i.e. that
\[
\phi \circ \eta^j=h^j\circ \phi
\]
We then only need to show that $\phi$ commutes with the differentials too. Pick $x\in E_{m+k}$. Then
\[
d_{m}(\phi_m(x\otimes \varepsilon_{i_1}\wedge \dots \wedge \varepsilon_{i_k}))=d_{m}(h^{i_1}_{m+1}\circ \dots \circ h^{i_k}_{m+k}(x))\underbracket{=}_{(\ref{commutation d and h})} 
\]
\[
\sum_{j=1}^k(-1)^{j+1}f_{i_j}h^{i_1}_{m+2}\circ \dots \circ \widehat{h^{i_j}} \circ \dots \circ h^{i_k}_{m+k}(x)+(-1)^kh^{i_1}_{m+2}\circ \dots \circ h^{i_k}_{m+k+1}\circ d_{m+k}(x)
\]
On the other hand, we have that
\[
\phi_{m+1}(\partial_m(x\otimes \varepsilon_{i_1}\wedge \dots \varepsilon_{i_k}))=
\]
\[
\phi_{m+1}\Bigl ((-1)^k d_{m+k}(x)\otimes \varepsilon_{i_1}\wedge \dots \wedge \varepsilon_{i_k}+\sum_{j=1}^n \bigl ((-1)^jf_{i_j}x \otimes \varepsilon_{i_1}\wedge \dots \wedge \widehat{\varepsilon_{i_j}} \wedge \dots \wedge \varepsilon_{i_k} \bigr ) \Bigr )
\]
\[
=(-1)^kh^{i_1}_{m+2}\circ \dots \circ h^{i_k}_{m+k+1}\circ d_{m+k}(x)+\sum_{j=1}^k(-1)^{j+1}f_{i_j}h^{i_1}_{m+2}\circ \dots \circ \widehat{h^{i_j}} \circ \dots \circ h^{i_k}_{m+k}(x)
\]
If $k=0$, then $\phi_m(x)=x$ and there is nothing to show.
\end{proof}
\begin{rmk}
As the source of $\phi$ is a perfect $K(B,\underline{f})$-dg module, it follows that $(E,d,\{h^i\}_{i=1,\dots,n})$ and $cone(\phi)$ are equivalent in $\sing{B,\underline{f}}$.
\end{rmk}
\begin{prop}
Assume that $(E,d,\{h^i\}_{i=1,\dots,n})$ as above is concentrated in degrees $[m',m]$, where $m-m'\geq n+1$ (i.e. the dg module is concentrated in at least $n+2$ degrees). Then $cone(\phi)$ is equivalent, in $\sing{B,\underline{f}}$, to a $K(B,\underline{f})$-dg module concentrated in degrees $[m',m-1]$
\end{prop}
\begin{proof}
We claim that we can exhibit $cone(\phi)$ as the cone of a cocycle morphism of $K(B,\underline{f})$-dg modules whose domain is 
\begin{equation}
-\bigl ( ( E_{m'}\xrightarrow{d_{m'}}\dots \xrightarrow{d_{m'+n}}E_{m'+n})\otimes_B K(B,\underline{f})\bigr )
\end{equation}
which is a perfect $K(B,\underline{f})$-dg module concentrated in degrees $[m'-n,m'+n]$. The $-$ above means that we change the sign of all the $\delta_i$'s and $\mu_i^j$'s. Notice that it is a $K(B,\underline{f})$-sub-dg module of $cone(\phi)$ and that $\delta$ and the $\mu^i$'s coincide with the ones induced by this inclusion.

Now consider the $K(B,\underline{f})$-sub-dg module of $cone(\phi)$, which in degree $s$ is the projective $B$-module
\begin{equation}
E_{s} \oplus \Bigl (\bigoplus_{\substack{j=0 \\\ j+s\geq m'+n}}^n E_{j+s+1}\otimes_B \bigwedge^{j}(B\varepsilon_1\oplus \dots \oplus B\varepsilon_n) \Bigr)\subseteq \bigl ( cone(\phi) \bigr )_s
\end{equation}
which we will refer to as $(F,\partial,\{\eta^i\}_{i=1,\dots n})$\footnote{Clearly, $\partial_.$ and $\{\eta_.^i\}_{i=1,\dots,n}$ are induced by $cone(\phi)$}. This is still a $K(B,\underline{f})$-dg module as $\partial_.$ and $\{\eta_.^i\}_{i=1,\dots,n}$ are well defined on it, i.e. $\partial_s(F_s)\subseteq F_{s+1}$ and $\eta_{s}^i(F_s)\subseteq F_{s-1}$.  Notice that, for $s\leq m'-1$, $F_s=0$ and, for $s\geq m'+n$, $F_s=cone(\phi)_s$. This means that $(F,\partial,\{\eta^i\}_{i=1,\dots n})$ is a $K(B,\underline{f})$-dg module concentrated in degrees $[m',m]$, as $cone(\phi)_s=0$ if $s>m$. Label $\iota : (F,\partial,\{\eta^i\}_{i=1,\dots n})\rightarrow cone(\phi)$ the canonical inclusion and $\pi : cone(\phi)\rightarrow (F,\partial,\{\eta^i\}_{i=1,\dots n})$ the canonical projection.

Notice that, for any $s$, we have that
\begin{equation}
F_s\oplus  \bigl (- ( E_{m'}\xrightarrow{d_{m'}}\dots \xrightarrow{d_{m'+n}}E_{m'+n}\bigr )\otimes_B K(B,\underline{f})\bigr )_{s+1}=
\end{equation}
\[
E_s \oplus \Bigl (\bigoplus_{\substack{j=0 \\\ j+s\geq m'+n}}^n E_{j+s+1}\otimes_B \bigwedge^{j}(B\varepsilon_1\oplus \dots \oplus B\varepsilon_n) \Bigr) \oplus \Bigl ( \bigoplus_{\substack{j=0 \\\ j+s+1\leq m'+n}}^{n} E_{j+s+1}\otimes_B \bigwedge^j ( B\varepsilon_1\oplus \dots \oplus B\varepsilon_n) \Bigr )
\]
\[
\simeq E_s \oplus \Bigl (\bigoplus_{j=0}^n E_{j+s+1}\otimes_B\bigwedge^{j}(B\varepsilon_1\oplus \dots \oplus B\varepsilon_n) \Bigr)=cone(\phi)_s
\]
Define 
\begin{equation}
\psi : -\bigl ( ( E_{m'}\xrightarrow{d_{m'}}\dots \xrightarrow{d_{m'+n}}E_{m'+n-1})\otimes_B K(B,\underline{f})\bigr ) \rightarrow (F,\partial, \{\eta^i\}_{i=1,\dots,n})
\end{equation}
in every degree as the composition
\[
\Bigl ( \bigoplus_{\substack{j=0 \\\ j+s\leq m'+n}} E_{s+j}\otimes_B\bigwedge^j ( B\varepsilon_1\oplus \dots \oplus B\varepsilon_n) \Bigr ) \subseteq cone(\phi)_{s-1}\xrightarrow{\partial_{s-1}}cone(\phi)_s\xrightarrow{\pi} F_s
\]
This is a cocycle morphism by construction. Notice that, as 
\begin{equation*}
    \eta_s^i(F_s)\subseteq F_{s-1}
\end{equation*}
and 
\begin{equation*}
\eta_s^i\bigl ( ( E_{m'}\xrightarrow{d_{m'}}\dots \xrightarrow{d_{m'+n}}E_{m'+n})\otimes_B K(B,\underline{f})\bigr )_s\subseteq \bigl ( ( E_{m'}\xrightarrow{d_{m'}}\dots \xrightarrow{d_{m'+n}}E_{m'+n})\otimes_B K(B,\underline{f})\bigr )_{s-1}
\end{equation*}
by Lemma \ref{cone K(B,f) dg mod}  we find that $cone(\psi)=cone(\phi)$.

To conclude, notice that $(F,\partial,\{\eta^i\}_{i=1,\dots,n})$ coincides, in degrees $m-1$ and $m$, with
\[
E_{m-1}\oplus E_m \xrightarrow{\begin{bmatrix} d_{m-1} 1\end{bmatrix}} E_m
\]
Therefore, $(F,\partial,\{\eta^i\}_{i=1,\dots,n})$ is quasi-isomorphic to 
\[
\begindc{\commdiag}[20]
\obj(0,0)[]{$F_{m'}$}
\obj(30,0)[]{$\dots$}
\obj(60,0)[]{$F_{m-2}$}
\obj(140,0)[]{$Ker(\begin{bmatrix} d_{m-1} 1\end{bmatrix})\simeq E_{m-1}$}
\mor(2,-2)(28,-2){$\partial_{m'}$}[\atright,\solidarrow]
\mor(28,2)(2,2){$\eta_{m'+1}^i$}[\atright,\solidarrow]
\mor(32,-2)(58,-2){$\partial_{m-3}$}[\atright,\solidarrow]
\mor(58,2)(32,2){$\eta_{m-2}^i$}[\atright,\solidarrow]
\mor(65,-2)(112,-2){$\tilde{\partial}_{m-2}$}[\atright,\solidarrow]
\mor(112,2)(65,2){$\tilde{\eta}_{m-1}^i$}[\atright,\solidarrow]
\enddc
\]
As $(F,\partial,\{\eta^i\}_{i=1,\dots,n})$ is equivalent to $(E,d\{h^i\}_{i=1,\dots,n})$ in $\sing{B,\underline{f}}$, we have proved the proposition. 
\end{proof}
\begin{rmk}{[B.~Keller]}
The previous proposition holds in a more general setting. More precisely: let $B$ a commutative ring, $n$ an integer and let $R$ be a (possibly non-commutative) dg algebra such that
\begin{itemize}
    \item $R_p=0$ if $p\neq [-n,0]$;
    \item $R_0=B$.
\end{itemize}
Let $m,m'\in \Z$ such that $m-m'\geq n+1$ and let $E$ be a right dg $R$-module such that
\begin{itemize}
    \item $E_p=0$ if $p\geq [m',m]$;
    \item $E_p$ is a projective right $B$-module of finite type, for all $p\in Z$.
\end{itemize}
Let $\textup{\textbf{Mod}}_R$ denote the dg-category of right dg $R$-modules and $\perf{R}$ the sub-dg category spanned by perfect right dg $R$-modules. Then $E$ is equivalent in $\textup{\textbf{Mod}}_R/\perf{R}$ to a right dg $R$-module $F$ verifying the following two conditions:
\begin{itemize}
    \item $F_p=0$ if $p\neq [m',m-1]$;
    \item $F_p$ is a projective right $B$-module of finite type, for all $p\in \Z$.
\end{itemize}
\end{rmk}
\begin{rmk}{[B.~Keller]} It seems that the content of the previous proposition is close to the "fundamental domain" theorem due to C.~Amiot (see \cite{am09}) and generalized by O.~Iyama and D.~Yang (see \cite{iy17}). The precise comparison will be investigated elsewhere.
\end{rmk}
Then the following structure theorem holds:
\begin{thm}\label{structure theorem}
Let $(Spec(B),\underline{f})$ be a $n$-dimensional affine Landau-Ginzburg model  over $S$. Then every object in the dg category of relative singularities $\sing{B,\underline{f}}$ is a retract of an object represented by a $K(B,\underline{f})$-dg module concentrated in $n+1$ degrees.
\end{thm}
\begin{proof}
Let $(E,d,\{h^i\}_{i=1,\dots,n})$ be an object in $\textup{Coh}^s(B,\underline{f})$ concentrated in degrees $[m',m]$. We produce an inductive argument on amplitude $a=m-m'+1$ of the interval where $(E,d,\{h^i\}_{i=1,\dots,n})$ is nonzero. If $m-m'\leq n$ there is nothing to prove. Otherwise, apply the previous proposition. Now, as the homotopy category of $\sing{B,\underline{f}}$ coincides with the idempotent completion of the Verdier quotient of the homotopy category of $\cohbp{Spec(K(B,\underline{f}))}{Spec(B)}$ by the homotopy category of $\perf{Spec(K(B,\underline{f}))}$, we conclude.
\end{proof}

\subsection{\textsc{Matrix factorizations}}
It is well known (see \cite{orl04}, \cite{bw12}, \cite{efpo15}, \cite{brtv}) that the dg category of relative singularities $\sing{B,f}$ associated to a $1$-dimensional affine flat Landau-Ginzburg model over a regular local ring is equivalent to the dg category of matrix factorizations $\mf{B}{f}$ introduced by Eisenbud (\cite{eis80}).  In this section we shall recall what matrix factorizations are.

\begin{context}
In this section we will always work in the context of $1$-dimensional LG models. Therefore, we will omit to say it explicitly.
\end{context}
Let $(Spec(B),f)$ be an affine LG model over $S$.
\begin{defn}
A matrix factorization over $(B,f)$ is the datum of a pair of projective $B$-modules of finite type $E_0$, $E_1$ together with $B$-linear morphisms $E_0\xrightarrow{p_0} E_{1}$ and $E_1\xrightarrow{p_1} E_0$ such that $p_1\circ p_0=f$ and $p_0\circ p_1=f$.
\end{defn}
We can naturally organize matrix factorizations in a $\mathbb{Z}/2\mathbb{Z}$-graded dg category, denoted $\mf{B}{f}$, as follows:
\begin{itemize}
\item the objects of $\mf{B}{f}$ are matrix factorizations over $(B,f)$
\item given two matrix factorizations $(E,p)$ and $(F,q)$ over $(B,f)$, we define the morphisms in degree $0$ (resp. $1$) $\hm^{0}\bigl ( (E,p),(F,q) \bigr )$ (resp. $\hm^{1}\bigl ( (E,p),(F,q) \bigr )$) as the $B$-module of pairs of $B$-linear morphisms $(\phi_0 : E_0 \rightarrow F_0, \phi_1 : E_1 \rightarrow F_1)$ (resp. $(\psi_0 : E_0 \rightarrow F_1, \psi_1 : E_0 \rightarrow F_1)$)
\item given a map $(\chi_0,\chi_1): (E,p)\rightarrow (F,q)$ of degree $i$ ($i=0,1$), we define $\delta ((\chi_0,\chi_1)):= q \circ \chi -(-1)^i \chi \circ p$
\item composition and identities are defined in the obvious way
\end{itemize}

Then we can view $\mf{B}{f}$ as an $A$-linear dg category by means of the structure morphism $A\rightarrow B$.
\begin{rmk}
Notice that since we are considering projective $B$-modules and $B$ is flat over $A$, $\mf{B}{f}$ is a locally flat $A$-linear dg category.
\end{rmk}
The homotopy category of $MF(B,f)$ has a triangulated structure: the suspension is defined as
\begin{equation}
\bigl ( \dbarrow{E_0}{E_1}{p_0}{p_1}\bigr )[1]= \dbarrow{E_1}{E_0}{-p_1}{-p_0}
\end{equation}
and the cone of a closed morphism $(\phi):(E,p)\rightarrow (F,q)$ is defined by
\begin{equation}
\dbarrow{F_0\oplus E_1}{F_1\oplus E_0}{\begin{bmatrix} q_0 & \phi_1\\ 0 & -p_1 \end{bmatrix}}{\begin{bmatrix} q_1 & \phi_0\\ 0 & -p_0 \end{bmatrix}}
\end{equation}
See \cite{orl04} for more details.
Moreover, $\mf{B}{f}$ has a symmetric monoidal structure, defined by
\begin{equation}
(E,p)\otimes (F,q)=\begindc{\commdiag}
\obj(0,1)[]{$(E_0\otimes_B F_0)\oplus (E_1\otimes_B F_1)$}
\obj(50,1)[]{$(E_0\otimes_B F_1)\oplus (E_1\otimes_B F_0)$}
\mor(18,-1)(32,-1){$p\otimes q$}[\atright,\solidarrow]
\mor(32,3)(18,3){$p\otimes q$}[\atright,\solidarrow]
\enddc
\end{equation}
As explained in \cite{brtv}, it is possible to define a lax monoidal $\oo$-functor
\begin{equation}
\mf{\bullet}{\bullet}^{\otimes}: \lgm{1}^{\text{aff,op},\boxplus}\rightarrow \textbf{dgcat}_A^{\text{idem},\otimes}
\end{equation}
It is then possible to extend it to $\lgm{1}^{op,\boxplus}$ by Kan extension. With a little abuse of notation, we still denote this extension by
\begin{equation}\label{mf functor}
    \mf{\bullet}{\bullet}^{\otimes}: \lgm{1}^{\text{op},\boxplus}\rightarrow \textbf{dgcat}_A^{\text{idem},\otimes}
\end{equation}
We refer to \cite{brtv} for more details.
\begin{rmk}
There exists a second definition of matrix factorizations for non-affine LG-models $(X,f)$, see \cite{bw12}, \cite{efi18}, \cite{orl12}. If $X$ is a separated scheme with enough vector bundles, the two definitions agree.
\end{rmk}
\begin{rmk}
Being a lax monoidal $\oo$-functor, $(\ref{mf functor})$ factors through the full subcategory $\textup{Mod}_{\mf{A}{0}}(\dgcatm)^{\otimes}$.
\end{rmk}

\subsection{\textsc{More on the structure of} Sing(B,f)}
As the Koszul algebra $K(B,f)$ is particularly simple, in the case $n=1$ it is possible to give a more detailed description of the objects of $\sing{B,f}$. This is what we will do in the following.
Our first remark concerns the periodicity of the dg category $\sing{B,f}$.
\begin{lmm}\label{lemma del salto della quaglia}
Let $\dbarrow{\underbracket{E}_{n-1}}{\underbracket{F}_{n}}{p}{q}$ be an object in $\textup{Coh}^s(B,f)$. Then it is equivalent to \\  $\dbarrow{\underbracket{F}_{n-2}}{\underbracket{E}_{n-1}}{q}{p}$ in $\sing{B,f}$.
\end{lmm}
\begin{proof}
Consider $\Bigl ( E\xrightarrow{p}F \Bigr )\otimes_B K(B,f) \in \textup{Perf}^s(B,f)$. This is the $K(B,f)$ dg module
\begin{equation}
\begindc{\commdiag}[10]
\obj(0,0)[]{$\underbracket{E}_{n-2}$}
\obj(80,0)[]{$\underbracket{E\oplus F}_{n-1}$}
\obj(160,0)[]{$\underbracket{F}_{n}$}
\mor(3,0)(71,0){$\begin{bmatrix}
-p\\
f
\end{bmatrix}$}[\atright,\solidarrow]
\mor(89,0)(157,0){$\begin{bmatrix}
f & p
\end{bmatrix}$}[\atright,\solidarrow]
\mor(71,4)(3,4){$\begin{bmatrix}
0 & 1
\end{bmatrix}$}[\atright,\solidarrow]
\mor(157,4)(89,4){$\begin{bmatrix}
1 \\
0
\end{bmatrix}$}[\atright,\solidarrow]
\enddc
\end{equation}
Then let $\phi$ be the following morphism of $K(B,f)$ dg modules:
\begin{equation}
\begindc{\commdiag}[30]
\obj(0,10)[1]{$E$}
\obj(30,10)[2]{$F\oplus E$}
\obj(60,10)[3]{$F$}
\obj(30,-10)[4]{$E$}
\obj(60,-10)[5]{$\underbracket{F}_{n}$}
\mor(3,9)(27,9){$\begin{bmatrix}
-p \\
f
\end{bmatrix}$}[\atright, \solidarrow]
\mor(33,9)(57,9){$\begin{bmatrix}
f & p
\end{bmatrix}$}[\atright, \solidarrow]
\mor(57,11)(33,11){$\begin{bmatrix}
1 \\
0
\end{bmatrix}$}[\atright, \solidarrow]
\mor(27,11)(3,11){$\begin{bmatrix}
0 & 1
\end{bmatrix}$}[\atright, \solidarrow]
\mor(33,-11)(57,-11){$p$}[\atright, \solidarrow]
\mor(57,-9)(33,-9){$q$}[\atright, \solidarrow]
\mor{2}{4}{$\begin{bmatrix}
q & 1
\end{bmatrix}$}[\atright, \solidarrow]
\mor{3}{5}{$1$}[\atright, \solidarrow]
\enddc
\end{equation}
This morphism exhibits an equivalence in $\sing{B,f}$ between $\dbarrow{E}{F}{p}{q}$ and $cone(\phi)$, which is
\begin{equation}
\begindc{\commdiag}[30]
\obj(0,0)[1]{$E$}
\obj(30,0)[2]{$F\oplus E$}
\obj(60,0)[3]{$F\oplus E$}
\obj(90,0)[4]{$\underbracket{F}_{n}$}
\mor(3,-1)(27,-1){$\begin{bmatrix}
p\\
-f
\end{bmatrix}$}[\atright, \solidarrow]
\mor(33,-1)(57,-1){$\begin{bmatrix}
-f & -p\\
q & 1
\end{bmatrix}$}[\atright, \solidarrow]
\mor(63,-1)(87,-1){$\begin{bmatrix}
1 & p
\end{bmatrix}$}[\atright, \solidarrow]
\mor(27,1)(3,1){$\begin{bmatrix}
0 & -1
\end{bmatrix}$}[\atright, \solidarrow]
\mor(57,1)(33,1){$\begin{bmatrix}
-1 & 0\\
0 & 0
\end{bmatrix}$}[\atright, \solidarrow]
\mor(87,1)(63,1){$\begin{bmatrix}
0 \\
q
\end{bmatrix}$}[\atright, \solidarrow]
\enddc
\end{equation}
and can be written as the cone of the following morphism of $K(B,f)$ dg modules
\begin{equation}
\begindc{\commdiag}[30]
\obj(0,20)[1]{$E$}
\obj(30,20)[2]{$E$}
\obj(0,0)[3]{$F$}
\obj(30,0)[4]{$F\oplus E$}
\obj(60,0)[5]{$\underbracket{F}_{n}$}
\mor(3,19)(27,19){$f$}[\atright, \solidarrow]
\mor(27,21)(3,21){$1$}[\atright, \solidarrow]
\mor{1}{3}{$p$}[\atright, \solidarrow]
\mor{2}{4}{$\begin{bmatrix}
-p \\
1
\end{bmatrix}$}[\atleft, \solidarrow]
\mor(3,-1)(27,-1){$\begin{bmatrix}
-f \\
q
\end{bmatrix}$}[\atright, \solidarrow]
\mor(33,-1)(57,-1){$\begin{bmatrix}
1 & p
\end{bmatrix}$}[\atright, \solidarrow]
\mor(27,1)(3,1){$\begin{bmatrix}
-1 & 0
\end{bmatrix}$}[\atright, \solidarrow]
\mor(57,1)(33,1){$\begin{bmatrix}
0 \\
q
\end{bmatrix}$}[\atright, \solidarrow]
\enddc
\end{equation}
Notice that the source of this morphism is $E\otimes_B K(B,f)$, where $E$ is a complex concentrated in degree $n-1$. In particular, as $E$ is a projective $B$-module, it is a perfect $K(B,f)$ dg module. Therefore, in $\sing{B,f}$, the target of this morphism is equivalent to $\dbarrow{E}{F}{p}{q}$. Then consider the following morphism of $K(B,f)$ dg modules:
\begin{equation}
\begindc{\commdiag}[30]
\obj(0,10)[1]{$F$}
\obj(30,10)[2]{$F\oplus E$}
\obj(60,10)[3]{$F$}
\obj(0,-10)[4]{$F$}
\obj(30,-10)[5]{$\underbracket{E}_{n-1}$}

\mor(3,9)(27,9){$\begin{bmatrix}
-f \\
q
\end{bmatrix}$}[\atright, \solidarrow]
\mor(33,9)(57,9){$\begin{bmatrix}
1 & p
\end{bmatrix}$}[\atright, \solidarrow]
\mor(27,11)(3,11){$\begin{bmatrix}
-1 & 0
\end{bmatrix}$}[\atright, \solidarrow]
\mor(57,11)(33,11){$\begin{bmatrix}
0 \\
q
\end{bmatrix}$}[\atright, \solidarrow]
\mor(3,-11)(27,-11){$q$}[\atright, \solidarrow]
\mor(27,-9)(3,-9){$p$}[\atright, \solidarrow]
\mor{4}{1}{$1$}[\atright, \solidarrow]
\mor(30,-7)(30,9){$\begin{bmatrix}
-p\\
1
\end{bmatrix}$}[\atright, \solidarrow]
\enddc
\end{equation}
It is not hard to verify that this is a quasi-isomorphism. Following the chain of equivalences in $\sing{B,f}$ we get that 
\[
\dbarrow{\underbracket{E}_{n-1}}{\underbracket{F}_{n}}{p}{q}\simeq \dbarrow{\underbracket{F}_{n-2}}{\underbracket{E}_{n-1}}{q}{p}
\]
\end{proof}

\begin{cor}\label{suspension in Sing}
Let $\dbarrow{\underbracket{E}_{n-1}}{\underbracket{F}_{n}}{p}{q}$ be in $\textup{Coh}^s(B,f)$. Then 
\begin{equation}
\bigl ( \dbarrow{\underbracket{E}_{n-1}}{\underbracket{F}_{n}}{p}{q} \bigr )[1]\simeq \bigl ( \dbarrow{\underbracket{F}_{n-1}}{\underbracket{E}_{n}}{-q}{-p} \bigr )
\end{equation}
in $\sing{B,f}$.
\end{cor}
\begin{proof}
In $\textup{Coh}^s(B,f)$, we know that $\bigl ( \dbarrow{\underbracket{E}_{n-1}}{\underbracket{F}_{n}}{p}{q} \bigr )[1]$ is equivalent to \\$\bigl ( \dbarrow{\underbracket{E}_{n-2}}{\underbracket{F}_{n-1}}{-p}{-q} \bigr )$. Then, by Lemma \ref{lemma del salto della quaglia} we get 
\[
\bigl ( \dbarrow{\underbracket{E}_{n-1}}{\underbracket{F}_{n}}{p}{q} \bigr )[1]\simeq \bigl ( \dbarrow{\underbracket{E}_{n-2}}{\underbracket{F}_{n-1}}{-p}{-q} \bigr ) \simeq \bigl ( \dbarrow{\underbracket{
F}_{n-1}}{\underbracket{E}_{n}}{-q}{-p} \bigr )
\]
\end{proof}
We will now provide an explicit description of the image of an object via the quotient functor 
\[
\textup{Coh}^s(B,f)\rightarrow \sing{B,f}
\]
\begin{thm}\label{structure theorem n=1}
Let
\begin{equation}
\begindc{\commdiag}[20]
\obj(-10,0)[]{$(E,d,h)= 0$}
\obj(30,0)[]{$E_m$}
\obj(60,0)[]{$E_{m+1}$}
\obj(90,0)[]{$\dots$}
\obj(120,0)[]{$E_{m'-1}$}
\obj(150,0)[]{$E_{m'}$}
\obj(180,0)[]{$0$}
\mor(3,0)(27,0){$$}
\mor(33,-1)(57,-1){$d_{m}$}[\atright, \solidarrow]
\mor(57,1)(33,1){$h_{m+1}$}[\atright, \solidarrow]
\mor(63,-1)(87,-1){$d_{m+1}$}[\atright, \solidarrow]
\mor(87,1)(63,1){$$}[\atright, \solidarrow]
\mor(93,-1)(117,-1){$$}[\atright, \solidarrow]
\mor(117,1)(93,1){$h_{-1}$}[\atright, \solidarrow]
\mor(123,-1)(147,-1){$d_{-1}$}[\atright, \solidarrow]
\mor(147,1)(123,1){$h_0$}[\atright, \solidarrow]
\mor(153,0)(177,0){$$}
\enddc
\end{equation}
be an object in $\textup{Coh}^s(B,f)$. Then the following equivalence holds in $\sing{B,f}$:
\begin{equation}
(E,d,h)\simeq 
\begindc{\commdiag}[25]
\obj(0,0)[]{$\underbracket{\bigoplus_{i\in \mathbb{Z}}E_{2i-1}}_{-1}$}
\obj(40,0)[]{$\underbracket{\bigoplus_{i\in \mathbb{Z}}E_{2i}}_{0}$}
\mor(8,-1)(32,-1){$d+h$}[\atright, \solidarrow]
\mor(32,1)(8,1){$d+h$}[\atright, \solidarrow]
\enddc
\end{equation}
Moreover, it is natural in $(E,d,h)$,
\end{thm}
\begin{proof}
The first part of the proof is the same as the one of Theorem \ref{structure theorem}, but we rewrite it in an explicit manner for the reader's convenience. Moreover, we will assume that $m=-2n+1$ for some $n>0$ (if $m=-2n+2$, just put $E_{-2n+1}=0$) and that $m'=0$. It is clear that this does not compromise the generality of the proof.

Consider the perfect $K(B,f)$-dg module $(E,d)\otimes_B K(B,f)$\footnote{recall that $(E,d)$ is a perfect $B$-dg module} and the following morphism 
\begin{equation}
\phi:(E,d)\otimes_B K(B,f)\rightarrow (E,d,h)
\end{equation}
of $K(B,f)$-dg modules:
\begin{equation}
\begindc{\commdiag}[20]
\obj(-10,30)[1]{$\underbracket{E_{-2n+1}}_{-2n}$}
\obj(40,30)[2]{$\underbracket{E_{-2n+2}\oplus E_{-2n+1}}_{-2n+1}$}
\obj(100,30)[3]{$\underbracket{E_{-2n+3}\oplus E_{-2n+2}}_{-2n+2}$}
\obj(140,30)[4]{$\dots$}
\obj(170,30)[5]{$\underbracket{E_{0}\oplus E_{-1}}_{-1}$}
\obj(200,30)[6]{$\underbracket{E_{0}}_{0}$}
\mor(-5,29)(20,29){$\small{\begin{bmatrix}
-d_{-2n+1} \\
f
\end{bmatrix}}$}[\atright, \solidarrow]
\mor(20,31)(-5,31){$\small{\begin{bmatrix}
0 & 1
\end{bmatrix}}$}[\atright, \solidarrow]
\mor(60,29)(80,29){$\small{\begin{bmatrix}
-d_{-2n-2} & 0 \\
f & d_{-2n+1}
\end{bmatrix}}$}[\atright, \solidarrow]
\mor(80,31)(60,31){$\small{\begin{bmatrix}
0 & 1 \\
0 & 0
\end{bmatrix}}$}[\atright, \solidarrow]
\mor(123,29)(140,29){$\small{\begin{bmatrix}
-d_{-2n+3} & 0 \\
f & d_{-2n+2}
\end{bmatrix}}$}[\atright, \solidarrow]
\mor(140,31)(123,31){$\small{\begin{bmatrix}
0 & 1 \\
0 & 0
\end{bmatrix}}$}[\atright, \solidarrow]
\mor(140,29)(160,29){$$}
\mor(160,31)(140,31){$$}
\mor(180,29)(197,29){$\small{\begin{bmatrix}
f & d_{-1}
\end{bmatrix}}$}[\atright, \solidarrow]
\mor(197,31)(180,31){$\small{\begin{bmatrix}
1 \\
0
\end{bmatrix}}$}[\atright, \solidarrow]
\obj(40,-10)[7]{$\underbracket{E_{-2n+1}}_{-2n+1}$}
\obj(100,-10)[8]{$\underbracket{E_{-2n+2}}_{-2n+2}$}
\obj(140,-10)[9]{$\dots$}
\obj(170,-10)[10]{$\underbracket{ E_{-1}}_{-1}$}
\obj(200,-10)[11]{$\underbracket{E_{0}}_{0}$}
\mor{2}{7}{$\small{\begin{bmatrix}
h_{-2n+2} & 1
\end{bmatrix}}$}
\mor{3}{8}{$\small{\begin{bmatrix}
h_{-2n+3} & 1
\end{bmatrix}}$}
\mor{5}{10}{$\small{\begin{bmatrix}
h_{0} & 1
\end{bmatrix}}$}
\mor{6}{11}{$1$}
\mor(45,-11)(95,-11){$d_{-2n+1}$}[\atright, \solidarrow]
\mor(95,-9)(45,-9){$h_{-2n+2}$}[\atright, \solidarrow]
\mor(105,-11)(125,-11){$$}
\mor(125,-9)(105,-9){$$}
\mor(145,-11)(165,-11){$$}
\mor(165,-9)(145,-9){$$}
\mor(175,-11)(195,-11){$d_{-1}$}[\atright, \solidarrow]
\mor(195,-9)(175,-9){$h_{0}$}[\atright, \solidarrow]
\enddc
\end{equation}
Then, $Cone(\phi)$ is equivalent to $(E,d,h)$ in $\sing{B,f}$. $Cone(\phi)$ is the $K(B,f)$-dg module
\begin{equation}
\begindc{\commdiag}[18]
\obj(-25,0)[]{$\small{\underbracket{E_{-2n+1}}_{-2n-1}}$}
\obj(30,0)[]{$\small{\underbracket{E_{-2n+2}\oplus E_{-2n+1}}_{-2n}}$}
\obj(120,0)[]{$\small{\underbracket{E_{-2n+3}\oplus E_{-2n+2}\oplus E_{-2n+1}}_{-2n+1}}$}
\mor(-17,1)(12,1){$\small{\begin{bmatrix}
d_{\footnotesize{-2n+1}} \\
-f
\end{bmatrix}}$}[\atright, \solidarrow]
\mor(12,5)(-17,5){$\small{\begin{bmatrix}
0 & -1
\end{bmatrix}}$}[\atright, \solidarrow]
\mor(48,1)(85,1){$\small{\begin{bmatrix}
d_{\footnotesize{-2n+2}} & 0 \\
-f & -d_{-2n+1} \\
h_{-2n+2} & 1
\end{bmatrix}}$}[\atright, \solidarrow]
\mor(85,5)(48,5){$\small{\begin{bmatrix}
0 & -1 & 0 \\
0 & 0 & 0
\end{bmatrix}}$}[\atright, \solidarrow]
\mor(170,1)(190,1){$\small{\begin{bmatrix}
d_{\footnotesize{-2n+3}} & 0 & 0 \\
-f & -d_{-2n+2} & 0 \\
h_{-2n+3} & 1 & d_{-2n+1}
\end{bmatrix}}$}[\atright, \solidarrow]
\mor(190,5)(170,5){$\small{\begin{bmatrix}
0 & -1 & 0 \\
0 & 0 & 0 \\
0 & 0 & h_{-2n+2}
\end{bmatrix}}$}[\atright, \solidarrow]
\obj(0,-50)[]{$\small{\underbracket{E_{-2n+4}\oplus E_{-2n+3}\oplus E_{-2n+2}}_{-2n+2}}$}
\obj(54,-50)[]{$\dots$}
\obj(90,-50)[]{$\small{\underbracket{E_{0}\oplus E_{-1}\oplus E_{-2}}_{-2}}$}
\obj(150,-50)[]{$\small{\underbracket{E_0\oplus E_{-1}}_{-1}}$}
\obj(180,-50)[]{$\small{\underbracket{E_0}_{0}}$}
\mor(35,-49)(53,-49){$$}[\atright, \solidarrow]
\mor(53,-45)(35,-45){$$}[\atright, \solidarrow]
\mor(52,-49)(70,-49){$$}[\atright, \solidarrow]
\mor(70,-45)(52,-45){$$}[\atright, \solidarrow]
\mor(108,-49)(140,-49){$\small{\begin{bmatrix}
-f & -d_{-1} & 0 \\
h_0 & 1 & d_{-2}
\end{bmatrix}}$}[\atright, \solidarrow]
\mor(140,-45)(108,-45){$\small{\begin{bmatrix}
-1 & 0 \\
0 & 0 \\
0 & h_{-1}
\end{bmatrix}}$}[\atright, \solidarrow]
\mor(160,-49)(178,-49){$\small{\begin{bmatrix}
1 & d_{-1}
\end{bmatrix}}$}[\atright, \solidarrow]
\mor(178,-45)(160,-45){$\small{\begin{bmatrix}
0 \\
h_0
\end{bmatrix}}$}[\atright, \solidarrow]
\enddc
\end{equation}
which can be seen as the cone of the following morphism (call it $\varphi$)
\begin{equation}
\begindc{\commdiag}[18]
\obj(0,0)[1]{$\small{\underbracket{E_{-2n+1}}_{-2n}}$}
\obj(50,0)[2]{$\small{\underbracket{E_{-2n+2}\oplus E_{-2n+1}}_{-2n+1}}$}
\obj(130,0)[3]{$\small{\underbracket{E_{-2n+2}}_{-2n+2}}$}
\mor(8,1)(32,1){$\small{\begin{bmatrix}
-d_{-2n+1} \\
f
\end{bmatrix}}$}[\atright, \solidarrow]
\mor(32,5)(8,5){$\small{\begin{bmatrix}
0 & 1
\end{bmatrix}}$}[\atright, \solidarrow]
\mor(68,1)(122,1){$\small{\begin{bmatrix}
f & d_{-2n+1}
\end{bmatrix}}$}[\atright, \solidarrow]
\mor(122,5)(68,5){$\small{\begin{bmatrix}
1\\
0
\end{bmatrix}}$}[\atright, \solidarrow]
\obj(50,-50)[4]{$\small{\underbracket{E_{-2n+3}\oplus E_{-2n+1}}_{-2n+1}}$}
\obj(130,-50)[5]{$\small{\underbracket{E_{-2n+4}\oplus E_{-2n+3}\oplus E_{-2n+2}}_{-2n+2}}$}
\obj(180,-50)[]{$\small{\dots}$}
\obj(200,-50)[]{$\small{\underbracket{E_0}_{0}}$}
\mor(68,-49)(97,-49){$\small{\begin{bmatrix}
d_{-2n+3} & 0 \\
-f & 0 \\
h_{-2n+3} & d_{-2n+1}
\end{bmatrix}}$}[\atright, \solidarrow]
\mor(97,-45)(68,-45){$\small{\begin{bmatrix}
0 & -1 & 0 \\
0 & 0 & h_{-2n+2}
\end{bmatrix}}$}[\atright, \solidarrow]
\mor(163,-49)(178,-49){$$}[\atright, \solidarrow]
\mor(178,-45)(163,-45){$$}[\atright, \solidarrow]
\mor(182,-49)(198,-49){$$}[\atright, \solidarrow]
\mor(198,-45)(182,-45){$$}[\atright, \solidarrow]
\mor{2}{4}{$\small{\begin{bmatrix}
d_{-2n+2} & 0 \\
h_{-2n+1} & 1
\end{bmatrix}}$}[\atright, \solidarrow]
\mor{3}{5}{$\small{\begin{bmatrix}
0 \\
-d_{-2n+2} \\
1
\end{bmatrix}}$}
\enddc
\end{equation}
As the source of this morphism is $\bigl ( \underbracket{E_{-2n+1}}_{-2n+1}\xrightarrow{d_{-2n+1}} \underbracket{E_{-2n+2}}_{-2n+2} \bigr )\otimes_B K(B,f)$, it is a perfect $K(B,f)$-dg module. Therefore, in $\sing{B,f}$ we have that
\[
(E,d,h)\simeq cone(\phi) = cone(\varphi) \simeq target(\varphi)
\]
The cohomology groups in degree $-1$ and $0$ of $target(\varphi)$ vanish. Therefore, we have found that in $\sing{B,f}$ $(E,d,h)$ is equivalent to 
\begin{equation}\label{[-2n+1,-2]}
\begindc{\commdiag}[18]
\obj(50,0)[4]{$\small{\underbracket{E_{-2n+3}\oplus E_{-2n+1}}_{-2n+1}}$}
\obj(130,0)[5]{$\small{\underbracket{E_{-2n+4}\oplus E_{-2n+3}\oplus E_{-2n+2}}_{-2n+2}}$}
\obj(180,0)[]{$\small{\dots}$}
\obj(230,0)[]{$\underbracket{K}_{-2}$}
\mor(68,1)(97,1){$\small{\begin{bmatrix}
d_{-2n+3} & 0 \\
-f & 0 \\
h_{-2n+3} & d_{-2n+1}
\end{bmatrix}}$}[\atright, \solidarrow]
\mor(97,5)(68,5){$\small{\begin{bmatrix}
0 & -1 & 0 \\
0 & 0 & h_{-2n+2}
\end{bmatrix}}$}[\atright, \solidarrow]
\mor(163,1)(178,1){$$}[\atright, \solidarrow]
\mor(178,5)(163,5){$$}[\atright, \solidarrow]
\mor(182,1)(220,1){$\small{\begin{bmatrix}
d_{-1} & 0 & 0 \\
-f & -d_{-2} & 0 \\
h_{-1} & 1 & d_{-3}
\end{bmatrix}}$}[\atright, \solidarrow]
\mor(220,5)(182,5){$\small{\begin{bmatrix}
0 & -1 & 0 \\
0 & 0 & 0 \\
0 & 0 & h_{-2}
\end{bmatrix}}$}[\atright, \solidarrow]
\enddc
\end{equation}
where
\[
K=Ker\Bigl (\begin{bmatrix}
-f & -d_{-1} & 0 \\
h_0 & 1 & d_{-2}
\end{bmatrix}\Bigr )
\]

This is still an element in $\textup{Coh}^s(B,f)$. Indeed, from the short exact sequence of $B$-modules
\[
0 \rightarrow Ker(\begin{bmatrix} 1 & d_{-1} \end{bmatrix}) \rightarrow E_0\oplus E_{-1} \xrightarrow{\begin{bmatrix} 1 & d_{-1} \end{bmatrix}} E_0 \rightarrow 0
\]
since $E_0$ and $E_{-1}$ are $B$-projective, we conclude that $Ker(\begin{bmatrix} 1 & d_{-1} \end{bmatrix})$ is $B$-projective too. As the complex $cone(\varphi)$ is exact in degree $-1$, we also have the following short exact sequence of $B$-modules:
\[
0 \rightarrow K \rightarrow E_{0} \oplus E_{-1} \oplus E_{-2} \xrightarrow{\begin{bmatrix} -f & -d_{-1} & 0 \\ h_0 & 1 & d_{-2} \end{bmatrix}}  \underbracket{Im\Bigl ( \begin{bmatrix} -f & -d_{-1} & 0 \\ h_0 & 1 & d_{-2} \end{bmatrix} \Bigr )}_{= Ker \Bigl ( \begin{bmatrix} f & d_{-1} \end{bmatrix} \Bigr )} \rightarrow 0
\]
As $E_0$, $E_{-1}$, $E_{-2}$ and $Ker \Bigl ( \begin{bmatrix} f & d_{-1} \end{bmatrix} \Bigr )$ are projective $B$-modules, we conclude. 

Notice that we have found, in $\sing{B,f}$, and equivalence between $(E,d,h)$ (which is concentrated in degrees $[-2n+1,0]$) and an object represented by a complex concentrated in degrees $[-2n+1,-2]$. Therefore, by induction, we have proved that the image of $(E,d,h)$ is equivalent, in $\sing{B,f}$,to a $K(B,f)$-dg module concentrated in degrees $0$ and $-1$ (i.e. by a matrix factorization). Nevertheless, we can do better than this. Indeed, notice  that the $K(B,f)$-dg module (\ref{[-2n+1,-2]}) can be written as the cone of the following morphism of $K(B,f)$-dg modules:
\begin{equation}
\begindc{\commdiag}[18]
\obj(50,50)[1]{$\small{\underbracket{E_{-2n+3}}_{-2n+2}}$}
\obj(130,50)[2]{$\small{\underbracket{E_{-2n+3}}_{-2n+3}}$}
\mor(55,51)(125,51){$\small{f}$}[\atright, \solidarrow]
\mor(125,55)(55,55){$\small{1}$}[\atright, \solidarrow]
\obj(-5,0)[3]{$\small{\underbracket{E_{-2n+1}}_{-2n+1}}$}
\obj(50,0)[4]{$\small{\underbracket{E_{-2n+4}\oplus E_{-2n+2}}_{-2n+2}}$}
\obj(130,0)[5]{$\small{\underbracket{E_{-2n+5}\oplus E_{-2n+4}\oplus E_{-2n+3}}_{-2n+3}}$}
\obj(177,0)[]{$\small{\dots}$}
\obj(200,0)[]{$\underbracket{K}_{-2}$}
\mor(5,1)(30,1){$\small{\begin{bmatrix}
0 \\
d_{-2n+1}
\end{bmatrix}}$}[\atright, \solidarrow]
\mor(30,5)(5,5){$\small{\begin{bmatrix}
0 & h_{-2n+2}
\end{bmatrix}}$}[\atright, \solidarrow]
\mor(69,1)(95,1){$\small{\begin{bmatrix}
d_{-2n+4} & 0 \\
-f & 0 \\
h_{-2n+4} & d_{-2n+2}
\end{bmatrix}}$}[\atright, \solidarrow]
\mor(95,5)(69,5){$\small{\begin{bmatrix}
0 & -1 & 0 \\
0 & 0 & h_{-2n+3}
\end{bmatrix}}$}[\atright, \solidarrow]
\mor(162,1)(178,1){$$}[\atright, \solidarrow]
\mor(178,5)(162,5){$$}[\atright, \solidarrow]
\mor(179,1)(195,1){$$}[\atright, \solidarrow]
\mor(195,5)(179,5){$$}[\atright, \solidarrow]
\mor{1}{4}{$\small{\begin{bmatrix}
d_{-2n+3} \\
h_{-2n+3}
\end{bmatrix}}$}[\atright, \solidarrow]
\mor{2}{5}{$\small{\begin{bmatrix}
0 \\
-d_{-2n+3} \\
1
\end{bmatrix}}$}
\enddc
\end{equation}
As the source of this morphism is $\underbracket{E_{-2n+3}}_{-2n+3}\otimes_B K(B,f)$, and $E_{-2n+3}$ is a perfect $B$-module, this morphism provides an equivalence between $(E,d,h)$ and the target in $\sing{B,f}$. Moreover, we can iterate this procedure: the target of this morphism can be written as the cone of the following morphism:
\begin{equation}
\begindc{\commdiag}[18]
\obj(50,50)[1]{$\small{\underbracket{E_{-2n+4}}_{-2n+3}}$}
\obj(130,50)[2]{$\small{\underbracket{E_{-2n+4}}_{-2n+4}}$}
\mor(55,51)(125,51){$\small{f}$}[\atright, \solidarrow]
\mor(125,55)(55,55){$\small{1}$}[\atright, \solidarrow]
\obj(-35,0)[3]{$\small{\underbracket{E_{-2n+1}}_{-2n+1}}$}
\mor(-28,1)(-7,1){$\small{d_{-2n+1}}$}[\atright, \solidarrow]
\mor(-7,5)(-28,5){$\small{h_{-2n+2}}$}[\atright, \solidarrow]
\obj(0,0)[3]{$\small{\underbracket{E_{-2n+2}}_{-2n+2}}$}
\obj(50,0)[4]{$\small{\underbracket{E_{-2n+5}\oplus E_{-2n+3}}_{-2n+3}}$}
\obj(130,0)[5]{$\small{\underbracket{E_{-2n+6}\oplus E_{-2n+5}\oplus E_{-2n+4}}_{-2n+4}}$}
\obj(177,00)[]{$\small{\dots}$}
\obj(195,0)[]{$\underbracket{K}_{-2}$}
\mor(6,1)(32,1){$\small{\begin{bmatrix}
0 \\
d_{-2n+2}
\end{bmatrix}}$}[\atright, \solidarrow]
\mor(32,5)(6,5){$\small{\begin{bmatrix}
0 & h_{-2n+3}
\end{bmatrix}}$}[\atright, \solidarrow]
\mor(69,1)(95,1){$\small{\begin{bmatrix}
d_{-2n+5} & 0 \\
-f & 0 \\
h_{-2n+5} & d_{-2n+3}
\end{bmatrix}}$}[\atright, \solidarrow]
\mor(95,5)(69,5){$\small{\begin{bmatrix}
0 & -1 & 0 \\
0 & 0 & h_{-2n+4}
\end{bmatrix}}$}[\atright, \solidarrow]
\mor(162,1)(178,1){$$}[\atright, \solidarrow]
\mor(178,5)(162,5){$$}[\atright, \solidarrow]
\mor(179,1)(195,1){$$}[\atright, \solidarrow]
\mor(195,5)(179,5){$$}[\atright, \solidarrow]
\mor{1}{4}{$\small{\begin{bmatrix}
d_{-2n+4} \\
h_{-2n+4}
\end{bmatrix}}$}[\atright, \solidarrow]
\mor{2}{5}{$\small{\begin{bmatrix}
0 \\
-d_{-2n+4} \\
1
\end{bmatrix}}$}
\enddc
\end{equation}
Once again, as the source of this morphism of $K(B,f)$-dg modules is perfect, we obtain an equivalence between $(E,d,h)$ and the target of the morphism in $\sing{B,f}$. Proceeding this way, we obtain a chain of equivalences between our initial $K(B,f)$-dg module and the following:
\begin{equation}\label{ci siamo quasi}
\begindc{\commdiag}[20]
\obj(25,0)[]{$\small{\underbracket{E_{-2n+1}}_{-2n+1}}$}
\obj(60,0)[]{$\small{\underbracket{E_{-2n+2}}_{-2n+2}}$}
\obj(90,0)[]{$\small{\dots}$}
\obj(120,0)[]{$\small{\underbracket{E_{-4}}_{-4}}$}
\obj(150,0)[]{$\small{\underbracket{E_{-1}\oplus E_{-3}}_{-3}}$}
\obj(195,0)[]{$\underbracket{K}_{-2}$}
\mor(32,1)(52,1){$\small{d_{-2n+1}}$}[\atright, \solidarrow]
\mor(52,5)(32,5){$\small{h_{-2n+2}}$}[\atright, \solidarrow]
\mor(68,1)(88,1){$\small{d_{-2n+2}}$}[\atright, \solidarrow]
\mor(88,5)(68,5){$\small{
h_{-2n+3}}$}[\atright, \solidarrow]
\mor(95,1)(115,1){$\small{d_{-5}}$}[\atright, \solidarrow]
\mor(115,5)(95,5){$\small{h_{-4}}$}[\atright, \solidarrow]
\mor(122,1)(142,1){$\small{\begin{bmatrix}
0 \\
d_{-4}
\end{bmatrix}}$}[\atright, \solidarrow]
\mor(142,5)(122,5){$\small{\begin{bmatrix}
0 & h_{-3}
\end{bmatrix}}$}[\atright, \solidarrow]
\mor(158,1)(185,1){$\small{\begin{bmatrix}
d_{-1} & 0 \\
-f & 0 \\
h_{-1} & d_{-3}
\end{bmatrix}}$}[\atright, \solidarrow]
\mor(185,5)(158,5){$\small{\begin{bmatrix}
0 & -1 & 0 \\
0 & 0 & h_{-2}
\end{bmatrix}}$}[\atright, \solidarrow]
\enddc
\end{equation}

By an induction argument, this $K(B,f)$-dg module is equivalent, in $\sing{B,f}$, to 
\begin{equation}
    \begindc{\commdiag}[20]
    \obj(0,0)[1]{$\small{\bigoplus_{i\in \mathbb{Z}}E_{2i-1}}$}
    \obj(150,0)[2]{$\small{E_{-2n+2}\oplus E_{-2n+4} \oplus \dots \oplus E_{-4}\oplus K}$}
    \mor(10,-1)(110,-1){$\small{\begin{bmatrix}
d_{-2n+1} & h_{-2n+3} & \dots & 0 & 0 \\
0 & d_{-2n+3} & \dots & 0 & 0 \\
 & & \vdots & & \\
 0 & 0 & \dots & h_{-3} & 0 \\
 0 & 0 & \dots & 0 & d_{-1} \\
 0 & 0 & \dots & 0 & -f \\
 0 & 0 & \dots & d_{-3} & h_{-1}
\end{bmatrix}}$}[\atright, \solidarrow]
\mor(110,1)(10,1){$\small{\begin{bmatrix}
h_{-2n+2} & 0 & \dots & 0 & 0 & 0 \\
d_{-2n+2} & h_{-2n+4} & \dots & 0 & 0 & 0 \\
 & & \vdots & &  &\\
 0 & 0 & \dots & 0 & 0 & h_{-2} \\
 0 & 0 & \dots & 0 & -1 & 0
\end{bmatrix}}$}[\atright, \solidarrow]
    \enddc
\end{equation}
We can finally consider the following morphism of $K(B,f)$-dg modules concentrated in degrees $-1$ and $0$
\begin{equation}\label{final step}
\begindc{\commdiag}[20]
\obj(-3,63)[1]{$\small{\bigoplus_{i\in \mathbb{Z}}E_{2i-1}}$}
\obj(-3,0)[2]{$\small{\bigoplus_{i\in \mathbb{Z}}E_{2i}}$}
\obj(120,63)[3]{$\small{\bigoplus_{i\in \mathbb{Z}}E_{2i-1}}$}
\obj(120,0)[4]{$\small{E_{-2n+2}\oplus E_{-2n+4} \oplus \dots \oplus E_{-4}\oplus K}$}
\mor(-2,55)(-2,5){$\small{d+h}$}[\atright, \solidarrow]
\mor(2,5)(2,55){$\small{d+h}$}[\atright, \solidarrow]
\mor(118,55)(118,5){$\small{\begin{bmatrix}
d_{-2n+1} & h_{-2n+3} & \dots & 0 & 0 \\
0 & d_{-2n+3} & \dots & 0 & 0 \\
 & & \vdots & & \\
 0 & 0 & \dots & h_{-3} & 0 \\
 0 & 0 & \dots & 0 & d_{-1} \\
 0 & 0 & \dots & 0 & -f \\
 0 & 0 & \dots & d_{-3} & h_{-1}
\end{bmatrix}}$}[\atright, \solidarrow]
\mor(122,5)(122,55){$\small{\begin{bmatrix}
h_{-2n+2} & 0 & \dots & 0 & 0 & 0 \\
d_{-2n+2} & h_{-2n+4} & \dots & 0 & 0 & 0 \\
 & & \vdots & &  &\\
 0 & 0 & \dots & 0 & 0 & h_{-2} \\
 0 & 0 & \dots & 0 & -1 & 0
\end{bmatrix}}$}[\atright, \solidarrow]
\mor{2}{4}{$\small{\begin{bmatrix}
1 & 0 & \dots & 0 & 0 & 0 \\
0 & 1 & \dots & 0 & 0 & 0 \\
 & & \vdots &  &  & \\
 0 & 0 & \dots & 1 & 0 & 0 \\
 0 & 0 & \dots & 0 & 0 & 1 \\
 0 & 0 & \dots & 0 & -d_{-2} & -h_0 \\
  0 & 0 & \dots & 0 & 1 & 0 \\
\end{bmatrix}}$}[\atright,\solidarrow]
\mor{1}{3}{$id$}[\atleft, \solidarrow]
\enddc
\end{equation}
It is not hard to check that morphism (\ref{final step}) is a quasi-isomorphism. Notice that the target of (\ref{final step}) is equivalent in $\sing{B,f}$ to the $K(B,f)$-dg module $(E,d,h)$.

Also notice that since all the passages above are functorial, the equivalence is natural in $(E,d,h)$. In particular, a morphism of $K(B,f)$-dg modules $\phi : (E,d,h)\rightarrow (E',d',h')$ corresponds, under this equivalence, to
\[
\begindc{\commdiag}
\obj(0,30)[1]{$\small{\bigoplus_{i\in \mathbb{Z}}E_{2i-1}}$}
\obj(40,30)[2]{$\small{\bigoplus_{i\in \mathbb{Z}}E_{2i}}$}
\obj(0,0)[3]{$\small{\bigoplus_{i\in \mathbb{Z}}E'_{2i-1}}$}
\obj(40,0)[4]{$\small{\bigoplus_{i\in \mathbb{Z}}E'_{2i}}$}
\mor(10,30)(30,30){$d+h$}[\atright,\solidarrow]
\mor(30,35)(10,35){$d+h$}[\atright,\solidarrow]
\mor(10,0)(30,0){$d'+h'$}[\atright,\solidarrow]
\mor(30,5)(10,5){$d'+h'$}[\atright,\solidarrow]
\mor{1}{3}{$\oplus \phi_{2i-1}$}
\mor{2}{4}{$\oplus \phi_{2i}$}
\enddc
\]
\end{proof}
\begin{rmk}
The algorithm we have provided actually puts the final $K(B,f)$-dg module
\[
\begindc{\commdiag}[20]
\obj(0,100)[1]{$\small{E_{-2n+1}\oplus E_{-2n+3}\oplus \dots \oplus E_{-3}\oplus E_{-1}}$}
\obj(120,100)[2]{$\small{E_{-2n+2}\oplus E_{-2n+4} \oplus \dots \oplus E_{-2}\oplus E_{0}}$}
\mor(50,100)(70,100){$d+h$}[\atright, \solidarrow]
\mor(70,105)(50,105){$d+h$}[\atright, \solidarrow]
\enddc
\]
in degrees $-2n+1$ and $-2n+2$. However, thanks to Lemma \ref{lemma del salto della quaglia}, this is equivalent in $\sing{B,f}$ to the same dg module concentrated in degrees $-1$ and $0$
\end{rmk}
\begin{cor}\label{cones Sing(B,f)}
Let $\phi_. : \dbarrow{E_{-1}}{E_0}{d}{h}\rightarrow \dbarrow{E'_{-1}}{E'_0}{d'}{h'}$ be a closed morphism in $\textup{Coh}^s(B,f)$. Then $cone(\phi_.)$ is equivalent to $\dbarrow{\small{E'_{-1}\oplus E_0}}{\small{E'_0\oplus E_{-1}}}{\begin{bmatrix} d' & \phi_0 \\ 0 & -h\end{bmatrix}}{\begin{bmatrix} h' & \phi_{-1}\\ 0 & -d\end{bmatrix}}$ in $\sing{B,f}$.
\end{cor}
\begin{proof}
This is a straightforward consequence of the computation of $Cone(\phi)$ in $\textup{Coh}^s(B,f)$ and of the previous theorem. 
\end{proof}
\begin{cor}\label{Orlov equivalence}
The lax monoidal $\oo$-natural transformation
\begin{equation}
Orl^{-1,\otimes} : \sing{\bullet,\bullet} \rightarrow \mf{\bullet}{\bullet} : \lgm{1}^{op,\boxplus}\rightarrow \dgcatmo
\end{equation}
constructed in \cite[\S2.4]{brtv} defines a lax-monoidal $\oo$-natural equivalence.
\end{cor}
\begin{proof}
By Kan extension and descent, it is sufficient to consider the affine case. Let $(Spec(B),f)\in \lgm{1}^{\textup{aff},op}$.
As the dg categories $\sing{B,f}$ and $\mf{B,f}$ are triangulated, it is sufficient to show that the induced functor
\[
[Orl^{-1}] : [\sing{B,f}] \rightarrow [\mf{B}{f}]
\]
\[
(E,d,h)\mapsto \dbarrow{\bigoplus_{i\in \mathbb{Z}}E_{2i-1}}{\bigoplus_{i\in \mathbb{Z}}E_{2i}}{d+h}{d+h}
\]
is an equivalence. Consider
\[
Orl : [\mf{B}{f}] \rightarrow [\sing{B,f}]
\] 
\[
\dbarrow{E}{F}{p}{q} \mapsto \dbarrow{\underbracket{E}_{-1}}{\underbracket{F}_{0}}{p}{q}
\]
This is an exact functor between triangulated categories by Corollary \ref{suspension in Sing} and by Corollary \ref{cones Sing(B,f)}.
It is clear that $[Orl^{-1}]\circ Orl$ is the identity functor. By Theorem \ref{structure theorem n=1}, $Orl\circ [Orl^{-1}]$ is equivalent to the identity functor too.
\end{proof}
\begin{rmk}
Notice that $Orl$ is a derived version of the "Cok" functor introduced in \cite{orl04}. Indeed, when $f$ is flat, the $K(B,f)$-dg module $coker(p)$ concentrated in degree $0$ is quasi-isomorphic to $\dbarrow{\underbracket{E}_{-1}}{\underbracket{F}_{0}}{p}{q}$.
\end{rmk}
\begin{rmk}
In \cite{efpo15}, the authors also introduced a coherent version of $\mf{B}{f}$. When $f$ is flat, they proved it to be equivalent to another category of singularities, defined as the Verdier quotient
\begin{equation}
\sing{B,f}_{coh}=\cohb{B/f}/\textup{\textbf{E}}
\end{equation}
where $\textup{\textbf{E}}$ is the thick subcategory of $\cohb{B/f}$ generated by the image of the pullback $\iota^*: \cohb{B}\rightarrow \cohb{B/f}$. 

Our proof of Theorem \ref{structure theorem n=1} also tells us that, for any $f$, all objects in this triangulated category can be represented by $K(B,f)$-dg modules concentrated in degrees $[-1,0]$. This can be used to show that the equivalence proven in \cite{efpo15} holds for any potential $f$, provided that we consider the derived fiber instead of $B/f$.
\end{rmk}

\section{\textsc{A new dg category of relative singularities}}
One of the most important features of $\sing{X,f}\simeq \mf{X}{f}$ is the fact that it is a $2$-periodic dg category. This is not true anymore in the higher dimensional case. Indeed, in the monopotential case, $2$-periodicity in $\sing{X,f}$ ($\mf{X}{f}$ is $2$-periodic by definition) is given by the follwoing exact triangles ( for $\Ec \in \cohbp{X_0}{X}$)
\begin{equation}
    \mathfrak{i}^*\mathfrak{i}_*\Ec\rightarrow \Ec \rightarrow \Ec[2]
\end{equation}
but in the pluripotential case the cofiber of $\mathfrak{i}_0^*\mathfrak{i}_{0*}\Ec \rightarrow \Ec$ can not be identified with a shift of $\Ec$. However, there is a canonical way to define a $2$-periodic dg category out of $\sing{X,\underline{f}}$. 
\subsection{\textsc{Forcing \texorpdfstring{$2$}{2}-periodicity}}
Let $(X,\underline{f})$ be a $n$-dimensional LG model. For each $k=1,\dots,n$, consider the $(n-1)$-dimensional LG model $(X,\underline{f}_k)$, where $\underline{f}_k=(f_1,\dots,f_{k-1},f_{k+1},\dots,n)$. Label $X_k$ the derived zero locus of $\underline{f}_k$:
\begin{equation}
    \begindc{\commdiag}[18]
    \obj(-30,10)[1]{$X_k:=X\times^h_{\aff{n-1}{S}}S$}
    \obj(30,10)[2]{$X$}
    \obj(-30,-10)[3]{$S$}
    \obj(30,-10)[4]{$\aff{n-1}{S}$}
    \mor{1}{2}{$$}
    \mor{1}{3}{$$}
    \mor{2}{4}{$\underline{f}_k$}
    \mor{3}{4}{$0$}
    \enddc
\end{equation}
We can also write $X_0=X\times^h_{\aff{n}{S}}S$ as the derived zero locus of $f_k:X_k\rightarrow \aff{1}{S}$. In particular, we get $n$ one-codimensionsional lci closed embeddings
\begin{equation}
    \mathfrak{i}_{0k}:X_0\rightarrow X_k
\end{equation}
\begin{lmm}
Each embedding $\mathfrak{i}_{0k}:X_0\rightarrow X_k$ induces a pair of adjoint dg functors
\begin{equation}
    \mathfrak{i}_{0k}^*:\sing{X,\underline{f}_k}\rightarrow \sing{X,\underline{f}} \hspace{1cm} \mathfrak{i}_{0k*}:\sing{X,\underline{f}}\rightarrow \sing{X,\underline{f}_k}
\end{equation}
\end{lmm}
\begin{proof}
Recall that we can write 
\begin{equation*}
    \sing{X,\underline{f}}\simeq \cohbp{X_0}{X}/\perf{X_0}
\end{equation*}
\begin{equation*}
    \sing{X,\underline{f}_k}\simeq \cohbp{X_k}{X}/\perf{X_k}
\end{equation*} 
The statement follows easily from the derived base change theorem applied to the homotopy pullback
\begin{equation*}
    \begindc{\commdiag}[18]
    \obj(-30,15)[1]{$X_0$}
    \obj(30,15)[2]{$X_k$}
    \obj(-30,-15)[3]{$Y_k$}
    \obj(30,-15)[4]{$X$}
    \mor{1}{2}{$\mathfrak{i}_{0k}$}
    \mor{1}{3}{$\mathfrak{j}_{0k}$}
    \mor{2}{4}{$\mathfrak{i}_k$}
    \mor{3}{4}{$\mathfrak{j}_k$}
    \enddc
\end{equation*}
where $Y_k$ is the derived zero locus of $f_k$ and from the fact that all these morphisms are lci.
\end{proof}
Then, for each $\Ec \in \sing{X,\underline{f}}$, we have an exact triangle
\begin{equation}
    \mathfrak{i}_{0k}^*\mathfrak{i}_{0k*}\Ec \rightarrow \Ec \xrightarrow{\chi_k} \Ec[2]
\end{equation}
\begin{rmk}
It follows from the functoriality of the counit of the adjunction $(\mathfrak{i}_{0k}^*,\mathfrak{i}_{0k*})$ and from that of the cone construction (since we are working with dg categories), that the $\chi_k$'s define natural transformation of $\oo$-functors
\begin{equation}
    \chi_k:id\rightarrow id[2]:\sing{X,\underline{f}}\rightarrow \sing{X,\underline{f}}
\end{equation}
\end{rmk}
These natural transformation $\chi_k$ should correspond to Eisenbud's operators (\cite{eis80}, \cite{bw15}). Localizing $\sing{X,\underline{f}}$ with respect to the class of morphisms $\{ \chi_k:\Ec \rightarrow \Ec[2] \}$ we obtain a $2$-periodic dg category $\sing{X,\underline{f}}[\chi_k^{-1}]$. Also notice that each $\chi_k$ imposes $2$-periodicity. Therefore, we have $n-1$ non trivial natural equivalences 
\begin{equation}
    \chi_k^{-1}\circ \chi_l : id\rightarrow id : \sing{X,\underline{f}}[\chi_k^{-1}]\rightarrow \sing{X,\underline{f}}[\chi_k^{-1}]
\end{equation}
\subsection{\textsc{Connection to other dg categories}}
It is a theorem due to D.~Orlov and J.~Burke-M.~Walker that, if $\underline{f}$ is a regular sequence in $X=Spec(B)$, there is an equivalence of dg categories
\begin{equation}
    \sing{B,\underline{f}}\simeq \sing{\proj{n-1}{B},W_{\underline{f}}}
\end{equation}
where $W_{\underline{f}}=f_1\cdot T_1+\dots +f_n\cdot T_n \in \Oc(1)(\proj{n-1}{B})$.
\begin{rmk}
We will discuss this equivalence more in detail at the end of the next chapter.
\end{rmk}
Moreover, J.~Burke and M.~Walker prove (see \cite[\S3, \S4]{bw15}) that the Eisenbud operators on $\sing{B,\underline{f}}$ correspond, under the above mentioned equivalence, to multiplication by the parameters $T_i \in \Oc(1)(\proj{n-1}{B})$ on $\sing{\proj{n-1}{B},W_{\underline{f}}}$. Therefore, the dg category $\sing{X,\underline{f}}[\chi_k^{-1}]$ should correspond to $\sing{\proj{n-1}{B},W_{\underline{f}}}[T_k^{-1}]$. We expect the latter dg category to be equivalent to $\sing{\gm{B}^{n-1},W_{\underline{f}|}}\simeq \mf{\gm{B}^{n-1}}{W_{\underline{f}|}}$. Moreover, fix the coordinates on $\gm{B}^{n-1}$: $t_i:=T_i/T_n$, for $i=1,\dots,n-1$. Then
\begin{equation}
    W_{\underline{f}|}:=f_1\cdot t_1+\dots +f_{n-1}\cdot t_{n-1}+f_n
\end{equation}
Notice that we have an equivalence in $\lgm{n}^{\boxplus}$
\begin{equation}
    (\gm{B}^{n-1},f_1\cdot t_1+\dots +f_{n-1}\cdot t_{n-1}+f_n)=(\gm{B},f_1\cdot t_1)\boxplus \dots \boxplus (\gm{B},f_{n-1}\cdot t_{n-1})\boxplus (\gm{B},f_n)
\end{equation}
In the special case we are working over a field $k$ of characteristic $0$, the Thom-Sebastiani theorem for matrix factorizations (\cite[Theorem 3.2.1.3]{pr11}) gives us the following equivalence
\begin{equation}
    \mf{\gm{B}^{n-1}}{W_{\underline{f}|}}_{V(f_1\cdot t_1)\times \dots \times V(f_{n-1}\cdot t_{n-1})\times V(f_n)}\simeq 
\end{equation}
\begin{equation*}
\mf{\gm{B}}{f_1\cdot t_1} \otimes_{k\drl \beta \drr}\dots \otimes_{k\drl \beta \drr}\mf{\gm{B}}{f_{n-1}\cdot t_{n-1}} \otimes_{k\drl \beta \drr} \mf{B}{f_n}
\end{equation*}

%% file: chapters/chapter03.tex
In this chapter we will study the motivic and $\ell$-adic realization of dg categories of singularities of a twisted LG model whose underlying scheme is regular. We will then use the formula obtained in this way together with Theorem \ref{reduction of codimension theorem} to give a formula for the $\ell$-adic realization of the dg category of singularities of the special fiber of a scheme over a regular local ring of dimension $n$.

\section{\textsc{What is in this chapter}}
We start this chapter with a quick review of the main theorem of \cite{brtv}, that served both as a model and as a motivation for the investigations presented later. The main purpose of the above-mentioned paper is to identify a classical object of singularity theory, namely the $\ell$-adic sheaf of (inertia invariant) vanishing cycles, with the $\ell$-adic cohomology of a non-commutative space, the dg category of singularities of the special fiber. The connection between these two objects has been known for a while, thanks to works of T.~Dyckerhoff (\cite{dy11}), A.~Prygel (\cite{pr11}), A.~Efimov (\cite{efi18}) and many others. The main result of A.~Blanc, M.~Robalo, B.~Toën and G.~Vezzosi's paper (\cite[Theorem 4.39]{brtv}) reads as follows
\begin{theorem-non}\label{main theorem brtv}
Let $p:X\rightarrow S$ be a proper, flat, regular scheme over an excellent strictly henselian trait $S$. There is an equivalence of $i_{\sigma}^*\Rlv_S(\sing{S,0})\simeq \Ql{,\sigma}(\beta)\otimes_{\Ql{,\sigma}} \Ql{,\sigma}^{\textup{h}I}$-modules
\begin{equation}
    i_{\sigma}^*\Rlv_S(\sing{X_{\sigma}})\simeq p_*\Phi_{p}(\Ql{,X}(\beta))^{\textup{h}I}[-1]
\end{equation}
Here $i_{\sigma}:\sigma \hookrightarrow S$ is the embedding of the closed point in $S$, $\Phi_p(\Ql{,X}(\beta))$ is the $\ell$-adic sheaf of vanishing cycles associated to $\Ql{,X}(\beta)$ and $I$ is the inertia group (as $S$ is strictly henselian, it coincides with the absolute Galois group of the open point in $S$).
\end{theorem-non}
It might be useful to consider the following mind map:
\begin{equation*}
    \begindc{\commdiag}[18]
    \obj(0,80)[1]{$p:X\rightarrow S \text{ (initial data)}$}
    \obj(-70,60)[2]{$(X,\pi \circ p:X\rightarrow \aff{1}{S})\in \lgm{1}$}
    \obj(-70,50)[3]{$(\text{for } \pi \text{ a uniformizer of } S)$}
    \obj(-70,30)[4]{$\sing{X,\pi \circ p}\simeq \mf{X}{\pi \circ p}$}
    \obj(-70,20)[5]{$\simeq \sing{X_{\sigma}} \text{ (as } X \text{ is regular)}$}
    \obj(0,0)[6]{$i_{\sigma}^*\Rlv_S(\sing{X_{\sigma}})\simeq p_*\Phi_{p}(\Ql{,X}(\beta))^{\textup{h}I}[-1]$}
    \obj(70,60)[7]{$\Phi_p(\Ql{,X}(\beta))$}
    \obj(70,50)[8]{$\text{(vanishing cycles of } X)$}
    \obj(70,30)[9]{$\Phi_p(\Ql{,X}(\beta))^{\textup{h}I}$}
    \obj(70,20)[10]{$\text{(inertia invariant vanishing cycles)}$}
    \mor(-20,75)(-50,65){$$}
    \mor{3}{4}{$$}
    \mor(-50,15)(-20,5){$$}
    \mor(20,75)(50,65){$$}
    \mor{8}{9}{$$}
    \mor(50,15)(20,5){$$}
    \obj(-70,100)[]{$\text{NON COMMUTATIVE SIDE}$}
    \obj(70,100)[]{$\text{VANISHING CYCLES SIDE}$}
    \enddc
\end{equation*}
We can ask whether the above diagram makes sense in more general situations. 
\begin{itemize}
    \item one can start with the datum of a proper, flat and regular scheme over an excellent local regular ring of dimension $n$. One recovers the case treated in \cite{brtv} when $n=1$
\end{itemize}

It is immediate to observe that the left hand side of the diagram makes sense without any change:
\begin{itemize}
    \item assume that we are given a proper, flat morphism $p:X\rightarrow S$, with $X$ regular and $S$ local, regular of dimension $n$. Let $\underline{\pi}=(\pi_1,\dots,\pi_n)$ be a collection of generators of the closed point of $S$. Then we can consider the morphism $\underline{\pi}\circ p:X\rightarrow \aff{n}{S}$ and its fiber $X_0$ along the origin $S\rightarrow \aff{n}{S}$. It makes perfectly sense to consider $\Rlv_S(\sing{X_0})$
\end{itemize}
It comes up that this generalization is related to the following one
\begin{itemize}
    \item one can consider pairs $(X,s_X)$ where $X$ is a regular scheme over $S$ and $s_X:X\rightarrow Tot(\Lc_S)$ is a morphism towards the total space of a line bundle over $S$. One recovers the situation pictured above when $\Lc_S=\Oc_S$ is the trivial line bundle. In this situation, we want to compute $\Rlv_S(\sing{X_0})$, where $X_0$ is the fiber of $s_X:X\rightarrow Tot(\Lc_S)$ along the zero section $S\rightarrow Tot(\Lc_S)$
\end{itemize}. 

One can view the former generalization as a particular case of the latter thanks to a theorem of D.~Orlov (\cite{orl06}) and J.~Burke-M.~Walker (\cite{bw15}), which tells us that the dg category of singularities of $(X,\underline{\pi}\circ p)$ is equivalent to the dg category of singularities of $(\proj{n-1}{X},(\pi_1\circ p)\cdot T_1+\dots + (\pi_n\circ p)\cdot T_n \in \Oc(1)(\proj{n-1}{X}))$

Thus, we only need to find the appropriate generalization of (inertia invariant) vanishing cycles. The first thing we can think of are vanishing cycles over general bases, developed by G.~Laumon in \cite{lau82} following ideas of P.~Deligne. However, it seems that this is not the right point of view for our purposes. Instead, we will pursue the following analogy
\begin{equation*}
    \begindc{\commdiag}[18]
    \obj(-80,60)[]{$\text{TOPOLOGY}$}
    \obj(0,60)[]{$\text{VANISHING CYCLES}$}
    \obj(80,60)[]{$\text{OUR SETTING}$}
    \obj(-80,40)[]{$\mathbb{D}\text{ unital}$}
    \obj(-80,30)[]{$\text{open disk}$}
    \obj(0,40)[]{$S \text{ strictly}$}
    \obj(0,30)[]{$ \text{henselian trait}$}
    \obj(80,40)[]{$Tot(\Lc)\text{ total space}$}
    \obj(80,30)[]{$\text{of a line bundle}$}
    \obj(-80,10)[]{$0\hookrightarrow \mathbb{D}$}    \obj(-80,0)[]{$\text{origin}$}
    \obj(0,10)[]{$\sigma \hookrightarrow S$}
    \obj(0,0)[]{$\text{special point}$}
    \obj(80,10)[]{$S\hookrightarrow Tot(\Lc)$}
    \obj(80,0)[]{$\text{zero section}$}
    \obj(-80,-20)[]{$\mathbb{D}^*=\mathbb{D}-\{0\}\hookrightarrow \mathbb{D}$}    \obj(-80,-30)[]{$\text{punctured disk}$}
    \obj(0,-20)[]{$\eta \hookrightarrow S$}
    \obj(0,-30)[]{$\text{generic point}$}
    \obj(80,-20)[]{$\mathcal{U}=Tot(\Lc)-S \hookrightarrow Tot(\Lc)$}
    \obj(80,-30)[]{$\text{open complementary}$}
    \obj(-80,-50)[]{$\widetilde{\mathbb{D}^*} \rightarrow \mathbb{D}^*$}    \obj(-80,-60)[]{$\text{universal cover}$}
     \obj(-80,-70)[]{$\text{of the punctured disk}$}
    \obj(0,-50)[]{$\bar{\eta}\rightarrow \eta$}
    \obj(0,-60)[]{$\text{separable closure}$}
    \obj(0,-70)[]{$\text{of the generic point}$}
    \obj(80,-60)[]{$?$}
    \enddc
\end{equation*}
As we will only need to define the analogous of inertia-invariant vanishing cycles, we will not face in this thesis the problem of filling the empty spot in the mental map above. Nevertheless, we will come back to this matter at the end of the chapter, presenting a strategy to complete the picture.

We will define an appropriate generalization of $\Phi_p(\Ql{,X}(\beta))^{\textup{h}I}$ and to prove a generalization of the formula stated in Theorem \ref{main theorem brtv}. Our theorem will then look as follows
\begin{theorem-non}{\textup{\ref{main theorem}}}
Let $X$ be a regular scheme and let $s_X$ be a global section of a line bundle $\Lc_X$. Denote $X_0$ the zero locus of $s_X$. Then
\begin{equation*}
    \Rlv_{X_0}(\sing{X_0})\simeq \Phi^{\textup{mi}}_{(X,s_X)}(\Ql{}(\beta))[-1]
\end{equation*}
Here $\Phi^{\textup{mi}}_{(X,s_X)}(\Ql{}(\beta))$ is what we call the monodromy-invariant vanishing cycles $\ell$-adic sheaf (see Definition \ref{mododromy-invariant vanishing cycles def}). It coincides with inertia invariant vanishing cycles when the line bundle is trivial.
\end{theorem-non}
Using this theorem combined with the above mentioned result of D.~Orlov and J.~Burke-M.~Walker (Theorem \ref{reduction of codimension theorem}), we will deduce the following formula for the situation in which one starts with a regular ring $B$ and a regular sequence $\underline{f}$
\begin{theorem-non}{\textup{\ref{formula for (B,f1,fn)}}}
Let $B$ a regular ring and $\underline{f}$ a regular sequence. Then the following equivalence holds
\begin{equation*}
    \Rlv_B(\sing{B/\underline{f}})\simeq p_*i_*\Phi^{\textup{mi}}_{(\proj{n-1}{B},W_{\underline{f}})}(\Ql{}(\beta))[-1]
\end{equation*}
where $p:\proj{n-1}{B}\rightarrow B$, $i:V(W_{\underline{f}})\rightarrow \proj{n-1}{B}$ and $W_{\underline{f}}=f_1\cdot T_1+\dots +f_n\cdot T_n \in \Oc(1)(\proj{n-1}{B})$.
\end{theorem-non}

We will end this chapter with some remarks on the following two problems:
\begin{itemize}
    \item It seems possible to define a formalism of vanishing cycles in twisted situations, i.e. in the situation in which we have a morphism $s_X:X\rightarrow V(\Lc_X)$ for $\Lc_X \in \textup{\textbf{Pic}}(X)$. This is all about completing the empty slot in the mind map above. One should then be able to find $\Phi^{\textup{mi}}_{(X,s_X)}(\Ql{}(\beta))$ via a procedure that corresponds to taking homotopy fixed points in the usual situation
    \item We will comment the regularity hypothesis that appears both in A.~Blanc-M.~Robalo-B.~Toën-G.~Vezzosi's theorem and in the generalization we provide
\end{itemize}
\section{\textsc{Motivic realization of twisted LG models}}
\begin{notation}\label{context non affine base}
Let $S$ be a noetherian (not necessarely affine) regular scheme. We will label $\sch$ the category of $S$-schemes of finite type.
\end{notation}
\subsection{\textsc{The category of twisted LG models}}
Consider the category $\sch$ and let $\Lc_S$ be a line bundle on $S$. Then we can consider the category of Landau-Ginzburg models over $(S,\Lc_S)$, defined as follows:
\begin{itemize}
    \item objects consists of pairs $(X,s_X)$, where $p:X\rightarrow S$ is a (flat) $S$-scheme and $s_X$ is a section of $\Lc_X:=p^*\Lc_S$.
    \item morphisms $(X,s_X)\xrightarrow{f}(Y,s_Y)$ consists of morphisms $f:X\rightarrow Y$ in $\sch$ such that $s_X=f^*s_Y$.
    \item composition and identity morphisms are clear.
\end{itemize}
We will denote this category by $\tlgm{S}$.
\begin{rmk}\label{tlgm trivial line bundle}
If $\Lc_S$ is the trivial line bundle, then $\tlgm{S}$ coincides with the category of usual Landau-Ginzburg models over $S$. Indeed, in this case $Tot(\Lc_S)=Spec_{\Oc_S}(\Oc_S[t])$. Therefore, for any $X\in \sch$, a section of $\Lc_X$ consists of a morphism $\Oc_X[t]\rightarrow \Oc_X$, i.e. of a global section of $\Oc_X$.
\end{rmk}
\begin{construction}
It is possible to endow $\tlgm{S}$ with a symmetric monoidal structure 
\begin{equation}
    \boxplus: \tlgm{S}\times \tlgm{S}\rightarrow \tlgm{S}
\end{equation}
\[
\bigl ((X,s_X),(Y,s_Y) \bigr )\mapsto (X,s_X)\boxplus (Y,s_Y)=(X\times_S Y,s_X\boxplus s_Y)
\]
where $s_x\boxplus s_Y=p_X^*s_X+p_Y^*s_Y$. This is clearly associative and the unit is given by $(S,0)$, where $0$ stands for the zero section of $\Lc_S$.
\end{construction}
\begin{rmk}
If $\Lc_S$ is the trivial line bundle, then $\boxplus$ coincides with the tensor product on the category of LG modules over $S$.
\end{rmk}
We will now exhibit twisted LG models as a fibered category.
\begin{defn}
Let $\fclgm$ be the category defined as follows:
\begin{itemize}
    \item objects are triplets $(f:Y\rightarrow X,\Lc_X,s_Y)$ where $f$ is a flat morphism between $S$-schemes, $\Lc_X$ is a line bundle on $X$ and $s_Y$ is a global section of $f^*\Lc_X$.
    \item given two objects $(f_i:Y_i\rightarrow X_i,\Lc_{X_i},s_{Y_i})$ ($i=1,2$), a morphism from the first to the second is the datum of a commutative diagram
    \begin{equation}
        \begindc{\commdiag}[16]
        \obj(-15,15)[1]{$Y_1$}
        \obj(15,15)[2]{$X_1$}
        \obj(-15,-15)[3]{$Y_2$}
        \obj(15,-15)[4]{$X_2$}
        \mor{1}{2}{$f_1$}
        \mor{1}{3}{$g_Y$}
        \mor{2}{4}{$g_X$}
        \mor{3}{4}{$f_2$}
        \enddc
    \end{equation}
    and of an isomorphism $\alpha: g_X^*\Lc_{X_2}\rightarrow \Lc_{X_1}$ such that $s_{Y_1}$ corresponds to $s_{Y_2}$ under the isomorphism
    \begin{equation}
        g_{Y}^*f_2^*\Lc{X_2}\simeq f_1^*g_{X}^*\Lc_{X_2}\simeq f_1^*\Lc_{X_1}
    \end{equation}
    By abuse of notation, we will say in the future that $g_Y^*(s_{Y_2})=s_{Y_1}$ if this condition is satisfied. We will denote such a morphism by $(g,\alpha)$.
    \item Composition and identities are defined in an obvious way.
\end{itemize}
We will refer to this category as the category of generalized Landau-Ginzburg models.
\end{defn}
Notice that there is a functor
\begin{equation}\label{fibration LG}
    \pi:\fclgm \rightarrow \int_{X\in \sch} \bgm{S}(X)
\end{equation}
defined as
\begin{equation*}
    \begindc{\commdiag}[20]
    \obj(0,15)[1]{$(f_1:Y_1\rightarrow X_1,\Lc_{X_1},s_{Y_1})$}
    \obj(60,15)[2]{$(X_1,\Lc_{X_1})$}
    \obj(0,-15)[3]{$(f_2:Y_2\rightarrow X_2,\Lc_{X_2},s_{Y_2})$}
    \obj(60,-15)[4]{$(X_2,\Lc_{X_2})$}
    \mor{1}{3}{$(g,\alpha)$}
    \mor{2}{4}{$(g_X,\alpha)$}
    \obj(38,0)[5]{$\mapsto$}
    \enddc
\end{equation*}
\begin{lmm}
The functor $(\ref{fibration LG})$ exhibits $\fclgm$ as a fibered category over $\int_{X\in \sch} \bgm{S}(X)$.
\end{lmm}
\begin{proof}
Consider a map $(g_X,\alpha):(X_1,\Lc_{X_1})\rightarrow (X_2,\Lc_{X_2})$ in $\int_{X\in \sch} \bgm{S}(X)$ and let $(f_2:Y_2\rightarrow X_2,\Lc_{X_2},s_{Y_2})$ be an object of $\fclgm$ over $(X_2,\Lc_{X_2})$. Consider the morphism 
\[
(g,\alpha):(f_1:Y_1:=X_1\times_{X_2}Y_2\rightarrow X_1,\Lc_{X_1},s_{X_1})\rightarrow (f_2:Y_2\rightarrow X_2,\Lc_{X_2},s_{Y_2})
\]
where $g_Y$ is the projection $Y_1\rightarrow X_1$ (which is flat as it is the pullback of a flat morphism), $f_1$ is the projection $Y_1\rightarrow Y_2$ and $s_{Y_1}$ is $g_Y^*(s_{Y_2})$. It is clear that it is a morphism of $\fclgm$ over $(g_X,\alpha)$. We need to show that it is cartesian. Consider
\begin{equation}
    \begindc{\commdiag}[20]
    \obj(0,15)[1]{$(f_1:Y_1\rightarrow X_1,\Lc_{X_1},s_{Y_1})$}
    \obj(0,-15)[2]{$(X_1,\Lc_{X_1})$}
    \obj(80,15)[3]{$(f_2:Y_2\rightarrow X_2,\Lc_{X_2},s_{Y_2})$}
    \obj(80,-15)[4]{$(X_2,\Lc_{X_2})$}
    \mor{1}{3}{$(g,\alpha)$}[\atright,\solidarrow]
    \mor{2}{4}{$(g_X,\alpha)$}
    \mor{1}{2}{$$}
    
    \mor{3}{4}{$$}
    \obj(-55,45)[5]{$(h:Z\rightarrow W,\Lc_{W},s_{Z})$}
    
    \obj(-55,0)[6]{$(W,\Lc_{W})$}
    
    \mor{5}{6}{$$}
    
    \mor{6}{2}{$(p,\beta)$}
    \mor{5}{3}{$(q,\gamma)$}
    \mor(-40,40)(-5,20){$((r,p),\beta)$}[\atright,\dashArrow]
    \enddc
\end{equation}
Then the universal property of $Y_1$ gives us an unique morphism $r:Z\rightarrow Y_1$ such that the compositions with $f_1$ and $g_Y$ are $p\circ h$ and $q$ respectively. We just need to show that $r^*(s_{Y_1})=s_Z$. But this is clear since $s_{Y_1}=g_Y^*(s_{Y_2})$, $q_Y^*(s_{Y_2})=s_Z$ and $g_Y\circ r=q_Y$
\end{proof}
\begin{rmk}
Let $(X,\Lc_X)$ be an object of $\int_{X\in \sch} \bgm{S}(X)$. Then the fiber of $(X,\Lc_X)$ along $\pi$ is $\tlgm{X}$.
\end{rmk}
\begin{defn}
We say that a collection of maps $\{(g_i,\alpha_i):(U_i,\Lc_{U_i})\rightarrow (X,\Lc_{X})\}_{i\in I}$ in $\int_{X\in \sch} \bgm{S}(X)$ is a Zariski covering if $\{g_i:U_i\rightarrow X\}_{i\in I}$ is so. They clearly define a pre-topology on $\int_{X\in \sch} \bgm{S}(X)$\footnote{Notice that pullbacks exists in $\int_{X\in \sch} \bgm{S}(X)$: the pullback of $(Y_1,\Lc_{Y_1})\xrightarrow{(f_1,\alpha_1)}(X,\Lc_X)\xleftarrow{(f_2,\alpha_2)}(Y_2,\Lc_{Y_2})$ is $(Y_1\times_X Y_2,\Lc_{Y_1\times_X Y_2})$ with the projections defined in an obvious way}. We will refer to the corresponding topology by Zarisky topology on $\int_{X\in \sch} \bgm{S}(X)$. 
\end{defn}
\begin{lmm}
$\fclgm$ is a stack over $\int_{X\in \sch} \bgm{S}(X)$ endowed with the Zariski topology.
\end{lmm}
\begin{proof}
This is a simple consequence of the fact that a morphism of schemes is uniquely determined by its restriction to a Zariski covering and that line bundles are sheaves.
\end{proof}
We conclude this section with the following observation:
\begin{lmm}\label{zariski descent tlgmm}
Let $\{(g_i,\alpha_i):(U_i,\Lc_{U_i})\rightarrow (X,\Lc_{X})\}_{i\in I}$ be a Zariski covering in $\int_{X\in \sch} \bgm{S}(X)$. Then the canonical functor
\begin{equation}\label{compatibility limits with tlgmo}
    \tlgm{X}^{\boxplus}\rightarrow \varprojlim \tlgm{U_i}^{\boxplus}
\end{equation}
is a symmetric monoidal equivalence.
\end{lmm}
\begin{proof}
Consider the functors
\begin{equation}
    \tlgm{X}\rightarrow \tlgm{U_i}
\end{equation}
\begin{equation*}
    (f:Y\rightarrow X,s_Y)\mapsto (f_i:Y\times_X U_i\rightarrow U_i,s_{Y|U_i})
\end{equation*}
It is easy to see that they respect the tensor structure. Therefore, we get the desired symmetric monoidal functor $(\ref{compatibility limits with tlgmo})$ in the (big) category of symmetric monoidal categories and symmetric monoidal functors. Then, in order to prove that is a symmetric monoidal equivalence, it suffices to show that the underlying functor 
\begin{equation}
     Forget(\tlgm{X}^{\boxplus}\rightarrow \varprojlim \tlgm{U_i}^{\boxplus})\simeq (\tlgm{X}\rightarrow Forget(\varprojlim \tlgm{U_i}^{\boxplus}))
\end{equation}
As the functor which forgets the symmetric monoidal structure of a category is a right adjoint, it preserves limits. Thus
\begin{equation}
    Forget(\varprojlim \tlgm{U_i}^{\boxplus})\simeq \varprojlim \tlgm{U_i}
\end{equation}
The assertion now follows from the previous lemma.
\end{proof}
\subsection{\textsc{The dg category of singularities of a twisted LG model}}\label{The dg category of singularities of a twisted LG model}
Let $(X,s_X)$ be a generalized $(S,\Lc_S)$ LG-model. The section $s_X$ defines a closed sub-scheme of $X$. Since we are not assuming that the section is regular, some derived structure might appear. More precisely: let $\Oc_X\rightarrow \Lc_X$ be the morphism of $\Oc_X$-modules associated to $s_X$. Then, taking duals, it defines a morphism $\Lc_X^{\vee}\rightarrow \Oc_X$. If we apply the relative spectrum functor, we get a morphisms $Spec_{\Oc_X}(\Oc_X)=X\rightarrow Spec_{\Oc_X}(Sym_{\Oc_X}(\Lc_X^{\vee}))=V(\Lc_X)$, i.e. a section of the vector bundle associated to $\Lc_X$. By abuse of notation, we will label this morphism by $s_X$. We can also consider the zero section $X\rightarrow V(\Lc_X)$, which is the morphism associated to $0\in \Lc_X(X)$. Then the derived pullback diagram
\begin{equation}\label{pullback square derived zero locus}
    \begindc{\commdiag}[20]
    \obj(0,15)[1]{$X_0$}
    \obj(30,15)[2]{$X$}
    \obj(0,-15)[3]{$X$}
    \obj(30,-15)[4]{$V(\Lc_X)$}
    \mor{1}{2}{$\mathfrak{i}$}
    \mor{1}{3}{$$}
    \mor{2}{4}{$s_X$}
    \mor{3}{4}{$0$}
    \enddc
\end{equation}

gives us the desired derived closed subscheme of $X$.
\begin{rmk}
Notice that the zero section $X\xrightarrow{0}V(\Lc_X)$ is an lci morphism (Zariski locally, it is just the zero section a scheme $Y$ in the affine line $\aff{1}{Y}$). It is well known that this class of morphisms is stable under derived pullbacks. In particular, $\mathfrak{i}:X_0\rightarrow X$ is a derived lci morphism.
\end{rmk}
\begin{rmk}
There is a truncation morphism $\mathfrak{t}:\pi_0(X_0)\rightarrow X$, where $\pi_0(X_0)$ is the classical scheme associated to $X_0$. Whenever $s_X$ is a flat morphism, the truncation morphism is an equivalence in the $\oo$-category of derived schemes.
\end{rmk}

\begin{rmk}
Of major importance for our purposes is that $\mathfrak{i}$ is an lci morphism. Indeed, by \cite{to12}, if $f:Y\rightarrow Z$ is an lci morphism of derived schemes and $E\in \perf{Y}$, then $f_*E\in \perf{Z}$. In particular, we get that $\perf{X_0}$ is a full subcategory of $\cohbp{X_0}{X}$.
\end{rmk}

The dg category of absolute singularities of a derived scheme is a non-commutative invariant (in the sense of Kontsevich) which captures the singularities of the scheme. For the reader's convenience, let us recall the definition of this crucial character, that we have already introduced in the previous chapter:
\begin{defn}
Let $Z$ be a derived scheme of finite type over $S$ whose structure sheaf is cohomologically bounded. Then the dg category of absolute singularities is the dg quotient in $\dgcatm$
\begin{equation}
    \sing{Z}:=\cohb{Z}/\perf{Z}
\end{equation}
\end{defn}
\begin{rmk}
It is a classical theorem due to M.~Auslander-D.A.~Buchsbaum (\cite[Theorem 4.1]{ab56}) and J.P.~Serre (\cite[Théorème 3]{se55}) that a noetherian ring $A$ is regular if and only if it has finite global dimension. This means that the category of perfect complexes $\perf{A}$ coincides with that of coherent bounded complexes $\cohb{A}$, i.e. $\sing{A}\simeq 0$. This explains the name of this object.
\end{rmk}
\begin{rmk}
Notice that the hypothesis on $Z$ are crucial, as in general the structure sheaf of a derived scheme might not be cohomologically bounded. As an example of a derived scheme that doesn't sit in this class of objects, consider 
\[
Spec(\C[x]\otimes^{\mathbb{L}}_{\C[x,y]/(xy)}\C[y])\simeq Spec(\C[x])\times^h_{Spec(\C[x,y]/(xy))}Spec(\C[y])
\]
By tensoring the projective $\C[x,y]/(xy)$-resolution of $\C[x  ]$
\[
\dots \xrightarrow{y} \C[x,y]/(xy)\xrightarrow{x}\C[x,y]/(xy)\xrightarrow{y}\C[x,y]/(xy)\rightarrow 0
\]
with $\C[y]$ over $\C[x,y]/(xy)$ we get that this derived scheme is the spectrum of the cdga
\[
\dots \xrightarrow{y} \C[y]\xrightarrow{0}\C[y]\xrightarrow{y}\C[y]\rightarrow 0
\]
which has nontrivial cohomology in every even negative degree.
\end{rmk}
Since $\mathfrak{i}$ is an lci morphism of derived schemes, by \cite{to12} and by \cite{gr17} the pushforward induces a morphism 
\begin{equation}\label{push sing}
\mathfrak{i}_*: \sing{X_0}\rightarrow \sing{X}
\end{equation}
\begin{defn}
Let $(X,s_X)$ be a twisted LG model over $(S,\Lc_S)$. Its dg category of singularities is defined as the fiber in $\dgcatm$ of $(\ref{push sing})$:
\begin{equation}
    \sing{X,s_X}:=fiber\bigl(\mathfrak{i}_*:\sing{X_0}\rightarrow \sing{X} \bigr)
\end{equation}
\end{defn}
\begin{rmk}
It is a consequence of the functoriality properties of the pullback and that of taking quotients and fibers in $\dgcatm$ that the assignments $(X,s_X)\mapsto \sing{X,s_X}$ can be organized in an $\oo$-functor
\begin{equation}\label{functor sing tlgm}
    \sing{\bullet,\bullet}:\tlgm{S}^{\textup{op}}\rightarrow \dgcatm
\end{equation}
\end{rmk}
We will need to endow the $\oo$-functor $(\ref{functor sing tlgm})$ with the structure of a lax monoidal $\oo$-functor. In order to do so, we will introduce a strict model for $\cohbp{X_0}{X}$, following the strategy exploited in \cite{brtv}.
\begin{construction}\label{strict model for the derived zero locus} Let $(X,s_X)$ be a twisted LG model over $S$. The identity morphism on $\Lc^{\vee}_X$ induces a $Sym_{\Oc_X}(\Lc^{\vee}_X)$-projective resolution of $\Oc_X$, namely
\begin{equation}
    \Lc^{\vee}_X\otimes_{\Oc_X}Sym_{\Oc_X}(\Lc^{\vee}_X)\rightarrow Sym_{\Oc_X}(\Lc^{\vee}_X)
\end{equation}
Tensoring up along the morphism $Sym_{\Oc_X}(\Lc^{\vee}_X)\rightarrow \Oc_X$ induced by $s_X$ we obtain that the cdga associated to $X_0$ is the spectrum of $\Lc^{\vee}_X\xrightarrow{s_X} \Oc_X$.
\end{construction}
\begin{ex}\label{derived zero locus zero section}
Consider the zero section of $\Lc_X$. Then the cdga associated to the structure sheaf of $X\times^h_{V(\Lc_X)}X$ is the Koszul algebra (in degrees $[-1,0]$)
\begin{equation}
    \Lc^{\vee}_X\xrightarrow{0}\Oc_X
\end{equation}
\end{ex}
\begin{construction}
Let $S=Spec(A)$ be a noetherian, regular affine scheme and let $L$ be a projective $A$-moddule of rank $1$. Let $\tlgm{S}^{\textup{aff}}$ be the category of affine twisted LG models over $(S,L)$, i.e. of pairs $(Spec(B),s)$ where $B$ is a flat $A$-algebra of finite type and $s\in L_B=L\otimes_A B$. For any object $(X,s)=(Spec(B),s)$, the derived zero locus of $s$ is given by the spectrum of the Koszul algebra $K(B,s)$. The associated cdga is
\begin{equation*}
    K(B,s)=L_B^{\vee}\xrightarrow{s}B
\end{equation*}
concentrated in degrees $-1$ and $0$, where we label $s:L_B^{\vee}\rightarrow B$ the morphism induced by $s:B\rightarrow L_B$. Then, similarly to \cite[Remark 2.30]{brtv}, there is an equivalence between cofibrant $K(B,s)$-dg modules (denoted $\widehat{K(B,s)}$) and $\qcoh{X_0}$. Under this equivalence, $\cohb{X_0}$ corresponds to the full sub-dg category of $\widehat{K(B,s)}$ spanned by those cohomologically bounded dg modules with coherent cohomology (over $B/s=coker(L_B^{\vee}\xrightarrow{s}B)$) and $\perf{X_0}$ corresponds to homotopically finitely presented $K(B,s)$-dg modules. Moreover, the pushforward
\begin{equation*}
\mathfrak{i}_*:\qcoh{X_0}\rightarrow \qcoh{X}
\end{equation*}
corresponds to the forgetful functor
\begin{equation*}
    \widehat{K(B,s)}\rightarrow \widehat{B}
\end{equation*}

Define $\cohs{B,s}$ as the dg category of $K(B,s)$-dg modules whose underlying $B$-dg module is perfect. Since $X$ is an affine scheme, every perfect $B$-dg module is equivalent to a strictly perfect one. Therefore, we can take $\cohs{B,s}$ as the dg category of $K(B,s)$-dg moodules whose underlying $B$-module is strictly perfect. More explicitly, such an object corresponds to a triplet $(E,d,h)$ where $(E,d)$ is a strictly perfect $B$-dg module and 
\begin{equation*}
    h:E\rightarrow E\otimes_B L_B[1]
\end{equation*}
is a morphism such that 
\begin{equation*}
    (h\otimes id_{L_B}[1])\circ h=0 \hspace{0.5cm} 
    [d,h]=h\circ d+(d\otimes id[1])\circ h=id \otimes s
\end{equation*}
Notice that $\cohs{B,s}$ is a locally flat $A$-linear dg category. This follows immediately from the fact that the underlying complexes of $B$-modules are strictly perfect and from the fact that $B$ is a flat $A$-algebra.

\end{construction}

This strict dg category is analogous to the one that has been introduced in \cite{brtv}. It gives us a strict model for the $\oo$-category $\cohbp{X_0}{X}$. Indeed, the following result holds
\begin{lmm}\textup{\cite[Lemma 2.33]{brtv}}\\
Let $\cohs{B,s}^{\textup{acy}}$ be the full sub-category of $\cohs{B,s}$ spanned by acyclic dg modules. Then the cofibrant replacement induces an equivalence of dg categories
\begin{equation}
    \cohs{B,s}[\textup{q.iso}^{-1}]\simeq \cohs{B,s}/\cohs{B,s}^{\textup{acy}}\simeq \cohbp{X_0}{X}
\end{equation}
As a consequence, if we label $\textup{Perf}^s(B,s)$ the full subcategory of $\cohs{B,s}$ spanned by perfect $K(B,s)$-dg modules, there are equivalances of dg categories
\begin{equation}
   \cohs{B,s}/\textup{Perf}^s(B,s)\simeq \cohbp{X_0}{X}/\perf{X_0}\simeq \sing{X,s}
   \end{equation}
\end{lmm}

\begin{construction}
We can endow the assignment $(B,s)\rightarrow \cohs{B,s}$ with the structure of a pseudo-functor: let $f:(B,s)\rightarrow (C,t)$ be a morphism of affine twisted LG models, i.e. a morphism of $A$-algebras $f:B\rightarrow C$ such that the induced morphism $id\otimes f:L_B\rightarrow L_C$ sends $s$ to $t$. Then we can define the dg functor
\begin{equation}
    -\otimes_B C :\cohs{B,s}\rightarrow \cohs{C,t}
\end{equation}
\begin{equation*}
    (E,d,h)\mapsto (E\otimes_B C,d\otimes id, h\otimes id)
\end{equation*}
Indeed, it is clear that if $(E,d,h)$ is a $K(B,s)$-dg module whose underling $B$-module is perfect, then the $K(C,t)$-dg module $(E\otimes_BC,d\otimes id,h\otimes id)$ is levelwise $C$-projective and strictly bounded.

This yields a pseudo-functor
\begin{equation}\label{cohs}
    \cohs{\bullet,\bullet}:\textup{LG}^{\textup{aff,op}}_{(S,L)}\rightarrow \dgcatlf
\end{equation}
Moreover, we can endow $(B,s)\mapsto \cohs{B,s}$ with a weakly associative and weakly unital lax monoidal structure. For any pair of twisted affine LG models $(B,s)$, $(C,t)$, we need to construct a morphism
\begin{equation}
    \cohs{B,s}\otimes \cohs{C,t}\rightarrow \cohs{B\otimes_AC,s\boxplus t}
\end{equation}
For a pair $(B,s)$, let $Z(s)$ denote $Spec(K(B,s))$. Moreover, let $V(L)=Spec(Sym(L^{\vee}))$. Then consider the following diagram
\begin{equation}
    \begindc{\commdiag}[16]
    
    \obj(0,30)[1]{$Z(s)\times^h_{S}Z(t)$}
    \obj(60,30)[2]{$Z(s\boxplus t)$}
    \obj(0,0)[3]{$S$}
    \obj(60,0)[4]{$V(L)$}
    \obj(120,30)[5]{$Spec(B\otimes_A C)$}
    \obj(120,0)[6]{$V(L)\times_{S}V(L)$}
    \obj(60,-30)[7]{$S$}
    \obj(120,-30)[8]{$V(L)$}
    \mor{1}{2}{$\phi$}
    \mor{1}{3}{$$}
    \mor{2}{4}{$$}
    \mor{3}{4}{$\textup{zero}$}
    \mor{2}{5}{$$}
    \mor{4}{6}{$(id,-id)$}
    \mor{5}{6}{$s\times t$}
    \mor{4}{7}{$$}
    \mor{6}{8}{$+$}
    \mor{7}{8}{$\textup{zero}$}
    \enddc
\end{equation}

where $\phi$ is the morphism corresponding to
\begin{equation}
    \begindc{\commdiag}[16]
    \obj(-100,20)[a]{$K(B,s)\otimes_A K(C,t)$}
    \obj(-50,20)[1]{$L_{B\otimes_AC}^{\vee}$}
    \obj(0,20)[2]{$L_{B\otimes_AC}^{\vee}\oplus L_{B\otimes_AC}^{\vee}$}
    \obj(50,20)[3]{$B\otimes_AC$}
    \obj(-100,-20)[b]{$K(B\otimes_AC,s\boxplus t)$}
    \obj(0,-20)[4]{$L_{B\otimes_AC}^{\vee}$}
    \obj(50,-20)[5]{$B\otimes_AC$}
    \mor{1}{2}{$$}
    \mor{2}{3}{$$}
    \mor{4}{5}{$$}
    \mor{4}{2}{$\begin{bmatrix}
    1\\ 1
    \end{bmatrix}$}
    \mor{5}{3}{$1$}
    \mor{b}{a}{$\phi$}
    \enddc
\end{equation}
\end{construction}

Consider the two projections
\begin{equation}
    \begindc{\commdiag}[20]
    \obj(0,0)[1]{$Z(s)\times^h_{S}Z(t)$}
    \obj(-30,0)[2]{$Z(s)$}
    \obj(30,0)[3]{$Z(t)$}
    \mor{1}{2}{$p_s$}[\atright,\solidarrow]
    \mor{1}{3}{$p_t$}
    \enddc
\end{equation}
Given $(E,d,h)\in \cohs{B,s}$ and $(E',d',h')\in \cohs{C,t}$, define 
\begin{equation}\label{multiplication lax monoidal structure cohs}
    (E,d,h)\boxtimes (E',d',h')=\phi_*\bigl ( p_s^*(E,d,h)\otimes p_t^*(E',d',h') \bigr )
\end{equation}
This is a strictly perfect complex of $B\otimes_{A}C$-modules: its underlying complex of $A$-modules is $(E,d)\otimes_A(E',d')$, which immediately implies that it is strictly bounded. In order to see that each component $(E\otimes_A E')_n=\oplus_{i=0 }^nE_i\otimes_AE'_{n-i}$ is a projective $B\otimes_A C$-module, the same argument given in \cite[Construction 2.35]{brtv} applies. The morphism
\begin{equation*}
    E\otimes_AE'\rightarrow E\otimes_AE'\otimes_{B\otimes_AC}L_{B\otimes_AC}[1] 
\end{equation*}
is $h\otimes 1+ 1\otimes h'$. In particular, $\phi_*\bigl ( p_s^*(E,d,h)\otimes p_t^*(E',d',h') \bigr )$ lives in $\cohs{B\otimes_AC,s\boxplus t}$.
The map 
\begin{equation}\label{unit lax monoidal structure cohs}
    \underline{A}\rightarrow \cohs{A,0}
\end{equation}
where $\underline{A}$ is the unit in $\textup{dgCat}^{\textup{lf}}_{A}$, i.e. the dg category with only one object whose endomorphism algebra is $A$, is defined by the assignment
\begin{equation*}
    \bullet \mapsto A
\end{equation*}
It is clear that (\ref{multiplication lax monoidal structure cohs}) and (\ref{unit lax monoidal structure cohs}) satisfy the associativity and unity axioms, i.e. they enrich the pseudo-functor (\ref{cohs}) with a lax monoidal structure
\begin{equation}\label{cohs lax monoidal}
    \cohs{\bullet,\bullet}^{\otimes}:\textup{LG}_{(S,L)}^{\textup{aff,op},\boxplus}\rightarrow \textup{dgCat}_{A}^{\textup{lf},\otimes}
\end{equation}

As in \cite[Construction 2.34, Construction 2.37]{brtv}, consider $\textup{Pairs-dgCat}^{\textup{lf}}_{A}$, the category of pairs $(T,F)$, where $T\in \textup{dgCat}^{\textup{lf}}_{A}$ and $F$ is a class of morphisms in $T$. A morphisms $(T,F)\rightarrow (T',F')$ is a dg functor $T\rightarrow T'$ sending $F$ in $F'$. There is an obvious functor $\textup{Pairs-dgCat}^{\textup{lf}}_{A}\rightarrow \textup{dgCat}^{\textup{lf}}_{A}$ defined as $(T,F)\mapsto T$. We say that a morphism $(T,F)\rightarrow (T',F')$ in $\textup{Pairs-dgCat}^{\textup{lf}}_{A}$ is a DK-equivalence if the dg functor $T\rightarrow T'$ is a DK-equivalence in $\textup{dgCat}^{\textup{lf}}_{A}$. Denote this class of morphisms in $\textup{Pairs-dgCat}^{\textup{lf}}_{A}$ by $W_{DK}$.

Also notice that the symmetric monoidal structure on $\textup{dgCat}^{\textup{lf}}_{A}$ induces a symmetric monoidal structure on $\textup{Pairs-dgCat}^{\textup{lf}}_{A}$
\begin{equation}
    (T,F)\otimes (T',F'):= (T\otimes T',F\otimes F')
\end{equation}
It is clearly associative and unital, where the unit is $(\underline{A},\{id\})$. We will denote this symmetric monoidal structure by $\textup{Pairs-dgCat}^{\textup{lf},\otimes}_{A}$. Note that this symmetric monoidal structure is compatible to the class of morphisms in $W_{DK}$, as we are working with locally-flat dg categories.

Then  as in \cite[Construction 2.34, Construction 2.37]{brtv}, we can construct a symmetric monoidal $\oo$-functor
\begin{equation}
    \textup{loc}_{dg}^{\otimes}: \textup{Pairs-dgCat}^{\textup{lf},\otimes}_{A}[W_{DK}^{-1}]\rightarrow \textup{dgCat}^{\textup{lf},\otimes}_{A}[W_{DK}^{-1}]\simeq \dgcatdko
\end{equation}
which, on objects, is defined by sending $(T,F)$ to $T[F^{-1}]_{dg}$, the dg localization of $T$ with respect to $F$. Let $\{q.iso\}$ be the class of quasi-isomorphisms in $\textup{Coh}^s(R_i,s_{R_i})$. Consider the functor
\begin{equation}
    \textup{LG}^{\textup{aff,op},\boxplus}_{A}\rightarrow \textup{Pairs-dgCat}^{\textup{lf},\otimes}_{A}
\end{equation}
\begin{equation*}
    (B,s)\mapsto \bigl (\textup{Coh}^{s}(B,s), \{q.iso\} \bigr )
\end{equation*}

and compose it with the functor
\begin{equation}
    \textup{Pairs-dgCat}^{\textup{lf},\otimes}_{A}\xrightarrow{\textup{loc}} \textup{Pairs-dgCat}^{\textup{lf},\otimes}_{A}[W_{DK}^{-1}]\xrightarrow{\textup{loc}_{dg}}\textup{dgCat}^{\textup{lf},\otimes}_{A}[W_{DK}^{-1}]\simeq \textup{\textbf{dgCat}}^{\otimes}_{A}\xrightarrow{\textup{loc}} \textup{\textbf{dgCat}}^{\textup{idm},\otimes}_{A}
\end{equation}
we get, after a suitable monoidal left Kan extension, the desired lax monoidal $\oo$-functor
\begin{equation}\label{cohbp monoidal}
    \cohbp{\bullet}{\bullet}^{\otimes}:\textup{LG}^{\textup{op},\boxplus}_{A}\rightarrow \textup{\textbf{dgCat}}^{\textup{idm},\otimes}_{A}
\end{equation}

\begin{lmm}\label{RHom(A,A)=Sym_A(L[-2])}
Let $S=Spec(A)$ be a regular noetherian ring. Let $L$ be a line bundle over $S$. The following equivalence holds in $\textup{CAlg}(\dgcatmo)$
\begin{equation}\label{equivalence cohb(S_0) perf(R[u])}
    \cohbp{S_0}{S}^{\otimes}\simeq \cohb{S_0}^{\otimes}\xrightarrow{\sim} \perf{Sym_{\Oc_S}(\Lc[-2])}^{\otimes}
\end{equation}
\end{lmm}
\begin{proof}
The first equivalence is an immediate consequence of the regularity hypothesis on $S$.
As we have remarked in the Example (\ref{derived zero locus zero section}), the cdga assocaiated via the Dold-Kan correspondence to $S_0$ is $L^{\vee}\xrightarrow{0} A$. By definition, $A$ be a compact generator of $\perf{S}$. Also by definition, $L^{\vee}\xrightarrow{0}A$ is a compact generator of $\perf{S_0}$. It is easy to see that $A$, with the trivial $(L^{\vee}\xrightarrow{0}A)$-dg module structure, is a generator of $\cohb{S_0}$: if $M$ is a nonzero object in $\cohb{S_0}$, then there exists $0\neq m \in \textup{H}^i(M)$ for some $i$. We can moreover assume that $i=0$. Then $m$ induces a non-zero morphism $A\rightarrow M$.
In particular, we get the equivalence $\perf{\rhom_{(L^{\vee}\xrightarrow{0}A)}(A,A)}\simeq \cohb{S_0}$.
The endomorphism dg algebra $\rhom_{(L^{\vee}\xrightarrow{0}A)}(A,A)$ can be computed explicitly by means of the following $L^{\vee}\xrightarrow{0}A$-resolution of $A$ 
\begin{equation}
    \dots \underbracket{L^{\vee \otimes 2}}_{-3}\xrightarrow{0}\underbracket{L^{\vee}}_{-2}\xrightarrow{1}\underbracket{L^{\vee}}_{-1}\xrightarrow{0}\underbracket{A}_{0}
\end{equation}
Applying $\hm_{A}(-,A)$ we obtain
\begin{equation}
\underbracket{A}_{0}\xrightarrow{0}\underbracket{\hm_{A}(L^{\vee},A)}_{1}\xrightarrow{1}\underbracket{\hm_{A}(L^{\vee},A)}_{2}\xrightarrow{0}\underbracket{\hm_{A}(L^{\vee \otimes 2},A)}_{3}\dots
\end{equation}
\begin{equation*}
\simeq
    \underbracket{A}_{0}\xrightarrow{0}\underbracket{L}_{1}\xrightarrow{1}\underbracket{L}_{2}\xrightarrow{0}\underbracket{L^{\otimes2}}_{3}\dots
\end{equation*}
which is quasi isomorphic to $A$. However, when we ask for $(L^{\vee}\xrightarrow{0}A)$-linearity, the copies of $L^{\otimes n}$ in odd degree disappear, as the local generators $\varepsilon$ in degree $-1$ of $(L^{\vee}\xrightarrow{0}A)$ act via the identity on $(L^{\vee}\xrightarrow{0}A)$. Therefore we find that $\rhom_{(L^{\vee}\xrightarrow{0} A)}(A,A)$ is quasi-isomorphic to
\begin{equation}
    \underbracket{A}_{0}\rightarrow \underbracket{0}_{1}\rightarrow \underbracket{L}_{2}\rightarrow \dots \simeq Sym_{A}(L[-2])
\end{equation}
This shows that there exists an equivalence
\begin{equation}\label{equivalence modules}
    \rhom_{(L^{\vee}\xrightarrow{0}A)}(A,A)\simeq Sym_{A}(L[-2]) 
\end{equation}
as $(L^{\vee}\xrightarrow{0}A)$-dg modules. Notice that both these objects carry a canonical algebra structure. 
In order to conclude that the two commutative algebra structures coincide, we consider the dg functor
\begin{equation}
    \rhom_{(L^{\vee}\xrightarrow{0}A)}(A,-): \cohs{S,0}\rightarrow Sym_{A}(L[-2])-dgmod
\end{equation}

Notice that $Sym_{A}(L[-2])$ is a strict cdga, seen as a commutative algebra object in $\qcoh{A}$.
Similarly to \cite[Lemma 2.39]{brtv}, for $(E,d,h)$ an object in $\cohs{A,0}$, the $Sym_A(L[-2])$-dg module $\rhom_{(L^{\vee}\xrightarrow {0}A)}(A,(E,d,h))$ can be computed (in degrees $i$ and $i+1$) as
\begin{equation}
    \bigoplus_{n\geq 0}E_{i-2n}\otimes_A L^{\otimes n} \xrightarrow{d+h} \bigoplus_{n\geq 0}E_{i+1-2n}\otimes_A L^{\otimes n}
\end{equation}
The same arguments given in \textit{loc.cit.} hold mutatis mutandis in our situation and therefore we obtain a symmetric monoidal functor
\begin{equation}
      \cohs{S,0}^{\otimes}\rightarrow Sym_{A}(L[-2])-dgmod^{\otimes}
\end{equation}
which preserves quasi-isomorphisms. If we localize both the l.h.s. and the r.h.s. we thus obtain a symmetric monoidal $\oo$-functor
\begin{equation}
    \cohb{S_0}^{\otimes}\rightarrow \qcoh{Sym_A(L[-2])}^{\otimes}
\end{equation}
from which one recovers the equivalence of symmetric monoidal dg categories
\begin{equation}
    \cohb{S_0}^{\otimes}\simeq \perf{Sym_A(L[-2])}
\end{equation}
\end{proof}
\begin{cor}
Let $S=Spec(A)$ be a regular noetherian affine scheme and let $L$ be a line bundle over $S$. The lax monoidal $\oo$-functors $(\ref{cohbp monoidal})$ factors as
\begin{equation}\label{cohbp monoidal modules 2}
    \cohbp{\bullet}{\bullet}^{\otimes}:\textup{LG}_{(S,L)}^{\textup{op},\boxplus}\rightarrow \md_{\perf{Sym_{A}(L[-2])}}(\dgcatm)^{\otimes}
\end{equation}
\end{cor}
Then, for a noetherian regular scheme $S$ with a line bundle $\Lc_S$, we obtain a lax monoidal $\oo$-functor
\begin{equation}\label{cohbp non affine}
    \cohbp{\bullet}{\bullet}^{\otimes}:\tlgm{S}^{\textup{op},\otimes}\rightarrow \dgcatmo
\end{equation}
as the limit
\begin{equation}
    \varprojlim_{(Spec(A),\Lc_{S|Spec(A)})}\bigl ( \cohbp{\bullet}{\bullet}^{\otimes}:\textup{LG}_{Spec(A),\Lc_{S|Spec(A)}}^{\textup{op},\boxplus}\rightarrow \textup{\textbf{dgCat}}_A^{\textup{idm},\otimes}\bigr )
\end{equation}
where $Spec(A)\rightarrow S$ is a Zariski open subscheme. We have used Lemma \ref{zariski descent tlgmm} and the equivalence
\begin{equation}
    \varprojlim_{Spec(A)\rightarrow S}\textup{\textbf{dgCat}}_A^{\textup{idm},\otimes}\simeq \dgcatmo
\end{equation}
\begin{rmk}\label{actions}
The monoidal structure on $\cohbp{\bullet}{\bullet}^{\otimes}$ implies that each dg category $\cohbp{X_0}{X}$ is endowed with an action of $\cohb{S_0}$ (recall that we are assuming that $S$ is regular). Similarly to \cite[Remark 2.38]{brtv} we can describe this action. Consider the diagram
\begin{equation}
    \begindc{\commdiag}[20]
    \obj(20,30)[1]{$X_0\times^h_S S_0$}
    \obj(0,15)[2]{$X_0$}
    \obj(70,30)[3]{$X_0$}
    \obj(50,15)[4]{$X$}
    \obj(20,0)[5]{$S_0$}
    \obj(0,-15)[6]{$S$}
    \obj(70,0)[7]{$S$}
    \obj(50,-15)[8]{$V(\Lc_S)$}
    \mor{1}{2}{$p_{X_0}$}[\atright,\solidarrow]
    \mor{1}{3}{$a$}
    \mor{2}{4}{$\mathfrak{i}$}
    \mor{3}{4}{$\mathfrak{i}$}
    \mor{5}{6}{$$}
    \mor{5}{7}{$$}
    \mor{6}{8}{$\textup{zero}$}
    \mor{7}{8}{$\textup{zero}$}
    \mor{1}{5}{$$}
    \obj(14,9)[a]{$p_{S_0}$}
    \mor{2}{6}{$$}
    \mor{3}{7}{$$}
    \mor{4}{8}{$$}
    \enddc
\end{equation}
Notice that $X\times^{h}_S S_0\simeq X_0\times^h_X X_0$. Given $\Ec \in \cohbp{X_0}{X}$ and $\Fc \in \cohb{S_0}$, then $\Ec \boxtimes \Fc= a_*(p_{X_0}^* \Ec \otimes p_{S_0}^*\Fc)$. In particular, when $\Fc=\Lc_S^{\vee}\xrightarrow{0}\Oc_S=\Oc_{S_0}$, then
\begin{equation}
    \Ec \boxtimes (\Lc_S^{\vee}\xrightarrow{0}\Oc_S)=a_*(p_{X_0}^* \Ec\otimes p_{S_0}^*\Oc_{S_0})\simeq a_*(p_{X_0}^* \Ec\otimes a^*\Oc_{X_0})
\end{equation}
\begin{equation*}
    \underbracket{\simeq}_{\textup{proj. form.}}a_*p_{X_0}^*\Ec \underbracket{\simeq}_{\textup{der. prop. base change}}\mathfrak{i}^*\mathfrak{i}_*\Ec
\end{equation*}
and when $\Fc=\Oc_S=t_*\Oc_S$, where $t:S=\pi_0(S_0)\rightarrow S_0$ is the truncation morphism, the (homotopy) cartesian square
\begin{equation}
    \begindc{\commdiag}[20]
    \obj(0,15)[1]{$X_0\times_{S}S$}
    \obj(50,15)[2]{$X_0\times_S^hS_0$}
    \obj(0,-15)[3]{$S$}
    \obj(50,-15)[4]{$S_0$}
    \mor{1}{2}{$id\times t $}
    \mor{1}{3}{$q$}
    \mor{2}{4}{$p_{S_0}$}
    \mor{3}{4}{$t$}
    \enddc
\end{equation}
implies that
\begin{equation}
    \Ec\boxtimes \Oc_S=a_*(p_{X_0}^*\Ec\otimes p_{S_0}^*t_*\Oc_S)\underbracket{\simeq}_{\textup{der. prop. base change}} a_*(p_{X_0}^*\Ec\otimes (id\times t)_*\Oc_{X_0})
\end{equation}
\begin{equation*}
    \underbracket{\simeq}_{\textup{proj. form.}}a_*(id\times t)_*((id\times t)^*p_{X_0}^*\Ec)\simeq \Ec
\end{equation*}
Finally, we can consider $\Lc^{\vee}_S$ endowed with the trivial $\Lc^{\vee}_S\xrightarrow{0}\Oc_S$ dg module structure, namely $t_*\Lc^{\vee}_S$. Then, for any $\Ec \in \cohbp{X_0}{X}$, we have
\begin{equation}
    \Ec \boxtimes t_*\Lc^{\vee}_S=a_*(p_{X_0}^*\Ec\otimes p_{S_0}^*t_*\Lc_S^{\vee})\underbracket{\simeq}_{\text{der. prop. base change}} a_*(p_{X_0}^*\Ec\otimes (id\times t)_*q^*\Lc_S^{\vee})
\end{equation}
Notice that $q^*\Lc^{\vee}_S\simeq \Lc^{\vee}_X\otimes_{\Oc_X}\Oc_{X_0}$, thus we can continue
\begin{equation}
    \simeq  a_*(p_{X_0}^*\Ec\otimes (id\times t)_*(\Lc^{\vee}_X\otimes_{\Oc_X}\Oc_{X_0}))\simeq \Ec\otimes_{\Oc_{X_0}}\Oc_{X_0}\otimes_{\Oc_X}\Lc^{\vee}_X\simeq \Ec \otimes_{\Oc_X}\Lc_X^{\vee}
\end{equation}
\end{rmk}
We will construct, using $\cohbp{\bullet}{\bullet}^{\otimes}$, a lax monoidal $\oo$-functor $\sing{\bullet,\bullet}^{\otimes}:\textup{LG}_{(S,\Lc_S)}^{\textup{op},\boxplus}\rightarrow \dgcatmo$, following \cite{pr11} and \cite{brtv}. 
\begin{construction}
As in Lemma \ref{RHom(A,A)=Sym_A(L[-2])}, consider the strict cdga $Sym_{\Oc_S}(\Lc_S[-2])$ as a commutative algebra object in $\qcoh{S}$. Consider the $Sym_{\Oc_S}(\Lc_S[-2])$-algebra
\begin{equation}
    \Rc:=Sym_{\Oc_S}(\Lc_S[-2])[\nu^{-1}]
\end{equation}
where $\nu$ is the generator in degree $2$. More explicitly, 
\begin{equation}
    \Rc= \dots \rightarrow \underbracket{\Lc_S^{\vee,\otimes 2}}_{-4} \rightarrow \underbracket{\Lc_S^{\vee}}_{-2}\rightarrow 0\rightarrow \underbracket{\Oc_S}_{0}\rightarrow 0 \rightarrow \underbracket{\Lc_S}_{2}\rightarrow 0 \rightarrow \underbracket{\Lc_S^{\otimes 2}}_{4} \rightarrow \dots
\end{equation}
Then we have a symmetric monoidal $\oo$-functor 
\begin{equation}\label{invert u}
\begindc{\commdiag}[18]
\obj(-85,0)[1]{$\md_{\perf{Sym_{\Oc_S}(\Lc_S[-2])}}(\dgcatm)^{\otimes}$}
\obj(85,0)[2]{$ \md_{\perf{\Rc}}(\dgcatm)^{\otimes}$}
\mor{1}{2}{$-\otimes_{\perf{Sym_{\Oc_S}(\Lc_S[-2])}}\perf{\Rc}$}   
\enddc
\end{equation}
Composing it with $(\ref{cohbp monoidal modules 2})$ we obtain a lax monoidal $\oo$-functor
\begin{equation}\label{cohbp otimes perf R[u,u^{-1}]}
    \tlgm{S}^{\textup{op},\boxplus}\rightarrow \md_{\perf{\Rc}}(\dgcatm)^{\otimes}
\end{equation}
which, at the level of objects, is defined by
\begin{equation*}
    (X,s_X)\mapsto \cohbp{X_0}{X}\otimes_{\perf{Sym_{\Oc_S}(\Lc[-2])}}\perf{\Rc}
\end{equation*}
\end{construction}
\begin{rmk}
Let $\mathcal{U}=\{ U_i= Spec(A_i) \}$ be a Zariski affine covering of $S$, such that $\Lc_{S|U_i}$ is equivalent to $A_i$. Then the restriction of $\Rc$ to $U_i$ is equivalent to the dg algebra $A_i[u,u^{-1}]$, where $u$ sits in cohomological degree $2$.
\end{rmk}
\begin{lmm} Let $\mathcal{U}=\{U_i=Spec(B_i)\}$ be an affine open covering which trivializes $\Lc_S$ and $\Lc_S^{\vee}$.

Under the equivalence $(\ref{equivalence cohb(S_0) perf(R[u])})$, the full sub-category $\perf{S_0}$ of $\cohb{S_0}$ corresponds to the full sub-category  $\perf{Sym_{\Oc_S}(\Lc_S[-2])}^{\textup{tors}}$ of $\perf{Sym_{\Oc_S}(\Lc_S[-2])}$ spanned by perfect $Sym_{\Oc_S}(\Lc_S[-2])$-modules whose restriction to each $U_i$ is a $u$-torsion perfect $B_i[u]$-dg module.
\end{lmm}
\begin{proof}
Let $\Mc\in \perf{S_0}$. By definitions this means that, for every Zariski open $f:Spec(A)\rightarrow S_0$, the pullback $f^*\Mc\in \perf{A}$. Consider our covering $\mathcal{U}$. Then, for every $i\in I$, $\Mc_{|U_i}$ is a perfect $B_i$-module. By the same argument given above, we get equivalences $\cohb{U_{i,0}}^{\otimes}\simeq \perf{B_i[u]}^{\otimes}$ and the argument provided in \cite[Proposition 2.43]{brtv} guarantees that $\perf{U_{i,0}}$ corresponds to $\perf{B_i[u]}^{u-\textup{tors}}$ under this equivalence. This shows that every perfect complex is locally of $u$-torsion. The converse follows by the local characterization of perfect complexes and by the characterization given in \cite[Proposition 2.43]{brtv}.
\end{proof}
\begin{lmm}
The dg quotient $\cohb{S_0}\rightarrow \sing{S_0}$ corresponds to the base-change 
\begin{equation}
  -\otimes_{Sym_{\Oc_S}(\Lc_S[-2])}\Rc : \perf{Sym_{\Oc_S}(\Lc_S[-2])}\rightarrow \perf{\Rc}
\end{equation}
under the equivalence $(\ref{equivalence cohb(S_0) perf(R[u])})$
\end{lmm}
\begin{proof}
This follows by the fact that the categories involved satisfy Zariski descent, by the previous lemmas and by \cite[Proposition 2.43]{brtv}.
\end{proof}
\begin{lmm}
There is an equivalence of dg categories
\begin{equation}
    \cohbp{X_0}{X}\otimes_{\perf{Sym_{\Oc_S}(\Lc_S[-2])}}\perf{Sym_{\Oc_S}(\Lc_S[-2])}^{\textup{tors}}\simeq \perf{X_0}
\end{equation}
\end{lmm}
\begin{proof}
By definition, $\cohbp{X_0}{X}\otimes_{\perf{Sym_{\Oc_S}(\Lc_S[-2])}}\perf{Sym_{\Oc_S}(\Lc_S[-2])}^{\textup{tors}}$ is the subcategory spanned by locally $u$-torsion modules in $\cohbp{X_0}{X}$. We will show that they coincide with the subcategory of perfect complexes, using \cite[Proposition 2.45]{brtv}. Suppose that $\Mc$ is in $\perf{X_0}$. Then, for every affine open covering $j:Spec(A)\rightarrow X_0$ which trivializes $\Lc_{X_0}$, $j^*\Mc$ is perfect, and thus it is a $u$-torsion perfect $A[u]$ module by \textit{loc. cit}. Conversely, if $\Mc$ is locally $u$-torsion, each restriction $j^*\Mc$ is $u$-torsion, i.e. perfect. This proves the lemma.
\end{proof}
\begin{cor}
Let $(X,s_X)$ be a twisted LG model over $(S,\Lc_S)$. Then the exact sequence in $\dgcatm$
\begin{equation}
    \perf{X_0}\rightarrow \cohbp{X_0}{X}\rightarrow \sing{X,s_X}
\end{equation}
is equivalent to
\begin{equation}
    \cohbp{X_0}{X}\otimes_{\perf{Sym_{\Oc_S}(\Lc_S[-2])}}\perf{Sym_{\Oc_S}(\Lc_S[-2])}^{\textup{tors}}\rightarrow \cohbp{X_0}{X}\rightarrow
    \end{equation}
\begin{equation*}
    \rightarrow \cohbp{X_0}{X}\otimes_{\perf{Sym_{\Oc_S}(\Lc_S[-2])}}\perf{\Rc}
\end{equation*}
\end{cor}
\begin{proof}
This is an obvious consequence of the previous lemmas and of the fact that the $\oo$-functor
\[
T\mapsto \cohbp{X_0}{X}\otimes_{\perf{Sym_{\Oc_S}(\Lc_S[-2])}}T
\]preserves exact sequences in $\md_{\perf{Sym_{\Oc_S}(\Lc_S[-2])}}(\dgcatm)$.
\end{proof}
\begin{construction}
Consider the following lax monoidal $\oo$-functor
\begin{equation}
\begindc{\commdiag}[18]
\obj(-80,10)[1]{$\tlgm{S}^{\textup{op},\boxplus}$}
\obj(80,10)[2]{$\md_{\perf{Sym_{\Oc_S}(\Lc_S[-2])}}(\dgcatm)^{\otimes}$}
\obj(80,-10)[3]{$\md_{\perf{\Rc}}(\dgcatm)^{\otimes}$}

\mor{1}{2}{$\cohbp{\bullet}{\bullet}^{\otimes}$}
\mor{2}{3}{$-\otimes_{Sym_{\Oc_S}(\Lc_S[-2])}\Rc$}
\enddc    
\end{equation}
The previous corollary means that its underlying $\oo$-functor, composed with the forgetful functor
\begin{equation}
    \md_{\perf{\Rc}}(\dgcatm)\rightarrow \dgcatm
\end{equation}
is defined on objects by the assignment
\begin{equation*}
    (X,s)\mapsto \sing{X,s}
\end{equation*}
\end{construction}
\subsection{\textsc{The motivic realization of} Sing(X,s)}
Recall from Section \ref{the bridge between motives and non commutative motives} that there is a motivic realization lax monoidal $\oo$-functor 
\begin{equation}
    \Mv_S:\dgcatmo\rightarrow \md_{\bu_S}(\sh)^{\otimes}
\end{equation}
with the following properties
\begin{itemize}
    \item $\Mv_S$ commutes with filtered colimits
    \item For every dg category $T$, $\Mv_S(T)$ is the spectrum
    \begin{equation*}
        Y\in \sm \mapsto \textup{HK}(\perf{Y}\otimes_S T)
    \end{equation*}
    \item $\Mv_S$ sends exact sequences of dg categories to fiber-cofiber sequences in $\md_{\bu_S}(\shm)$ 
\end{itemize}
Our main scope in this section is to study the motivic realization of the dg category $\sing{X,s_X}$ associated to $(X,s_X)\in \tlgm{S}$, under the assumption that $X$ is regular.
The first important fact is the following one
\begin{prop}{\textup{\cite[Proposition 3.13]{brtv}}}
Let $p:X\rightarrow S$ be in $\sch$. Then $\Mc^{\vee}_X(\perf{X})\simeq \bu_X$.
\end{prop}
\begin{proof}
Consider the construction given in Section \ref{the bridge between motives and non commutative motives} with $S$ replaced by $X$:
\begin{equation}
    \textup{\textbf{Sm}}^{\times}_X \xrightarrow{\perf{\bullet}} \textup{\textbf{dgCat}}^{\textup{idm,op},\otimes}\xrightarrow{\iota} \textup{\textbf{SH}}^{\textup{nc,op},\otimes}_X \xrightarrow{\rhom(-,1_X^{nc})} \textup{\textbf{SH}}^{\textup{nc},\otimes}_X \xrightarrow{\Mc_X} \textup{\textbf{SH}}_X^{\otimes}
\end{equation}
We wish to show that $\Mc_X (\rhom(\iota \circ \perf{X},1^{nc}_X))\simeq \bu_X$. Notice that since the $\oo$-functor above is right lax monoindal and $\perf{X}$ is a commutative algebra in $\textup{\textbf{dgCat}}^{\textup{idm,op},\otimes}$, there is a canonical morphism of commutative algebras 
\begin{equation}
    \bu_X\rightarrow \Mc_X (\rhom(\iota \circ \perf{X},1^{nc}_X))
\end{equation} 
It suffices to show that they represent the same $\oo$-functor. Let $\mathcal{Y}$ be an object in $\textup{\textbf{SH}}_X^{\otimes}$. Then
\begin{equation}
    \Map_{\textup{\textbf{SH}}_X}(\mathcal{Y},\Mc_X (\rhom(\iota \circ \perf{X},1^{nc}_X)))\simeq \Map_{\textup{\textbf{SH}}^{\textup{nc}}_X}(\textup{R}_{\textup{Perf}}(\mathcal{Y}),\rhom(\iota \circ \perf{X},1^{nc}_X))
\end{equation}
\begin{equation*}
    \simeq \Map_{\textup{\textbf{SH}}^{\textup{nc}}_X}(\textup{R}_{\textup{Perf}}(\mathcal{Y})\otimes \textup{R}_{\textup{Perf}}(\Sigma^{\oo}_+X),1^{nc}_X)\simeq \Map_{\textup{\textbf{SH}}_X^{nc}}(\textup{R}_{\textup{Perf}}(\mathcal{Y}\otimes \Sigma^{\oo}_+X),1^{nc}_X)
\end{equation*}
\begin{equation*}
    \simeq \Map_{\textup{\textbf{SH}}_X}(\textup{R}_{\textup{Perf}}(\mathcal{Y}),1_X^{nc})\simeq \Map_{\textup{\textbf{SH}}_X}(\mathcal{Y},\Mc_X(1_X^{nc}))\simeq \Map_{\textup{\textbf{SH}}_X}(\mathcal{Y},\bu_X)
\end{equation*}
where we have used that $\iota \circ \perf{\bullet}\simeq \textup{R}_{\textup{Perf}}\circ \Sigma^{\oo}_+$, which are all symmetric monoidal functors, $\Sigma^{\oo}_+(X)$ is the unit in $\textup{\textbf{SH}}_X^{\otimes}$, that $\textup{R}_{\textup{Perf}}$ is left adjoint to $\Mc_X$ and Robalo's result $\Mc_X(1^{nc}_X)\simeq \bu_X$ (\cite[Theorem 1.8]{ro15}). The fact that the equivalence holds in $\textup{CAlg}(\textup{\textbf{SH}}_X^{\otimes})$ follows immediately from the conservativity of the forgetful functor
\begin{equation}
   \textup{CAlg}(\textup{\textbf{SH}}_X^{\otimes})\rightarrow \textup{\textbf{SH}}_X^{\otimes}
\end{equation}
\end{proof}
\begin{rmk}\label{multiplication by perfect complex in HK}
Notice that since the equivalence $\bu_X\rightarrow \Mc^{\vee}_X(\perf{X})$ holds in $\textup{CAlg}(\textup{\textbf{SH}}_X)$, if we consider a perfect complex $\Ec$ scheme $X$, then the image of
\begin{equation}
    \perf{X}\xrightarrow{-\otimes \Ec}\perf{X}
\end{equation}
along $\Mc^{\vee}_X$ coincides with multiplication by the class $[\Ec]\in \textup{HK}_0(X)$ in the commutative algebra $\bu_X$, which we denote $m_{\Ec}$. In a formula : $\Mc^{\vee}_X(-\otimes \Ec)\simeq m_{\Ec}\in \Map_{\bu_X}(\bu_X,\bu_X)$.
Moreover, for any exact triangle $\Ec'\rightarrow \Ec \rightarrow \Ec''$ we get $[\Ec]=[\Ec']+[\Ec'']\in \textup{HK}_0(X)$. Then $m_{\Ec}=m_{\Ec'}+m_{\Ec''}$. In particular, $m_{\Ec[1]}=-m_{\Ec}$.
\end{rmk}
The next step will then to understand the motivic realization of the category of coherent bounded complexes $\cohb{Z}$ of an $S$-scheme $Z$, at least when it is possible to regard it as a closed subscheme of a regular $S$-scheme $X$.
\begin{prop}{\textup{\cite[Proposition 3.17]{brtv}}}\label{motivic realization cohb}
Let $p:X\rightarrow S$ be a regular $S$-scheme of finite type. Consider a closed subscheme $i: Z\rightarrow X$ and label $j:U\rightarrow X$ the embedding of the complementary open subscheme. Then there is an equivalence 
\begin{equation}
    \Mc_X(\cohb{Z})\simeq i_*i^!\bu_X
\end{equation}
\end{prop}
\begin{proof}
Consider the exact sequence of dg categories 
\begin{equation}
    \cohb{X}_Z \rightarrow \cohb{X}\xrightarrow{j^*} \cohb{U}
\end{equation}
where $\cohb{X}_Z$ denotes the subcategory of objects in $\cohb{X}$ whose support is in $Z$. The regularity hypothesis imposed on $X$ implies that $\cohb{X}\simeq \perf{X}$ and $\cohb{U}\simeq \perf{U}$. If we apply $\Mc_X^{\vee}$ and use the previous proposition, we obtain
\begin{equation}
    \Mc^{\vee}_X(\cohb{X}_Z)\rightarrow \bu_X \xrightarrow{\Mc^{\vee}_X(j^*)} \Mc^{\vee}_X(\perf{U})
\end{equation}
which is a fiber-cofiber sequence in $\textup{\textbf{SH}}_X$ (since $\Mv_X$ sends exact triangles of dg categories to fiber-cofiber sequences). Moreover, $\Mc_X^{\vee}$ is compatible with pushforwards and $\Mc^{\vee}_X(j^*)\sim j^*$. As the spectrum $\bu$ of non-connective homotopy-invariant K-theory is compatible with pullbacks (see \cite{cd19}), the previous fiber-cofiber sequence is nothing but
\begin{equation}
    \Mc^{\vee}_X(\cohb{X}_Z)\rightarrow \bu_X \rightarrow j_*j^*\bu_X
\end{equation}
where the map on the right is induced by the unit of the adjunction $(j^*,j_*)$. In particular, we get a canonical equivalence
\begin{equation}
    \Mc^{\vee}_X(\cohb{X}_Z)\simeq i_*i^!\bu_X
\end{equation}
and therefore we are left to show that  we have an equivalence \begin{equation}
    \Mc_X^{\vee}(\cohb{Z})\simeq \Mc_X^{\vee}(\cohb{X}_Z)
\end{equation} 
Notice that there is a canonical morphism $\Mc^{\vee}_X(i_*):\Mc^{\vee}_X(\cohb{Z}\rightarrow \Mc^{\vee}_X(\cohb{X}_Z))$.The collection of objects $\Sigma^{\oo}_+Y\otimes \bu_X$, where $Y\in \textup{\textbf{Sm}}_X$, forms a family of compact generators of $\md_{\bu_X}(\textup{\textbf{SH}}_X)$. As $\Mc^{\vee}_X$ commutes with filtered colimits, it suffices to show that the morphism
\begin{equation}
    \Map_{\bu_X}(\Sigma^{\oo}_+Y\otimes \bu_X,\Mc^{\vee}_X(\cohb{Z}))\rightarrow \Map_{\bu_X}(\Sigma^{\oo}_+Y\otimes \bu_X,\Mc^{\vee}_X(\cohb{X}_Z))
\end{equation}
is an equivalence of spectra. Notice that
\begin{equation}
     \Map_{\bu_X}(\Sigma^{\oo}_+Y\otimes \bu_X,\Mc^{\vee}_X(\cohb{Z}))\simeq  \Map_{\textup{\textbf{SH}}_X}(\Sigma^{\oo}_+Y,\Mc^{\vee}_X(\cohb{Z}))
\end{equation}
\begin{equation*}
    \simeq \Map_{\textup{\textbf{SH}}^{nc}_X}(\textup{R}_{\textup{Perf}}\circ \Sigma^{\oo}_+Y,\rhom(\cohb{Z},1^{nc}_X))\simeq \Map_{\textup{\textbf{SH}}^{nc}_X}(\perf{Y}\otimes \cohb{Z},1^{nc}_X)
\end{equation*}
By \cite{pr11}, $\perf{Y}\otimes \cohb{Z}\simeq \cohb{Y}\otimes \cohb{Z}\simeq \cohb{Y\times_X Z}$, and therefore the spectrum above coincides with $\textup{HK}(\cohb{Y\times_X Z})$, which by the homotopy-invariance of G-theory and by the theorem of the Heart (\cite[Corollary 6.4.1]{ba15}) coincides with the $G$-theory of $Y\times_X Z$. In the same way, we obtain that $\Map_{\textup{\textbf{SH}}_X}(\Sigma^{\oo}_+Y,\Mc^{\vee}_X(\cohb{X}_Z))$ coincides with the $G$-theory spectrum of the abelian category $\textup{\textbf{Coh}}(X)_Z$. The claim now follows from Quillen's dévissage.
\end{proof}
As a final observation, we remark that the assignemnt $Z\mapsto \Mc^{\vee}_S(\cohb{Z})$ is insensible to (derived) thickenings.
\begin{lmm}{\textup{\cite[Proposition 3.24]{brtv}}}
Let $Z$ be a derived scheme of finite type over $S$ and let $t:\pi_0(Z)\rightarrow Z$ be the canonical closed embedding of the underlying scheme of $Z$. Then
\begin{equation}
    \Mc^{\vee}_S(\cohb{Z})\simeq \Mc^{\vee}_S(\cohb{\pi_0(Z)})
\end{equation}
\end{lmm}
\begin{proof}
By the proof of the previous proposition, $\Mc^{\vee}_S(\cohb{Z})$ represents the spectra valued sheaf $Y\mapsto \textup{HK}(\cohb{Y\times_S Z})$. Similarly, $\Mc^{\vee}_S(\cohb{\pi_0(Z)})$ represents the spectra-valued sheaf $Y\mapsto \textup{HK}(\cohb{Y\times_S \pi_0(Z)})$. Notice that $\pi_0(Y\times_S Z)=\pi_0(Y)\times_{\pi_0(S)}\pi_0(Z)=Y\times_S \pi_0(Z)$ and that the heart of $\cohb{Y\times_S Z}$ is equivalent to the heart of $\cohb{Y\times_S \pi_0(Z)}$. Then the theorem of the Heart and the computation above allows us to conclude.
\end{proof}
\begin{rmk}
Obviously, if $Z$ is a derived scheme of finite type over an $S$-scheme of finite type $X$, the same result holds if we consider $\Mc^{\vee}_X$.
\end{rmk}
Now consider a twisted LG model $(X,s_X)$ over $(S,\Lc_S)$ and assume that $X$ is a regular scheme. As above, let $X_0\xrightarrow{\mathfrak{i}}X$ be the derived zero section of $s_X$ in $X$ and let $j:X_{\mathcal{U}}=X-X_0\rightarrow X$ be the corresponding open embedding. In this case, we have that $\sing{X,s_X}\simeq \sing{X_0}$. Consider the diagram in $\textup{\textbf{dgCat}}_X^{\textup{idm}}$ and its image in $\textup{\textbf{SH}}_X$
\begin{equation}
    \begindc{\commdiag}[13]
    \obj(-80,30)[1]{$\perf{X_0}$}
    \obj(-30,30)[2]{$\cohb{X_0}$}
    \obj(20,30)[3]{$\sing{X_0}$}
    \obj(-30,0)[4]{$\perf{X}$}
    \obj(-30,-30)[5]{$\perf{X_{\mathcal{U}}}$}
    \mor{1}{2}{$$}[\atright,\injectionarrow]
    \mor{2}{3}{$$}
    \mor{2}{4}{$\mathfrak{i}_*$}
    \mor{4}{5}{$j^*$}
    \mor{1}{4}{$\mathfrak{i}_*$}[\atright,\solidarrow]
    \obj(80,30)[1a]{$\Mc^{\vee}_X(\perf{X_0})$}
    \obj(150,30)[2a]{$\Mc^{\vee}_X(\cohb{X_0})$}
    \obj(220,30)[3a]{$\Mc^{\vee}_X(\sing{X_0})$}
    \obj(150,0)[4a]{$\Mc^{\vee}_X(\perf{X})$}
    \obj(150,-30)[5a]{$\Mc^{\vee}_X(\perf{X_{\mathcal{U}}})$}
    \mor{1a}{2a}{$$}
    \mor{2a}{3a}{$$}
    \mor{2a}{4a}{$\Mc^{\vee}_X(\mathfrak{i}_*)$}
    \mor{4a}{5a}{$j^*$}
    \mor{1a}{4a}{$\Mc^{\vee}(\mathfrak{i}_*)$}[\atright,\solidarrow]
    \obj(50,0)[6]{$\rightsquigarrow$}
    \obj(50,-10)[7]{$\Mv_X$}
    \enddc
\end{equation}
By the previous results, the diagram on the right can be rewritten and completed as follows
\begin{equation}\label{fundamental diagram motivic realization sing}
    \begindc{\commdiag}[18]
    \obj(80,30)[1a]{$\Mc^{\vee}_X(\perf{X_0})$}
    \obj(160,30)[2a]{$i_*i^!\bu_X$}
    \obj(240,30)[3a]{$\Mc^{\vee}_X(\sing{X_0})$}
    \obj(160,0)[4a]{$\bu_X$}
    \obj(160,-30)[5a]{$j_*j^*\bu_X$}
    \obj(80,0)[6a]{$j_!j^*\bu_X$}
    \obj(240,0)[7a]{$i_*i^*\bu_X$}
    \mor{1a}{2a}{$$}
    \mor{2a}{3a}{$$}
    \mor{2a}{4a}{$\Mc^{\vee}_X(\mathfrak{i}_*)$}
    \mor{4a}{5a}{$\text{counit }(j^*,j_*)$}
    \mor{1a}{4a}{$\Mc^{\vee}(\mathfrak{i}_*)$}
    \mor{6a}{4a}{}
    \mor{4a}{7a}{}
    \mor{6a}{5a}{$\alpha$}
    \enddc
\end{equation}
where $i=\mathfrak{i}\circ t: \pi_0(X_0)\rightarrow X$ is the closed embedding of the underlying scheme of $X_0$ in $X$.
\begin{rmk}
The morphism $\bu_X\rightarrow i_*i^*\bu_X$ can be factored as
\begin{equation}
    \bu_X\simeq \Mc^{\vee}_X(\perf{X})\xrightarrow{\Mc^{\vee}_X(\mathfrak{i}^*)} \Mc^{\vee}_X(\perf{X_0}) \xrightarrow{\Mc^{\vee}_X(t^*)} \Mc^{\vee}_X(\perf{\pi_0(X_0)})\simeq i_*i^*\bu_X
\end{equation}
\end{rmk}
\begin{lmm}
The endomorphism $\Mc^{\vee}_X(\mathfrak{i}^*)\circ \Mc^{\vee}_X(\mathfrak{i}_*):\Mc^{\vee}_X(\perf{X_0})\rightarrow \Mc^{\vee}_X(\perf{X_0})$ is homotopic to $1-m_{\Lc^{\vee}_{X_0}}$, where $m_{\Lc^{\vee}_{X_0}}$ is the autoequivalence induced by $-\otimes \Lc^{\vee}_{X_0}:\perf{X_0}\rightarrow \perf{X_0}$.
\end{lmm}
\begin{proof}
Indeed, by Remark \ref{actions}, the endomorphism $\mathfrak{i}^*\mathfrak{i}_*$ is equivalent to $-\boxtimes (\Lc^{\vee}_S\xrightarrow{0}\Oc_S)$, the identity is equivalent to $-\boxtimes \Oc_S$ and $-\boxtimes \Lc^{\vee}_S$ corresponds to $-\otimes \Lc^{\vee}_{X_0}:\perf{X_0}\rightarrow \perf{X_0}$. Then, considering the cofiber sequence of dg functors
\begin{equation}
    \begindc{\commdiag}[15]
    \obj(-30,15)[1]{$-\otimes_{\Oc_S}\Lc_S^{\vee}$}
    \obj(30,15)[2]{$-\otimes_{\Oc_S}\Oc_S \simeq id$}
    \obj(-30,-15)[3]{$0$}
    \obj(30,-15)[4]{$\mathfrak{i}^*\mathfrak{i}_*$}
    \mor{1}{2}{$0$}
    \mor{1}{3}{$$}
    \mor{2}{4}{$$}
    \mor{3}{4}{$$}
    \enddc
\end{equation}

we see that $\mathfrak{i}^*\mathfrak{i}_*$ is equivalent to the dg functor
\begin{equation}
    \perf{X_0}\rightarrow \perf{X_0}
\end{equation}
\begin{equation*}
    \Ec\mapsto \Ec \oplus \Ec\otimes_{\Oc_{X_0}}\Lc^{\vee}_{X_0}[1]
\end{equation*}
By Remark \ref{multiplication by perfect complex in HK} we conclude that
\begin{equation}
    \Mv_X(\mathfrak{i}^*\mathfrak{i}_*)\simeq m_{\Oc_{X_0}\oplus \Lc^{\vee}_{X_0}[1]}\simeq 1-m_{\Lc^{\vee}_{X_0}}
\end{equation}
\end{proof}
\begin{cor}\label{fundamental fiber cofiber sing}
Let $(X,s_X)\in \tlgm{S}$ and assume that $X$ is regular. Then there is a fiber-cofiber sequence in $\md_{\bu_X}(\textup{\textbf{SH}}_X)$
\begin{equation}
\begindc{\commdiag}[15]
   \obj(-70,15)[1]{$\Mc^{\vee}_X(\sing{X,s_X})$}
   \obj(80,15)[2]{$cofiber(\Mc^{\vee}_X(t^*)\circ (1-m_{\Lc^{\vee}_{X_0}}):\Mc^{\vee}_X(\perf{X_0})\rightarrow i_*i^*\bu_X)$}
   \obj(80,-15)[3]{$cofiber(\alpha:j_!j^*\bu_X\rightarrow j_*j^*\bu_X)$}
   \mor{1}{2}{$$}
   \mor{2}{3}{$$}
\enddc
\end{equation}
In particular, if we apply the $\oo$-functor $i^*$, we get the following fiber-cofiber sequence:
\begin{equation}
    i^*\Mc^{\vee}_X(\sing{X,s_X})\rightarrow i^*cofiber(\Mc^{\vee}_X(\perf{X_0})\xrightarrow{\Mc^{\vee}_X(t^*)\circ (1-m_{\Lc^{\vee}_{X_0}})} i_*i^*\bu_X)\rightarrow i^*j_*j^*\bu_X
\end{equation}
\end{cor}
\begin{proof}
The second statement is an immediate consequence of the first and of the equivalence $i^*j_!\simeq 0$. The first fiber-cofiber sequence can be obtained by applying the octahedreon axiom to the triangle
\begin{equation*}
    \begindc{\commdiag}[15]
    \obj(-30,15)[1]{$\Mv_X(\perf{X_0})$}
    \obj(30,15)[2]{$i_*i^!\bu_X$}
    \obj(30,-15)[3]{$i_*i^*\bu_X$}
    \mor{1}{2}{$$}
    \mor{2}{3}{$$}
    \mor{1}{3}{$\Mv_X(i^*)\circ \Mv_X(\mathfrak{i}_*)$}[\atright,\solidarrow]
    \enddc
\end{equation*}
which appears in diagram (\ref{fundamental diagram motivic realization sing}) and by the fact that
\begin{equation*}
    cofib(\alpha)\simeq cofib(i_*i^!\bu_X\rightarrow i_*i^*\bu_X)
\end{equation*}
(see Lemma \ref{general lemma stable category}).
\end{proof}
Let $\Mv_{S,\Q}$ denote the composition
\begin{equation}
    \dgcatmo \rightarrow \md_{\bu_{S}}(\sh)^{\otimes}\xrightarrow{-\otimes \textup{H}\Q} \md_{\bu_{S,\Q}}(\sh)^{\otimes}
\end{equation}
We shall now study the commutative algebra object 
\begin{equation}
\Mc^{\vee}_{S,\Q}(\sing{S,0})\underbracket{\simeq}_{S \text{ regular}}\Mc^{\vee}_{S,\Q}(\sing{S_0})
\end{equation}
This is a particularly important object for our purposes as, for every $(X,s_X)\in \tlgm{S}$, the $\Q$-linear motivic realization of $\sing{X,s_X}$ lies in $\md_{\Mc^{\vee}_{S,\Q}(\sing{S_0})}(\sh)$.
\begin{lmm}\label{motivic realization sing(S,0)}
There are equivalences in $\textup{CAlg}(\sh)$:
\begin{equation}
    \Mc^{\vee}_{S,\Q}(\sing{S,0})\simeq cofiber(\bu_{S,\Q}\xrightarrow{1-m_{\Lc^{\vee}_S}}\bu_{S,\Q})\simeq i_0^*j_{0*}\bu_{\mathcal{U},\Q}
\end{equation}
Here $i_0:S\rightarrow V(\Lc_S)$ is the zero section, $j_0:\mathcal{U}=V(\Lc_S)-S\rightarrow V(\Lc_S)$ the open complementary and 
\begin{equation}
    cofiber(\bu_{S,\Q}\xrightarrow{1-m_{\Lc^{\vee}_S}}\bu_{S,\Q}):=\bu_{S,\Q}\otimes_{\textup{H}\Q[t]}\textup{H}\Q \in \textup{CAlg}(\sh)
\end{equation}
where we consider the morphisms $\textup{H}\Q[t]\rightarrow \textup{H}\Q$ induced by $t\mapsto 0$ and $\textup{H}\Q[t]\rightarrow \bu_{S,\Q}$ induced by $t\mapsto 1-m_{\Lc_S^{\vee}}$.
\end{lmm}
\begin{proof}
Notice that the underlying object of the commutative algebra $\bu_{S,\Q}\otimes_{\textup{H}\Q[t]}\textup{H}\Q$ is $cofiber(\bu_{S,\Q}\xrightarrow{1-m_{\Lc^{\vee}_S}}\bu_{S,\Q})$. This follows easily from the fact that 
\begin{equation}
    \textup{H}\Q[t]\xrightarrow{t}\textup{H}\Q[t]
\end{equation}
is a free resolution of $\textup{H}\Q$.
Consider the diagram
\begin{equation}
    \begindc{\commdiag}[20]
    \obj(-40,10)[1]{$S_0$}
    \obj(0,10)[2]{$S$}
    \obj(50,10)[3]{$\emptyset$}
    \obj(-40,-10)[4]{$S$}
    \obj(0,-10)[5]{$V:=V(\Lc_S)$}
    \obj(50,-10)[6]{$\mathcal{U}=V-S$}
    \mor{1}{2}{$\mathfrak{i}$}
    \mor{3}{2}{$j$}
    \mor{4}{5}{$i_0$}
    \mor{6}{5}{$j_0$}
    \mor{1}{4}{$$}
    \mor{2}{5}{$i_0$}
    \mor{3}{6}{$$}
    \enddc
\end{equation}
where both squares are (homotopy) cartesian. Notice that in this case $\pi_0(S_0)=S$ and therefore $i=\mathfrak{i}\circ t=id$. Then, the fiber-cofiber sequence of the previous corollary gives us an equivalence
\begin{equation}
    \Mc^{\vee}_{S,\Q}(\sing{S_0})\simeq cofiber(\Mc^{\vee}_{S,\Q}(\perf{S_0}\xrightarrow{\Mc^{\vee}_{S,\Q}(t^*)\circ (1-m_{\Lc^{\vee}_S})}\bu_{S,\Q}))
\end{equation}
Since the square
\begin{equation}
    \begindc{\commdiag}[18]
    \obj(-40,15)[1]{$\Mc^{\vee}_{S,\Q}(\perf{S_0})$}
    \obj(40,15)[2]{$\Mc^{\vee}_{S,\Q}(\perf{S_0})$}
    \obj(-40,-15)[3]{$\Mc^{\vee}_{S,\Q}(\perf{S})$}
    \obj(40,-15)[4]{$\Mc^{\vee}_{S,\Q}(\perf{S})$}
    \mor{1}{2}{$1-m_{\Lc^{\vee}_{S_0}}$}
    \mor{1}{3}{$\Mc^{\vee}_S(t^*)$}
    \mor{2}{4}{$\Mc^{\vee}_S(t^*)$}
    \mor{3}{4}{$1-m_{\Lc^{\vee}_{S}}$}
    \enddc
\end{equation}
is commutative up to coherent homotopy, to get the first desired equivalence
\begin{equation}
    \Mv_{S,\Q}(\sing{S,0})\simeq \Mv_{S,\Q}(\sing{S_0})\simeq cofiber(\Mc^{\vee}_{S,\Q}(\perf{S_0}\xrightarrow{\Mc^{\vee}_{S,\Q}(t^*)\circ (1-m_{\Lc^{\vee}_{S}})}\bu_{S,\Q}))
\end{equation}
it suffices to show that $\Mc^{\vee}_{S,\Q}(t^*)$ is an equivalence. This is true as $\mathfrak{i}\circ t=id$ and $t\circ \mathfrak{i}$ is homotopic to the identity (see \cite[Remark 3.31]{brtv}).

To get the second equivalence, consider: 
\begin{equation}
    i_0^*j_{0*}\bu_{\mathcal{U},\Q}\simeq cofiber(i^*_0i_{0*}i^!_0\bu_{V,\Q}\simeq i_0^!\bu_{V,\Q}\xrightarrow{c} \bu_{S,\Q}\simeq i^*_0\bu_{V,\Q})
\end{equation}
where the map $c$ is the one induced by the counit $i_{0*}i_0^!\bu_V\rightarrow \bu_V$. Notice that $i_0:S\rightarrow V$ is a closed embedding between regular schemes. In particular, absolute purity holds. It follows from \cite[Remark 13.5.5]{cd19} that the composition
\begin{equation}
    \bu_{V,\Q}\rightarrow i_{0*}i^*_0\bu_{V,\Q}\xrightarrow{\eta_{i_0}}\bu_{V,\Q}
\end{equation}
corresponds to $1-m_{\Lc^{\vee}_S}$, as the conormal sheaf of $i_0:S\rightarrow V$ is $\Lc_{S}^{\vee}$. If we apply $i^!$ we obtain
\begin{equation}
   i^!(1-m_{\Lc^{\vee}_S}):  i^!\bu_{V,\Q}\xrightarrow{c}\bu_S\underbracket{\xrightarrow{i^!\eta_{i_0}'}}_{\simeq \text{ abs. pur.}} i^!\bu_{V,\Q} 
\end{equation}
which under the equivalence $\bu_{S,\Q}\simeq i^!\bu_{V,\Q}$ corresponds to $1-m_{\Lc^{\vee}_S}$.
\end{proof}
\begin{construction}
Let $(X,s_X)\in \tlgm{S}$. Consider the morphism
\begin{equation}
    \bu_X\xrightarrow{1-m_{\Lc^{\vee}_X}} \bu_X
\end{equation} 
Since $j:X_{\mathcal{U}}\rightarrow X$ is the complementary open subscheme to the zero section of $V(\Lc^{\vee}_X)$, it follows that $j^*\Lc^{\vee}_X\simeq \Oc_{X_{\mathcal{U}}}$. In particular, 
\begin{equation}
j^*\circ (1-m_{\Lc^{\vee}_X})\simeq 1-m_{\Lc^{\vee}_{X_{\mathcal{U}}}}\sim 0
\end{equation}
Then we obtain a morphism
\begin{equation}
   sp^{\textup{mot}}_{X}: cofiber(\bu_X\xrightarrow{1-m_{\Lc^{\vee}_X}}\bu_X)\rightarrow j_*j^*\bu_X
\end{equation}
\end{construction}
\begin{prop}
Let $(X,s_X)\in \tlgm{S}$ and assume that $X$ is regular and that $X\xrightarrow{s_X}V(\Lc_S)\xleftarrow{\text{zero}}S$ is Tor-independent (i.e. $ X\times_{V(\Lc_S)}S\simeq X\times^h_{V(\Lc_S)}S$). Then
\begin{equation}
    i^*\Mc^{\vee}_X(\sing{X,s_X})\simeq fiber(i^*sp^{\textup{mot}}_X)
\end{equation}
\begin{proof}
This follows immediately from the octahedron property applied to the triangle in the following diagram
\begin{equation}
    \begindc{\commdiag}[20]
    \obj(-30,15)[1]{$\bu_{X_0}$}
    \obj(0,15)[2]{$i^!\bu_X$}
    \obj(40,15)[3]{$i^*\Mc^{\vee}_X(\sing{X_0})$}
    \obj(0,0)[4]{$\bu_{X_0}$}
    \obj(0,-15)[5]{$i^*j_*j^*\bu_X$}
    \mor{1}{2}{$$}
    \mor{2}{3}{$$}
    \mor{2}{4}{$$}
    \mor{4}{5}{$$}
    \mor{1}{4}{$1-m_{\Lc^{\vee}_{X_0}}$}[\atright,\solidarrow]
    \enddc
\end{equation}
and from the compatibility of $1-m_{\Lc^{\vee}_X}$ with pullbacks.
\end{proof}
\end{prop}
\section{\textsc{Monodromy-invariant vanishing cycles}}\label{Monodromy-invariant vanishing cycles}
\subsection{\textsc{The formalism of vanishing cycles}}\label{the formalism of vanishing cycles}
We shall begin with a quick review of the formalism of vanishing cycles.
We hope that this will be useful to understand the analogy which led to the definition of monodromy-invariant vanishing cycles.
\begin{notation}\label{notation trait}
Throughout this section, $A$ will be an excellent henselian trait and $S$ will denote the associated affine scheme. Label $k$ its residue field and $K$ its fraction field. Let $\sigma$ be $Spec(k)$ and $\eta$ be $Spec(K)$. Fix algebraic closures $k^{alg}$ and $K^{alg}$ of $k$ and $K$ respectively. We will consider
\begin{itemize}
\item the maximal separable extension $k^{sep}$ of $k$ inside $k^{alg}$ and let $\bar{\sigma}=Spec(k^{sep})$
\item the maximal unramified extension $K^{unr}$ of $K$ in $K^{alg}$ and let $\eta^{unr}=Spec(K^{unr})$
\item the maximal tamely ramified extension $K^t$ of $K$ inside $K^{alg}$ and let $\eta^t=Spec(K^t)$
\item the maximal separable extension $K^{sep}$ of $K$ inside $K^{alg}$ and let $\bar{\eta}=Spec(K^{sep})$
\end{itemize}
Moreover, we will fix an uniformizer $\pi$. 
\end{notation}
\begin{rmk}
It is well known that there is an equivalence between the category of separable extensions of $k$ and that of unramified extensions of $K$. In particular, $\textup{Gal}(k^{sep}/k)\simeq \textup{Gal}(K^{unr}/K)$. This, together with the fundamental theorem of Galois theory, implies that there is an exact sequence of groups
\begin{equation}
    1\rightarrow \textup{Gal}(K^{sep}/K^{unr})\rightarrow \textup{Gal}(K^{sep}/K)\rightarrow \textup{Gal}(k^{sep}/k)\rightarrow 1
\end{equation}
The Galois group on the left is called the inertia group and it is usually denoted by $I$. The chain of extensions $K^{unr}\subseteq K^{t}\subseteq K^{sep}$ gives us the following decomposition of $I$:
\begin{equation}
    1\rightarrow \underbracket{\textup{Gal}(K^{sep}/K^{t})}_{=:I_w} \rightarrow I \rightarrow \underbracket{\textup{Gal}(K^{sep}/K)}_{=:I_t}\rightarrow 1
\end{equation}
The Galois group $I_w$ is called wild inertia group, while $I_t$ is called tame inertia group. 
See \cite{se62} for more details.
\end{rmk}
\begin{rmk}
Let $A^{sh}$ be a strict henselianisation of $A$ and let $\bar{S}=Spec(A^{sh})$. We then have the following picture
\begin{equation}
    \begindc{\commdiag}[20]
    \obj(0,15)[1]{$\bar{\sigma}$}
    \obj(30,15)[2]{$\bar{S}$}
    \obj(60,15)[3]{$\eta^{unr}$}
    \obj(0,-5)[4]{$\sigma$}
    \obj(30,-5)[5]{$S$}
    \obj(60,-5)[6]{$\eta$}
    \obj(90,15)[7]{$\eta^t$}
    \obj(120,15)[8]{$\bar{\eta}$}
    \mor{1}{2}{$\bar{i}$}[\atright,\injectionarrow]
    \mor{3}{2}{$j_{unr}$}[\atleft,\injectionarrow]
    \mor{4}{5}{$i$}[\atright,\injectionarrow]
    \mor{6}{5}{$j$}[\atright,\injectionarrow]
    \mor{1}{4}{$v_{\sigma}$}
    \mor{2}{5}{$v$}
    \mor{3}{6}{$v_{\eta}$}
    \mor{7}{3}{$$}
    \mor{8}{7}{$$}
    \cmor((90,18)(90,21)(87,24)(60,24)(34,24)(31,21)(31,19)) \pdown(60,22){$j_t$}
    \cmor((120,18)(120,24)(117,27)(74,27)(32,27)(29,24)(29,19)) \pdown(75,32){$\bar{j}$}
    \enddc
\end{equation}
where both squares are cartesian. In particular, notice that $(\sigma,\eta)$ and $(\bar{\sigma},\eta^{unr})$ form closed-open coverings of $S$ and $\bar{S}$ respectively.
\end{rmk}
Let $p:X\rightarrow S$ be of finite type and let $\bar{X}=X\times_S \bar{S}$. Consider the following diagram, cartesian over the base (this is also called Grothendieck's trick in the literature)
\begin{equation}
    \begindc{\commdiag}[20]
    \obj(0,10)[1]{$X_{\bar{\sigma}}$}
    \obj(30,10)[2]{$\bar{X}$}
    \obj(60,10)[3]{$X_{\eta^{unr}}$}
    \obj(90,10)[4]{$X_{\eta^t}$}
    \obj(120,10)[5]{$X_{\bar{\eta}}$}
    \obj(0,-10)[6]{$\bar{\sigma}$}
    \obj(30,-10)[7]{$\bar{S}$}
    \obj(60,-10)[8]{$\eta^{unr}$}    \obj(90,-10)[9]{$\eta^t$}
    \obj(120,-10)[10]{$\bar{\eta}$}
    \mor{1}{2}{$$}[\atright,\injectionarrow]
    \mor{3}{2}{$$}[\atright,\injectionarrow]
    \mor{4}{3}{$$}
    \mor{5}{4}{$$}
    \mor{6}{7}{$$}[\atright,\injectionarrow]
    \mor{8}{7}{$$}[\atright,\injectionarrow]
    \mor{9}{8}{$$}
    \mor{10}{9}{$$}
    \mor{1}{6}{$$}
    \mor{2}{7}{$$}
    \mor{3}{8}{$$}
    \mor{4}{9}{$$}
    \mor{5}{10}{$$}
    \enddc
\end{equation}

\begin{rmk}
The $\oo$-category $\shvl(\bar{X})$ is the recollement of $\shvl(X_{\bar{\sigma}})$ and $\shvl(X_{\eta^{unr}})$ in the sense of \cite[Appendix A.8]{ha}\footnote{Strictly speaking, one should consider the full subcategories of $\shvl(\bar{X})$ spanned by the images of $\bar{i}_*$ and $j_{\eta^{unr}}$ which are closed under equivalences.}. Indeed, $\shvl(\bar{X})$ is a stable $\oo$-category and therefore it admits finite limits. Moreover, both $\bar{i}_*:\shvl(X_{\bar{\sigma}})\rightarrow \shvl(\bar{X})$ and $j^{unr}_*: \shvl(X_{\eta^{unr}})\rightarrow \shvl(\bar{X})$ are fully faithful $\oo$-functors: these are immediate consequences of the proper and smooth base change applied to the cartesian squares
\begin{equation}
    \begindc{\commdiag}[15]
    \obj(-15,15)[1]{$X_{\bar{\sigma}}$}
    \obj(15,15)[2]{$X_{\bar{\sigma}}$}
    \obj(-15,-15)[3]{$X_{\bar{\sigma}}$}
    \obj(15,-15)[4]{$\bar{X}$}
    \mor{1}{2}{$id$}
    \mor{1}{3}{$id$}
    \mor{2}{4}{$\bar{i}$}
    \mor{3}{4}{$\bar{i}$}
    \obj(45,15)[5]{$X_{\bar{\sigma}}$}
    \obj(75,15)[6]{$X_{\bar{\sigma}}$}
    \obj(45,-15)[7]{$X_{\bar{\sigma}}$}
    \obj(75,-15)[8]{$\bar{X}$}
    \mor{5}{6}{$id$}
    \mor{5}{7}{$id$}
    \mor{6}{8}{$j_{\eta^{unr}}$}
    \mor{7}{8}{$j_{\eta^{unr}}$}
    \enddc
\end{equation}
It is known that the conditions of \cite[Definition A.8.1]{ha} are verified. In particular, using the equivalence $\shvl(X_{\eta^{unr}})\simeq \shvl(X_{\bar{\eta}})^I$ that we have discussed in Section \ref{l-adic sheaves}, we have that $\shvl(\bar{X})$ is the recollement of $\shvl(X_{\bar{\sigma}})$ and $\shvl(X_{\bar{\eta}})^{I}$.
\end{rmk}
\begin{rmk}
Notice that we have the following diagram where each arrow is an equivalence:
\begin{equation}
    \begindc{\commdiag}[15]
    \obj(-30,20)[1]{$\shvl(X_{\eta^{unr}})$}
    \obj(0,-20)[2]{$\shvl(X_{\eta^t})^{I_t}$}
    \obj(30,20)[3]{$\shvl(X_{\bar{\eta}})^{I}$}
    \mor{1}{2}{$j_{t}^*$}
    \mor{1}{3}{$\bar{j}^*$}
    \mor{2}{3}{$j_w^*$}
    \cmor((30,25)(30,30)(28,32)(0,32)(-28,32)(-30,30)(-30,25))
    \pdown(0,39){$\bar{j}_*(-)^{\textup{h}I}$}
    \mor(40,17)(13,-17){$j_{w*}(-)^{\textup{h}I_w}$}
    \mor(-13,-17)(-40,17){$j_{t*}(-)^{\textup{h}I_t}$}
    \enddc
\end{equation}
\end{rmk}
\begin{construction}
We will now sketch the construction of Deligne's topos in the $\oo$-categorical world. Notice that, by the recollement technique, the étale $\oo$-topos of $\bar{S}$ is the recollement of the étale $\oo$-topos of $\bar{\sigma}$ and that of $\eta^{unr}$. Moreover, notice that the étale $\oo$-topos of $\eta^{unr}$ is equivalent to the $\oo$-category of spaces with a continuous action of $I$.  For a scheme $Y$ over $\bar{\sigma}$, we consider the lax $\oo$-limit $\tilde{Y}_{et}\times_{\tilde{\bar{\sigma}}_{et}}\tilde{\bar{S}}_{et}$, which exists by \cite[pag. 614]{htt}. The decomposition of $\tilde{\bar{S}}_{et}$ gives us a digram
\begin{equation}
    \begindc{\commdiag}[20]
    \obj(-30,15)[1]{$\tilde{Y}_{et}$}
    \obj(0,15)[2]{$\tilde{Y}_{et}\times_{\tilde{\bar{\sigma}}_{et}}\tilde{\bar{S}}_{et}$}
    \obj(45,15)[3]{$\tilde{Y}_{et}\times_{\tilde{\bar{\sigma}}_{et}}\tilde{\eta^{unr}}_{et}$}
    \obj(-30,-15)[4]{$\tilde{\bar{\sigma}}_{et}$}
    \obj(0,-15)[5]{$\tilde{\bar{S}}_{et}$}
    \obj(45,-15)[6]{$\tilde{\eta^{unr}}_{et}$}
    \mor{1}{2}{$$}
    \mor{3}{2}{$$}
    \mor{1}{4}{$$}
    \mor{2}{5}{$$}
    \mor{3}{6}{$$}
    \mor{4}{5}{$$}
    \mor{6}{5}{$$}
    \enddc
\end{equation}

We will label the $\oo$-category of $\ell$-adic sheaves on $\tilde{Y}_{et}\times_{\tilde{\bar{\sigma}}_{et}}\tilde{\eta^{unr}}_{et}$ by $\shvl(Y)^I$, as it is the $\oo$-category of $\ell$-adic sheaves on $Y$ endowed with a continuous action of $I$. 
The $\oo$-category $\shvl(Y)^{I_t}$ of $\ell$-adic sheaves on $Y$ endowed with a continuous action of $I_t$ identifies with the full subcategory of $\shvl(Y)^I$ such that the induced action of $I_w$ is trivial.
\end{construction}
\begin{defn} Let $p:\bar{X}\rightarrow \bar{S}$ be an $\bar{S}$-scheme. 
Let $\Ec\in \shvl(\bar{X}_{\eta^{unr}})$. The $\ell$-adic sheaf of tame nearby cycles of $\Ec$ is defined as
\begin{equation}
\Psi_t(\Ec):=\bar{i}^*j_{t*}\Ec_{|X_{\eta_t}} \in \shvl(X_{\bar{\sigma}})^{I_t}
\end{equation}
Analogously, the $\ell$-adic sheaf of nearby cycles of $\Ec$ is defined as
\begin{equation}
    \Psi(\Ec):=\bar{i}^*\bar{j}_*\Ec_{|X_{\bar{\eta}}}\in \shvl(X_{\bar{\sigma}})^I
\end{equation}
\end{defn}
\begin{rmk}
With the same notation as above, the following equivalence holds in $\shvl(\bar{X}_{\bar{\sigma}})^{I_t}$:
\begin{equation}
    \Psi_t(\Ec)\simeq \Psi(\Ec)^{\textup{h}I_w}
\end{equation}
\end{rmk}
\begin{defn}
Let $p:\bar{X}\rightarrow \bar{S}$ be an $\bar{S}$-scheme. Let $\Fc \in \shvl(\bar{X})$. The unit of the adjuction $(j_t^*,j_{t*})$ induces a morphism in $\shvl(\bar{X}_{\bar{\sigma}})^{I_t}$
\begin{equation}
    \bar{i}^*\Fc \rightarrow \Psi_t(\Fc_{\bar{X}_{\eta^{unr}}})
\end{equation}
when we regard the object on the left endowed with the trivial $I_t$-action. The cofiber of this morphism is by definition the $\ell$-adic sheaf of tame vanishing cycles of $\Fc$, which we will denote by $\Phi_t(\Fc)$. In a similar way, we define the $\ell$-adic sheaf of vanishing cycles $\Phi(\Fc)$ as the cofiber of the morphism (called the specialization morphism in literature)
\begin{equation}
    \bar{i}^*\Fc \rightarrow \Psi(\Fc_{\bar{X}_{\eta^{unr}}})
\end{equation}
induced by the unit of the adjunction $(\bar{j}^*,\bar{j}_*)$.
\end{defn}
\subsection{\textsc{Monodromy-invariant vanishing cycles}}
\begin{context}
Assume that $S$ is strictly henselian, i.e. that $k$ is a separably closed field and that $p:X\rightarrow S$ is a proper morphism.
\end{context}
For our purposes, we are interested in the image of $\Phi(\Ql{,X})$ via the $\oo$-functor 
\begin{equation*}
    (-)^{\textup{h}I}:\shvl(X_{\sigma})^{I}\rightarrow \shvl(X_{\sigma})
\end{equation*}
It is then important to remark that it is possible to determine $p_{\sigma*}\Phi(\Ql{,X})^{\textup{h}I}$ without ever mentioning this $\oo$-functor. Notice that $\Phi(\Ql{,X})^{\textup{h}I}$ is the cofiber of the image of the specialization morphism via $(-)^{\textup{h}I}$:
\begin{equation}\label{fiber cofiber sphI}
    \Ql{,X_0}^{\textup{h}I}\xrightarrow{(sp)^{\textup{h}I}}\Psi(\Ql{,X_{\bar{\eta}}})^{\textup{h}I}\rightarrow \Phi(\Ql{,X})^{\textup{h}I}
\end{equation}
is a fiber-cofiber sequence in $\shvl(X_{\sigma})$.
\begin{rmk}
There are equivalences in $\textup{CAlg}(\shvl(\sigma))$
\begin{itemize}
    \item $\bigl( p_{\sigma*}\Ql{,X_{\sigma}} \bigr)^{\textup{h}I}\simeq p_{\sigma*}\Ql{,X_{\sigma}}\otimes_{\Ql{,\sigma}}\Ql{,\sigma}^{\textup{h}I}$ (see \cite[Proposition 4.32]{brtv})
    \item $\bigl( p_{\sigma*}i^*\bar{j}_*\Ql{,X_{\bar{\eta}}} \bigr)^{\textup{h}I}\simeq p_{\sigma*}i^*j_*\Ql{,X_{\eta}}$ (see \cite[Proposition 4.31]{brtv})
    \item $\Ql{,\sigma}^{\textup{h}I}\simeq i_0^*j_{0*}\Ql{,\eta}$ (see \cite[Lemma 4.34]{brtv})
\end{itemize}
\end{rmk}
Consider the object $p_{\sigma*}p_{\sigma}^*i_0^*j_{0*}\Ql{,\eta}$ in $\shvl(\sigma)$. The following chain of equivalences, together with the previous remark, identify it with $\bigl(p_{\sigma*}\Ql{,X_{\sigma}} \bigr)^{\textup{h}I}$
\begin{equation}
    p_{\sigma*}p_{\sigma}^*i_0^*j_{0*}\Ql{,\eta}\simeq p_{\sigma*}\bigl(\Ql{,X_{\sigma}}\otimes_{\Ql{,X_{\sigma}}}p_{\sigma}^*i_0^*j_{0*}\Ql{,\eta}\bigr)\underbracket{\simeq}_{\text{proj.form.}}p_{\sigma*}\Ql{,X_{\sigma}}\otimes_{\Ql{,\sigma}}i_0^*j_{0*}\Ql{,\eta}
\end{equation}
\begin{defn}
Consider the canonical morphism in $\shvl(X_{\sigma})$
\begin{equation}
    sp_p^{\textup{hm}}:p_{\sigma}^*i_0^*j_{0*}\Ql{,\eta}\simeq i^*p^*j_{0*}\Ql{,\eta}\rightarrow i^*j_*p_{\eta}^*\Ql{,\eta}\simeq i^*j_*\Ql{,X_{\eta}}
\end{equation}
induced by the base change natural transformation $p^*j_{0*}\rightarrow j_*p_{\eta}^*$. We will refer to it as the monodromy invariant specialization morphism. We will refer to the cofiber of this morphism in $\shvl(X_{\sigma})$ as monodromy invariant vanishing cycles, which we will denote $\Phi_{p}^{\textup{hm}}(\Ql{})$.
\end{defn}
In order to justify the choice of the name of this morphism, we shall prove the following:
\begin{lmm}\label{comparison monodromy invariant inertia invariant vanishing cycles}
The two arrows $p_{\sigma*}(sp)^{\textup{h}I}$ and $p_{\sigma}sp_p^{\textup{hm}}$ are homotopic under the equivalences $\bigl( p_{\sigma*}i^*\bar{j}_*\Ql{,X_{\bar{\eta}}} \bigr)^{\textup{h}I}\simeq p_{\sigma*}i^*j_*\Ql{,X_{\eta}}$ and $\bigl( p_{\sigma*}\Ql{,X_{\sigma}} \bigr)^{\textup{h}I}\simeq p_{\sigma*}\Ql{,X_{\sigma}}\otimes_{\Ql{,\sigma}}\Ql{,\sigma}^{\textup{h}I}\simeq p_{\sigma*}p_{\sigma}^*i_0^*j_{0*}\Ql{,\eta}$  
\end{lmm}
\begin{proof}
We shall view both arrows as maps $p_{\sigma*}\Ql{,X_{\sigma}}\otimes_{\Ql{,\sigma}}\Ql{,\sigma}^{\textup{h}I}\rightarrow p_{\sigma*}i^*j_*\Ql{,X_{\eta}}$

Since the tensor product defines a cocartesian symmetric monoidal structure on $\textup{CAlg}(\shvl(\sigma))$ (see \cite{ha}), it suffices to show that the maps of commutative algebras obtained from $p_{\sigma*}(sp)^{\textup{h}I}$ and $p_{\sigma*}sp_p^{\textup{hm}}$ by precomposition with $p_{\sigma*}\Ql{,X_{\sigma}}\rightarrow p_{\sigma*}\Ql{,X_{\sigma}}\otimes_{\Ql{,\sigma}}\Ql{,\sigma}^{\textup{h}I}$ and $\Ql{,\sigma}^{\textup{h}I}\rightarrow p_{\sigma*}\Ql{,X_{\sigma}}\otimes_{\Ql{,\sigma}}\Ql{,\sigma}^{\textup{h}I}$ respectively are equivalent. It is essentially the proof of the main theorem in $\cite{brtv}$ that the precompositions of $p_{\sigma*}(sp)^{\textup{h}I}$ with these two maps are induced by the canonical morphisms $1\rightarrow j_*j^*$ and $1\rightarrow p_{\eta*}p_{\eta}^*$. We are left to show that the same holds true if we consider $sp_p^{\textup{hm}}$ instead of $(sp)^{\textup{h}I}$. Recall that $sp_p^{\textup{hm}}$ is induced by the base change morphism $p^*j_{0*}\rightarrow j_*p_{\eta}^*$, i.e. the morphism which corresponds under adjunction to $j_{0*}\rightarrow j_{0*}p_{\eta*}p_{\eta}^*$, induced by the unit of the adjunction $(p_{\eta}^*,p_{\eta*})$. From this we obtain the following diagram
\begin{equation}
    \begindc{\commdiag}[18]
    \obj(-60,20)[1]{$i_0^*j_{0*}\Ql{,\eta}$}
    \obj(60,20)[2]{$i_0^*j_{0*}p_{\eta*}p_{\eta}^*\Ql{,\eta}\simeq p_{\sigma*}i^*j_*p_{\eta}^*\Ql{,\eta}$}
    \obj(-60,0)[3]{$p_{\sigma*}p_{\sigma}^*i_0^*j_{0*}\Ql{,\eta}$}
    \obj(60,0)[4]{$p_{\sigma*}p_{\sigma}^*p_{\sigma*}i^*j_*p_{\eta}^*\Ql{,\eta}$}
    \obj(60,-20)[6]{$p_{\sigma*}i^*j_*p_{\eta}^*\Ql{,\eta}$}
    \mor{1}{2}{$1\rightarrow p_{\eta*}p_{\eta}^*$}
    \mor{1}{3}{$1\rightarrow p_{\sigma*}p_{\sigma}^*$}
    \mor{2}{4}{$$}
    \mor{3}{4}{$$}
    \obj(-120,0)[5]{$p_{\sigma*}p_{\sigma}^*\Ql{,\sigma}$}
    \mor{5}{3}{$1\rightarrow j_{0}^*j_{0*}$}
    \mor{3}{6}{$p_{\sigma*}sp_p^{\textup{hm}}$}[\atright,\solidarrow]
    \mor{4}{6}{$$}
    \cmor((85,16)(88,15)(90,14)(95,0)(90,-14)(88,-15)(85,-16)) \pleft(100,0){$id$}
    \cmor((-120,-5)(-70,-16)(-50,-19)(-35,-20)(-20,-20)(0,-20)(40,-20)) \pright(0,-22){$$}
    \enddc
\end{equation}
In particular the composition $i_0^*j_{0*}\Ql{,\eta}\rightarrow p_{\sigma*}p_{\sigma}^*i_0^*j_{0*}\Ql{,\eta} \xrightarrow{p_{\sigma*}sp_p^{\textup{hm}}} p_{\sigma*}i^*j_*\Ql{,X_{\eta}}$ is homotopic to the map induced by $1\rightarrow p_{\eta*}p_{\eta}^*$. The composition $p_{\sigma*}p_{\sigma}^*\Ql{,\sigma}\rightarrow p_{\sigma*}p_{\sigma}^*i_0^*j_{0*}\Ql{,\eta} \rightarrow p_{\sigma*}i^*j_*\Ql{,X_{\eta}}$ is homotopic to the map $p_{\sigma*}\Ql{,X_{\sigma}}\rightarrow p_{\sigma*}i^*j_{*}\Ql{,X_{\eta}}$ induced by $1\rightarrow j_*j^*$ by the following commutative triangle
\begin{equation}
    \begindc{\commdiag}[18]
    \obj(-25,0)[1]{$p^*\Ql{,S}$}
    \obj(25,15)[2]{$p^*j_{0*}j_0^*\Ql{,S}$}
    \obj(25,-15)[3]{$j_*p_{\eta}^*j_0^*\Ql{,S}\simeq j_*j^*p^*\Ql{,S}$}
    \mor{1}{2}{$1\rightarrow j_{0*}j_0^*$}
    \mor{2}{3}{$$}
    \mor{1}{3}{$1\rightarrow j_*j^*$}[\atright,\solidarrow]
    \enddc
\end{equation}
\end{proof}

\begin{rmk}
Notice that what we said above is actually true for any $\ell$-adic sheaf in the image of $p^*:\shvl(S)\rightarrow \shvl(X)$.
\end{rmk}
The advantage of this reformulation is that it allows to define inertia-invariant vanishing cycles without ever mention the inertia group. We will adopt it for the situation we are interested in.

Assume that $S$ is a noetherian regular scheme.
\begin{defn}\label{mododromy-invariant vanishing cycles def}
Let $(X,s_X)$ be a twisted LG model over $(S,\Lc_S)$. Consider the diagram obtained from the zero section of the vector bundle associated to $\Lc_X^{\vee}$:
\begin{equation}
  \begindc{\commdiag}[18]
    \obj(0,10)[1]{$\pi_0(X_0)$}
    \obj(30,10)[2]{$X$}
    \obj(0,-10)[3]{$X$}
    \obj(30,-10)[4]{$V(\Lc_X)$}
    \obj(80,10)[5]{$X_{\mathcal{U}}$}
    \obj(80,-10)[6]{$\mathcal{U}_X:=V(\Lc_X)-X$}
    \mor{1}{2}{$i$}
    \mor{1}{3}{$s_0$}
    \mor{2}{4}{$s_X$}
    \mor{3}{4}{$i_0$}[\atright,\solidarrow]
    \mor{5}{2}{$j$}[\atright,\solidarrow]
    \mor{5}{6}{$s_{\mathcal{U}}$}
    \mor{6}{4}{$j_0$}
    \enddc
\end{equation}
Define the monodromy invariant specialization morphism associated to $(X,s_X)$ as the map in $\shvl(\pi_0(X_0))$
\begin{equation}
   sp_{(X,s_X)}^{\textup{mi}}: s_0^*i_0^*j_{0*}\Ql{,\mathcal{U}_X}\simeq i^*s_{X}^*j_{0*}\Ql{,\mathcal{U}_X}\rightarrow i^*j_*\Ql{,X_{\mathcal{U}}}\simeq i^*j_*s_{\mathcal{U}}^*\Ql{,\mathcal{U}_X}
\end{equation}
induced by the natural transformation $s_X^*j_{0*}\rightarrow j_*s_{\mathcal{U}}^*$. We will refer to the cofiber of $sp_{(X,s_X)}^{\textup{mi}}$ as monodromy invariant vanishing cycles of $(X,s_X)$, that we will denote $\Phi^{\textup{mi}}_{(X,s_X)}(\Ql{})$.
\end{defn}
\begin{prop}\label{cofiber first chern class}
Let $(X,s_X)$ be as above and assume $X$ regular. There is an equivalence
\begin{equation}
    i_0^*j_{0*}\Ql{,\mathcal{U}_X}\simeq cofib \bigl ( c_1(\Lc_X):\Ql{,X}(-1)[-2]\rightarrow \Ql{,X} \bigr )
\end{equation}
\end{prop}
\begin{proof}
Consider the diagram
\begin{equation*}
    X\xrightarrow{i_0}V\xleftarrow{j_0}\mathcal{U}
\end{equation*}
The $\ell$-adic sheaf $i_0^*j_{0*}\Ql{,\mathcal{U}}$ is the cofiber of the morphism $i_0^*\bigl( i_{0*}i_0^{!}\Ql{,V}\xrightarrow{\text{counit}} \Ql{,V} \bigr)$ in $\shvl(X)$. Recall that $i_0^*i_{0*}\simeq id$. Therefore, we can consider the following triangle
\begin{equation}
    \begindc{\commdiag}[18]
    \obj(-20,15)[1]{$i_0^!\Ql{,V}$}
    \obj(20,15)[2]{$\Ql{,X}$}
    \obj(0,-15)[3]{$\Ql{,X}(-1)[-2]$}
    \mor{1}{2}{$$}
    \mor{3}{1}{$\text{abs. pur.}$}
    \obj(-7,0)[4]{$\simeq$}
    \mor{3}{2}{$$}
    \enddc
\end{equation}
The absolute purity isomorphism is given by the class $cl(X)\in \textup{H}^2_X(V,\Ql{}(1))$, whose image in $\textup{H}^2(X,\Ql{}(1))$ is $c_1(\mathcal{N}_{X|V}^{\vee})$, the first Chern class of the conormal bundle (see \cite{fuj02}). The map $\Ql{,X}\rightarrow \Ql{,X}(1)[2]$ corresponding to this class is the image of the right vertical arrow in the diagram above via the $\oo$-functor $-(1)[2]$. It suffices to notice that $\mathcal{N}_{X|V}^{\vee}\simeq \Lc_X$ to conclude.
\end{proof}
\section{\textsc{The comparison theorem}}
\subsection{\textsc{\texorpdfstring{$2$}{2}-periodic \texorpdfstring{$\ell$}{l}-adic sheaves}}
We shall now approach the comparison between monodromy-invariant vanishing cycles and the $\ell$-adic realization of the dg category of singularity $\sing{X_0}$. In order to do so, we need to work with the category of $\Ql{,S}(\beta)$-modules in $\shvl(S)$. Recall that there is an adjucntion of $\oo$-functors
\begin{equation}
    \shvl(S)\rightleftarrows \md_{\Ql{,S}(\beta)}(\shvl(S))
\end{equation}
given by $-\otimes_{\Ql{,S}}\Ql{,S}(\beta)$ and by the forgetful functor.

Notice that $\md_{\Ql{,\bullet}}(\shvl(\bullet))$ defines a fibered category over the category of schemes, which satisfies Grothendieck's six functors formalism.
\begin{defn}\label{Phimi}
With the same notation as above let 
\begin{equation}\label{definition spmi beta}
    s_0^*i_0^*j_{0*}\Ql{,\mathcal{U}_X}(\beta)\simeq i^*s_X^*j_{0*}\Ql{,\mathcal{U}_X}(\beta)\rightarrow i^*j_*s_{\mathcal{U}}^*\Ql{,\mathcal{U}_X}(\beta)
\end{equation}
be the arrow in $\md_{\Ql{,X}(\beta)}(\shvl(X))$ induced by the base change map $s_X^*j_{0*}\rightarrow j_*s_{\mathcal{U}}^*$. Denote by $\Phi^{\textup{mi}}_{(X,s_X)}(\Ql{}(\beta))$ its cofiber.
\end{defn}

\begin{prop}\label{compatibility Phimi with beta}
The maps $(\ref{definition spmi beta})$ and $sp^{\textup{mi}}_{(X,s_X)}\otimes_{\Ql{,X}}\Ql{,X}(\beta)$ are homotopic. In particular,
there is an equivalence in $\md_{\Ql{,X}(\beta)}\bigl( \shvl(X) \bigr)$:
\begin{equation}
    \Phi^{\textup{mi}}_{(X,s_X)}(\Ql{})\otimes_{\Ql{,X}}\Ql{,X}(\beta)\simeq \Phi_{(X,s_X)}^{\textup{mi}}(\Ql{,X}(\beta))
\end{equation}
\end{prop}
\begin{proof}This is analogous to \cite[Proposition 4.28]{brtv}. As a first step, notice that both $*$-pullbacks and $*$-pushforwards commute with Tate and usual shifts. Since $-\otimes_{\Ql{}}\Ql{}(\beta)$ commutes with $*$-pullbacks, we have
\begin{equation}
    i^*s_X^*j_{0*}\Ql{,\mathcal{U}_X}\otimes_{\Ql{,X}}\Ql{,X}(\beta)\simeq i^*s_X^*(j_{0*}\Ql{,\mathcal{U}_X}\otimes_{\Ql{,V_X}}\Ql{,V_X}(\beta))
\end{equation}
\begin{equation}
     i^*j_{*}\Ql{,X_{\mathcal{U}}}\otimes_{\Ql{,X}}\Ql{,X}(\beta)\simeq i^*(j_{*}\Ql{,X_{\mathcal{U}}}\otimes_{\Ql{,V_X}}\Ql{,V_X}(\beta))
\end{equation}
Using the equivalence $\Ql{}(\beta)\simeq \bigoplus_{i\in \Z}\Ql{}(i)[2i]$, we see that
\begin{equation}
    j_{0*}\Ql{,\mathcal{U}_X}\otimes_{\Ql{,V_X}}\Ql{,V_X}(\beta)\simeq \bigoplus_{i\in \Z}(j_{0*}\Ql{,\mathcal{U}_X}(i)[2i])
\end{equation}
We shall show that the canonical map $\bigoplus_{i\in \Z}(j_{0*}\Ql{,\mathcal{U}_X}(i)[2i])\rightarrow j_{0*}(\bigoplus_{i\in \Z}\Ql{,\mathcal{U_X}}(i)[2i])$. This follows immediately from the fact that the $*$-pushforward commutes with filtered colimits and from the equivalence $\bigoplus_{i\in \Z} \Ql{,\mathcal{U}_X}(i)[2i]\simeq colim_{i\geq 0}\bigoplus_{k=-i}^{i}\Ql{,\mathcal{U}_X}(k)[2k]$. This shows that
\begin{equation}
i^*s_X^*j_{0*}\Ql{,\mathcal{U}_X}\otimes_{\Ql{,X}}\Ql{,X}(\beta)\simeq i^*s_X^*j_{0*}\Ql{,\mathcal{U}_X}(\beta)
\end{equation}
The same argument applies to show that 
\begin{equation}
i^*j_{*}\Ql{,X_{\mathcal{U}}}\otimes_{\Ql{,X}}\Ql{,X}(\beta)\simeq i^*j_{*}\Ql{,X_{\mathcal{U}}}(\beta)
\end{equation}
To show that the two maps are homotopic, it suffices to show that they are so before applying $i^*$. Notice that $(\Ql{,\mathcal{U}_X}\rightarrow s_{\mathcal{U}*}s_{\mathcal{U}}^*\Ql{,\mathcal{U}_X})\otimes_{\Ql{,\mathcal{U}_X}}\Ql{,\mathcal{U}_X}(\beta)$ is homotopic to $\Ql{,\mathcal{U}_X}(\beta)\rightarrow s_{\mathcal{U}*}s_{\mathcal{U}}^*\Ql{,\mathcal{U}_X}(\beta)$, as $-\otimes_{\Ql{,\mathcal{U}_X}}\Ql{,\mathcal{U}_X}(\beta)$ is compatible with the $*$-pullback and with the $!$-pushforward, which coincides with the $*$-pushforward for $s_{\mathcal{U}}$ as it is a closed morphism. For what we have said above,
\begin{equation}
    j_{0*}(\Ql{,\mathcal{U}_X}\rightarrow s_{\mathcal{U}*}s_{\mathcal{U}}^*\Ql{,\mathcal{U}_X})\otimes_{\Ql{,\mathcal{U}_X}}\Ql{,\mathcal{U}_X}(\beta)\simeq j_{0*}\Ql{,\mathcal{U}_X}(\beta)\rightarrow j_{0*}s_{\mathcal{U}*}s_{\mathcal{U}}^*\Ql{,\mathcal{U}_X}(\beta)
\end{equation}
and $-\otimes_{\Ql{}}\Ql{}(\beta)$ is compatible with $*$-pullbacks. The last assertion follows as $-\otimes_{\Ql{}}\Ql{}(\beta)$ is exact.
\end{proof}
\begin{cor}
Let $(X,s_X)\in \tlgm{S}$ and assume that $X$ is a regular scheme. There is an equivalence in $\md_{\Ql{,X}(\beta)}(\shvl(X))$
\begin{equation}
    i_{0}^*j_{0*}\Ql{,\mathcal{U}_X}(\beta)\simeq cofiber \bigl(c_1(\Lc_X)\otimes_{\Ql{,X}}\Ql{,X}(\beta):\Ql{,X}(\beta)\rightarrow \Ql{,X}(\beta) \bigr)
\end{equation}
\end{cor}
\begin{proof}
This follows immediately from the previous Proposition, Proposition \ref{cofiber first chern class} and from the Bott equivalence $\Ql{,X}(\beta)(-1)[-2]\simeq \Ql{,X}(\beta)$.
\end{proof}
\begin{rmk}
Notice that $i_0^*j_{0*}$ and $i^*j_*$ are lax monoidal functors. In particular, $s_0^*i_0^*j_{0*}\Ql{,\mathcal{U}_X}(\beta)$ and $i^*j_*\Ql{,X_{\mathcal{U}}}(\beta)$ have commutative algebra structures. The map $sp_{(X,s_X)}^{\textup{mi}}\otimes_{\Ql{,X_0}}\Ql{,X_0}(\beta)$ is a map of commutative algebras as $s_X^*j_{0*}\Ql{,X_{\mathcal{U}_X}}(\beta)\rightarrow j_*s_{\mathcal{U}}^*\Ql{,X_{\mathcal{U}_X}}(\beta)$ is so. In particular, the $\ell$-adic sheaf $\Phi^{\textup{mi}}_{(X,s_X)}(\Ql{}(\beta))$ lives in $\md_{i^*s_X^*j_{0*}\Ql{,\mathcal{U}_X}}(\beta)(\shvl(X_0))$.
\end{rmk}

\subsection{\textsc{The main theorem}}
\begin{prop}\label{compatibility chern characters}
Let $X$ be a regular scheme and let $\Lc \in Pic(X)$. Then
\begin{equation}
    cofib(1-\Rl_X(m_{\Lc}):\Rl_X(\bu_X)\rightarrow \Rl_X(\bu_X))
\end{equation}
\begin{equation*}
\simeq cofib(c_1(\Lc)\otimes_{\Ql{,X}}\Ql{,X}(\beta):\Ql{,X}(\beta)\rightarrow \Ql{,X}(\beta))
\end{equation*}
where $c_1(\Lc)$ is the first Chern class of $\Lc$.

Moreover, 
\begin{equation}
    cofib(c_1(\Lc)\otimes_{\Ql{,X}}\Ql{,X}(\beta))\simeq  cofib(c_1(\Lc^{\vee})\otimes_{\Ql{,X}}\Ql{,X}(\beta))
\end{equation}
\end{prop}
\begin{proof}
We start by noticing that the composition
\begin{equation}
    \textup{HK}_0(X)\otimes_{\Z} \Q \simeq H^{2*,*}_{\mathscr{M}}(X,\Z)\otimes_{\Z}\Q \simeq H^{2*,*}_{\mathscr{M}}(X,\Q)
\end{equation}
where $H^{2*,*}_{\mathscr{M}}$ denotes motivic cohomology, coincides with the Chern character (see \cite[\S 11.3.6]{cd19}) which is compatible with the $\ell$-adic Chern character (see \cite[\S 3.6 - \S 3.7]{brtv} and \cite{cd19}). Now, according to \cite[Proposition 12.2.9]{cd19} and \cite[Proposition 3.8]{de08}, $c_1(\Lc)$ is nilpotent and therefore the Chern character of $\Lc$ can be written as $1+c_1+\frac{1}{2}c_1(\Lc)^2+\dots+\frac{1}{m!}c_1(\Lc)^{m}$ for some $m\geq 1$. Then 
\begin{equation*}
    1-\Rl_X(m_{\Lc})=-c_1(\Lc)(1+\frac{1}{2}c_1(\Lc)+\dots+\frac{1}{m!}c_1(\Lc)^{m-1})
\end{equation*}
As $c_1(\Lc)$ is nilpotent, so is $\frac{1}{2}c_1(\Lc)+\dots+\frac{1}{m!}c_1(\Lc)^{m-1}$, whence $(1+\frac{1}{2}c_1(\Lc)+\dots+\frac{1}{m!}c_1(\Lc)^{m-1})$ is invertible. This shows
\begin{equation*}
    cofiber(1-\Rl_X(m_{\Lc}))\simeq cofiber(-c_1(\Lc)\otimes{\Ql{,X}}\Ql{,X}(\beta))\simeq cofiber(c_1(\Lc)\otimes{\Ql{,X}}\Ql{,X}(\beta))
\end{equation*}
The last equivalence follows from the fact that, by \cite[Proposition 12.2.9 and Remark 13.2.2]{cd19}, we have 
\begin{equation}
    0=c_1(\Lc \otimes \Lc^{\vee})\simeq F(c_1(\Lc),c_1(\Lc^{\vee}))=c_1(\Lc)+c_1(\Lc^{\vee})+\beta^{-1}\cdot c_1(\Lc)\cdot c_1(\Lc^{\vee})
\end{equation}
In particular, 
\begin{equation}
    c_1(\Lc^{\vee})\simeq -c_1(\Lc)(1+\beta^{-1}\cdot c_1(\Lc^{\vee}))
\end{equation}
and, as we have remarked above, $-(1+\beta^{-1}\cdot c_1(\Lc^{\vee}))$ is invertible as $c_1(\Lc^{\vee})$ is nilpotent.
\end{proof}
We are finally ready to state our main theorem:
\begin{thm}\label{main theorem}
Let $(X,s_X)$ be a twisted LG model over $(S,\Lc_S)$. Assume that $X$ is regular and that $0,s_X:X\rightarrow V(\Lc_X)$ are Tor-independent. The following equivalence holds in $\md_{i^*s_X^*j_{0*}\Ql{,\mathcal{U}_X}(\beta)}(\shvl(X_0))$:
\begin{equation}
    i^*\Rlv_{X}(\sing{X,s_X})\simeq \Phi^{\textup{mi}}_{(X,s_X)}(\Ql{}(\beta))[-1]
\end{equation}
where the $i^*s_X^*j_{0*}\Ql{,\mathcal{U}_X}(\beta)$-module structure on the l.h.s. is the one induced by the equivalence of commutative algebras
\begin{equation}
   i^*\Rlv_X(\sing{X,0})\simeq   i^*s_X^*j_{0*}\Ql{,\mathcal{U}_X}(\beta)
\end{equation}
\end{thm}
\begin{proof}
We start by noticing that $i=s_0:X_0\rightarrow X$. Indeed, $i_0 \circ s_0=s_X\circ i$ and both $i_0$ and $s_X$ are sections of the canonical morphism $V(\Lc_X)\rightarrow X$. In particular
\begin{equation}
     i^*\Rlv_X(\sing{X,0})\simeq  s_0^*\Rlv_X(\sing{X,0})\underbracket{\simeq}_{\textup{Lemma }\ref{motivic realization sing(S,0)}}
      i^*s_X^*j_{0*}\Ql{,\mathcal{U}_X}(\beta)
\end{equation}
Since $X$ is supposed to be regular, $\sing{X,s_X}\simeq \sing{X_0}$. By Propositions \ref{fundamental fiber cofiber sing} and \ref{compatibility chern characters} we have a fiber-cofiber sequence in $\md_{i^*s_X^*j_{0*}\Ql{,\mathcal{U}_X}(\beta)}(\shvl(X_0)$
\begin{equation}
    i^*\Rl_X(\sing{X_0})\rightarrow cofiber(c_1(\Lc_{X_0})\otimes_{\Ql{,X_0}}\Ql{,X_0}(\beta):\Ql{,X_0}(\beta)\rightarrow \Ql{,X_0}(\beta))\xrightarrow{f} i^*j_*\Ql{,X_{\mathcal{U}}}(\beta)
\end{equation}
On the other hand, by Proposition \ref{compatibility Phimi with beta}, we have another fiber-cofiber sequence 
\begin{equation}
   \Phi^{\textup{mi}}_{(X,s_X)}(\Ql{}(\beta))[-1]\rightarrow cofiber(c_1(\Lc_{X_0})\otimes_{\Ql{,X_0}}\Ql{,X_0}(\beta):\Ql{,X_0}(\beta)\rightarrow \Ql{,X_0}(\beta))\xrightarrow{g} i^*j_*\Ql{,X_{\mathcal{U}}}(\beta)
\end{equation}
Both the arrows on the right in these two fiber-cofiber sequences are defined by the universal property of the object in the middle. Therefore, if we consider the diagram
\begin{equation}
    \begindc{\commdiag}[18]
    \obj(-80,20)[1]{$\Ql{,X_0}(\beta)$}
    \obj(0,20)[2]{$\Ql{,X_0}(\beta)$}
    \obj(80,20)[3]{$cofiber(c_1(\Lc_{X_0}^{\vee})\otimes_{\Ql{,X_0}}\Ql{,X_0}(\beta))$}
    \obj(0,-20)[4]{$i^*j_*\Ql{,X_{\mathcal{U}}}(\beta)$}
    \mor{1}{2}{$c_1(\Lc_{X_0}^{\vee})\otimes_{\Ql{,X_0}}\Ql{,X_0}(\beta)$}
    \mor{2}{3}{$$}
    \mor{1}{4}{$\sim 0$}
    \mor(-3,17)(-3,-17){$\tilde{f}$}[\atright,\solidarrow]
    \mor(3,17)(3,-17){$\tilde{g}$}
    \mor(76,17)(4,-16){$f$}[\atright,\solidarrow]
    \mor(84,17)(12,-16){$g$}
    \enddc
\end{equation}

where the middle vertical arrows induce the morphism $f$ and $g$ respectively, it suffices to show $\tilde{f}\sim \tilde{g}$. By definition, $\tilde{f}$ is induced by the counit $\Ql{,X}(\beta)\rightarrow j_*j^*\Ql{,X}(\beta)$, while $\tilde{g}$ is induced by the base change morphism $b.c.:s_X^*j_{0*}\Ql{,\mathcal{U}_X}(\beta)\rightarrow j_*s_{\mathcal{U}}^*\Ql{,\mathcal{U}_X}(\beta)$:
\begin{equation}
    \tilde{g}:i^*s_X^*\Ql{,V_X}(\beta)\rightarrow i^*s_X^*j_{0*}j_0^*\Ql{,V_X}(\beta)\rightarrow i^*j_*s_{\mathcal{U}}^*j_{0}^*\Ql{,V_X}(\beta)
\end{equation}
It also suffices to show that the two arrows are homotopic before applying $i^*$. Consider the diagram
\begin{equation}
    \begindc{\commdiag}[18]
    \obj(-80,10)[1]{$s_X^*\Ql{,V_X}(\beta)$}
    \obj(0,10)[2]{$s_X^*j_{0*}j_0^*\Ql{,V_X}(\beta)$}
    \obj(80,10)[3]{$j_*s_{\mathcal{U}}^*j_0^*\Ql{,V_X}(\beta)\simeq j_*j^*s_X^*\Ql{,V_X}(\beta)$}
    \obj(0,-10)[4]{$s_X^*s_{X*}j_*s_{\mathcal{U}}^*j_0^*\Ql{,V_X}(\beta)$}
    \obj(0,-26)[5]{$s_X^*s_{X*}j_*j^*s_X^*\Ql{,V_X}(\beta)$}
    \obj(0,-18)[]{$\simeq$}
    \mor{1}{2}{$$}
    \mor{2}{3}{$b.c.$}
    \mor{2}{4}{$$}
    \mor(25,-18)(80,7){$s_X^*s_{X*}\rightarrow 1$}[\atright,\solidarrow]
    \mor(-80,7)(-25,-18){$$}
    \enddc
\end{equation}

If we analyse the commutative triangle on the left, we see that the objlique arrow corresponds to the unit of the adjunction $((s_X\circ j)^*,(s_X\circ j)_*)$. Indeed, by definition, the vertical arrow is $s_X^*j_{0*}(j_0^*\Ql{,V_X}(\beta)\xrightarrow{\text{unit }(s_{\mathcal{U}}^*,s_{\mathcal{U}*})} s_{\mathcal{U}*}s_{\mathcal{U}}^*j_0^*\Ql{,V_X}(\beta))$ and the horizontal arrow is the one induced by the unit of $(j_0^*,j_{0*})$. The composition is exactly the unit of the adjunction $((j_0\circ s_{\mathcal{U}})^*,(j_0\circ s_{\mathcal{U}})_*)$. Hence, the claim is proved as $j_0\circ s_{\mathcal{U}}\simeq s_X\circ j$. Now notice that $s_{X*}$ is conservative ($s_X$ is a closed morphism), i.e. the counit $s_X^*s_{X*}\rightarrow 1$ is a natural equivalence. If we compose the oblique arrow with it, we obtain the unit of $(j^*,j_*)$ evaluated in $s_X^*\Ql{,V_X}(\beta)$. The statement about the module structures is clear.
\end{proof}
\begin{rmk}
Notice that our main theorem provides a generalization of the formula proved in \cite{brtv}: assume that we are given a proper flat morphism $p:X\rightarrow S$ from a regular scheme to an excellent strictly henselian trait\footnote{We shall recycle the notation that we have introduced in $\S$ \ref{the formalism of vanishing cycles}}. Let $s_X:X\rightarrow \aff{1}{X}$ be the pullback of the section $S\rightarrow \aff{1}{S}$ given by the uniformizer $\pi$. Then we can consider the diagram
\begin{equation}
    \begindc{\commdiag}[18]
    \obj(-40,20)[1]{$\sigma$}
    \obj(0,20)[2]{$S$}
    \obj(40,20)[3]{$\eta$}
    \obj(-40,0)[4]{$X_0$}
    \obj(0,0)[5]{$X$}
    \obj(40,0)[6]{$X_{\eta}$}
    \obj(-40,-20)[7]{$X$}
    \obj(0,-20)[8]{$\aff{1}{X}$}
    \obj(40,-20)[9]{$\gm{X}$}
    \mor{1}{2}{$i_{\sigma}$}
    \mor{3}{2}{$j_{\eta}$}[\atright,\solidarrow]
    \mor{4}{1}{$p_{\sigma}$}
    \mor{5}{2}{$p$}
    \mor{6}{3}{$p_{\eta}$}[\atright,\solidarrow]
    \mor{4}{5}{$i$}
    \mor{6}{5}{$j$}[\atright,\solidarrow]
    \mor{4}{7}{$s_0$}[\atright,\solidarrow]
    \mor{5}{8}{$s_X$}
    \mor{6}{9}{$s_{\gm{X}}$}
    \mor{7}{8}{$i_0$}[\atright,\solidarrow]
    \mor{9}{8}{$j_0$}
    \enddc
\end{equation}
where all squares are cartesian. The theorem we have just proved tells us that
\begin{equation*}
    i^*\Rlv_X(\sing{X,s_X})\simeq \Phi_{(X,s_X)}^{\textup{mi}}(\Ql{}(\beta))[-1]
\end{equation*}
By the regularity assumption that we have made on $X$, $\sing{X,s_X}\simeq \sing{X_0}$. If we apply $p_{\sigma*}$ to the formula above we find:
\begin{equation}
    p_{\sigma*}i^*\Rlv_X(\sing{X_0})\simeq i_{\sigma}^*\Rlv_S(\sing{X_0})\simeq
\end{equation}
\begin{equation*}
    p_{\sigma*}\Phi_{(X,s_X)}^{\textup{mi}}(\Ql{}(\beta))[-1]\underbracket{\simeq}_{\text{Lemma }\ref{comparison monodromy invariant inertia invariant vanishing cycles}} \bigl( p_{\sigma*}\Phi_p(\Ql{,X}(\beta))[-1] \bigr)^{\textup{h}I}
\end{equation*}
which is exactly the content of \cite[Theorem 4.39]{brtv}. Also notice that in this case 
\begin{equation}
i^*s_X^*j_{0*}\Ql{,\gm{X}}(\beta)\simeq \Ql{,X_{0}}(\beta)\oplus \Ql{,X_{0}}(\beta)[1]
\end{equation}
as the first Chern class of the trivial line bundle is zero. This recovers the equivalence 
\begin{equation}
    (\Ql{,\sigma}(\beta))^{\textup{h}I}\simeq \Ql{,\sigma}(\beta)\oplus \Ql{,\sigma}(\beta)[1]
\end{equation}
proved in \textit{loc.cit.}
\end{rmk}
\section{\textsc{The \texorpdfstring{$\ell$}{l}-adic realization of the dg category of singularities of an \texorpdfstring{$n$}{n}-dimensional LG model}}

In this section we will explain the reason that led us to consider twisted LG model, as defined in $\S 3.1$.
\subsection{\textsc{Reduction of codimension}}
\begin{context}
Let $S=Spec(A)$ be a regular noetherian affine scheme of finite Krull dimension.
\end{context}
\begin{rmk}
We will use the same notation introduced in the previous chapter. In particular, $\lgm{n}^{\boxplus}$ denotes the symmetric monoidal category of $n$-dimensional LG models over $S$.
\end{rmk}
\begin{rmk}
This section is only concerned with dg categories of singularities and not with their motivic and $\ell$-adic realizations. Therefore, we could have inserted it in the previous chapter. However, since Theorem \ref{reduction of codimension theorem} is the reason why we have started to study the $\ell$-adic realization of dg categories of twisted LG models, for expository reasons we have preferred to include it in this chapter.
\end{rmk}
\begin{construction}
We will construct an (ordinary) symmetric monoidal functor
\begin{equation}
\Xi^{\boxplus}:\lgm{n}^{\boxplus}\rightarrow \tlgmproj{S}^{\boxplus}
\end{equation}
following the lead of \cite{bw15} and \cite{orl06}.

Let $(X,\underline{f})\in \lgm{n}$ and consider the $(n-1)$-dimensional projective space $\proj{n-1}{X}=\textup{Proj}_X(\Oc_X[T_1,\dots,T_n])$. We will write $\Oc(1)$ instead of $\Oc_{\proj{n-1}{X}}(1)$ to light the notation. Consider $W_{\underline{f}}=f_1T_1+\dots+f_nT_n\in \Gamma(\proj{n-1}{X},\Oc(1))$. The assignment 
\begin{equation*}
    (X,\underline{f})\mapsto (\proj{n-1}{X},W_{\underline{f}})
\end{equation*}
together with the obvious law for morphisms defines a functor
\begin{equation}
\Xi:    \lgm{n}\rightarrow \tlgmproj{S}
\end{equation}
It is immediate to observe that $(S,\underline{0})\mapsto (\proj{n-1}{S},0)$, i.e. the functor is compatible with the units of the two symmetric monoidal structures. It remains to show that
\begin{equation*}
    \Xi((X,\underline{f})\boxplus (Y,\underline{g}))\simeq \Xi((X,\underline{f}))\boxplus \Xi((Y,\underline{g}))
\end{equation*}
On the left and side we have
\begin{equation*}
    (\proj{n-1}{X\times_S Y},W_{\underline{f}\boxplus \underline{g}}=(f_1\boxplus g_1)\cdot T_1+\dots+(f_n\boxplus g_n)\cdot T_n)
\end{equation*}
while on the right hand side we have
\begin{equation*}
    (\proj{n-1}{X}\times_{\proj{n-1}{S}}\proj{n-1}{Y},W_{\underline{f}}\boxplus W_{\underline{g}})
\end{equation*}
Since $\proj{n-1}{X}\times_{\proj{n-1}{S}}\proj{n-1}{Y}\simeq (X\times_S \proj{n-1}{S})\times_{\proj{n-1}{S}}(Y\times_S \proj{n-1}{S})\simeq (X\times_S Y)\times_S \proj{n-1}{S}\simeq \proj{n-1}{X\times_S Y}$, it suffices to show that $W_{\underline{f}\boxplus \underline{g}}=W_{\underline{f}\boxplus W_{\underline{g}}}$. It is enough to do it locally ($\Oc(1)$ is a sheaf). We may consider the covering of $\proj{n-1}{X\times_S Y}$ consisting of open affine subschemes $Spec\bigl (B\otimes_A C[\frac{T_1}{T_j},\dots,\frac{T_n}{T_j}] \bigr )$, where $Spec(A)$ (resp. $Spec(B)$, resp. $Spec(C)$) is an open affine subscheme of $S$ (resp. $X$, resp. $Y$). The restriction to this open subset of $W_{\underline{f}\boxplus \underline{g}}$ is
\begin{equation*}
    (f_1\otimes1+1\otimes g_1)\cdot \frac{T_1}{T_j}+\dots +(f_n\otimes 1+1\otimes g_n)\cdot \frac{T_n}{T_j}
\end{equation*}
while that of $W_{\underline{f}}\boxplus W_{\underline{g}}$ is
\begin{equation*}
    \bigl(f_1\cdot \frac{T_1}{T_j}\otimes 1 + 1\otimes g_1\cdot \frac{T_1}{T_j}\bigr)+\dots+\bigl(f_n\cdot \frac{T_n}{T_j}\otimes1+1\otimes g_n\cdot \frac{T_n}{T_j}\bigr)
\end{equation*}
The claim follows and thus we obtain the desired symmetric monoidal functor.
\end{construction}
If we compose this symmetric monoidal functor with the lax monoidal $\oo$-functor defined in Section \ref{The dg category of singularities of a twisted LG model} we get
\begin{equation}\label{sing projectivization}
    \sing{\proj{n-1}{\bullet},W_{\bullet}}^{\otimes} : \lgm{n}^{\textup{op},\boxplus}\rightarrow \tlgmproj{S}^{\textup{op},\boxplus}\rightarrow \md_{\sing{\proj{n-1}{S},0}}\bigl(\dgcatm\bigr)^{\otimes}
\end{equation}
which, at the level of objects, corresponds to the assignment $(X,\underline{f})\mapsto \sing{\proj{n-1}{X},W_{\underline{f}}}$.

Recall that we define the dg category of singularities of an $n$-dimensional LG model $(X,\underline{f})$ in the following way: consider the derived zero locus of $\underline{f}$, defined as the homotopy pullback of $\underline{f}$ along the zero section 
\begin{equation}
    \begindc{\commdiag}[12]
    \obj(-20,20)[1]{$X_0$}
    \obj(20,20)[2]{$X$}
    \obj(-20,-20)[3]{$S$}
    \obj(20,-20)[4]{$\aff{n}{S}$}
    \mor{1}{2}{$i$}
    \mor{1}{3}{$$}
    \mor{2}{4}{$\underline{f}$}
    \mor{3}{4}{$\underline{0}$}
    \enddc
\end{equation}

As $\underline{0}:S\rightarrow \aff{n}{S}$ is a closed lci morphism, so is $i:X_0\rightarrow X$. Thus the pushforward induces a functor
\begin{equation}
    i_*:\sing{X_0}\rightarrow \sing{X}
\end{equation}
Then
\begin{equation}
    \sing{X,\underline{f}}:=Ker\bigl(i_*:\sing{X_0}\rightarrow \sing{X} \bigr)
\end{equation}
Also recall that in Section \ref{dg Category of singularities of an n-LG model over S} we defined a lax monoidal $\oo$-functor 
\begin{equation}
    \sing{\bullet,\bullet}^{\otimes}:\lgm{n}^{\textup{op},\boxplus}\rightarrow \md_{\sing{S,\underline{0}}}\bigl( \dgcatm \bigr)^{\otimes}
\end{equation}
At the level of objects, it is defined by $(X,\underline{f})\mapsto \sing{X,\underline{f}}$.

At this point, the reader might complain that there are too many $\textup{\textbf{Sing}}$'s involved, and this may cause confusion. As a partial justification to the choice of notation we made, let us show that these different $\oo$-functors are closely related.
\begin{construction}
Let $(X,\underline{f})\in \lgm{n}$ and consider the following picture, where both squares are homotopy cartesian. 
\begin{equation}
    \begindc{\commdiag}[18]
    \obj(0,20)[1]{$\proj{n-1}{X}$}
    \obj(30,20)[2]{$V(\Oc(1))$}
    \obj(-30,0)[3]{$\proj{n-1}{X_0}$}
    \obj(0,0)[4]{$V(W_{\underline{f}})$}
    \obj(30,0)[5]{$\proj{n-1}{X}$}
    \obj(-30,-20)[6]{$X_0$}
    \obj(30,-20)[7]{$X$}
    \mor{1}{2}{$W_{\underline{f}}$}
    \mor{4}{1}{$k$}
    \mor{5}{2}{$s_0$}[\atright,\solidarrow]
    \mor{4}{5}{$k$}
    \mor{3}{4}{$j$}
    \mor{3}{6}{$p_0$}
    \mor{5}{7}{$p$}
    \mor{6}{7}{$i$}
    \enddc
\end{equation}
To see that the two morphisms $V(W_{\underline{f}})\rightarrow \proj{n-1}{X}$ coincide, it suffices to notice that the square is commutative (in the $\oo$-categorical sense) and both $W_{\underline{f}}$ and $s_0$ are sections of the canonical morphism $V(\Oc(-1))\rightarrow \proj{n-1}{X}$.
By \cite[Chapter 4, Lemma 3.1.3, Lemma 5.1.4]{gr17}, $j_*p_0^*$ preserves complexes with coherent bounded cohomology. The derived proper base change equivalence $k_*j_*p_0^*\simeq p^*i_*$ implies that we have a dg functor
\begin{equation}\label{jp0}
    \widetilde{\Upsilon}_{(X,\underline{f})}:=j_*p_0^*:\cohbp{X_0}{X}\rightarrow \cohbp{V(W_{\underline{f}})}{\proj{n-1}{X}}
\end{equation}
Moreover, it is obviously functorial in $(X,\underline{f})$. It is also compatible with the tensor structures of $(X,\underline{f})\mapsto \cohbp{X_0}{X}$ and of $(X,\underline{f})\mapsto \cohbp{V(W_{\underline{f}})}{\proj{n-1}{X}}$. We can show this using explicit models. We can also restrict to the affine case. Let $(X=Spec(B),\underline{f})$ and $(Y=Spec(C),\underline{g})$ be affine $n$-dimensional LG models. Since $\cohbp{V(W_{\underline{f}})}{\proj{n-1}{X}}^{\otimes}\simeq \varprojlim_{i=1,\dots,n} \cohbp{V(W_{\underline{f}})\times_{\proj{n-1}{X}}D_{(T_i)}}{D_{(T_i)}}^{\otimes}$, it suffices to show that the claim holds true for
\begin{equation*}
    \widetilde{\Upsilon}_{(X,\underline{f}),i}:\cohbp{X_0}{X}\xrightarrow{\widetilde{\Upsilon}_{(X,\underline{f})}} \cohbp{V(W_{\underline{f}})}{\proj{n-1}{X}}\xrightarrow{u_i^*} \cohbp{V(W_{\underline{f}})\times_{\proj{n-1}{X}}D_{(T_i)}}{D_{(T_i)}}^{\otimes}
\end{equation*}
Without loss of generality, we shall assume that $i=n$. Let $\tilde{f}=f_1t_1+\dots+f_{n-1}t_{n-1}+f_n$, where $t_i$ are the parameters of $D_{(T_i)}\simeq \aff{n-1}{X}$. Recall that we have models $\cohs{B,\underline{f}}$ for $\cohbp{X_0}{X}$ and $\cohs{\aff{n-1}{X},\tilde{f}}$ for $\cohbp{V(\tilde{f})}{\aff{n-1}{X}}$. We refer to \cite{brtv} and \cite{pi19} for the definitions of these two objects. Then $\widetilde{\Upsilon}_{(X,\underline{f}),i}$ is induced by
\begin{equation}
    \cohs{X,\underline{f}}\rightarrow \cohs{\aff{n-1}{X},\tilde{f}}
\end{equation}
\begin{equation*}
    (E_j,d_j,\{h_j^s\}_{s=1,\dots,n})\mapsto (E_j[t_1,\dots,t_{n-1}],d_j,h^1t_1+\dots+h^{n-1}t_{n-1}+h^n)
\end{equation*}
We need to show that we have a commutative square
\begin{equation}
\begindc{\commdiag}[12]
\obj(-80,20)[1]{$\cohs{\aff{n-1}{B},\tilde{f}}\otimes \cohs{\aff{n-1}{C},\tilde{g}}$}
\obj(80,20)[2]{$\cohs{\aff{n-1}{B\otimes_A C},\tilde{f}\boxplus \tilde{g}}$}
\obj(-80,-20)[3]{$\cohs{B,\underline{f}}\otimes \cohs{C,\underline{g}}$}
\obj(80,-20)[4]{$\cohs{B\otimes_A C,\underline{f}\boxplus \underline{g}}$}
\mor{1}{2}{$$}
\mor{3}{1}{$$}
\mor{4}{2}{$$}
\mor{3}{4}{$$}
\enddc
\end{equation}
This is immediate after the definitions.

As a consequence, we obtain a lax monoidal $\oo$-natural transformation
\begin{equation}\label{widetilde Upsilon}
    \widetilde{\Upsilon}^{\otimes}:\cohbp{\bullet}{\bullet}\rightarrow \cohbp{V(W_{\bullet})}{\proj{n-1}{\bullet}}:\lgm{n}^{\textup{op},\otimes}\rightarrow \md_{\cohb{S_0}}(\dgcatm)^{\otimes}
\end{equation}
where $\cohbp{V(W_{\underline{\bullet}})}{\proj{n-1}{\bullet}}$ is defined by the composition of $\Sigma^{\boxplus}$,
$(\ref{cohbp non affine})$ (over the base $(\proj{n-1}{S},\Oc(1))$) and
\begin{equation*}
\md_{\cohb{(\proj{n-1}{S})_0}}(\dgcatm)^{\otimes}\rightarrow \md_{\cohb{S_0}}(\dgcatm)^{\otimes}
\end{equation*}
\end{construction}

\begin{lmm}
$j:\proj{n-1}{X_0}\rightarrow V(W_{\underline{f}})$ is an lci morphism of derived schemes.
\end{lmm}
\begin{proof}
Since $j$ is clearly of finite presentation, we will only need to show that the relative cotangent complex is of Tor amplitude $[-1,0]$. Recall that we have the fundamental fiber-cofiber sequence
\begin{equation}
    j^*\Lb_{V(W_{\underline{f}})/\proj{n-1}{X}}\rightarrow \Lb_{\proj{n-1}{X_0}/\proj{n-1}{X}}\rightarrow \Lb_{\proj{n-1}{X_0}/V(W_{\underline{f}})}
\end{equation}
at our disposal. Since $\proj{n-1}{X_0}\simeq S\times^h_{\aff{n}{S}}\proj{n-1}{X}$, $\Lb_{\proj{n-1}{X_0}/\proj{n-1}{X}}$ is equivalent to the pullback of $\Lb_{S/\aff{n}{S}}\simeq \Oc_S^{n}[1]$ along $\proj{n-1}{X_{0}}\rightarrow S$, i.e. $\Lb_{\proj{n-1}{X_0}/\proj{n-1}{X}}\simeq \Oc_{\proj{n-1}{X_0}}^n[1]$. On the other hand, $\Lb_{V(W_{\underline{f}})/\proj{n-1}{X}}\simeq k^*\Lb_{\proj{n-1}{X}/V(\Oc(-1))}\simeq k^*(\Oc(-1)[1])$.
Therefore, $j^*\Lb_{V(W_{\underline{f}})/\proj{n-1}{X}}\simeq \Oc_{\proj{n-1}{X_0}}\otimes_{\Oc_{\proj{n-1}{X}}}\Oc(-1)$. We just need to identify the morphism $j^*\Lb_{V(W_{\underline{f}})/\proj{n-1}{X}}\simeq \Oc_{\proj{n-1}{X_0}}\otimes_{\Oc_{\proj{n-1}{X}}}\Oc(-1)[1]\rightarrow \Oc_{\proj{n-1}{X_0}}^{n}[1]\simeq \Lb_{\proj{n-1}{X_0}}$. It coincides with the morphism $(T_1,\dots,T_n)$. In particular, the cofiber of this morphism, i.e. $\Lb_{\proj{n-1}{X_0}V(W_{\underline{f}})}$, is equivalent to $\Oc_{\proj{n-1}{X_0}}^{n-1}[1]$.
\end{proof}
\begin{construction}
Since $j$ is an lci morphism, the dg functor $(\ref{jp0})$ preserves perfect complexes. Therefore, for every $(X,\underline{f})\in \lgm{n}$, we have an induced dg functor
\begin{equation}
    \Upsilon_{(X,\underline{f})}:=j_*p_0^*:\sing{X,\underline{f}}\rightarrow \sing{\proj{n-1}{X},W_{\underline{f}}}
\end{equation}
Starting from $(\ref{widetilde Upsilon})$, by the usual standard arguments we thus obtain a lax monoidal $\oo$-natural transformation
\begin{equation}\label{reduction of codimension natural transformation}
    \Upsilon^{\otimes}:\sing{\bullet,\bullet}^{\otimes}\rightarrow \sing{\proj{n-1}{\bullet},W_{\bullet}}^{\otimes}:\lgm{n}^{\textup{op},\boxplus}\rightarrow \md_{\sing{S,\underline{0}}}(\dgcatm)^{\otimes}
\end{equation}
where $\sing{\proj{n-1}{\bullet},W_{\bullet}}^{\otimes}$ denotes the composition of ($\ref{sing projectivization}$) with
\begin{equation*}
    \md_{\sing{\proj{n-1}{S},0}}(\dgcatm)^{\otimes}\rightarrow \md_{\sing{S,\underline{0}}}(\dgcatm)^{\otimes}
\end{equation*}
\end{construction}
\begin{thm}{\textup{(\cite{orl06}, \cite{bw15})}}\label{reduction of codimension theorem}
The lax monoidal $\oo$-natural transformation $(\ref{reduction of codimension natural transformation})$ induces an equivalence
\begin{equation}
   \Upsilon_{(Spec(B),\underline{f})}: \sing{Spec(B),\underline{f}}\simeq \sing{\proj{n-1}{B},W_{\underline{f}}}
\end{equation}
of dg categories whenever $\underline{f}$ is a regular sequence. 
\end{thm}
\begin{proof}
This is an immediate consequence of \cite[Theorem 2.10]{bw15}. Indeed, as the dg categories are triangulated, it suffices to show that the statement is true on the induced functor of triangulated categories. This coincides by construction with that of \textit{loc. cit.}
\end{proof}
\begin{rmk}
We actually believe that the statement above remains true even if $\underline{f}$ is not assumed to be regular, as long as one considers the derived zero locus of $\underline{f}$.
\end{rmk}

\subsection{\textsc{The \texorpdfstring{$\ell$}{l}-adic realization of the dg Category of singularities of an \texorpdfstring{$n$}{n}-dimensional LG model}}
It is now easy to obtain the following computation
\begin{thm}\label{formula for (B,f1,fn)}
Let $(Spec(B),\underline{f})$ be an affine $n$-dimensional LG model. Assume that $B$ is a regular ring and that $\underline{f}$ is a regular sequence. The following equivalence holds in $\md_{\Rlv_B(\sing{Spec(B),\underline{0}})}(\shvl(Spec(B)))$
\begin{equation}
    \Rlv_B(\sing{Spec(B),\underline{f}})\simeq p_*i_*\Phi^{\textup{mi}}_{(\proj{n-1}{B},W_{\underline{f}})}(\Ql{}(\beta))[-1]
\end{equation}
\end{thm}
\begin{proof}
As $B$ is regular and $\underline{f}$ is a regular sequence, $\sing{Spec(B),\underline{f}}\simeq \sing{Spec(B/\underline{f})}$. Notice that in this situation $Spec(B/\underline{f})$ coincides with the derived zero locus of $\underline{f}$. By Theorem \ref{reduction of codimension theorem}, we have that $\sing{B,\underline{f}}\simeq \sing{V(W_{\underline{f}})}\simeq \sing{\proj{n-1}{B},W_{\underline{f}}}$. Then
\begin{equation}
    \Rlv_B(\sing{Spec(B),\underline{f}})\simeq p_*\Rlv_{\proj{n-1}{B}}(\sing{\proj{n-1}{B}},W_{\underline{f}})
\end{equation}
\begin{equation*}
\underbracket{\simeq}_{\textup{Theorem } \ref{main theorem}}p_*i_*\Phi^{\textup{mi}}_{(\proj{n-1}{B},W_{\underline{f}})}(\Ql{}(\beta))[-1]
\end{equation*}
where $p:\proj{n-1}{B}\rightarrow Spec(B)$ and $i:V(W_{\underline{f}})\rightarrow \proj{n-1}{B}$.
The fact that this equivalence holds in  $\md_{\Rlv_B(\sing{Spec(B),\underline{0}})}(\shvl(Spec(B)))$ follows immediately from the fact that $\Rlv_{\proj{n-1}{B}}(\sing{\proj{n-1}{B}},W_{\underline{f}})\simeq i_*\Phi^{\textup{mi}}_{(\proj{n-1}{B},W_{\underline{f}})}(\Ql{}(\beta))[-1]$ is an equivalence of $\Rlv_{\proj{n-1}{B}}(\sing{\proj{n-1}{B},0})$-modules and from the fact that \begin{equation}\Rlv_B(\sing{Spec(B),\underline{0}})\rightarrow p_*\Rlv_{\proj{n-1}{B}}(\sing{\proj{n-1}{B},0})
\end{equation}
is a morphism of commutative algebras. 
\end{proof}
\begin{cor}\label{answer intial question}
Assume that $S=Spec(A)$ is a noetherian regular local ring of dimension $n$ and let $\pi_1,\dots,\pi_n$be generators of the maximal ideal. Let $p:X=Spec(B)\rightarrow S=Spec(A)$ be a regular, flat, affine $S$-scheme of finite type. Let $\underline{\pi}:S\rightarrow \aff{n}{S}$ be the closed embedding associated to $\pi_1,\dots,\pi_n$. Then $\underline{\pi}\circ p$ is a regular sequence of $B$. Then the equivalence
\begin{equation}
    \Rlv_B(\sing{B,\underline{\pi}\circ p})\simeq q_*i_*\Phi^{\textup{mi}}_{(\proj{n-1}{B},W_{\underline{\pi}\circ p})}(\Ql{}(\beta))[-1]
\end{equation}
holds in $\md_{\Rlv_B(\sing{B,\underline{0}})}(\shvl(Spec(B)))$, where $q:\proj{n-1}{B}\rightarrow Spec(B)$ and $i:V(W_{\underline{\pi}\circ p})\rightarrow \proj{n-1}{B}$.
\end{cor}
\begin{rmk}
The corollary above generalizes to non-affine schemes as both members of the equivalence satisfy Zariski descent. This answers the question we asked ourselves at the beginning of this chapter.
\end{rmk}
\section{\textsc{Towards a vanishing cycles formalism over \texorpdfstring{$\ag{S}$}{ag}}}
The discussion in Section \ref{Monodromy-invariant vanishing cycles} suggests that it should be possible to construct a vanishing cycles formalism where the role of the base henselian trait is played by some more general geometric object. Moreover, the monodromy-invariant part of this construction should recover the sheaf of Definition \ref{Phimi}.
\subsection{\textsc{Tame vanishing cycles over \texorpdfstring{$\aff{1}{S}$}{A1}}}
It is well known that even if tame vanishing cycles where first defined for schemes over an henselian trait, what one actually uses is the geometry of $\aff{1}{S}$. Indeed, one can consider the following diagram
\begin{equation}
      \begindc{\commdiag}[18]
      \obj(-40,0)[1]{$S$}
      \obj(0,0)[2]{$\aff{1}{S}=Spec_S(\Oc_S[t])$}
      \obj(60,0)[3]{$\gm{S}$}
      \obj(120,0)[4]{$\varinjlim_{n\in \gm{S}(S)} \gm{S}[x]/(x^{n}-t)$}
      \mor{1}{2}{$i_0$}[\atleft,\injectionarrow]
      \mor{3}{2}{$j_0$}[\atright,\injectionarrow]
      \mor{4}{3}{$$}
      \enddc
\end{equation}

and use it to replace the one considered in \cite{sga7ii} when $S=Spec(A)$ is a strictly henselian trait
\begin{equation}
    \begindc{\commdiag}[15]
    \obj(0,0)[1]{$\sigma$}
    \obj(30,0)[2]{$S$}
    \obj(60,0)[3]{$\eta$}
    \obj(90,0)[7]{$\eta^t$}
    \mor{1}{2}{$$}[\atright,\injectionarrow]
    \mor{3}{2}{$$}[\atleft,\injectionarrow]
    \mor{7}{3}{$$}
    \cmor((90,4)(89,6)(87,7)(60,8)(33,7)(31,6)(30,4)) \pdown(60,14){$$}
    \enddc
\end{equation}

If we pullback the first diagram along the morphism $S\rightarrow \aff{1}{S}$ given by an uniformizer, we recover the second one. It is not surprising that in this way we are only able to recover the so called tame vanishing cycles. Indeed, the wild inertia group is of arithmetic nature.
One can define nearby cycles for schemes over $\aff{1}{S}$ in the usual way: for $f:X\rightarrow \aff{1}{S}$ consider the diagram
\begin{equation}
    \begindc{\commdiag}[18]
    \obj(-30,10)[1]{$X_0$}
    \obj(0,10)[2]{$X$}
    \obj(60,10)[3]{$X_{\mathcal{U}_n}$}
    \obj(150,10)[4]{$X_{\mathcal{U}_{\oo}}$}
    \obj(-30,-10)[5]{$S$}
    \obj(0,-10)[6]{$\aff{1}{S}$}
    \obj(60,-10)[7]{$\gm{S}[x]/(x^n-t)=\mathcal{U}_n$}
    \obj(150,-10)[8]{$\mathcal{U}_{\oo}=\varprojlim_{n\in \Oc_S^{\times}} \mathcal{U}_n$}
    \mor{1}{2}{$i$}
    \mor{3}{2}{$j_n$}[\atright,\solidarrow]
    \mor{4}{3}{$\phi_n$}[\atright,\solidarrow]
    \mor{1}{5}{$f_0$}[\atright,\solidarrow]
    \mor{2}{6}{$f$}[\atright,\solidarrow]
    \mor{3}{7}{$f_n$}
    \mor{4}{8}{$f_{\oo}$}
    \mor{5}{6}{$i_0$}[\atright,\solidarrow]
    \mor{7}{6}{$j_{0,n}$}
    \mor{8}{7}{$\phi_{0,n}$}
    \cmor((150,13)(150,17)(147,20)(75,20)(3,20)(0,17)(0,13)) \pdown(75,25){$j_{\oo}$}
    \cmor((150,-13)(150,-17)(147,-20)(75,-20)(3,-20)(0,-17)(0,-13)) \pup(75,-25){$j_{0,\oo}$}
    \enddc
\end{equation}
\begin{defn}
Let $\Fc \in \shvl(X)$. 
\begin{itemize}
    \item The $\ell$-adic sheaf of nearby cycles of $\Fc$ is defined as
\begin{equation}
    \Psi_f(\Fc):=i^*j_{\oo*}j_{\oo}^*\Fc \simeq \varinjlim_{n\in \Oc_S^{\times}} i^*j_{n*}j_n^*\Fc
\end{equation}
Notice that, since $\mu_{n,S}:=\Oc_S[x]/(x^n-1)$ naturally acts on $X_{\mathcal{U}_n}$, then $i^*j_{n*}j_n^*\Fc$ has a natural induced action. Then $\Psi_f(\Fc)$ has a natural action of $\varprojlim_{n\in \Oc_S^{\times}}\mu_{n,S}=:\mu_{\oo,S}$
\item There is a natural morphism $i^*\Fc \rightarrow \Psi_f(\Fc)$, that we can see as a $\mu_{\oo,S}$-equivariant morphism if we endow $i^*\Fc$ with the trivial action. The sheaf of vanishing cycles $\Phi_f(\Fc)$ of $\Fc$ is then defined as the cofiber of this morphism, with the induced $\mu_{\oo,S}$-action
\end{itemize}
\end{defn}
\subsection{\textsc{Tame vanishing cycles over} \texorpdfstring{$\ag{S}$}{ag}}
We can reproduce the situation described above for schemes over $\ag{S}$. In this case, the role of the zero section is played by $\bgm{S} \rightarrow \ag{S}$, and that of the open complementary by $S\simeq \gm{S}/\gm{S}\rightarrow \ag{S}$. We can consider elevation to the $n^{th}$ power in this case too:
\begin{equation}
    \Theta_n: \ag{S} \rightarrow \ag{S}
\end{equation}
which can be described as the functor
\begin{equation}
    (p:X\rightarrow S, \Lc_X, s_X)\mapsto (p:X\rightarrow S, \Lc_X^{\otimes n},s_X^{\otimes n})
\end{equation}
\begin{rmk}
Notice that the $\Theta_n$ are compatible, i.e.
\begin{equation}\label{compatibility Theta_n}
    \Theta_n\circ \Theta_m\simeq \Theta_{mn} \simeq \Theta_m\circ \Theta_n
\end{equation}
for any $n,m\in \N$
\end{rmk}
We can therefore consider the diagrams (for each $n\in \Oc_S^{\times}$)
\begin{equation}\label{base vanishing cycles ag}
    \begindc{\commdiag}[18]
    \obj(-70,-10)[1]{$\bgm{S}$}
    \obj(0,-10)[2]{$\ag{S}$}
    \obj(50,-10)[3]{$S$}
    \obj(-70,10)[4]{$\bgm{S} \times_{\ag{S}}\ag{S}$}
    \obj(0,10)[5]{$\ag{S}$}
    \obj(50,10)[6]{$S$}
    \obj(-140,10)[7]{$\bgm{S}$}
    \mor{1}{2}{$i_0$}[\atright,\solidarrow]
    \mor{3}{2}{$j_0$}
    \mor{4}{1}{$$}
    \mor{5}{2}{$\Theta_n$}
    \mor{6}{3}{$id$}
    \mor{4}{5}{$$}
    \mor{6}{5}{$j_{0,n}$}[\atright,\solidarrow]
    \mor{7}{4}{$t_{0,n}$}
    \mor{7}{1}{$\Theta_{0,n}$}[\atright,\solidarrow]
    \cmor((-140,14)(-140,18)(-138,20)(-70,20)(-2,20)(0,18)(0,14)) \pdown(-70,25){$i_{0,n}$}
    \enddc
\end{equation}
Notice that even if $t_{0,n}:\bgm{S} \rightarrow \bgm{S}\times_{\ag{S}}\ag{S}$ is not an equivalence, it shows $\bgm{S}\times_{\ag{S}}\ag{S}$ as a nilpotent thickening of $\bgm{S}$. More explicitly, objects of $\bgm{S}$ are pairs $(p:X\rightarrow S,\Lc_X)$, while those of $\bgm{S} \times_{\ag{S}}\ag{S}$ are triplets $(p:X\rightarrow S,\Lc_X,s_X)$, where $s_X$ is a $n$-torsion global section of $\Lc_X$. What is important for us is that pullback along $t_{0,n}$ induces an equivalence in étale cohomology. Moreover, $\Theta_n$ is a finite morphism of stacks. In particular, the base change theorem should be valid for the cartesian squares above. However, even if the $6$-functors formalism for the étale cohomology of stacks has been developed (\cite{lz15}, \cite{lz17i}, \cite{lz17ii}), the proper base change theorem has been proved in the case of representable morphisms. The morphisms $\Theta_n$ are not representable. For example, the pullback of $\Theta_n$ along the canonical atlas $S:\rightarrow \ag{S}$ is $\textup{B}\mu_n$.

For any $X$-point $(p:X\rightarrow S, \Lc_X,s_X)$ of $\ag{S}$ (e.g. for any twisted LG model over $S$), we can consider the following diagram, cartesian over (\ref{base vanishing cycles ag})
\begin{equation}
    \begindc{\commdiag}[18]
    \obj(-70,-10)[1]{$X_0$}
    \obj(0,-10)[2]{$X$}
    \obj(50,-10)[3]{$X_{\mathcal{U}}$}
    \obj(-70,10)[4]{$X_0 \times_{X}X_n$}
    \obj(0,10)[5]{$X_n$}
    \obj(50,10)[6]{$X_{\mathcal{U}}$}
    \obj(-140,10)[7]{$X_{0,n}$}
    \mor{1}{2}{$i$}[\atright,\solidarrow]
    \mor{3}{2}{$j$}
    \mor{4}{1}{$$}
    \mor{5}{2}{$\Omega_n$}
    \mor{6}{3}{$id$}
    \mor{4}{5}{$$}
    \mor{6}{5}{$j_{n}$}[\atright,\solidarrow]
    \mor{7}{4}{$t_n$}
    \mor{7}{1}{$\Omega_{0,n}$}[\atright,\solidarrow]
    \cmor((-140,14)(-140,18)(-138,20)(-70,20)(-2,20)(0,18)(0,14)) \pdown(-70,25){$i_{n}$}
    \enddc
\end{equation}
The same observations we made for the base diagram remain valid in this case too.
\begin{notation}
For a stack $Y$, and a torsion abelian group $\Lambda$ (we will only be interested in the case $\Lambda=\Z/\ell^n\Z$, $\ell \in \Oc_S^{\times}$) let $\shv(Y,\Lambda)$ the $\oo$-category of $\md_{\Lambda}$-valued étale sheaves on $Y$.
\end{notation}
We can then propose the following definition
\begin{defn}
Let $\Fc \in \shv(X,\Lambda)$. Then
\begin{equation}
    \Psi_{n}(\Fc):=i_n^*j_{n*}j^*\Fc\simeq i_n^*j_{n*}j_n^*\Omega_n^*\Fc \in \shv(X_{0,n},\Lambda)
\end{equation}
Moreover, let 
\begin{equation}
    sp_n:i_n^*\Omega_n^*\Fc\rightarrow \Psi_n(\Fc)
\end{equation}
be the morphism induced by the unit of $(j_n^*,j_{n*})$.
\end{defn}
We can then consider the images of these morphisms in $\shv(X_0,\Lambda)$
\begin{equation}
    \Omega_{0,n*}(sp_n:i_n^*\Omega_n^*\Fc\rightarrow \Psi_n(\Fc))
\end{equation}
Consider the diagram
\begin{equation}
    \begindc{\commdiag}[18]
    \obj(-50,30)[1]{$X_{0,nm}$}
    \obj(0,30)[2]{$X_0\times_X X_{nm}$}
    \obj(50,30)[3]{$X_{nm}$}
    \obj(100,30)[4]{$X_{\mathcal{U}}$}
    \obj(-50,0)[5]{$X_{0,n}$}
    \obj(0,0)[6]{$X_0\times_X X_{n}$}
    \obj(50,0)[7]{$X_{n}$}
    \obj(100,0)[8]{$X_{\mathcal{U}}$}
    \obj(0,-30)[9]{$X_{0}$}
    \obj(50,-30)[10]{$X$}
    \obj(100,-30)[11]{$X_{\mathcal{U}}$}
    \mor{1}{2}{$$}
    \mor{2}{3}{$$}
    \mor{4}{3}{$j_{nm}$}[\atright,\solidarrow]
    \mor{5}{6}{$$}
    \mor{6}{7}{$$}
    \mor{8}{7}{$j_n$}
    \mor{9}{10}{$i_0$}
    \mor{11}{10}{$j_0$}
    \mor{1}{5}{$\Omega_{0,m}^n$}[\atright,\solidarrow]
    \mor{2}{6}{$$}
    \mor{3}{7}{$\Omega_m^n$}
    \mor{4}{8}{$id$}
    \mor{5}{9}{$\Omega_{0,n}$}[\atright,\solidarrow]
    \mor{6}{9}{$$}
    \mor{7}{10}{$\Omega_n$}
    \mor{8}{11}{$id$}
    \cmor((-50,33)(-50,37)(-47,40)(0,40)(47,40)(50,37)(50,33)) \pdown(0,45){$i_{nm}$}
    \cmor((-60,30)(-70,30)(-75,25)(-75,0)(-70,-15)(-45,-30)(-7,-30)) \pright(-85,0){$\Omega_{0,nm}$}
    \cmor((-45,5)(-43,10)(-40,11)(0,11)(40,11)(43,10)(45,5)) \pdown(-10,15){$i_n$}
    \cmor((55,25)(60,24)(65,21)(66,0)(65,-21)(60,-24)(55,-25)) \pleft(73,10){$\Omega_{nm}$}
    \enddc
\end{equation}
Then we have
\begin{equation}
    \begindc{\commdiag}[18]
    \obj(-50,15)[1]{$sp_n:i_n^*\Omega_n^*\Fc$}
    \obj(50,15)[2]{$i_n^*j_{n*}j_n^*\Omega_n^*\Fc$}
    \obj(-60,-15)[3]{$\Omega_{0,m*}^n(sp_{nm}):\Omega_{0,m*}^ni_{nm}^*\Omega_{nm}^*\Fc$}
    \obj(50,-15)[4]{$\Omega_{0,m*}^ni_{nm}^*j_{nm*}j_{nm}^*\Omega_{nm}^*\Fc$}
    \mor{1}{2}{$$}
    \mor(-40,14)(-40,-14){$1\rightarrow \Omega_{m*}^n\Omega_m^{n*}$}
    \mor{2}{4}{$\simeq$}
    \mor{3}{4}{$$}
    \obj(-40,-25)[5]{$\simeq i_n^*\Omega_{m*}^n\Omega_m^{n*}\Omega_n^*\Fc$}
    \obj(50,-25)[6]{$\simeq i_n^*j_{n*}j_n^*\Omega_n^*\Fc$}
    \enddc
\end{equation}
We can then consider 
\begin{equation}
    \varinjlim_{n\in \N^{\times}}\Omega_{0,n*}i_n^*\Omega_n^*\Fc\rightarrow i_0^*j_{0*}j_0^*\Fc \in \shv(X_0,\Lambda)
\end{equation}
For $\Lambda=\Z/\ell^d\Z$ one can then consider the induced morphism on $\ell$-adic sheaves obtained by taking the limit over $d$ and then tensor with $\Ql{}$. It is expected that in this way one is able to recover monodromy invariant vanishing cycles. In order to explain way the first Chern class of the line bundle appears in the computation, it might be useful to look at the base diagram. The $\ell$-adic cohomology of $\bgm{S}$ is $\Ql{}[c_1]$, where $c_1$ is the universal first Chern class and lies in degree $2$. We expect that 
\begin{equation}
 \Ql{}\otimes_{\Z_{\ell}}\bigl (\varprojlim_{d}   \varinjlim_{n\in \N^{\times}}\Theta_{0,n*}i_n^*\Theta_n^* \Z/\ell^d\Z{}(\beta)\bigr )\simeq cofib(\Ql{}(\beta)\xrightarrow{c_1} \Ql{}(\beta))
\end{equation}
and one recovers $cofib(\Ql{,X_0}(\beta)\xrightarrow{c_1(\Lc_{X_0})}\Ql{,X_0}(\beta))$
by a formula of the kind
\begin{equation}
    \Ql{}[c_1]\simeq \textup{H}^*(\bgm{S},\Ql{})\rightarrow \textup{H}^*(X_0,\Ql{})
\end{equation}
\begin{equation*}
    c_1\mapsto c_1(\Lc_{X_0})
\end{equation*}
where $\Lc_{X_0}$ is the line bundle which determines the morphism $X_0\rightarrow \bgm{S}$

\section{\textsc{Some remarks on the regularity hypothesis}}\label{about the regularity hypothesis}
In the theorems about the $\ell$-adic realizaiton of the dg category of a (twisted, $n$-dimensional) LG model the regularity assumption on the ambient scheme is crucial. Indeed, we are not able to compute the motivic realization (and hence, the $\ell$-adic one) of $\cohbp{X_0}{X}$. However, this dg category sits in the following pullback diagram
\begin{equation}\label{pullback cohbp}
    \begindc{\commdiag}[18]
    \obj(-30,15)[1]{$\cohbp{X_0}{X}$}
    \obj(30,15)[2]{$\perf{X}_{X_0}$}
    \obj(-30,-15)[3]{$\cohb{X_0}$}
    \obj(30,-15)[4]{$\cohb{X}_{X_0}$}
    \mor{1}{2}{$\mathfrak{i}_*$}
    \mor{1}{3}{$$}[\atright, \injectionarrow]
    \mor{2}{4}{$$}[\atright, \injectionarrow]
    \mor{3}{4}{$\mathfrak{i}_*$}
    \enddc
\end{equation}
It is well known after Quillen's dévissage for G-theory that 
\begin{equation}
    \Mv(\cohb{X_0})\simeq \Mv(\cohb{X}_{X_0})
\end{equation}
Therefore, to say that $\Mv(\mathfrak{i}_*:\cohbp{X_0}{X}\rightarrow \perf{X}_{X_0})$ is an equivalence or to say that the image of square (\ref{pullback cohbp}) via $\Mv$ is still a pullback square are equivalent statements. We believe that this is the case, even though this matter will be investigated elsewhere. Notice that to know that such an equivalence holds true would allow to compute the motivic realization of $\cohbp{X_0}{X}$ (under the additional hypothesis $X_0\simeq \pi_0(X_0)$). Indeed, the localization sequence for K-theory would allow us to compute it as the fiber
\begin{equation}
    \Mv_X(\cohbp{X_0}{X})\simeq \Mv_X(\perf{X}_{X_0})\rightarrow \Mv_X(\perf{X})\rightarrow \Mv_X(\perf{X-X_0})
\end{equation}
In other words,
\begin{equation}
    \Mv_X(\cohbp{X_0}{X})\simeq i_*i^!\bu_X
\end{equation}
In this case, the proof of Theorem \ref{main theorem} would work without any changes.

%% file: chapters/chapter04.tex
In this chapter we will investigate the situation analogous to \cite{brtv}, but instead of considering an excellent strictly henselian trait as a base, we will consider a discrete valuation ring of rank $2$.
\section{\textsc{What is in this chapter}}
Always bearing \cite{brtv} in mind as a model, we consider the situation in which we are given a proper, flat regular scheme over a discrete valuation ring of rank $2$. In this case, we don't have just two varieties (the special and the generic fiber) that appear in the rank $1$ setting. Instead, there are $3$ varieties: the special fiber lying over the unique closed point of the base ring, the generic fiber lying over the unique open point of the base, and a third variety that lives in the middle. Another way to picture this situation is the following: a discrete valuation ring of rank $2$ is a way to glue together two rank $1$ discrete valuation rings such that the generic point of the first coincide with the closed point of the second. Then the variety in the middle is just the fiber over this special point.

We propose a definition of vanishing cycles in this setting. As the generic point of the base ring is $\C \drl x \drr \drl y \drr$, these $\ell$-adic sheaves come equipped with an action of its absolute Galois group, which is isomorphic to $\hat{\Z}\times \hat{\Z}$. We then dedicate a consistent part of this chapter to study this action, and to understand what happens if one applies the (homotopy) fixed points $\oo$-functor. We think that this has its own interest. 

In the last part of this section we investigate the connection between the vanishing cycles sheaves we have defined and the $\ell$-adic realization of certain dg categories of singularities. For the same reasons explained in Section \ref{about the regularity hypothesis} we will need to assume that the scheme is regular. Moreover, we will need to assume that the $\ell$-adic realization $\oo$-functor that we have considered so far is available in the non-noetherian situation too and that the results in \cite[\S B.4]{pr11} extend to our context. 

In the end we are able to prove
\begin{theorem-non}{\textup{\ref{theorem chpt 4}}}
Let $X$ be a regular scheme, proper and flat over $S=Spec(\C\dsl x \dsr + y\cdot \C\drl x \drr \dsl y \dsr)$. Let $Y$ be the fiber over $Spec(\C\dsl x \dsr + y\cdot \C\drl x \drr \dsl y \dsr/<y>)$. Then there is a fiber-cofiber sequence in $\shvl(\C)$
\begin{equation*}
    i_{0}^*\Rlv_{S}(\sing{Y})\rightarrow p_{\C*}\Phi^{(2)}(\Ql{}(\beta))[-1]^{\textup{h}Gal(\bar{\eta}/\eta)}\rightarrow cofib(\xi)
\end{equation*}
where $i_{01}:Spec(\C)\rightarrow S$, $\Phi^{(2)}(\Ql{}(\beta))$ is the $\ell$-adic sheaf of second order vanishing cycles (see Definition \ref{def iterated vanishing cycles}) and $\xi$ is an explicit morphism.
\end{theorem-non}
\section{\textsc{The base ring}}
\begin{rmk}
For more details on higher valuation rings, we refer the reader to \cite{mor}.
\end{rmk}
In this section we will briefly introduce the discrete valuation ring of dimension $2$ we will deal with.

Consider the valuation
\begin{equation}
v: \C \drl x \drr \drl y \drr^{\x}\rightarrow \Z \x \Z
\end{equation}
\[
\sum_{i=n}^{\oo} \bigl ( \sum_{j=m_i}^{\oo} a_{ij}x^j\bigr )y^i \mapsto (n,m_n)
\]
and put the lexicographical order on $\Z \x \Z$: $(a,b) < (a',b')$ if and only if $a<a'$ or $a=a'$ and $b<b'$. Then, we can consider the ring of integers of $v$:
\begin{equation}
\mathcal{O}_v := \{ a\in \C\drl x\drr \drl y \drr^{\x} : v(a)\geq (0,0) \} \cup \{ 0 \}
\end{equation}
\begin{notation}
From now on, set $R:= \mathcal{O}_v$.
\end{notation}
\begin{rmk}
It is not hard to see that $R= \C\dsl x \dsr + y\cdot \C\drl x \drr \dsl y \dsr \subseteq \C\drl x\drr \drl y \drr$.
\end{rmk}
It is known that $R$ is a regular non-Noetherian ring of Krull dimension $2$. Indeed, this ring has three prime ideals $m\supset p \supset (0)$, where
\begin{equation}
m=x\cdot R \hspace{0.5cm} p=y\cdot \C\drl x \drr \dsl y \dsr =<y\cdot x^i,i\in \Z>
\end{equation}
Then, it is possible to compute
\begin{equation}
R/m\simeq \C  \hspace{0.5cm} R/p= \C \dsl x \dsr
\end{equation}
Put $S:= Spec(R)$, $S_1:= Spec(R/p)$ and $\sigma:= Spec(R/m)$. The open subschemes $S-\sigma$ and $S-S_1$ coincide with $Spec (\C \drl x \drr \dsl y \dsr)$ and $Spec(\C \drl x \drr \drl y \drr)$ respectively.
Thus, we have the following diagram
\begin{equation}
\begindc{\commdiag}[20]
\obj(0,30)[1]{$\sigma$}
\obj(30,30)[2]{$S_1$}
\obj(70,30)[3]{$\varepsilon=Spec(\C \drl x \drr)$}
\obj(30,0)[4]{$S$}
\obj(70,0)[5]{$U_0 := S-\sigma$}
\obj(30,-30)[6]{$\eta := S-S_1$}
\mor{1}{2}{$i_{01}$}[\atleft,\injectionarrow]
\mor{3}{2}{$l$}[\atright, \injectionarrow]
\mor{1}{4}{$i_0$}[\atright, \injectionarrow]
\mor{2}{4}{$i_1$}[\atleft,\injectionarrow]
\mor{3}{5}{$k$}[\atleft,\injectionarrow]
\mor{6}{4}{$j_1$}[\atleft,\injectionarrow]
\mor(39,-27)(70,-3){$j_{10}$}[\atright, \injectionarrow]
\mor{5}{4}{$j_0$}[\atleft,\injectionarrow]
\enddc
\end{equation}
where the square is cartesian.
\section{\textsc{Iterated vanishing cycles}}

\subsection{\textsc{First and second order vanishing cycles}}
Consider the scheme $S$ we introduced in the previous section. This comes with a flag of closed sub-schemes $\sigma \subset S_1\subset S$ which we will use to construct an analogue to classical nearby cycles.
\begin{notation}
\begin{itemize}
\item 
We fix a separable closure $\overline{\C \drl x \drr}$ of $\C \drl x \drr$ and we denote its spectrum by $\bar{\varepsilon}$. This coincides with an algebraic closure of $\C \drl x \drr$ as the field is perfect
\item 
We fix an algebraic closure $\overline{\C \drl x \drr \drl y\drr}$ of $\C \drl x \drr \drl y\drr$ and we denote its spectrum by $\bar{\eta}$
\item We consider the maximal unramified extension $\C \drl x \drr \drl y\drr^{mu}$ of $\C \drl x \drr \drl y\drr$ inside $\overline{\C \drl x \drr \drl y\drr}$ and we denote its spectrum by $\eta^{mu}$
\end{itemize}
\end{notation}
Therefore, we have the following diagram
\begin{equation}\label{base diagram}
\begindc{\commdiag}[20]
\obj(0,30)[1]{$\sigma$}
\obj(30,30)[2]{$S_1$}
\obj(60,30)[3]{$\varepsilon$}
\obj(90,30)[4]{$\bar{\varepsilon}$}
\obj(30,0)[5]{$S$}
\obj(60,0)[6]{$U_0=S-\sigma$}
\obj(90,0)[7]{$\overline{U_0}$}
\obj(30,-30)[8]{$\eta$}
\obj(60,-30)[9]{$\eta^{mu}$}
\obj(90,-30)[10]{$\bar{\eta}$}
\mor{1}{2}{$i_{01}$}[\atleft,\injectionarrow]
\mor{1}{5}{$i_0$}[\atright,\injectionarrow]
\mor{2}{5}{$i_1$}[\atleft,\injectionarrow]
\mor{3}{6}{$k$}[\atleft,\injectionarrow]
\mor{4}{7}{$\bar{k}$}[\atleft,\injectionarrow]
\mor{3}{2}{$l$}[\atright,\injectionarrow]
\mor{6}{5}{$j_{0}$}[\atright,\injectionarrow]
\mor{8}{5}{$j_{1}$}[\atright,\injectionarrow]
\mor{8}{6}{$j_{10}$}[\atright,\injectionarrow]
\mor{9}{7}{$j_{10}^{mu}$}[\atright,\injectionarrow]
\mor{4}{3}{$$}
\mor{7}{6}{$$}
\mor{9}{8}{$$}
\mor{10}{9}{$$}
\mor{10}{7}{$\bar{j_{10}}$}[\atright,\solidarrow]
\cmor((90,34)(90,38)(88,40)(60,40)(32,40)(30,38)(30,34)) \pdown(60,43){$\bar{l}$}
\cmor((88,5)(88,8)(86,10)(60,10)(34,10)(32,8)(32,5)) \pdown(75,14){$\bar{j}_0$}
\cmor((90,-34)(90,-38)(88,-40)(66,-40)(55,-40)(44,-40)(33,-40)(22,-40)(21,-39)(20,-38)(20,-28)(20,-15)(20,-2)(21,-1)(22,0)(24,0)(26,0)) \pright(60,-45){$\bar{j}_1$}
\enddc
\end{equation}
Here $\overline{U}_0$ is a strict henselization of $U_0$. It is easy to see that $\overline{U}_0=Spec(\overline{\C \drl x \drr} \dsl y \dsr)$.

By the classical theory of henselian rings, we have the following diagram of Galois groups
\begin{equation}
\begindc{\commdiag}[9]
\obj(100,60)[1]{$1$}
\obj(-100,30)[2]{$1$}
\obj(0,30)[3]{$Gal(\bar{\varepsilon}/\varepsilon^{mu})$}
\obj(100,30)[4]{$Gal(\bar{\varepsilon}/\varepsilon)$}
\obj(200,30)[5]{$Gal(\bar{\sigma}/\sigma)$}
\obj(300,30)[6]{$1$}
\obj(100,0)[7]{$Gal(\bar{\eta}/\eta)$}
\obj(100,-30)[8]{$Gal(\bar{\eta}/\eta^{mu})$}
\obj(100,-60)[9]{$1$}
\mor{2}{3}{$$}
\mor{3}{4}{$$}
\mor{4}{5}{$$}
\mor{5}{6}{$$}
\mor{9}{8}{$$}
\mor{8}{7}{$$}
\mor{7}{4}{$$}
\mor{4}{1}{$$}
\enddc
\end{equation}
where the sequences are exact.
It is moreover possible to be more precise. Indeed, we know that every finite extension of $\C \drl x \drr$ is of the form $\C \drl x^{\frac{1}{n}}\drr$ for some $n\in \N$ and, similarly, every finite extension of $\C \drl x \drr \drl y \drr$ is of the form $\C \drl x^{\frac{1}{m}} \drr \drl y^{\frac{1}{n}} \drr$ for some $m,n \in \N$. Moreover, unramified extensions of $\C \drl x \drr \drl y \drr$ correspond to separable extensions of $\C \drl x \drr$, the correspondence being
\begin{equation}
\C \drl x^{\frac{1}{m}}\drr  \mapsto \C \drl x \drr \drl y \drr \otimes_{\C\drl x \drr} \C \drl x^{\frac{1}{m}} \drr \simeq \C \drl x^{\frac{1}{m}} \drr \drl y \drr
\end{equation}
In particular, we obtain that $\C \drl x \drr \drl y \drr ^{mu}\simeq \overline{\C \drl x \drr}\drl y \drr$. Let us now give a look at the Galois groups of these field extensions. Fix a separable extension $\C \drl x \drr \drl y \drr \subset \C \drl x^{\frac{1}{m}}\drr \drl y^{\frac{1}{n}} \drr$. It is not hard to see that an element in the Galois group of this field extension is given by the $\C \drl x \drr \drl y \drr$-automorphism of $\C \drl x \drr \drl y \drr \subset \C \drl x^{\frac{1}{m}}\drr \drl y^{\frac{1}{n}} \drr$ uniquely determined by the assigments
\[
x^{\frac{1}{m}}\mapsto \zeta_{m}x^{\frac{1}{m}} \hspace{1cm} y^{\frac{1}{n}}\mapsto \zeta_{n}x^{\frac{1}{n}}
\]
where $\zeta_m$ and $\zeta_n$ are an $m^{th}$-root and an $n^{th}$-root of the unity respectively. This implies that 
\begin{equation}
Gal\bigl ( \C \drl x^{\frac{1}{m}} \drr \drl y^{\frac{1}{n}} \drr/\C \drl x \drr \drl y \drr \bigr )\simeq \Z/m\Z \x \Z/n\Z
\end{equation}
and therefore that 
\begin{equation}
Gal\bigl ( \overline{\C \drl x \drr \drl y \drr}/\C \drl x \drr \drl y \drr \bigr )=\varprojlim_{m,n} Gal\bigl ( \C \drl x^{\frac{1}{m}} \drr \drl y^{\frac{1}{n}} \drr/\C \drl x \drr \drl y \drr \bigr ) \simeq \hat{\Z}\x \hat{\Z}
\end{equation}
Similarly,
\begin{equation}
Gal\bigl (\overline{ \C \drl x \drr} \drl y^{\frac{1}{n}} \drr/\overline{\C \drl x \drr} \drl y \drr \bigr )\simeq \Z/n\Z
\end{equation}
and
\begin{equation}
Gal\bigl ( \overline{\C \drl x \drr \drl y \drr}/\overline{\C \drl x \drr} \drl y \drr \bigr )=\varprojlim_n Gal\bigl (\overline{ \C \drl x \drr} \drl y^{\frac{1}{n}} \drr/\overline{\C \drl x \drr} \drl y \drr \bigr )\simeq \hat{\Z}
\end{equation}
Finally, it is easy to see that the canonical morphism
\[
Gal\bigl ( \overline{\C \drl x \drr \drl y \drr}/\overline{\C \drl x \drr} \drl y \drr \bigr ) \rightarrow Gal\bigl ( \overline{\C \drl x \drr \drl y \drr}/\C \drl x \drr \drl y \drr \bigr )
\]
corresponds, under these isomorphism to
\begin{equation}
\hat{\Z}\rightarrow \hat{\Z}\x \hat{\Z}
\end{equation}
\[
\alpha \mapsto (0,\alpha)
\]
In particular, the vertical short exact sequence of groups in the diagram above is split exact.

Consider an $S$-scheme $p: X\rightarrow S$ and the diagram
\begin{equation}
\begindc{\commdiag}[20]
\obj(0,30)[1]{$X_{\sigma}$}
\obj(30,30)[2]{$X_1$}
\obj(60,30)[3]{$X_{\varepsilon}$}
\obj(90,30)[4]{$X_{\bar{\varepsilon}}$}
\obj(30,0)[5]{$X$}
\obj(60,0)[6]{$X_{U_0}$}
\obj(90,0)[7]{$X_{\overline{U_0}}$}
\obj(30,-30)[8]{$X_{\eta}$}
\obj(60,-30)[9]{$X_{\eta^{mu}}$}
\obj(90,-30)[10]{$X_{\bar{\eta}}$}
\mor{1}{2}{$i_{01}$}[\atleft,\injectionarrow]
\mor{1}{5}{$i_0$}[\atright,\injectionarrow]
\mor{2}{5}{$i_1$}[\atleft,\injectionarrow]
\mor{3}{6}{$k$}[\atleft,\injectionarrow]
\mor{4}{7}{$\bar{k}$}[\atleft,\injectionarrow]
\mor{3}{2}{$l$}[\atright,\injectionarrow]
\mor{6}{5}{$j_{0}$}[\atright,\injectionarrow]
\mor{8}{5}{$j_{1}$}[\atright,\injectionarrow]
\mor{8}{6}{$j_{10}$}[\atright,\injectionarrow]
\mor{9}{7}{$j_{10}^{mu}$}[\atright,\injectionarrow]
\mor{4}{3}{$v_{\varepsilon}$}
\mor{7}{6}{$v_{U_0}$}
\mor{9}{8}{$v_{\eta}$}
\mor{10}{9}{$$}
\mor{10}{7}{$\bar{j_{10}}$}[\atright,\solidarrow]
\cmor((90,34)(90,38)(88,40)(60,40)(32,40)(30,38)(30,34)) \pdown(60,43){$\bar{l}$}
\cmor((88,5)(88,8)(86,10)(60,10)(34,10)(32,8)(32,5)) \pdown(75,14){$\bar{j}_0$}
\cmor((90,-34)(90,-38)(88,-40)(66,-40)(55,-40)(44,-40)(33,-40)(22,-40)(21,-39)(20,-38)(20,-28)(20,-15)(20,-2)(21,-1)(22,0)(24,0)(26,0)) \pright(60,-45){$\bar{j}_1$}
\enddc
\end{equation}
cartesian over (\ref{base diagram}).
\begin{notation}
We will add subscrits when we want to consider the morphism induced by $p$ under pullback, e.g. we will write $p_{\sigma}:X_{\sigma}\rightarrow \sigma$.
\end{notation}
\begin{defn}
Fix an object $\Fc$ in $\shvl(X)$.
\begin{itemize}
\item The first order nearby cycles of $\Fc$ relative to $p$ is the object
\begin{equation}
\Psi^{(1)}_p(\Fc):=i_{01}^*\bar{l}_*\Fc_{|X_{\bar{\varepsilon}}} \in \shvl(X_{\sigma})^{Gal(\bar{\eta}/\eta)}
\end{equation}
\item The second order nearby cycles of $\Fc$ relative to $p$ is the object 
\begin{equation}
\Psi^{(2)}_p(\Fc):= i_{01}^*\bar{l}_*\bar{k}^*\bar{j}_{10*}\Fc_{|X_{\bar{\eta}}}  \in \shvl(X_{\sigma})^{Gal(\bar{\eta}/\eta)}
\end{equation}
\end{itemize}
Here  $ \shvl(X_{\sigma})^{Gal(\bar{\eta}/\eta)}$ is the stable $\oo$-category of $\ell$-adic sheaves on $X_{\sigma}$ equipped with an equivariant structure with respect to the continuous action of $Gal(\bar{\eta}/\eta)\simeq \hat{\Z} \x \hat{\Z}$ on $X_{\sigma}$ (i.e. the trivial action).
\end{defn}
\begin{rmk}
The first order nearby cycle object of an $\ell$-adic sheaf $\Fc$ is intimately related to the classical notion of nearby cycle object of $\Fc_{|X_{1}}$ relative to $p_1 : X_{1}\rightarrow S_1$. The latter is defined by
\begin{equation}
\Psi_{p_1}(\Fc_{|X_{1}})= i_{01}^*\bar{l}_*\Fc_{|X_{\bar{\varepsilon}}} \in \shvl(X_{\sigma})^{Gal(\bar{\varepsilon}/\varepsilon)}
\end{equation}
(see for example \cite[Definition 4.5]{brtv}). There is an adjunction
\begin{equation}
\begindc{\commdiag}[20]
\obj(-15,0)[1]{$\shvl(X_{\sigma})^{Gal(\bar{\varepsilon}/\varepsilon)}$} 
\obj(55,0)[2]{$\shvl(X_{\sigma})^{Gal(\bar{\eta}/\eta)}$}
\mor(3,3)(35,3){$trivial$}
\mor(35,-3)(3,-3){$\textup{h}Gal(\bar{\eta}/\eta^{mu})$}
\enddc
\end{equation}  
where the left adjoint sends an object in
$\shvl(X_{\sigma})^{Gal(\bar{\varepsilon}/\varepsilon)}$ to
the same sheaf, with the action of
$Gal(\bar{\eta}/\eta)\simeq Gal(\bar{\eta}/\eta^{mu})\x Gal(\bar{\varepsilon}/\varepsilon)$ 
where the first component acts trivially and the action of
the second component agrees with the starting one. Then, the 
first order nearby cycle object of $\Fc\in \shvl(X)$ relative to $p$ coincides with the image of $\Psi_{p_1}(\Fc_{|X_1})$ under the left adjoint functor.
\end{rmk}
Notice that we always have an equivariant morphism
\begin{equation}
\Psi^{(1)}_p(\Fc)=i_{01}^*\bar{l}_* \Fc_{|X_{\bar{\varepsilon}}}\simeq i_{01}^*\bar{l}_*\bar{k}^*\bar{j}_0^*(\Fc)\rightarrow i_{01}^*\bar{l}_*\bar{k}^*\bar{j}_{10*}\bar{j}_{10}^*\bar{j}_0^*(\Fc)\simeq \Psi^{(2)}_p(\Fc)
\end{equation}
induced by the counit of the adjuction $(\bar{j}_{10}^*,\bar{j}_{10*})$

Also notice that the image of the specialization morphism
\begin{equation}
\Fc_{|X_{\sigma}}\rightarrow \Psi_{p_1}(\Fc_{|X_{1}}) \in \shvl(X_{\sigma})^{Gal(\bar{\varepsilon}/\varepsilon)}
\end{equation}
induces a $Gal(\bar{\eta}/\eta)$-equivarant morphism 
\begin{equation}
\Fc_{|X_{\sigma}}\rightarrow \Psi^{(1)}_p(\Fc)
\end{equation}
via the $\oo$-functor considered above. Therefore, composing these two morphisms, we end up with a $Gal(\bar{\eta}/\eta)$-equivariant triangle
\begin{equation}\label{triangle nearby cycles}
\begindc{\commdiag}[20]
\obj(0,15)[1]{$\Fc_{|X_{\sigma}}$}
\obj(50,15)[2]{$\Psi^{(1)}_p(\Fc)$}
\obj(25,-15)[3]{$\Psi^{(2)}_p(\Fc)$}
\mor{1}{2}{$$}
\mor{1}{3}{$$}
\mor{2}{3}{$$}
\enddc
\end{equation}
\begin{rmk}

The $\oo$-category $\shvl(S)$ is the recollement of $\shv(S_1)$ and $\shvl(\eta)$. The former is itself the recollement of $\shvl(\sigma)$ and $\shvl(\varepsilon)$. Moreover, using the equivalences
\begin{equation}
\shvl(\varepsilon)\simeq \shvl(\bar{\varepsilon})^{Gal(\bar{\varepsilon}/\varepsilon)} \hspace{1cm} \shvl(\eta)\simeq \shvl(\bar{\eta})^{Gal(\bar{\eta}/\eta)}
\end{equation}
we obtain a description of $\shvl(S)$ where objects consist of triangles
\begin{equation}
\begindc{\commdiag}[20]
\obj(0,30)[1]{$\underbracket{\Fc_0}_{\in \shvl(\sigma)}$}
\obj(60,30)[2]{$\underbracket{\Fc_{\varepsilon}}_{\in \shvl(\bar{\varepsilon})^{Gal(\bar{\varepsilon}/\varepsilon)}}$}
\obj(30,0)[3]{$\underbracket{\Fc_{\eta}}_{\in \shvl(\bar{\eta})^{Gal(\bar{\eta}/\eta)}}$}
\mor{1}{2}{$$}
\mor{2}{3}{$$}
\mor{1}{3}{$$}
\enddc
\end{equation}
where the morphisms are equivariant.

If we consider a $\C$-scheme of finite type $X_0$, then we can consider the product of its étale topos with the étale topos of $S$, $\tilde{X}_0 \x \tilde{S}$. If we consider the stable $\oo$-category of $\ell$-adic sheaves on it, $\shvl(\tilde{X}_0 \x \tilde{S})$, we have a description of its objects as triangles
\begin{equation}
\begindc{\commdiag}[20]
\obj(0,30)[1]{$\underbracket{\Fc_0}_{\in \shvl(X_0)}$}
\obj(60,30)[2]{$\underbracket{\Fc_{\varepsilon}}_{\in \shvl(X_0)^{Gal(\bar{\varepsilon}/\varepsilon)}}$}
\obj(30,0)[3]{$\underbracket{\Fc_{\eta}}_{\in \shvl(X_0)^{Gal(\bar{\eta}/\eta)}}$}
\mor{1}{2}{$$}
\mor{2}{3}{$$}
\mor{1}{3}{$$}
\enddc 
\end{equation}
where the morphisms are again equivariant.


Then we can view the assignment $\Fc\in \shvl(X)\mapsto$ triangle (\ref{triangle nearby cycles}) as part of an $\oo$-functor
\begin{equation}
\shvl(X)\rightarrow \shvl(\tilde{X_{\sigma}}\x \tilde{S})
\end{equation}
\end{rmk}

\begin{defn}\label{def iterated vanishing cycles}
We define the following $\ell$-adic sheaves in $\shvl(X_{\sigma})^{Gal(\bar{\eta}/\eta)}$
\begin{itemize}
    \item $\Phi^{(1)}_p(\Fc):=cofib( \Fc_{X_{\sigma}}\rightarrow \Psi^{(1)}_p(\Fc))$ (first order vanishing cycles)
    \item $\Phi^{(2)}_p(\Fc):=cofib( \Fc_{X_{\sigma}}\rightarrow \Psi^{(2)}_p(\Fc))$ (second order vanishing cycles)
\end{itemize}
\end{defn}
\subsection{\textsc{The action of the torus \texorpdfstring{$\eta$}{eta}}}
\begin{defn}
\begin{itemize}
\item We will refer to 
$\Hl(\varepsilon):=i_{01}^*l_*\mathbb{Q}_{\ell,\varepsilon}$ as the cohomology of the punctured disk
\item we will refer to $\Hl(\eta):=i_{01}^*l_*k^*j_{10*}\mathbb{Q}_{\ell,\eta} \in \shvl(\sigma)$ as the cohomology of the torus
\end{itemize}
\end{defn}
\begin{lmm}

The objects $\Hl(\varepsilon)$ and $\Hl(\eta)$ are commutative algebras in $\shvl(\sigma)$ in a canonical way.\\

Moreover, there are equivalences
\begin{equation}\label{cohomology punctured disk}
\Hl(\varepsilon)\simeq \mathbb{Q}_{\ell,\sigma}\oplus \mathbb{Q}_{\ell,\sigma}
[-1]
\end{equation}
\begin{equation}\label{cohomology torus}
\Hl(\eta)\simeq  \mathbb{Q}_{\ell,\sigma}\oplus \mathbb{Q}_{\ell,\sigma}
[-1] \oplus  \mathbb{Q}_{\ell,\sigma}
[-1]\oplus \mathbb{Q}_{\ell,\sigma}
[-2]
\end{equation}
The algebra structures on the right hand sides are the usual ones (i.e. $Sym_{\Ql{,\sigma}}(\Ql{,\sigma}[-1])$ and $Sym_{\Ql{,\sigma}}(\Ql{,\sigma}^{\oplus 2}[-1])$)
\end{lmm}
\begin{proof}

This is \cite[Lemma 4.16]{brtv}. For what concerns the commutative algebra structures, this is a consequence of the fact that both $i_{01}^*l_*$ and $k^*j_{10_*}$ are lax monoidal $\oo$-functors. For the equivalences in the statement, the first one follows from the computation of \cite[Lemma 4.16]{brtv}. The second one, follows by a double application of $loc.cit$\footnote{we ignore Tate shifts as we are working over $\C$}:
\begin{equation}
\Hl(\eta)\simeq i_{01}^*l_*\bigl ( \mathbb{Q}_{\ell,\varepsilon}\oplus \mathbb{Q}_{\ell,\varepsilon}
[-1]\bigr )
\end{equation}
\begin{equation*}\simeq  \mathbb{Q}_{\ell,\sigma}\oplus \mathbb{Q}_{\ell,\sigma}
[-1] \oplus  \mathbb{Q}_{\ell,\sigma}
[-1]\oplus \mathbb{Q}_{\ell,\sigma}
[-2]
\end{equation*}
where the last equivalence follows from the compatibility of $i_{01}^*l_*$ with Grothendieck's six functors formalism.
The last assertion is a consequence of the fact that, since we are working over a field of characteristic zero, there is an unique cdga whose cohomology algebra is free.
\end{proof}

For a proper, flat $S$-scheme of finite presentation $p: X\rightarrow S$, we can define
\begin{equation}
\Hl(X_{\varepsilon}):= p_{\sigma *}i_{01}^*l_*\mathbb{Q}_{\ell,X_{\varepsilon}} \in \shvl(\sigma)
\end{equation}
\begin{equation}
\Hl(X_{\eta}):= p_{\sigma *}i_{01}^*l_*k^*j_{10 *}\mathbb{Q}_{\ell,X_{\eta}} \in \shvl(\sigma)
\end{equation}
Notice that proper base change implies that the object $\Hl(X_{\varepsilon})$ is endowed with a natural map of algebras
\begin{equation}
\Hl(\varepsilon)\rightarrow \Hl(X_{\varepsilon})
\end{equation}
and the object $\Hl(X_{\eta})$ is endowed with a natural map of algebras
\begin{equation}
\Hl(\eta)\rightarrow \Hl(X_{\eta})
\end{equation}

Moreover, both $\Hl(X_{\varepsilon})$ and $\Hl(X_{\eta})$ are also algebras over $p_{\sigma *}\mathbb{Q}_{\ell,X_{\sigma}}$. This gives us two morphisms of algebras

\begin{equation}
\Hl(\varepsilon)\otimes_{\mathbb{Q}_{\ell,\sigma}} p_{\sigma *}\mathbb{Q}_{\ell,X_{\sigma}} \rightarrow \Hl(X_{\varepsilon})
\end{equation}

\begin{equation}
\Hl(\eta)\otimes_{\mathbb{Q}_{\ell,\sigma}} p_{\sigma *}\mathbb{Q}_{\ell,X_{\sigma}} \rightarrow \Hl(X_{\eta})
\end{equation}

\subsection{\textsc{Homotopy fixed points}}
In this section we study the $Gal(\bar{\eta}/\eta)$-homotopy fixed points of the $\ell$-adic sheaves defined above. Indeed, the nearby cycles we have defined naturally live in the symmetric monoidal stable presentable $\oo$-category $\shvl(\sigma)^{Gal(\bar{\eta}/\eta)}$ of $\ell$-adic complexes with a continuous action of the Galois group $Gal(\bar{\eta}/\eta)$. First of all, notice that we have the following adjunctions
\begin{equation}
\begindc{\commdiag}[20]
\obj(0,15)[1]{$\shvl(\sigma)$}
\obj(80,15)[2]{$\shvl(\sigma)^{Gal(\bar{\eta}/\eta)}$}
\obj(40,-15)[3]{$\shvl(\sigma)^{Gal(\bar{\varepsilon}/\varepsilon)}$}
\mor(8,17)(62,17){$trivial$}
\mor(62,13)(8,13){$(-)^{\textup{h}Gal(\bar{\eta}/\eta)}$}
\mor(3,10)(32,-12){$trivial$}
\mor(25,-12)(-4,10){$(-)^{\textup{h}Gal(\bar{\varepsilon}/\varepsilon)}$}
\mor(48,-12)(72,10){$trivial$}
\mor(79,10)(55,-12){$(-)^{\textup{h}Gal(\bar{\eta}/\eta^{mu})}$}
\enddc
\end{equation}
In the diagram above the functors \textit{trivial}, that send an $\ell$-adic sheaf to itself considered with the trivial action of the relevant profinite group, are left adjoints. Since it is obvious that the composition of the two oblique \textit{trivial} functors gives the third one, by the uniqueness (up to equivalence) of the adjoint functors we get a natural equivalence of $\oo$-functors
\begin{equation}\label{equiv hgal 1}
\bigl ( (-)^{\textup{h}Gal(\bar{\eta}/\eta^{mu})}\bigr )^{\textup{h}Gal(\bar{\varepsilon}/\varepsilon)}\simeq (-)^{\textup{h}Gal(\bar{\eta}/\eta)} : \shvl(\sigma)^{Gal(\bar{\eta}/\eta)}\rightarrow \shvl(\sigma)
\end{equation}
Also notice that we have a second factorization of the $\oo$-functor $\textup{h}Gal(\bar{\eta}/\eta)$, given by the following digram:
\begin{equation}
\begindc{\commdiag}[20]
\obj(0,15)[1]{$\shvl(\sigma)$}
\obj(80,15)[2]{$\shvl(\sigma)^{Gal(\bar{\eta}/\eta)}$}
\obj(40,-15)[3]{$\shvl(\sigma)^{Gal(\bar{\eta}/\eta^{mu})}$}
\mor(8,17)(62,17){$trivial$}
\mor(62,13)(8,13){$(-)^{\textup{h}Gal(\bar{\eta}/\eta)}$}
\mor(3,10)(32,-12){$trivial$}
\mor(25,-12)(-4,10){$(-)^{\textup{h}Gal(\bar{\eta}/\eta^{mu})}$}
\mor(48,-12)(72,10){$trivial$}
\mor(79,10)(55,-12){$(-)^{\textup{h}Gal(\bar{\varepsilon}/\varepsilon)}$}
\enddc
\end{equation}
This gives us the natural equivalence of $\oo$-functors
\begin{equation}\label{equiv hgal 2}
\bigl ( (-)^{\textup{h}Gal(\bar{\varepsilon}/\varepsilon))}\bigr )^{\textup{h}Gal(\bar{\eta}/\eta^{mu})}\simeq (-)^{\textup{h}Gal(\bar{\eta}/\eta)} : \shvl(\sigma)^{Gal(\bar{\eta}/\eta)}\rightarrow \shvl(\sigma)
\end{equation}
Some care needs to be taken when we speak of $\shvl(\sigma)^{Gal(\bar{\eta}/\eta^{mu})}$. One way to define this object is the following: the natural inclusion $\C \subset \overline{\C\drl x \drr}$ induces a morphism $\bar{\varepsilon}\xrightarrow{\phi} \sigma$ and the pullback along it gives us an equivalence at the level of the $\oo$-categories of $\mathbb{Q}_{\ell}$-adic sheaves $\shvl(\sigma)\xrightarrow{\simeq}\shvl(\bar{\varepsilon})$. Then we can consider 
\begin{equation}\label{commutative square shvl sigma varepsilon}
\begindc{\commdiag}[20]
\obj(0,15)[1]{$\shvl(\sigma)$}
\obj(80,15)[2]{$\shvl(\bar{\varepsilon})$}
\obj(0,-15)[3]{$\shvl(\sigma)^{Gal(\bar{\eta}/\eta^{mu})}$}
\obj(80,-15)[4]{$\shvl(\bar{\varepsilon})^{Gal(\bar{\eta}/\eta^{mu})}$}
\mor{1}{2}{$\phi^* \simeq$}
\mor(76,10)(76,-10){$trivial$}[\atright,\solidarrow]
\mor(84,-10)(84,10){$(-)^{\textup{h}Gal(\bar{\eta}/\eta^{mu})}$}[\atright,\solidarrow]
\mor(-4,10)(-4,-10){$trivial$}[\atright,\solidarrow]
\mor(4,-10)(4,10){$(-)^{\textup{h}Gal(\bar{\eta}/\eta^{mu})}$}[\atright,\solidarrow]
\mor{3}{4}{$\phi^* \simeq$}
\enddc
\end{equation}
\begin{lmm}
The $\oo$-functors
\begin{equation}
(-)^{\textup{h}Gal(\bar{\eta}/\eta)}: \shvl(\sigma)^{Gal(\bar{\eta}/\eta)}\rightarrow \shvl(\sigma) 
\end{equation}
\begin{equation}
(-)^{\textup{h}Gal(\bar{\eta}/\eta^{mu})}: \shvl(\sigma)^{Gal(\bar{\mu}/\eta^{mu})}\rightarrow \shvl(\sigma)
\end{equation}
\begin{equation}
(-)^{\textup{h}Gal(\bar{\eta}/\eta^{mu})}: \shvl(\sigma)^{Gal(\bar{\eta}/\eta)}\rightarrow \shvl(\sigma)^{Gal(\bar{\varepsilon}/\varepsilon)}
\end{equation}
\begin{equation}
(-)^{\textup{h}Gal(\bar{\varepsilon}/\varepsilon)}: \shvl(\sigma)^{Gal(\bar{\eta}/\eta)}\rightarrow \shvl(\sigma)^{Gal(\bar{\eta}/\eta^{mu})}
\end{equation}
\begin{equation}
(-)^{\textup{h}Gal(\bar{\varepsilon}/\varepsilon)}: \shvl(\sigma)^{Gal(\bar{\varepsilon}/\varepsilon)}\rightarrow \shvl(\sigma)
\end{equation}
are lax-monoidal.
\end{lmm}
\begin{proof}
This is a formal consequence of the fact that all these $\oo$-functors are right adjoints to symmetric monoidal $\oo$-functors (i.e. the \textit{trivial} functors).
\end{proof}
We are now ready to compute the image of $\Ql{,\sigma}$ via the different homotopy invariant $\oo$-functors.
\begin{lmm}\label{homotopy invariant units}
We have the following equivalences of commutative algebra objects in $\shvl(\sigma)$:
\begin{equation}\label{first equivalence calg}
\Ql{,\sigma}^{\textup{h}Gal(\bar{\varepsilon}/\varepsilon)}\simeq \Hl(\varepsilon)
\end{equation}
\begin{equation}\label{second equivalence calg}
\Ql{,\sigma}^{\textup{h}Gal(\bar{\eta}/\eta^{mu})}\simeq \phi_*\bar{k}^*j_{10*}^{mu}\Ql{,\eta^{mu}}
\end{equation}
\begin{equation}\label{third equivalence calg}
\Ql{,\sigma}^{\textup{h}Gal(\bar{\eta}/\eta)}\simeq \Hl(\eta)
\end{equation}
where the commutative algebra structures on the left hand sides is the one transferred by the lax-monoidality of the $\oo$-functors $(-)^{\textup{h}Gal(\bullet)}$. Moreover,
at the level of underlying objects we have equivalences
\begin{equation}\label{first equivalence obj}
\Ql{,\sigma}^{\textup{h}Gal(\bar{\varepsilon}/\varepsilon)}\simeq \Ql{,\sigma}\oplus \Ql{,\sigma}
[-1]
\end{equation}
\begin{equation}\label{second equivalence obj}
\Ql{,\sigma}^{\textup{h}Gal(\bar{\eta}/\eta^{mu})}\simeq \Ql{,\sigma}\oplus \Ql{,\sigma}
[-1]
\end{equation}
\begin{equation}\label{third equivalence obj}
\Ql{,\sigma}^{\textup{h}Gal(\bar{\eta}/\eta)}\simeq \Ql{,\sigma}\oplus \Ql{,\sigma}
[-1]\oplus \Ql{,\sigma}
[-1]\oplus \Ql{,\sigma}
[-2]
\end{equation}   
\end{lmm}
\begin{proof}
The proof of the equivalences $(\ref{first equivalence calg})$ and $(\ref{first equivalence obj})$ is \cite[Lemma 4.34]{brtv}. We rewrite it for the reader's convenience:
the specialization map 
\begin{equation}
\Ql{,\sigma}\rightarrow i_{01}^*\bar{l}_*\Ql{,\bar{\varepsilon}}
\end{equation}
is an equivariant equivalence in $\shvl(\sigma)^{Gal({\bar{\varepsilon}}/\varepsilon)}$ as $S_1$ is smooth over itself via the identity map (this can be computed explicitly, see \cite[footnote 32]{brtv} ). Therefore, taking its image via $(-)^{\textup{h}Gal(\bar{\varepsilon}/\varepsilon)}$ we get the following equivalences in $\shvl(\sigma)$
\begin{equation}
\Ql{,\sigma}^{\textup{h}Gal(\bar{\varepsilon}/\varepsilon)}\simeq \bigl ( i_{01}^*\bar{l}_*\Ql{,\bar{\varepsilon}}\bigr )^{\textup{h}Gal(\bar{\varepsilon}/\varepsilon)}\simeq i_{01}^*l_*\Ql{,\varepsilon}=: \Hl(\varepsilon)
\end{equation}
where the last equivalence is given by the argument in \cite[Proposition 4.31]{brtv}. Then equivalence $(\ref{first equivalence obj})$ follows from equivalence (\ref{cohomology punctured disk}).

Let us now consider the equivalences $(\ref{second equivalence calg})$ and $(\ref{second equivalence obj})$. By the commutativity (up to coherent homotopy) of the square $(\ref{commutative square shvl sigma varepsilon})$ and by the fact that $\phi^*$ is symmetric monoidal we have a chain of equivalences
\begin{equation}
\phi^*\bigl ( (\Ql{,\sigma})^{\textup{h}Gal(\bar{\eta}/\eta^{mu})}\bigr )\simeq \bigl ( \phi^*\Ql{,\sigma} \bigr )^{\textup{h}Gal(\bar{\eta}/\eta^{mu})}\simeq (\Ql{,\bar{\varepsilon}})^{\textup{h}Gal(\bar{\eta}/\eta^{mu})}
\end{equation}
By the same argument above applied to the strictly local henselian trait $\bar{U}_0=Spec(\overline{\C\drl x \drr}\dsl y \dsr)$ the last term is equivalent to $\bar{k}^*j_{10*}^{mu}\Ql{,\eta^{mu}}$. Then, $(\ref{second equivalence calg})$ follows from the observation that the inverse to $\phi^*$ is the derived pushforward $\phi_*$. At the level of objects, we know that 
\begin{equation}
\bar{k}^*j_{10*}^{mu}\Ql{,\eta^{mu}}\simeq \Ql{,\bar{\varepsilon}}\oplus \Ql{,\bar{\varepsilon}}
[-1]
\end{equation}
and $(\ref{second equivalence obj})$ follows from the fact that $\phi_*(\Ql{,\bar{,\varepsilon}})\simeq \Ql{,\sigma}$ and from the compatibility of $\phi_*$ with the usual and the Tate shift.

Let us now focus on the last two equivalences. First of all, notice that we have an equivalence of $\oo$-functors
\begin{equation}
i_{01}^*l_* (-)^{\textup{h}Gal(\bar{\eta}/\eta^{mu})}\simeq \bigl ( i_{01}^*l_* (-)\bigr )^{\textup{h}Gal(\bar{\eta}/\eta^{mu})} : \shvl(\sigma)^{Gal(\bar{\eta}/\eta^{mu})}\rightarrow \shvl(\sigma)
\end{equation}
Indeed, the continuity of the action of $Gal(\bar{\eta}/\eta^{mu})$ implies that 
\begin{equation}
l_*(-)^{\textup{h}Gal(\bar{\eta}/\eta^{mu})}\simeq \varinjlim_n l_*(-)^{\textup{h}Gal(\eta^{mu}(y^{\frac{1}{n}})/\eta^{mu})}
\end{equation}
Now, $i_{01}^*$ commutes with colimits and with finite limits and since $(-)^{\textup{h}Gal(\eta^{mu}(y^{\frac{1}{n}})/\eta^{mu})}$ are finite limits ($\eta^{mu}\subseteq \eta^{mu}(y^{\frac{1}{n}})$ are finite Galois extensions) we get the desired equivalence. Therefore, we obtain
\begin{equation}
\bigl ( i_{01}^*l_* (\Ql{,\varepsilon})\bigr )^{\textup{h}Gal(\bar{\eta}/\eta^{mu})}\simeq i_{01}^*l_* (\Ql{,\varepsilon})^{\textup{h}Gal(\bar{\eta}/\eta^{mu})}
\end{equation}
Recall that the morphism $v_{\varepsilon}: \bar{\varepsilon}\rightarrow \varepsilon$ induces an equivalence of symmetric monoidal stable presentable $\oo$-categories
\begin{equation}\label{equivalence v varepsilon}
\begindc{\commdiag}[20]
\obj(0,0)[1]{$\shvl(\varepsilon)$}
\obj(80,0)[2]{$\shvl(\bar{\varepsilon})^{Gal(\bar{\varepsilon}/\varepsilon)}$}
\mor(8,3)(62,3){$v_{\varepsilon}^*$}
\mor(62,-3)(8,-3){$v_{\varepsilon*}(-)^{\textup{h}Gal(\bar{\varepsilon}/\varepsilon)}\simeq \bigl (v_{\varepsilon*}(-)\bigr )^{\textup{h}Gal(\bar{\varepsilon}/\varepsilon)}$}
\enddc
\end{equation}
Thus, we can continue our chain of equivalences
\begin{equation}
i_{01}^*l_* (\Ql{,\varepsilon})^{\textup{h}Gal(\bar{\eta}/\eta^{mu})}\simeq i_{01}^*l_*\bigl ((v_{\varepsilon*}\Ql{,\bar{\varepsilon}} )^{\textup{h}Gal(\bar{\varepsilon}/\varepsilon)}\bigr )^{\textup{h}Gal(\bar{\eta}/\eta^{mu})}\simeq  i_{01}^*l_*\bigl ((v_{\varepsilon*}\Ql{,\bar{\varepsilon}} )^{\textup{h}Gal(\bar{\eta}/\eta^{mu})}\bigr )^{\textup{h}Gal(\bar{\varepsilon}/\varepsilon)}
\end{equation}
where the last equivalence follows from $(\ref{equiv hgal 1})$ and $(\ref{equiv hgal 2})$. Since taking homotopy fixed points commutes with pushforwards, we also get
\begin{equation}
i_{01}^*l_*\bigl ((v_{\varepsilon*}\Ql{,\bar{\varepsilon}} )^{\textup{h}Gal(\bar{\eta}/\eta^{mu})}\bigr )^{\textup{h}Gal(\bar{\varepsilon}/\varepsilon)}\simeq i_{01}^*l_*v_{\varepsilon*}\bigl (\Ql{,\bar{\varepsilon}} ^{\textup{h}Gal(\bar{\eta}/\eta^{mu})}\bigr )^{\textup{h}Gal(\bar{\varepsilon}/\varepsilon)}
\end{equation}
\[
\simeq i_{01}^*l_*v_{\varepsilon*}\bigl ( \bar{k}^*j_{10*}^{mu}\Ql{,\eta^{mu}} \bigr )^{\textup{h}Gal(\bar{\varepsilon}/\varepsilon)}
\]
where the last equivalence follows from \cite[Lemma 4.34]{brtv}.

Now notice that $\overline{U}_0=\varprojlim(Spec(\C\drl x^{\frac{1}{n}}\drr \dsl y \dsr))$ and $\overline{U}_0\xrightarrow{v_{U_0}} U_0$ is the limit of the projective system of morphisms $Spec(\C\drl x^{\frac{1}{n}}\drr \dsl y \dsr)\rightarrow U_0$, which are smooth. Moreover, the morphisms $Spec(\C\drl x^{\frac{1}{n}}\drr \dsl y \dsr)\rightarrow Spec(\C\drl x^{\frac{1}{n'}}\drr \dsl y \dsr)$ are affine. It is then a consequence of the smooth base change theorem (see the version in \cite[Theorem 1.1.5]{lan}) that
\begin{equation}
v_{U_0}^*j_{10*}\simeq j_{10*}^{mu}v_{\eta}^*
\end{equation}
where $v_{\eta}: \eta^{mu}\rightarrow \eta$. Then we can continue 
\begin{equation}
 i_{01}^*l_*v_{\varepsilon*}\bigl ( \bar{k}^*j_{10*}^{mu}\Ql{,\eta^{mu}} \bigr )^{\textup{h}Gal(\bar{\varepsilon}/\varepsilon)}\simeq  i_{01}^*l_*v_{\varepsilon*}\bigl ( \bar{k}^*v_{U_0}^*j_{10*}\Ql{,\eta} \bigr )^{\textup{h}Gal(\bar{\varepsilon}/\varepsilon)}\simeq i_{01}^*l_*v_{\varepsilon*}\bigl (v_{\varepsilon}^*k^*j_{10*}\Ql{,\eta} \bigr )^{\textup{h}Gal(\bar{\varepsilon}/\varepsilon)}
\end{equation}
Since $v_{\varepsilon*}v_{\varepsilon}^*(-)^{\textup{h}Gal(\bar{\varepsilon}/\varepsilon)}\simeq id_{\varepsilon}$, we get
\begin{equation}
i_{01}^*l_*\Ql{,\varepsilon}^{\textup{h}Gal(\bar{\eta}/\eta)}\simeq i_{01}^*l_*k^*j_{10*}\Ql{,\eta}
\end{equation}
We now get $(\ref{third equivalence calg})$ as follows
\begin{equation}
\Ql{,\sigma}^{\textup{h}Gal(\bar{\eta}/\eta)}\simeq (\Ql{,\sigma}^{\textup{h}Gal(\bar{\varepsilon}/\varepsilon)})^{\textup{h}Gal(\bar{\eta}/\eta^{mu})}\simeq (i_{01}^*l_*\Ql{,\varepsilon})^{\textup{h}Gal(\bar{\eta}/\eta^{mu})}\simeq i_{01}^*l_*k^*j_{10*}
\Ql{,\eta}
\end{equation}
Finally, $(\ref{third equivalence obj})$ follows from $(\ref{first equivalence obj})$, $(\ref{second equivalence obj})$ and the compatibility of $(-)^{\textup{h}Gal(\bar{\eta}/\eta^{mu})}$ with direct sums, Tate shifts and usual shifts.
\end{proof}
Notice that by the unit of the adjuctions $(trivial, (-)^{\textup{h}Gal(\bar{\varepsilon}/\varepsilon)})$ and $(trivial,(-)^{\textup{h}Gal(\bar{\eta}/\eta^{mu})})$ we get the following two maps of commutative algebra objects in $\shvl(\sigma)$
\begin{equation}\label{first map calg}
\Hl(\varepsilon)\simeq \Ql{\sigma}^{\textup{h}Gal(\bar{\varepsilon}/\varepsilon)}\rightarrow \bigl ( \Ql{\sigma}^{\textup{h}Gal(\bar{\varepsilon}/\varepsilon)} \bigr )^{\textup{h}Gal(\bar{\eta}/\eta^{mu})}\simeq \Hl(\eta)
\end{equation}
\begin{equation}\label{second map calg}
\Ql{\sigma}^{\textup{h}Gal(\bar{\eta}/\eta^{mu})}\rightarrow  \bigl ( \Ql{\sigma}^{\textup{h}Gal(\bar{\eta}/\eta^{mu})} \bigr )^{\textup{h}Gal(\bar{\varepsilon}/\varepsilon)}\simeq \Hl(\eta)
\end{equation}
Moreover, these are compatible in the sense that the square
\begin{equation}
\begindc{\commdiag}[20]
\obj(0,15)[1]{$\Ql{\sigma}$}
\obj(80,15)[2]{$\Ql{\sigma}^{\textup{h}Gal(\bar{\varepsilon}/\varepsilon)}=:\Hl(\varepsilon)$}
\obj(0,-15)[3]{$\Ql{\sigma}^{\textup{h}Gal(\bar{\eta}/\eta^{mu})}$}
\obj(80,-15)[4]{$\Hl(\eta)$}
\mor{1}{2}{$(trivial,(-)^{\textup{h}Gal(\bar{\varepsilon}/\varepsilon)})$}
\mor{1}{3}{$ (trivial,(-)^{\textup{h}Gal(\bar{\eta}/\eta^{mu})})$}[\atright,\solidarrow]
\mor{2}{4}{$(\ref{first map calg})$}
\mor{3}{4}{$(\ref{second map calg})$}
\enddc
\end{equation}
where the left vertical arrow (resp. the top arrow) is induced by the unit of the adjunction $(trivial,(-)^{\textup{h}Gal(\bar{\eta}/\eta^{mu})})$ (resp. $(trivial,(-)^{\textup{h}Gal(\bar{\varepsilon}/\varepsilon)})$), is commutative (up to coherent homotopy) since
\begin{equation}
\bigl ( (-)^{\textup{h}Gal(\bar{\varepsilon}/\varepsilon)}\bigr )^{\textup{h}Gal(\bar{\eta}/\eta^{mu})}\simeq \bigl ( (-)^{\textup{h}Gal(\bar{\eta}/\eta^{mu})}\bigr )^{\textup{h}Gal(\bar{\varepsilon}/\varepsilon)}
\end{equation}
Hence, by the universal product of the tensor product \cite[Proposition 3.2.4.7]{ha} we get a map of commutative algebras
\begin{equation}\label{equivalence calg cohomology torus}
\Ql{\sigma}^{\textup{h}Gal(\bar{\varepsilon}/\varepsilon)}\otimes_{\Ql{\sigma}}\Ql{\sigma}^{\textup{h}Gal(\bar{\eta}/\eta^{mu})}\rightarrow \Hl(\eta)
\end{equation}
\begin{lmm}\label{lemma equivalence calg cohomology torus}
The map $(\ref{equivalence calg cohomology torus})$ is an equivalence of commutative algebra objects in $\shvl(\sigma)$.
\end{lmm}
\begin{proof}
Consider the projection $q: \textup{B}Gal(\bar{\varepsilon}/\varepsilon)\rightarrow \sigma$. The adjunction $(q^*,q_*)$ identifies with
\begin{equation}
\begindc{\commdiag}[20]
\obj(0,0)[1]{$\shvl(\sigma)$}
\obj(80,0)[2]{$\shvl(\sigma)^{Gal(\bar{\varepsilon}/\varepsilon)}$}
\mor(8,3)(62,3){$trivial$}
\mor(62,-3)(8,-3){$(-)^{\textup{h}Gal(\bar{\varepsilon}/\varepsilon)}$}
\enddc
\end{equation}
We know that $\Hl(\eta)\simeq \bigl ( \Ql{,\sigma}^{\textup{h}Gal(\bar{\eta}/\eta^{mu})} \bigr )^{\textup{h}Gal(\bar{\varepsilon}/\varepsilon)}$. The action of $Gal(\bar{\varepsilon}/\varepsilon)$ on $\Ql{,\sigma}^{\textup{h}Gal(\bar{\eta}/\eta^{mu})}$ is the trivial one. Therefore, by the projection formula applied to $q$ (\cite[Corollary 0.1.3]{lz17i}) and by the equivalence of commutative algebras $q^*\Ql{,\sigma}\simeq \Ql{,\textup{B}Gal(\bar{\varepsilon}/\varepsilon)}$ we get
\begin{equation}
\Hl(\eta)\simeq q_*q^*(\Ql{,\sigma}^{\textup{h}Gal(\bar{\eta}/\eta^{mu})})\simeq
q_*\bigl ( q^*\Ql{,\sigma}^{\textup{h}Gal(\bar{\eta}/\eta^{mu})}\otimes_{q^*\Ql{,\sigma}}q^*\Ql{,\sigma}\bigr )
\end{equation}
\begin{equation*}
\simeq \Ql{,\sigma}^{\textup{h}Gal(\bar{\eta}/\eta^{mu})}\otimes_{\Ql{,\sigma
}} \Ql{,\sigma}^{\textup{h}Gal(\bar{\varepsilon}/\varepsilon)}
\end{equation*}
It is then clear that $(\ref{equivalence calg cohomology torus})$ is an equivalence.
\end{proof}

\section{\textsc{Comparison with the \texorpdfstring{$\ell$}{l}-adic realization of dg categories of singularities}}
\subsection{\textsc{\texorpdfstring{$\rbul{}$}{RBUl}-valued cohomology}}
In this section we introduce $\rbul{}$-valued cohomology. We begin with the following definition
\begin{defn}
\begin{itemize}
\item The $\rbul{}$-valued cohomology of the punctured disk is defined as 
\begin{equation}
\Hr(\varepsilon):= i_{01}^*l_*\rbul{\varepsilon}
\end{equation}
in $\textup{Mod}_{\rbul{\sigma}}(\shvl(\sigma))$
\item The $\rbul{}$-valued cohomology of the torus is defined as 
\begin{equation}
\Hr(\eta):= i_{01}^*l_*k^*j_{10*}\rbul{\eta}
\end{equation}
in $\textup{Mod}_{\rbul{\sigma}}(\shvl(\sigma))$
\end{itemize}
\end{defn}
Notice that the lax monoidal structure of the $\oo$-functors $i_{01}^*l_*$ and $k^*j_{10*}$ give a commutative algebra structure to both $\Hr(\varepsilon)$ and $\Hr(\eta)$. We shall now give an alternative description of these two objects, following the lead of \cite{brtv} (as usual).

Let us start with the $\rbul{}$-valued cohomology of the punctured disk. The canonical morphism of commutative algebras in $\shvl(\varepsilon)$
\begin{equation}
\Ql{,\varepsilon} \rightarrow \rbul{\varepsilon}\simeq \bigoplus_{i\in \Z}\Ql{,\varepsilon}(i)[2i]
\end{equation}
induces, via $i_{01}^*l_*$, the following map of algebras in $\shvl(\sigma)$
\begin{equation}
\Hl(\varepsilon):=i_{01}^*l_*\Ql{\varepsilon}\rightarrow i_{01}^*l_*\rbul{\varepsilon}=: \Hr(\varepsilon) 
\end{equation}
On the other hand, the $\oo$-functor $\Rl$ is compatible with the $6$-functors formalism, and therefore the canonical map of algebras $\textup{B}\mathbb{U}_{\sigma}\rightarrow i_{01}^*l_*\textup{B}\mathbb{U}_{\varepsilon}$ induces
\begin{equation}
\rbul{\sigma}\rightarrow \Hr(\varepsilon)
\end{equation}
which gives to the latter a $\rbul{\sigma}$-algebra structure. It is then possible to consider the canonical map of commutative algebras
\begin{equation}
\Hl(\varepsilon)\otimes_{\Ql{,\sigma}}\rbul{\sigma}\rightarrow \Hr(\varepsilon)
\end{equation}
This is an equivalence in $\textup{CAlg}(\shvl(\sigma))$. Indeed, it is sufficient to verify that it is so at the level of the underlying objects:
\begin{equation}
\Hl(\varepsilon)\otimes_{\Ql{\sigma}}\rbul{\sigma}=i_{01}^*l_*
\Ql{,\varepsilon}\otimes_{\Ql{,\sigma}}\bigl ( \bigoplus_{i\in \Z}
\Ql{,\sigma}(i)[2i] \bigr )
\end{equation}
\begin{equation*}
\simeq \bigoplus_{i\in \Z}
\bigl ( i_{01}^*l_*\Ql{,\varepsilon}\otimes_{\Ql{,\sigma}} \Ql{,\sigma}(i)[2i] \bigr )
\simeq \bigoplus_{i\in \Z} i_{01}^*l_*\Ql{,\varepsilon}(i)[2i]\simeq i_{01}^*l_*\bigl (\bigoplus_{i\in \Z}\Ql{,\varepsilon}(i)[2i]\bigr )\simeq \Hr(\varepsilon)
\end{equation*}
where we used that $-\otimes_{\Ql{,\sigma}}-$ commutes with colimits separately in each variable, $i_{01}^*$ commutes with colimits (it is a left adjoint) and the same holds true for $l_*$ (see \cite[Example 9.4.6]{ro14}).

We can produce a similar argument for the $\rbul{}$-valued cohomology of the torus. The image via the lax monoidal $\oo$-functor $i_{01}^*l_*k^*j_{10*}$ of the canonical map of commutative algebras in $\shvl(\eta)$
\begin{equation}
\Ql{,\eta}\rightarrow \rbul{\eta}\simeq \bigoplus_{i\in \Z} \Ql{,\eta}(i)[2i]
\end{equation}
gives a $\Hl(\eta)$-algebra structure to $\Hr(\eta)$:
\begin{equation}
\Hl(\eta)\rightarrow \Hr(\eta)
\end{equation}
Moreover, the canonical morphism $\textup{B}\mathbb{U}_{\sigma}\rightarrow i_{01}^*l_*k^*j_{10*}\textup{B}\mathbb{U}_{\eta}$ provides a $\rbul{\sigma}$-algebra structure to $\Hr(\eta)$ by taking its image via $\Rl$. The universal property of the tensor product (\cite[Proposition 3.2.4.7]{ha}) gives us a canonical map in $\textup{CAlg}(\shvl(\sigma))$:
\begin{equation}
\Hl(\eta)\otimes_{\Ql{,\sigma}}\rbul{\sigma}\rightarrow \rbul{\eta}
\end{equation}
As for the previous case, in order to show that this map is an equivalence it is sufficient to consider the underlying objects. The computation we did before does the job in this case too:
\begin{equation}
\Hl(\eta)\otimes_{\Ql{,\sigma}}\rbul{\sigma}\simeq i_{01}^*l_*k^*j_{10*}\Ql{,\eta}\otimes_{\Ql{,\sigma}}\bigl ( \bigoplus_{i\in \Z}\Ql{,\sigma}(i)[2i]\bigr )
\end{equation}
\begin{equation*}
\simeq \bigoplus_{i\in \Z} \bigl ( i_{01}^*l_*k^*j_{10*}\Ql{,\eta}\otimes_{\Ql{,\sigma}} \Ql{,\sigma}(i)[2i]\bigr )\simeq  \bigoplus_{i\in \Z} i_{01}^*l_*k^*j_{10*}\Ql{,\eta}(i)[2i]
\end{equation*}
\begin{equation*}
\simeq i_{01}^*l_*k^*j_{10*} \bigl (\bigoplus_{i\in \Z} \Ql{,\eta}(i)[2i]\bigr )\simeq \Hr(\eta)
\end{equation*}

Let us summarize, for future reference, what we have shown so far:
\begin{lmm}
The following equivalences hold in $\textup{CAlg}(\shvl(\sigma))$:
\begin{enumerate}
\item $\Hl(\varepsilon)\otimes_{\Ql{,\sigma}} \rbul{\sigma}\simeq \Hr(\varepsilon)$
\item $\Hl(\eta)\otimes_{\Ql{,\sigma}} \rbul{\sigma}\simeq \Hr(\eta)$
\item $\Hr(\varepsilon)\otimes_{\Ql{,\sigma}}\Ql{,\sigma}^{\textup{h}Gal(\bar{\eta}/\eta^{mu})}\simeq \Hr(\eta)$
\end{enumerate}
\end{lmm}
\begin{proof}
We only need to show the last equivalence, which is a consequence of the first two combined with Lemma \ref{lemma equivalence calg cohomology torus}:
\begin{equation}
\Hr(\eta)\simeq \Hl(\eta)\otimes_{\Ql{,\sigma}} \rbul{\sigma} \simeq \Ql{,\sigma}^{\textup{h}Gal(\bar{\eta}/\eta^{mu})}\otimes_{\Ql{,\sigma}}\Hl(\varepsilon)\otimes_{\Ql{,\sigma}}\rbul{\sigma}
\end{equation}
\[
\simeq \Ql{,\sigma}^{\textup{h}Gal(\bar{\eta}/\eta^{mu})}\otimes_{\Ql{,\sigma}}\Hr(\varepsilon)
\]
\end{proof}
Let us go one step further and consider the analogous of the $\rbul{}$-valued cohomology of the punctured disk and of the $\rbul{}$-valued cohomology of the torus for an arbitrary proper, flat, finitely presented $S$-scheme $p:X\rightarrow S$.
\begin{defn}
For any $p:X\rightarrow S$ as above, define
\begin{equation}
\Hr(X_{\varepsilon}):=p_{\sigma*}i_{01}^*l_*\rbul{X_{\varepsilon}}\simeq i_{01}^*l_*p_{\varepsilon*}\rbul{X_{\varepsilon}}
\end{equation}
\begin{equation}
\Hr(X_{\eta}):=p_{\sigma*}i_{01}^*l_*k^*j_{10*}\rbul{X_{\eta}}\simeq i_{01}^*l_*k^*j_{10*}p_{\eta*}\rbul{X_{\eta}}
\end{equation}
\end{defn}
Note that we have the following maps of commutative algebras in $\shvl(\sigma)$
\begin{equation}\label{algebra structures}
\begindc{\commdiag}[17]
\obj(0,80)[1]{$\Hl(\varepsilon)$}
\obj(80,80)[2]{$\Hr(X_{\varepsilon})$}
\obj(40,40)[3]{$\Hr(\varepsilon)$}
\mor{1}{2}{$$}
\mor{1}{3}{$$}
\mor{3}{2}{$$}
\obj(120,80)[4]{$\Hl(\eta)$}
\obj(200,80)[5]{$\Hr(X_{\eta})$}
\obj(160,40)[6]{$\Hr(\eta)$}
\mor{4}{5}{$$}
\mor{4}{6}{$$}
\mor{6}{5}{$$}
\obj(60,20)[7]{$p_{\sigma*}\rbul{X_{\sigma}}$}
\obj(140,20)[8]{$\Hr(X_{\varepsilon})$}
\obj(100,-20)[9]{$\Hr(X_{\eta})$}
\mor{7}{8}{$$}
\mor{7}{9}{$$}
\mor{8}{9}{$$}
\enddc
\end{equation}

which, using \cite[Proposition 3.2.4.7]{ha} once again, give us maps 
\begin{equation}\label{first can map}
\Hr(\varepsilon)\otimes_{\rbul{\sigma}} p_{\sigma*}\rbul{X_{\sigma}}\simeq \Hl(\varepsilon)\otimes_{\Ql{\sigma}}p_{\sigma*}\rbul{X_{\sigma}}\rightarrow \Hr(X_{\varepsilon})
\end{equation}
\begin{equation}\label{second can map}
\Hr(\eta)\otimes_{\rbul{\sigma}} p_{\sigma*}\rbul{X_{\sigma}}\simeq \Hl(\eta)\otimes_{\Ql{,\sigma}}p_{\sigma*}\rbul{X_{\sigma}}\rightarrow \Hr(X_{\eta})
\end{equation}

\begin{rmk}
Note that $\Hr(X_{\eta})$ carries a structure of $\Ql{\sigma}^{\textup{h}Gal(\bar{\eta}/\eta^{mu})}$-algebra and of $\Hl(\varepsilon)$-algebra, by restriction of scalars.
\end{rmk}
In particular, starting from the third triangle in diagram $(\ref{algebra structures})$, we get the following one:
\begin{equation}
\begindc{\commdiag}[20]
\obj(0,20)[1]{$p_{\sigma*}\rbul{X_{\sigma}}\otimes_{\Ql{,\sigma}}\Ql{,\sigma}^{\textup{h}Gal(\bar{\eta}/\eta^{mu})}$}
\obj(120,20)[2]{$\Hr(X_{\varepsilon})\otimes_{\Ql{,\sigma}}\Ql{,\sigma}^{\textup{h}Gal(\bar{\eta}/\eta^{mu})}$}
\obj(60,-20)[3]{$\Hr(X_{\eta})$}
\mor{1}{2}{$-\otimes id_{\Ql{\sigma}^{\textup{h}Gal(\bar{\eta}/\eta^{mu})}}$}
\mor(0,15)(50,-15){$$}
\mor(120,15)(70,-15){$$}
\enddc
\end{equation}
Finally, using the fact that both $\Hr(X_{\varepsilon})$ and $\Hr(X_{\eta})$ are $\Hl(\varepsilon)$-algebras, we get
\begin{equation}\label{fundamental triangle 1}
\begindc{\commdiag}[17]
\obj(0,40)[1]{$\bigl (p_{\sigma*}\rbul{X_{\sigma}}\otimes_{\Ql{,\sigma}}\Hl(\varepsilon)\bigr )\otimes_{\Ql{,\sigma}}\Ql{,\sigma}^{\textup{h}Gal(\bar{\eta}/\eta^{mu})}$}
\obj(0,30)[a]{$\simeq p_{\sigma*}\rbul{X_{\sigma}}\otimes_{\Ql{,\sigma}}\Hl(\eta)$}
\obj(160,40)[2]{$\Hr(X_{\varepsilon})\otimes_{\Ql{,\sigma}}\Ql{,\sigma}^{\textup{h}Gal(\bar{\eta}/\eta^{mu})}$}
\obj(80,0)[3]{$\Hr(X_{\eta})$}
\mor{1}{2}{$(\ref{first can map})\otimes id_{\Ql{\sigma}^{\textup{h}Gal(\bar{\eta}/\eta^{mu})}}$}
\mor(15,25)(70,5){$(\ref{second can map})$}[\atright,\solidarrow]
\mor(145,25)(90,5){$$}
\enddc
\end{equation}

We shall now provide an equivalence between diagram $(\ref{fundamental triangle 1})$ and the image of the canonical triangle
\begin{equation} 
\begindc{\commdiag}[20]
\obj(60,20)[1]{$p_{\sigma*}\rbul{X_{\sigma}}$}
\obj(140,20)[2]{$i_{01}^*\bar{l}_*\rbul{X_{\bar{\varepsilon}}}$}
\obj(100,-20)[3]{$ i_{01}^*\bar{l}_*\bar{k}^*\bar{j}_{10*}\rbul{X_{\bar{\eta}}}$}
\mor{1}{2}{$sp_1$}
\mor{2}{3}{$sp_2$}
\mor{1}{3}{$(sp_2\circ sp_1)$}[\atright,\solidarrow]
\enddc
\end{equation}
via the $\oo$-functor $(-)^{\textup{h}Gal(\bar{\eta}/\eta)}$
where $sp_1$ and $sp_2$ are the specialization morphisms
\begin{equation}
sp_1: \rbul{X_{\sigma}}\rightarrow i_{01}^*\bar{l}_*\rbul{X_{\bar{\varepsilon}}}
\end{equation}
\begin{equation}
sp_2:i_{01}^*\bar{l}_*\rbul{X_{\bar{\varepsilon}}}\rightarrow i_{01}^*\bar{l}_*\bar{k}^*\bar{j}_{10*}\rbul{X_{\bar{\eta}}}
\end{equation}
As a preliminary step, we show that the vertices of the two triangles coincide.
\begin{prop}\label{gal invariant interesting stuff}
There are canonical equivalences in $\textup{CAlg}(\shvl(\sigma))$
\begin{equation}\label{rbul sigma hgal eta}
\bigl ( p_{\sigma*}\rbul{X_{\sigma}}\bigr )^{\textup{h}Gal(\bar{\eta}/\eta)}\simeq p_{\sigma*}\rbul{X_{\sigma}}\otimes_{\Ql{,\sigma}}\Ql{,\sigma}^{\textup{h}Gal(\bar{\eta}/\eta)}
\end{equation}
\begin{equation}\label{rbul varepsilon hgal eta}
\bigl ( p_{\sigma*}i_{01}^*\bar{l}_*\rbul{X_{\bar{\varepsilon}}}\bigr )^{\textup{h}Gal(\bar{\eta}/\eta)}\simeq \Hr(X_{\varepsilon})\otimes_{\Ql{,\sigma}}\Ql{,\sigma}^{\textup{h}Gal(\bar{\eta}/\eta^{mu})}
\end{equation}
\begin{equation}\label{rbul eta hgal eta}
\bigl ( p_{\sigma*}i_{01}^*\bar{l}_*\bar{k}^*\bar{j}_{10*}\rbul{X_{\bar{\eta}}}\bigr )^{\textup{h}Gal(\bar{\eta}/\eta)}\simeq \Hr(X_{\eta})
\end{equation}
\end{prop}
\begin{proof}
\begin{enumerate}
\item ($\ref{rbul sigma hgal eta}$): Consider the projection $q:\textup{B}Gal(\bar{\eta}/\eta)\rightarrow \sigma$. The adjunction $(q^*,q_*)$ identifies with
\begin{equation}
\begindc{\commdiag}[20]
\obj(0,0)[1]{$\shvl(\sigma)$}
\obj(80,0)[2]{$\shvl(\sigma)^{Gal(\bar{\eta}/\eta)}$}
\mor(8,3)(62,3){$trivial$}
\mor(62,-3)(8,-3){$(-)^{\textup{h}Gal(\bar{\eta}/\eta)}$}
\enddc
\end{equation}
As $p_{\sigma*}\rbul{X_{\sigma}}$ is endowed with the trivial $Gal(\bar{\eta}/\eta)$-action, the result follows from the projection formula \cite{lz17i}:
\begin{equation}
p_{\sigma*}\rbul{X_{\sigma}}^{\textup{h}Gal(\bar{\eta}/\eta)}\simeq q_*q^*p_{\sigma*}\rbul{X_{\sigma}}
\end{equation}
\begin{equation*}
\simeq q_*\bigl ( q^*p_{\sigma*}\rbul{X_{\sigma}}\otimes_{\Ql{,\textup{B}Gal(\bar{\eta}/\eta)}}\Ql{,\textup{B}Gal(\bar{\eta}/\eta)}\bigr )
\simeq p_{\sigma*}\rbul{X_{\sigma}}\otimes_{\Ql{,\sigma}} \Ql{,\sigma}^{\textup{h}Gal(\bar{\eta}/\eta)}
\end{equation*}
\item $(\ref{rbul varepsilon hgal eta})$: By the natural equivalence $(\ref{equiv hgal 2})$ we have
\begin{equation}
\bigl ( p_{\sigma*}i_{01}^*\bar{l}_*\rbul{X_{\bar{\varepsilon}}} \bigr )^{\textup{h}Gal(\bar{\eta}/\eta)}\simeq \Bigl (\bigl ( p_{\sigma*}i_{01}^*\bar{l}_*\rbul{X_{\bar{\varepsilon}}} \bigr )^{\textup{h}Gal(\bar{\varepsilon}/\varepsilon)}\Bigr )^{\textup{h}Gal(\bar{\eta}/\eta^{mu})}
\end{equation}
and by \cite[Proposition 4.32]{brtv} this is equivalent to
\[
\bigl ( p_{\sigma*}i_{01}^*l_*\rbul{X_{\varepsilon}}\bigr )^{\textup{h}Gal(\bar{\eta}/\eta^{mu})}
\]
The result then follows once again by the projection formula, applied to the projection $\textup{B}Gal(\bar{\eta}/\eta^{mu})\rightarrow \sigma$.
\item $(\ref{rbul eta hgal eta})$: By the natural equivalence $(\ref{equiv hgal 1})$ we have that
\begin{equation}
\bigl ( p_{\sigma*}i_{01}^*\bar{l}_*\bar{k}^*\bar{j}_{10*}\rbul{X_{\bar{\eta}}}\bigr )^{\textup{h}Gal(\bar{\eta}/\eta)}
\end{equation}
\begin{equation*} 
\simeq \Bigl ( \bigl ( p_{\sigma*}i_{01}^*\bar{l}_*\bar{k}^*\bar{j}_{10*}\rbul{X_{\bar{\eta}}} \bigr )^{\textup{h}Gal(\bar{\eta}/\eta^{mu})}\Bigr )^{\textup{h}Gal(\bar{\varepsilon}/\varepsilon)}
\end{equation*}
and by the same argument used in Lemma \ref{homotopy invariant units} and by \cite[Proposition 4.32]{brtv} this is also equivalent to
\begin{equation}\label{chain 1}
 \Bigl (  p_{\sigma*}i_{01}^*\bar{l}_*\bigl (\bar{k}^*\bar{j}_{10*}\rbul{X_{\bar{\eta}}} \bigr )^{\textup{h}Gal(\bar{\eta}/\eta^{mu})}\Bigr )^{\textup{h}Gal(\bar{\varepsilon}/\varepsilon)} 
 \end{equation}
 \begin{equation*}
 \simeq  \Bigl (  p_{\sigma*}i_{01}^*\bar{l}_*\bigl (\bar{k}^*j_{10*}^{mu}\rbul{X_{\eta}^{mu}} \bigr )\Bigr )^{\textup{h}Gal(\bar{\varepsilon}/\varepsilon)}
\end{equation*}
Note that by the equivalence
\begin{equation}
\begindc{\commdiag}[20]
\obj(0,0)[1]{$\shvl(X_{\varepsilon})$}
\obj(80,0)[2]{$\shvl(X_{\bar{\varepsilon}})^{Gal(\bar{\varepsilon}/\varepsilon)}$}
\mor(10,3)(60,3){$v_{\varepsilon}^*$}
\mor(60,-3)(10,-3){$v_{\varepsilon*}(-)^{\textup{h}Gal(\bar{\varepsilon}/\varepsilon)}\simeq (v_{\varepsilon*}(-))^{\textup{h}Gal(\bar{\varepsilon}/\varepsilon)}$}
\enddc
\end{equation}
we find that $id_{X_{\varepsilon}}\simeq (v_{\varepsilon*}v_{\varepsilon}^*(-))^{\textup{h}Gal(\bar{\varepsilon}/\varepsilon)}$ and therefore
\begin{equation}
k^*j_{10*}\simeq id_{X_{\varepsilon}}k^*j_{10*}\simeq (v_{\varepsilon*}v_{\varepsilon}^*k^*j_{10*})^{\textup{h}Gal(\bar{\varepsilon}/\varepsilon)}
\end{equation}
\begin{equation*}
\simeq (v_{\varepsilon*}\bar{k}^*v_{U_0}^*j_{10*})^{\textup{h}Gal(\bar{\varepsilon}/\varepsilon)}\simeq (v_{\varepsilon*}\bar{k}^*j_{10*}^{mu}v_{\eta}^*)^{\textup{h}Gal(\bar{\varepsilon}/\varepsilon)}
\end{equation*}
where the last equivalence follows from the smooth base change theorem (\cite[Theorem 1.1.5]{lan}). Then the chain $(\ref{chain 1})$ continues
\begin{equation}
\Bigl (  p_{\sigma*}i_{01}^*\bar{l}_*\bigl (\bar{k}^*j_{10*}^{mu}\rbul{X_{\eta}^{mu}} \bigr )\Bigr )^{\textup{h}Gal(\bar{\varepsilon}/\varepsilon)}
\end{equation}
\begin{equation*}
\simeq p_{\sigma*}i_{01}^*l_* \Bigl ( v_{\varepsilon*} \bar{k}^*j_{10*}^{mu}v_{\eta}^*\rbul{X_{\eta}}\Bigr )^{\textup{h}Gal(\bar{\varepsilon}/\varepsilon)}
\simeq p_{\sigma*}i_{01}^*l_*k^*j_{10*}\rbul{X_{\eta}}
\end{equation*}
\end{enumerate}
\end{proof}
\begin{prop}
The following square is commutative (up to coherent homotopy)
\begin{equation}
\begindc{\commdiag}[20]
\obj(0,20)[1]{$p_{\sigma*}\rbul{X_{\sigma}}\otimes_{\Ql{,\sigma}}\Hl(\eta)$}
\obj(100,20)[2]{$\bigl (p_{\sigma*}\rbul{X_{\sigma}} \bigr )^{\textup{h}Gal(\bar{\eta}/\eta)}$}
\obj(0,-20)[3]{$\Hr(X_{\varepsilon})\otimes_{\Ql{\sigma}}\Ql{,\sigma}^{\textup{h}Gal(\bar{\eta}/\eta^{mu})}$}
\obj(100,-20)[4]{$\bigl ( p_{\sigma*}i_{01}^*\bar{l}_*\rbul{X_{\bar{\varepsilon}}}
\bigr )^{\textup{h}Gal(\bar{\eta}/\eta)}$}
\mor{1}{2}{$(\ref{rbul sigma hgal eta})$}
\mor{1}{3}{$(\ref{first can map})\otimes id$}
\mor{2}{4}{$(sp_1)^{\textup{h}Gal(\bar{\eta}/\eta)}$}
\mor{3}{4}{$(\ref{rbul varepsilon hgal eta})$}
\enddc
\end{equation}
\end{prop}
\begin{proof}
By the universal property of the tensor product (\cite[Proposition 3.2.4.7]{ha}), we need to show that the two maps 
\[
\Hl(\eta)\rightarrow \bigl ( p_{\sigma*}i_{01}^*\bar{l}_*\rbul{X_{\bar{\varepsilon}}}\bigr )^{\textup{h}Gal(\bar{\eta}/\eta)}
\] 
and the two maps 
\[
p_{\sigma*}\rbul{X_{\sigma}}\rightarrow \bigl ( p_{\sigma*}i_{01}^*\bar{l}_*\rbul{X_{\bar{\varepsilon}}}\bigr )^{\textup{h}Gal(\bar{\eta}/\eta)}
\]
obtained from the precompositions with 
\[
\Hl(\eta)\xrightarrow{1\otimes id} p_{\sigma*}\rbul{X_{\sigma}}\otimes_{\Ql{,\sigma}}\Hl(\eta)
\] 
and 
\[
p_{\sigma*}\rbul{X_{\sigma}}\xrightarrow{id\otimes 1}p_{\sigma*}\rbul{X_{\sigma}}\otimes_{\Ql{\sigma}}\Hl(\eta)
\]
respectively, are equivalent. Consider the diagram
\begin{equation*}
\begindc{\commdiag}[17]
\obj(-70,120)[1]{$\Hl(\varepsilon)\otimes_{\Ql{,\sigma}}\Ql{,\sigma}^{\textup{h}Gal(\bar{\eta}/\eta^{mu})}\simeq \Hl(\eta)$}
\obj(70,120)[2]{$\bigl ( p_{\sigma*}\rbul{X_{\sigma}} \bigr )^{\textup{h}Gal(\bar{\eta}/\eta)}$}
\obj(-70,70)[3]{$\bigl ( p_{\sigma*}i_{01}^*\bar{l}_*\rbul{X_{\bar{\varepsilon}}} \bigr )^{\textup{h}Gal(\bar{\varepsilon}/\varepsilon)}\otimes_{\Ql{,\sigma}}\Ql{,\sigma}^{\textup{h}Gal(\bar{\eta}/\eta^{mu})}$}
\obj(70,70)[4]{$\hspace{0.6cm} \bigl ( p_{\sigma*}i_{01}^*\bar{l}_*\rbul{X_{\bar{\varepsilon}}}\bigr )^{\textup{h}Gal(\bar{\eta}/\eta)}$}
\mor{1}{2}{$$}
\mor{1}{3}{$(i_{01}^*\bar{l}_*\Ql{\bar{\varepsilon}})^{\textup{h}Gal(\bar{\varepsilon}/\varepsilon)}\otimes id$}[\atright,\solidarrow]
\mor{2}{4}{$(sp_1)^{\textup{h}Gal(\bar{\eta}/\eta)}$}
\mor{3}{4}{$\simeq$}
\obj(-40,30)[5]{$ (i_{01}^*\bar{l}_*\Ql{,\bar{\varepsilon}})^{\textup{h}Gal(\bar{\varepsilon}/\varepsilon)}$}
\obj(-40,20)[5b]{$\simeq \Ql{,\sigma}^{\textup{h}Gal(\bar{\varepsilon}/\varepsilon)}$}
\obj(40,30)[6]{$\bigl ( p_{\sigma*}\rbul{X_{\sigma}} \bigr )^{\textup{h}Gal(\bar{\varepsilon}/\varepsilon)}$}
\obj(-40,-20)[7]{$\Hr(X_{\varepsilon})$}
\obj(40,-20)[8]{$\bigl ( p_{\sigma*}i_{01}^*\bar{l}_*\rbul{X_{\bar{\varepsilon}}}\bigr )^{\textup{h}Gal(\bar{\varepsilon}/\varepsilon)}$}
\mor{5}{6}{$$}
\mor{5b}{7}{$$}
\mor{6}{8}{$(sp_1)^{\textup{h}Gal(\bar{\varepsilon}/\varepsilon)}$}[\atright,\solidarrow]
\mor{7}{8}{$\simeq$}
\mor(-40,34)(-50,66){$$}[\atleft,\solidline]
\mor(-53,74)(-66,116){$$}[\atleft,\solidarrow]
\mor(40,34)(50,66){$$}[\atleft,\solidline]
\mor(53,74)(66,116){$$}[\atleft,\solidarrow]
\mor(-42,-16)(-52,21){$$}[\atleft,\solidline]
\mor(-56,33)(-66,66){$$}[\atleft,\solidarrow]
\mor(42,-16)(54,26){$$}[\atleft,\solidline]
\mor(56,33)(66,66){$$}[\atleft,\solidarrow]
\enddc
\end{equation*}
where the left diagonal arrows are the ones induced by the equivalences 
\begin{equation*}
    \Ql{,\sigma}^{\textup{h}Gal(\bar{\varepsilon}/\varepsilon)}\xrightarrow{\sim}\Hl(\varepsilon) \hspace{0.5cm} \Hr(X_{\varepsilon})\simeq \bigl ( p_{\sigma*}i_{01}^*\bar{l}_*\rbul{X_{\bar{\varepsilon}}} \bigr )^{\textup{h}Gal(\bar{\varepsilon}/\varepsilon)}
\end{equation*} 
and by the unit $\Ql{,\sigma}\rightarrow \Ql{,\sigma}^{\textup{h}Gal(\bar{\eta}/\eta^{mu})}$, while the right diagonal arrows are induced by the counit of the adjunction $(trivial,(-)^{\textup{h}Gal(\bar{\eta}/\eta^{mu})})$.
Then the front square, the left and the right diagonal squares commute. Moreover, the map $\Hl(\varepsilon)\xrightarrow{id\otimes 1} \Hl(\varepsilon)\otimes_{\Ql{,\sigma}}\Ql{,\sigma}^{\textup{h}Gal(\bar{\eta}/\eta^{mu})}\simeq \Hl(\eta)$ identifies with the counit of the adjuction $(trivial, (-)^{\textup{h}Gal(\bar{\eta}/\eta^{mu})})$. Therefore, the two maps 
\[
\Hl(\varepsilon)\rightarrow \bigl ( p_{\sigma*}i_{01}^*\bar{l}_*\rbul{X_{\bar{\varepsilon}}}\bigr )^{\textup{h}Gal(\bar{\eta}/\eta)}
\]
coincide (up to homotopy). On the other hand. the two maps 
\[
\Ql{,\sigma}^{\textup{h}Gal(\bar{\eta}/\eta^{mu})}\rightarrow \bigl ( p_{\sigma*}i_{01}^*\bar{l}_*\rbul{X_{\bar{\varepsilon}}}\bigr )^{\textup{h}Gal(\bar{\eta}/\eta)}
\]
both coincide (up to homotopy) with the morphism that endows the latter with a $\Ql{\sigma}^{\textup{h}Gal(\bar{\eta}/\eta^{mu})}$-algebra structure.
By the universal property of the tensor product, this shows that the two maps 
\[
\Hl(\eta)\rightarrow \bigl ( p_{\sigma*}i_{01}^*\bar{l}_*\rbul{X_{\bar{\varepsilon}}}\bigr )^{\textup{h}Gal(\bar{\eta}/\eta)}
\]
are equivalent (up to coherent homotopy).

It remains to show that the two maps 
\[
p_{\sigma*}\rbul{X_{\sigma}}\rightarrow \bigl ( p_{\sigma*}i_{01}^*\bar{l}_*\rbul{X_{\bar{\varepsilon}}}\bigr )^{\textup{h}Gal(\bar{\eta}/\eta)}
\]
also coincide. In order to do so, one can look at the following diagram
\begin{equation*}
\begindc{\commdiag}[18]
\obj(0,120)[1]{$p_{\sigma*}\rbul{X_{\sigma}}$}
\obj(145,120)[2]{$p_{\sigma*}\rbul{X_{\sigma}}^{\textup{h}Gal(\bar{\eta}/\eta)}$}
\obj(0,80)[3]{$p_{\sigma*}\rbul{X_{\sigma}}^{\textup{h}Gal(\bar{\varepsilon}/\varepsilon)}\otimes_{\Ql{\sigma}}\Ql{\sigma}^{\textup{h}Gal(\bar{\eta}/\eta^{mu})}$}
\obj(103,80)[4]{$p_{\sigma*}\rbul{X_{\sigma}}^{\textup{h}Gal(\bar{\varepsilon}/\varepsilon)}$}
\obj(103,40)[5]{$\bigl ( p_{\sigma*}i_{01}^*\bar{l}_*\rbul{X_{\bar{\varepsilon}}} \bigr )^{\textup{h}Gal(\bar{\varepsilon}/\varepsilon)}$}
\obj(0,0)[6]{$\Hr(X_{\varepsilon})\otimes_{\Ql{\sigma}}\Ql{\sigma}^{\textup{h}Gal(\bar{\eta}/\eta^{mu})}$}
\obj(145,0)[7]{$\bigl ( p_{\sigma*}i_{01}^*\bar{l}_*\rbul{X_{\bar{\varepsilon}}} \bigr )^{\textup{h}Gal(\bar{\eta}/\eta)}$}
\mor{1}{2}{$\footnotesize{\text{counit } (\footnotesize{trivial},(-)^{\textup{h}Gal(\bar{\eta}/\eta)})}$}
\mor{1}{3}{$
\begin{minipage}[c]{3.3cm}
\centering
  \footnotesize{\text{counit }}
\begin{math}
  (\footnotesize{trivial},(-)^{\textup{h}Gal(\bar{\varepsilon}/\varepsilon)})
\end{math}
\end{minipage} \hspace{0.1cm}
\otimes 1$}[\atright,\solidarrow]

\mor{1}{4}{$\begin{minipage}[c]{3.3cm}
\centering
  \footnotesize{\text{counit }}
\begin{math}
  (\footnotesize{trivial},(-)^{\textup{h}Gal(\bar{\varepsilon}/\varepsilon)})
\end{math}
\end{minipage} $}
\mor(108,85)(140,115){$\begin{minipage}[c]{3.1cm}
\centering
  \footnotesize{\text{counit }}
\begin{math}
  (\footnotesize{trivial},(-)^{\textup{h}Gal(\bar{\eta}/\eta^{mu})})
\end{math}
\end{minipage}$}[\atright,\solidarrow]
\mor{3}{6}{$(sp_1)^{\textup{h}Gal(\bar{\varepsilon}/\varepsilon)}\otimes id$}[\atright,\solidarrow]
\mor{4}{3}{$id\otimes 1$}[\atright,\solidarrow]
\mor{4}{5}{$(sp_1)^{\textup{h}Gal(\bar{\varepsilon}/\varepsilon)}$}[\atright,\solidarrow]
\mor{5}{6}{$id\otimes 1$}[\atright,\solidarrow]
\mor{6}{7}{$\simeq$}
\mor(145,115)(145,95){$$}[\atright,\solidline]
\mor(145,85)(145,5){$(sp_1)^{\textup{h}Gal(\bar{\eta}/\eta)}$}[\atleft,\solidarrow]
\mor(108,35)(140,5){$\begin{minipage}[c]{3.3cm}
\centering
  \footnotesize{\text{counit }}
\begin{math}
  (\footnotesize{trivial},(-)^{\textup{h}Gal(\bar{\eta}/\eta^{mu})})
\end{math}
\end{minipage}$}[\atright,\solidarrow]
\enddc
\end{equation*}
\end{proof}
\begin{prop}
The following square is commutative (up to coherent homotopy)
\begin{equation}
\begindc{\commdiag}[18]
\obj(0,20)[1]{$\Hr(X_{\varepsilon})\otimes_{\Ql{,\sigma}}\Ql{,\sigma}^{\textup{h}Gal(\bar{\eta}/\eta^{mu})}$}
\obj(120,20)[2]{$\bigl ( p_{\sigma*}i_{01}^*\bar{l}_*\rbul{X_{\bar{\varepsilon}}}\bigr )^{\textup{h}Gal(\bar{\eta}/\eta^{mu})}$}
\obj(0,-20)[3]{$\Hr(X_{\eta})$}
\obj(120,-20)[4]{$\bigl ( p_{\sigma*}i_{01}^*\bar{l}_*\bar{k}^*\bar{j}_{10*}\rbul{X_{\bar{\eta}}}\bigr )^{\textup{h}Gal(\bar{\eta}/\eta)}$}
\mor{1}{2}{$(\ref{rbul varepsilon hgal eta})$}
\mor{1}{3}{$can$}
\mor{2}{4}{$(sp_2)^{\textup{h}Gal(\bar{\eta}/\eta)}$}
\mor{3}{4}{$(\ref{rbul eta hgal eta})$}
\enddc
\end{equation}
\end{prop}
\begin{proof}

Indeed, one can look at the following diagram
\begin{equation*}
\begindc{\commdiag}[18]
\obj(0,90)[1]{$\Hr(X_{\varepsilon})\otimes_{\Ql{,\sigma}}\Ql{,\sigma}^{\textup{h}Gal(\bar{\eta}/\eta^{mu})}$}
\obj(160,90)[2]{$\bigl ( p_{\sigma*}i_{01}^*l_*\rbul{X_{\varepsilon}}\bigr )^{\textup{h}Gal(\bar{\eta}/\eta^{mu})}$}
\obj(0,-20)[3]{$\Hr(X_{\eta})$}
\obj(160,-20)[4]{$\bigl ( p_{\sigma*}i_{01}^*\bar{l}_*\bar{k}^*\bar{j}_{10*}\rbul{X_{\bar{\eta}}}\bigr )^{\textup{h}Gal(\bar{\eta}/\eta)}$}
\obj(105,50)[5]{$\Hr(X_{\varepsilon})$}
\obj(105,15)[6]{$\bigl ( p_{\sigma*}i_{01}^*\bar{l}_*\bar{k}^*\bar{j}_{10*}\rbul{X_{\bar{\eta}}} \bigr )^{\textup{h}Gal(\bar{\varepsilon}/\varepsilon)}$}
\mor{1}{2}{$(\ref{rbul varepsilon hgal eta})$}
\mor{1}{3}{$can$}
\mor(160,87)(160,67){$$}[\atright,\solidline]
\mor(160,60)(160,-15){$(sp_2)^{\textup{h}Gal(\bar{\eta}/\eta)}$}
\mor{3}{4}{$(\ref{rbul eta hgal eta})$}[\atright,\solidarrow]
\mor{5}{1}{$id\otimes 1$}
\mor(105,55)(155,85){$\begin{minipage}[c]{3.3cm}
\centering
  \footnotesize{\text{counit }}
\begin{math}
  (\footnotesize{trivial},(-)^{\textup{h}Gal(\bar{\eta}/\eta^{mu})})
\end{math}
\end{minipage}$}[\atright,\solidarrow]
\mor{5}{3}{$
\begin{minipage}[c]{2.3cm}
\centering
\footnotesize{\text{counit } }
\begin{math}
(j_{10}^*,j_{10*})
\end{math}
\end{minipage}$}[\atright,\solidarrow]
\mor{5}{6}{$(sp_2)^{\textup{h}Gal(\bar{\varepsilon}/\varepsilon)}$}
\mor{4}{6}{$
\begin{minipage}[c]{3.3cm}
\centering
\footnotesize{\text{counit } }
\begin{math}
(\footnotesize{trivial}, (-)^{\textup{h}Gal(\bar{\eta}/\eta^{mu})})
\end{math}
\end{minipage}$}[\atleft,\solidarrow]
\enddc
\end{equation*}
and see that the two maps $\Hr(X_{\varepsilon})\rightarrow \bigl ( p_{\sigma*}i_{01}^*\bar{l}_*\bar{k}^*\bar{j}_{10*}\rbul{X_{\bar{\eta}}}\bigr )^{\textup{h}Gal(\bar{\eta}/\eta)}$ are equivalent. The two maps $\Ql{\sigma}^{\textup{h}Gal(\bar{\eta}/\eta^{mu})}\rightarrow \bigl ( p_{\sigma*}i_{01}^*\bar{l}_*\bar{k}^*\bar{j}_{10*}\rbul{X_{\bar{\eta}}}\bigr )^{\textup{h}Gal(\bar{\eta}/\eta)}$ coincide as both are equivalent to the morphism that endows $\bigl ( p_{\sigma*}i_{01}^*\bar{l}_*\bar{k}^*\bar{j}_{10*}\rbul{X_{\bar{\eta}}}\bigr )^{\textup{h}Gal(\bar{\eta}/\eta)}$ with the given $\Ql{\sigma}^{\textup{h}Gal(\bar{\eta}/\eta^{mu})}$-algebra structure. Then the lemma follows by the universal property of the tensor product.
\end{proof}
It then follows immediately that
\begin{cor}\label{fundamental triangles rbul cohomology}
The two triangles 
\[
\begindc{\commdiag}[20]
\obj(0,40)[1]{$p_{\sigma*}\rbul{X_{\sigma}}\otimes_{\Ql{,\sigma}}\Hl(\eta)$}
\obj(140,40)[2]{$\Hr(X_{\varepsilon})\otimes_{\Ql{,\sigma}}\Ql{,\sigma}^{\textup{h}Gal(\bar{\eta}/\eta^{mu})}$}
\obj(70,10)[3]{$\Hr(X_{\eta})$}
\mor(0,35)(65,15){$$}
\mor{1}{2}{$$}
\mor(135,35)(75,15){$$}
\obj(0,-20)[4]{$\bigl ( p_{\sigma*}\rbul{X_{\sigma}}\bigr )^{\textup{h}Gal(\bar{\eta}/\eta)}$}
\obj(140,-20)[5]{$\bigl ( p_{\sigma*}i_{01}^*\bar{l}_*\rbul{X_{\bar{\varepsilon}}}\bigr )^{\textup{h}Gal(\bar{\eta}/\eta)}$}
\obj(70,-50)[6]{$\bigl ( p_{\sigma*}i_{01}^*\bar{l}_*\bar{k}^*\bar{j}_{10*}\rbul{X_{\bar{\eta}}} \bigr )^{\textup{h}Gal(\bar{\eta}/\eta)}$}
\mor{4}{5}{$(sp_1)^{\textup{h}Gal(\bar{\eta}/\eta)}$}
\mor(135,-25)(75,-45){$(sp_2)^{\textup{h}Gal(\bar{\eta}/\eta)}$}
\mor(0,-25)(65,-45){$(sp_2\circ sp_1)^{\textup{h}Gal(\bar{\eta}/\eta)}$}[\atright,\solidarrow]
\enddc
\]
are equivalent.
\end{cor}
\subsection{\textsc{\texorpdfstring{$\ell$}{l}-adic realization of dg categories of relative singularities}}

In this section we will try to highlight the connection between the $\rbul{}$-valued cohomology we studied in the previous section, in particular the triangle of Corollary \ref{fundamental triangles rbul cohomology} with the $\ell$-adic realization of certain dg categories of singularities. 
\begin{caution}\label{caution}
In this section we work under the hypothesis that an $\ell$-adic realization functor exists in our context (i.e. over the base $S=Spec(\C \dsl x \dsr + y \cdot \C \drl x\drr \dsl y \dsr)$, that is not noetherian!) and moreover that it is compatible with the formalism of $6$-operations and Tate shifts as the one which is available in the noetherian case. If the author didn't misunderstand them, the experts believe that this is possible: the idea is to use the continuity properties of $\md_{\MB{}}(\textup{\textbf{SH}})$ and $\shv(\bullet,\Z_{\ell}/n\Z_{\ell})$ as showed in \cite{cd16} and \cite{cd19} to extend the $\ell$-adic realization functor to the non-noetherian world.
Here is a sketch of the construction:
\begin{itemize}
    \item For $X$ a noetherian scheme, $\ell$ invertible in every residue field of $X$ and $n\geq 1$, consider $\shv(X,\Z/\ell^{n}\Z)$, the $\oo$-category of $\md_{\Z/\ell^{n}\Z}$ étale hypersheaves. One then considers the sub $\oo$-category of locally constructible objects (\cite{cd16}), denoted $\shv_{lc}(X,\Z/\ell^{n}\Z)$
    \item If $X$ is an affine scheme, we can write it as a filtered limit of noetherian affine schemes (with affine transiction maps), $X=\varprojlim_{i\in I}X_i$. Then we define
    \begin{equation}
        \shv_{lc}(X,\Z/\ell^{n}\Z):=\varinjlim_{i\in I}\shv_{lc}(X_i,\Z/\ell^{n}\Z)
    \end{equation}
    We extend the above definition to non-affine schemes by Zariski descent. This definition is compatible with the previous one (this is the continuity theorem proved in \cite{cd16})
    \item Define
    \begin{equation}
        \shvl(X):=\textup{Ind}(\varprojlim_{n}\shv_{lc}(X,\Z/\ell^{n}\Z))\otimes \Q
    \end{equation}
    \item Since 
    \begin{equation}
        \textup{\textbf{SH}}_X \simeq \varinjlim_{i\in I}\textup{\textbf{SH}}_{X_i}
    \end{equation}
    and $\MB{}$ is compatible with filtered limits, we get an $\ell$-adic realization
    \begin{equation}
        \md_{\MB{}}(\textup{\textbf{SH}}_X)\rightarrow \shvl(X)
    \end{equation}
\end{itemize}
Special thanks to D.~C.~Cisinski, who kindly and patiently explained this to me. 

We will also need to assume that the results of \cite[\S B.4]{pr11} extend to our context. In particular, we need to assume that if $X$ is a regular, finitely presented, proper and flat $S$-scheme, $Y$ is smooth over $S$ and $Z=X\times_SSpec(R/y\cdot R)$, then
\begin{equation}
    \cohb{Z}\otimes_S \cohb{Y}\simeq \cohb{Z\times_S Y}
\end{equation}
and
\begin{equation}
    \cohb{X}_Z\otimes_S\cohb{Y}\simeq \cohb{X\times_SY}_{Z\times_SY}
\end{equation}
\end{caution}
\begin{context}
Let $p:X\rightarrow S$ be a flat, proper, finitely presented morphism and suppose that $X$ is regular, i.e. that $\perf{X}\simeq \cohb{X}$. Notice that, even if $X$ is not noetherian, it is coherent ($S$ is a coherent ring and $p$ is a finitely generated morphism) and it is possible to consider $\cohb{X}$.
\end{context}
\begin{rmk}
Notice that it follows immediately from \cite{brtv} and the previous sections that the cofibers of $(sp_1)^{\textup{h}Gal(\bar{\eta}/\eta)}$ and $(sp_2)^{\textup{h}Gal(\bar{\eta}/\eta)}$ are related to dg categories of relative singularities. Indeed, being $X$ regular, if we consider $p_{|X_{U_0}}:X_{U_0}\rightarrow U_0$ we are exactly in the situation studied in \textit{loc. cit.}. In particular, we have an equivalence of $\Ql{,\bar{\varepsilon}}^{\textup{h}Gal(\bar{\eta}/\eta^{mu})}$-modules
\begin{equation}
    \bar{k}^*\mathcal{R}^{\ell,\vee}_{\bar{U}_0}(\sing{X_{\bar{\varepsilon}}})\simeq (p_{X_{\bar{\varepsilon}}*}\Phi_{p_{\bar{U}_0}}(\Ql{,X_{\bar{U}_0}}(\beta)))^{\textup{h}Gal(\bar{\eta}/\eta^{mu})}[-1]
\end{equation}
Applying $v_{\bar{\varepsilon}*}(-)^{\textup{h}Gal(\bar{\varepsilon}/\varepsilon)}$ we obtain
\begin{equation}
    k^*\mathcal{R}^{\ell,\vee}_{U_0}(\sing{X_{\varepsilon}})\simeq (p_{X_{\varepsilon}*}\Phi_{p_{U_0}}(\Ql{,X_{U_0}}(\beta)))^{\textup{h}Gal(\bar{\eta}/\eta)}[-1]
\end{equation}
and therefore 
\begin{equation}
    i_{01}^*l_*\mathcal{R}^{\ell,\vee}_{\varepsilon}(\sing{X_{\varepsilon}})\simeq fib((sp_2)^{\textup{h}Gal(\bar{\eta}/\eta)})
\end{equation}
in $\md_{\Hr(\eta)}(\shvl(\sigma))$.

On the other hand, if we consider $p_1:X_1\rightarrow S_1$, we can't apply the results of \cite{brtv} as $X_1$ is not regular. Anyway, it can be easily seen that the only obstruction to do so is that we don't know whether the canonical morphism
\begin{equation}\label{conjectured devissage K theory}
    \mathcal{M}^{\vee}_{S_{1}}(\cohbp{X_0}{X_1})\rightarrow \mathcal{M}^{\vee}_{S_{1}}(\perf{X_1}_{X_0})
\end{equation}
is an equivalence in $\textup{\textbf{SH}}_{S_1}$. See also Section \ref{about the regularity hypothesis}.

However, it is possible to recover the $Gal(\bar{\varepsilon}/\varepsilon)$-invariant part of vanishing cycles if we consider the pushout of $\mathcal{R}^{\ell}_{S_1}(\perf{X_0})\rightarrow \mathcal{R}^{\ell}_{S_1}(\perf{X_1}_{X_0})$, which we shall label $C$. More precisely, there is an equivalence of $\Hr(\varepsilon)$-modules in $\shvl(\sigma)$
\begin{equation}
    i_{01}^*C\simeq (p_{\sigma*}(\Phi_{p_1}(\Ql{,X_1}(\beta)))[-1])^{\textup{h}Gal(\bar{\varepsilon}/\varepsilon)}
\end{equation}
Then, applying $(-)^{\textup{h}Gal(\bar{\eta}/\eta^{mu})}$ we get
\begin{equation}
    i_{01}^*C^{\textup{h}Gal(\bar{\eta}/\eta^{mu})}\simeq fiber((sp_1)^{\textup{h}Gal(\bar{\eta}/\eta)})
\end{equation}
in $\md_{\Hr(\eta)}(\shvl(\sigma))$.

If the map (\ref{conjectured devissage K theory}) is an equivalence, the left hand side can be identified with $i_{01}^*\mathcal{R}^{\ell}_{S_1}(\sing{X_1,x\circ p_1})^{\textup{h}Gal(\bar{\eta}/\eta^{mu})}$.
\end{rmk}

As a consequence of the previous remark, the only piece of (\ref{fundamental triangles rbul cohomology}) which remains to be understood in terms of non-commutative spaces is 
\begin{equation*}
   (sp_2\circ sp_1)^{\textup{h}Gal(\bar{\eta}/\eta)}:(p_{\sigma*}\rbul{X_{\sigma}})^{\textup{h}Gal(\bar{\eta}/\eta)}\rightarrow (\Psi^{(2)}_p(\rbul{X_{\bar{\eta}}}))^{\textup{h}Gal(\bar{\eta}/\eta)} 
\end{equation*}This is what we will try to do in the remaining part of this work.
\begin{construction}
Let $Y:=X\times_{S}Spec(R/y \cdot R)$. Consider the exact sequence in $\dgcatm$
\begin{equation}
    \perf{Y}\rightarrow \cohb{Y}\rightarrow \sing{Y}
\end{equation}
Since exact sequences in $\dgcatm$ are sent to fiber-cofiber sequences in $\textup{\textbf{Sh}}_{S}$ by the lax monoidal $\oo$-functor $\mathcal{R}^{\ell,\vee}_{S}$, we obtain
\begin{equation}\label{fiber cofiber Rl sing(X_1)}
    \mathcal{R}^{\ell,\vee}_{S}(\perf{Y})\rightarrow \mathcal{R}^{\ell,\vee}_{S}(\cohb{Y})\rightarrow \mathcal{R}^{\ell,\vee}_{S}(\sing{Y})
\end{equation}

Since $X$ is a regular scheme, $\perf{X}_{Y}$ coincides with $\cohb{X}_Y$ (as $\perf{X}\simeq \cohb{X}$ and $\perf{X_{\eta}}\simeq \cohb{X_{\eta}}$). We are assuming that the results of \cite[\S B.4]{pr11} extend to the case of coherent rings and therefore, as in Proposition \ref{motivic realization cohb}, by the theorem of the Heart (see \cite{ba15}) we get
\begin{equation}
 \mathcal{R}^{\ell,\vee}_{S}(\perf{X}_{Y}) \simeq \mathcal{R}^{\ell,\vee}_{S}(\cohb{X}_{Y}) \simeq \mathcal{R}^{\ell,\vee}_{S}(\cohb{Y})
\end{equation}
Indeed, $\textup{\textbf{Coh}}(Y)\rightarrow \textup{\textbf{Coh}}(X)_Y$ satisfies dévissage: as $\eta \simeq R[y^{-1}]$, we can describe $\textup{\textbf{Coh}}(X)_Y$ as those coherent sheaves killed by some power of $y$. Notice that $R/y \cdot R$ is a coherent ring (it is a quotient of a coherent ring by a finitely generated ideal) and $Y$ is finitely presented over $R/y \cdot R$, thus we are allowed to consider $\cohb{Y}$.

In particular, applying $i_{0}^*:\textup{\textbf{SH}}_{S}\rightarrow \textup{\textbf{SH}}_{\sigma}$ we can include the fiber-cofiber sequence (\ref{fiber cofiber Rl sing(X_1)}) in the following diagram
\begin{equation}
    \begindc{\commdiag}[10]
    \obj(-120,40)[1a]{$i_{0}^*\mathcal{R}^{\ell,\vee}_{S}(\perf{Y})$}
    \obj(0,40)[2a]{$i_{0}^*\mathcal{R}^{\ell,\vee}_{S}(\cohb{Y})$}
    \obj(120,40)[3a]{$i_{0}^*\mathcal{R}^{\ell,\vee}_{S}(\sing{Y})$}
    \obj(0,0)[4a]{$i_{0}^*\mathcal{R}^{\ell,\vee}_{S}(\perf{X})$}
    \obj(0,-40)[5a]{$i_0^* \mathcal{R}^{\ell,\vee}_{S}(\perf{X_{\eta}})$}
    \mor{1a}{2a}{$$}
    \mor{2a}{3a}{$$}
    \mor{2a}{4a}{$$}
    \mor{4a}{5a}{$$}
    \mor(-100,30)(-30,2){$\sim 0$}[\atright,\solidarrow]
    \enddc
\end{equation}
which is equivalent to
\begin{equation}
    \begindc{\commdiag}[10]
    \obj(-120,40)[1]{$p_{\sigma*}\Ql{,X_{\sigma}}(\beta)$}
    \obj(0,40)[2]{$p_{\sigma*}i_{01}^*i_1^!\Ql{,X}(\beta)$}
    \obj(120,40)[3]{$i_{0}^*\mathcal{R}^{\ell,\vee}_{S}(\sing{Y})$}
    \obj(0,0)[4]{$p_{\sigma*}\Ql{,X_{\sigma}}(\beta)$}
    \obj(0,-40)[5]{$p_{\sigma*}i_0^*j_{1*}\Ql{,X_{\eta}}(\beta)$}
    \mor{1}{2}{$$}
    \mor{2}{3}{$$}
    \mor{2}{4}{$$}
    \mor{4}{5}{$$}
    \mor(-100,30)(-20,2){$\sim 0$}[\atright,\solidarrow]
    \enddc
\end{equation}
If we apply the octahedron axiom (see \cite{ha}) to the triangle in the diagram above 
\begin{equation}\label{fundamental fiber cofiber sequence nc spaces step 1}
    i_{0}^*\mathcal{R}^{\ell,\vee}_{S}(\sing{Y})\rightarrow p_{\sigma*}\Ql{,X_{\sigma}}(\beta)\oplus p_{\sigma*}\Ql{,X_{\sigma}}(\beta)[1]\rightarrow p_{\sigma*}i_0^*j_{1*}\Ql{,X_{\eta}}(\beta)
\end{equation}
\end{construction}
Then consider the pullback square
\begin{equation}
    \begindc{\commdiag}[15]
    \obj(-15,15)[1]{$X_1$}
    \obj(15,15)[2]{$X_{\varepsilon}$}
    \obj(-15,-15)[3]{$X$}
    \obj(15,-15)[4]{$X_{U_0}$}
    \mor{2}{1}{$l$}
    \mor{1}{3}{$i_1$}
    \mor{4}{3}{$j_0$}
    \mor{2}{4}{$k$}
    \enddc
\end{equation}
and the morphism of $\oo$-functors
\begin{equation}\label{bcm}
    i_1^*j_{0*}\rightarrow l_*k^*
\end{equation}
which corresponds to
\begin{equation}
    j_{0*}\rightarrow i_{1*}l_*k^*\simeq j_{0*}k_*k^*
\end{equation}
whose fiber is $j_{0*}j_{10!}j_{10}^*$: for every object $E$ in $\textup{\textbf{SH}}_{U_0}$, there is a fiber-cofiber sequence
\begin{equation}\label{fiber cofiber sequence bcm}
    j_{0*}j_{10!}j_{10}^*E\rightarrow  j_{0*}E\rightarrow i_{1*}l_*k^*E\simeq j_{0*}k_*k^*E
\end{equation}
\begin{lmm}
There is a fiber cofiber sequcence in $\shvl(X_{\sigma})$
\begin{equation}
    i_0^*j_{0*}j_{10!}j_{10}^*j_{10*}\Ql{,X_{\eta}}(\beta)\rightarrow i_0^*j_{1*}\Ql{,X_{\eta}}(\beta)\rightarrow i_{01}^*l_*k^*j_{10*}\Ql{,X_{\eta}}
\end{equation}
\end{lmm}
\begin{proof}
In the fiber-cofiber sequence (\ref{fiber cofiber sequence bcm}), consider $E=j_{10*}\Ql{,X_{\eta}}(\beta)$ and apply $i_0^*$. What we get is
\begin{equation}
    i_0^*j_{0*}j_{10!}j_{10}^*j_{10*}\Ql{,X_{\eta}}(\beta)\rightarrow  i_0^*j_{0*}j_{10*}\Ql{,X_{\eta}}(\beta)\rightarrow  i_0^*i_{1*}l_*k^*j_{10*}\Ql{,X_{\eta}}(\beta)
\end{equation}
To conclude, it suffices to notice that
\begin{enumerate}
    \item $j_{0*}j_{10*}\simeq j_{1*}$
    \item $i_0^*\simeq i_{01}^*i_1^*$ and $i_1^*i_{1*}\simeq id$
\end{enumerate}
\end{proof}
\begin{construction}
Consider the following diagram
\begin{equation}
    \begindc{\commdiag}[18]
    \obj(-80,0)[1]{$i_{0}^*\mathcal{R}^{\ell,\vee}_{S}(\sing{Y})$}
    \obj(0,0)[2]{$p_{\sigma*}\Ql{,X_{\sigma}}(\beta)\oplus p_{\sigma*}\Ql{,X_{\sigma}}(\beta)[1]$}
    \obj(80,0)[3]{$p_{\sigma*}i_0^*j_{1*}\Ql{,X_{\eta}}(\beta)$}
    \obj(80,30)[4]{$p_{\sigma*}i_0^*j_{0*}j_{10!}j_{10}^*j_{10*}\Ql{,X_{\eta}}(\beta)$}
    \obj(80,-30)[5]{$p_{\sigma*}i_{01}^*l_*k^*j_{10*}\Ql{,X_{\eta}}(\beta)$}
    \mor{1}{2}{$$}
    \mor{2}{3}{$$}
    \mor{4}{3}{$$}
    \mor{3}{5}{$$}
    \mor{2}{5}{$\chi$}[\atright,\solidarrow]
    \enddc
\end{equation}

where the row is the fiber-cofiber sequence (\ref{fundamental fiber cofiber sequence nc spaces step 1}) and the coloumn on the right is the fiber-cofiber sequence of the previous lemma.

Then we can apply the octahedron property to the triangle on the right and obtain the following fiber cofiber sequence
\begin{equation}\label{fundamental fiber cofiber sequence nc spaces step 2}
    i_{0}^*\mathcal{R}^{\ell,\vee}_{S}(\sing{Y})\rightarrow fiber(\chi) \rightarrow p_{\sigma*}i_0^*j_{0*}j_{10!}j_{10}^*j_{10*}\Ql{,X_{\eta}}(\beta) 
\end{equation}
\end{construction}
Then we will need to understand the morphism $\chi$ from the vanishing cycles' theory point of view. By Lemma \ref{lemma equivalence calg cohomology torus} and Proposition \ref{gal invariant interesting stuff} , the following square is an homotopy pushout in $\textup{CAlg}(\shvl(\sigma))$:
\begin{equation}\label{diagram definition lambda}
    \begindc{\commdiag}[18]
    \obj(-50,15)[1]{$p_{\sigma*}\Ql{,X_{\sigma}}(\beta)$}
    \obj(40,15)[2]{$p_{\sigma*}\Ql{,X_{\sigma}}(\beta)\otimes_{\Ql{,\sigma}}\Ql{,\sigma}^{\textup{h}Gal(\bar{\varepsilon}/\varepsilon)}$}
    \obj(-50,-15)[3]{$p_{\sigma*}\Ql{,X_{\sigma}}(\beta)\otimes_{\Ql{,\sigma}}\Ql{,\sigma}^{\textup{h}Gal(\bar{\eta}/\eta^{mu})}$}
    \obj(40,-15)[4]{$p_{\sigma*}\Ql{,X_{\sigma}}(\beta)\otimes_{\Ql{,\sigma}}\Ql{,\sigma}^{\textup{h}Gal(\bar{\eta}/\eta)}$}
    \mor{1}{2}{$$}
    \mor{2}{4}{$$}
    \mor{3}{4}{$$}
    \mor{1}{3}{$$}
    \obj(100,-40)[5]{$\Hr(X_{\eta})$}
    \mor(50,-20)(80,-35){$$}
    \mor(-40,-20)(70,-40){$\lambda$}[\atright,\solidarrow]
    \cmor((77,15)(92,15)(100,12)(100,-25)(100,-30)) \pdown(95,-10){$$}
    \enddc
\end{equation}
and the map $p_{\sigma*}\Ql{,X_{\sigma}}(\beta)\otimes_{\Ql{,\sigma}}\Ql{,\sigma}^{\textup{h}Gal(\bar{\eta}/\eta)}\rightarrow p_{\sigma*}\Psi^{(2)}_p(\Ql{,X}(\beta))^{\textup{h}Gal(\bar{\eta}/\eta)}$ identifies with $(sp_2\circ sp_1)^{\textup{h}Gal(\bar{\eta}/\eta)}:(p_{\sigma*}\Ql{,X_{\sigma}})^{\textup{h}Gal{\bar{\eta}/\eta}}\rightarrow p_{\sigma*}\Psi^{(2)}_p(\Ql{,X}(\beta))^{\textup{h}Gal(\bar{\eta}/\eta)}$. As a next step, we will compare the two maps $\chi$ and $\lambda$. As a first observation, notice that the domain and the codomain of the two morphisms coincide (see Lemma \ref{homotopy invariant units}).
\begin{prop}
The two maps 
\begin{equation*}
\chi,\lambda : p_{\sigma*}\Ql{,X_{\sigma}}(\beta)\otimes_{\Ql{,\sigma}}\Ql{,\sigma}^{\textup{h}Gal(\bar{\eta}/\eta^{mu})}\rightarrow p_{\sigma*}i_{01}^*l_*k^*j_{10*}\Ql{,X_{\eta}}(\beta)=\Hr(X_{\eta})
\end{equation*}
are homotopic.
\end{prop}
\begin{proof}
The identification of the source and target of the two maps follows from the results of the previous sections. Then, by the universal property of the tensor product and since $\Ql{,\sigma}$ is the initial object in $\textup{CAlg}(\shvl(\sigma))$ (therefore, for any commutative algebra $\mathcal{C}\in \shvl(\sigma)$, $\Map_{\textup{CAlg}(\shvl(\sigma))}(\Ql{,\sigma},\mathcal{C})\simeq *$), it suffices to show that the two pairs of maps
\begin{equation}
    p_{\sigma*}\Ql{,X_{\sigma}}(\beta)\rightarrow p_{\sigma*}\Ql{,X_{\sigma}}(\beta)\otimes_{\Ql{,\sigma}}\Ql{,\sigma}^{\textup{h}Gal(\bar{\eta}/\eta^{mu})} \rightrightarrows \Hr(X_\eta)
\end{equation}

\begin{equation}
    \Ql{,\sigma}^{\textup{h}Gal(\bar{\eta}/\eta^{mu})}\rightarrow p_{\sigma*}\Ql{,X_{\sigma}}(\beta)\otimes_{,\Ql{\sigma}}\Ql{,\sigma}^{\textup{h}Gal(\bar{\eta}/\eta^{mu})}\rightrightarrows \Hr(X_{\eta})
\end{equation}
obtained by $\chi$ and $\lambda$ coincide . For what concerns the second pair of morphisms, it is immediate to see that they both coincide with the map which endows $\Hr(X_{\eta})$ with the structure of a $\Ql{,\sigma}^{\textup{h}Gal(\bar{\eta}/\eta^{mu})}$-algebra. For the first pair, notice that 
\begin{equation}
p_{\sigma*}\Ql{,X_{\sigma}}(\beta)\rightarrow p_{\sigma*}\Ql{,X_{\sigma}}(\beta)\otimes_{\Ql{,\sigma}}\Ql{,\sigma}^{\textup{h}Gal(\bar{\eta}/\eta^{mu})}\simeq p_{\sigma*}\Ql{,X_{\sigma}}(\beta)\oplus p_{\sigma*}\Ql{,X_{\sigma}}(\beta)[1]
\end{equation}
coincides with the inclusion of the first factor. The composition of this morphism with $\chi$ decomposes as
\begin{equation}
    \begindc{\commdiag}[18]
    \obj(-70,0)[1]{$p_{\sigma*}\Ql{,X_{\sigma}}(\beta)$}
    \obj(0,0)[2]{$p_{\sigma*}\Ql{,X_{\sigma}}(\beta)\oplus p_{\sigma*}\Ql{,X_{\sigma}}(\beta)[1]$}
    \obj(70,0)[3]{$p_{\sigma*}i_0^*j_{1*}\Ql{,X_{\eta}}(\beta)$}
    \obj(140,0)[4]{$\Hr(X_{\eta})$}
    \mor{1}{2}{$\begin{bmatrix}1\\0\end{bmatrix}$}
    \mor{2}{3}{$$}
    \mor{3}{4}{$i_1^*j_{0*}\rightarrow l_*k^*$}
    \cmor((-70,-4)(-70,-9)(-67,-12)(0,-12)(67,-12)(70,-9)(70,-4)) \pup(0,-17){$id\rightarrow j_{1*}j_1^*$}
    \cmor((0,4)(0,9)(3,12)(75,12)(137,12)(140,9)(140,5)) \pdown(75,17){$\chi$}
    \enddc
\end{equation}

that is equivalent to the morphism induced by $id\rightarrow i_{01}^*l_*k^*j_{10*}$. From Corollary \ref{fundamental triangles rbul cohomology}, it is easy to see that it is the one that defines
\begin{equation}
    \begindc{\commdiag}[18]
    \obj(-20,-40)[1]{$p_{\sigma*}\Ql{X_{,\sigma}}(\beta)$}
    \obj(-20,0)[2]{$p_{\sigma*}\Ql{,X_{\sigma}}(\beta)\otimes_{\Ql{\sigma}}\Ql{,\sigma}^{\textup{h}Gal(\bar{\eta}/\eta^{mu})}$}
    \obj(70,0)[3]{$p_{\sigma*}\Ql{,X_{\sigma}}(\beta)\otimes_{\Ql{,\sigma}}\Ql{,\sigma}^{\textup{h}Gal(\bar{\eta}/\eta)}$}
    \obj(170,0)[4]{$\Hr(X_{\eta})$}
    \mor{1}{2}{$\begin{bmatrix}1\\0\end{bmatrix}$}
    \mor{2}{3}{$$}
    \mor{3}{4}{$(sp_2\circ sp_1)^{\textup{h}Gal(\bar{\eta}/\eta)}$}[\atright,\solidarrow]
    \cmor((-20,4)(-20,9)(-17,12)(85,12)(167,12)(170,9)(170,6)) \pdown(75,17){$\lambda$}
    \enddc
\end{equation}
\end{proof}

The following lemma is a well known general fact about stable $\oo$-categories. We state it and prove it for the reader's convenience.
\begin{lmm}\label{general lemma stable category}
Let $\mathcal{C}$ be a stable $\oo$-category. Suppose we are given the following diagram 
\begin{equation}
    \begindc{\commdiag}[18]
    \obj(-30,0)[1]{$A$}
    \obj(0,0)[2]{$C$}
    \obj(30,0)[3]{$B$}
    \obj(0,30)[4]{$B'$}
    \obj(0,-30)[5]{$A'$}
    \mor{1}{2}{$u$}
    \mor{2}{3}{$v$}[\atright,\solidarrow]
    \mor{4}{2}{$v'$}[\atright,\solidarrow]
    \mor{2}{5}{$u'$}
    \mor{4}{3}{$v\circ v'$}
    \mor{1}{5}{$u'\circ u$}[\atright,\solidarrow]
    \enddc
\end{equation}
such that the column and the row are fiber cofiber sequences. Then $cofiber(u'\circ u)\simeq cofiber(v\circ v')$
\end{lmm}
\begin{proof}
The octahedron property (see \cite{ha}) applied to the lower triangle gives provides us with the following commutative diagram
\begin{equation}
    \begindc{\commdiag}[25]
    \obj(-45,15)[1]{$A$}
    \obj(-30,0)[2]{$C$}
    \obj(-15,15)[3]{$A'$}
    \mor{1}{2}{$u$}[\atright,\solidarrow]
    \mor{2}{3}{$u'$}
    \mor{1}{3}{$u'\circ u$}
    \obj(-15,-15)[4]{$B$}
    \obj(0,0)[5]{$cofiber(u'\circ u)$}
    \mor{2}{4}{$$}
    \mor{4}{5}{$$}
    \mor{3}{5}{$$}
    \obj(15,15)[6]{$B'[1]$}
    \obj(15,-15)[7]{$A[1]$}
    \mor{3}{6}{$$}
    \mor{5}{6}{$$}
    \mor{4}{7}{$$}
    \mor{5}{7}{$$}
    \obj(30,0)[8]{$C[1]$}
    \obj(45,15)[9]{$B[1]$}
    \mor{7}{8}{$u[1]$}[\atright,\solidarrow]
    \mor{8}{9}{$v[1]$}[\atright,\solidarrow]
    \mor{6}{8}{$v'[1]$}[\atright,\solidarrow]
    \mor{6}{9}{$(v\circ v')[1]$}
    \enddc
\end{equation}

where $cofiber(u'\circ u)\rightarrow B'[1]\xrightarrow{(v\circ v')[1]} B[1]$ is a fiber-cofiber sequence. Then the claim follows immediately.
\end{proof}

We are finally ready to state the main result of this section.
\begin{thm}\label{theorem chpt 4}
Let
\begin{equation}
    \xi:fiber(p_{\sigma*}\Ql{,X_{\sigma}}(\beta)^{\textup{h}Gal(\bar{\eta}/\eta^{mu})}\rightarrow p_{\sigma*}\Ql{,X_{\sigma}}(\beta)^{\textup{h}Gal(\bar{\eta}/\eta)})\rightarrow p_{\sigma*}i_0^*j_{0*}j_{10!}j_{10}^*j_{10*}\Ql{,X_{\eta}(\beta)}
\end{equation}
There exists a fiber-cofiber sequence in $\shvl(\sigma)$
\begin{equation}
\begindc{\commdiag}[18]
\obj(-80,0)[1]{$i_{0}^*\mathcal{R}^{\ell,\vee}_{S}(\sing{Y})$}
\obj(0,0)[2]{$(p_{\sigma*}\Phi^{(2)}_p(\Ql{,X}(\beta))[-1])^{\textup{h}Gal(\bar{\eta}/\eta)}$} 

\obj(70,0)[4]{$cofib(\xi)$}
\mor{1}{2}{$$}

\mor{2}{4}{$$}

\enddc
\end{equation}
\end{thm}
\begin{proof}
Consider the fiber-cofiber sequence (\ref{fundamental fiber cofiber sequence nc spaces step 2}) and the one that one obtains by the octahedron property applied to the lower triangle in diagram (\ref{diagram definition lambda}):
\begin{equation}
    fiber(p_{\sigma*}\Ql{,X_{\sigma}}(\beta)^{\textup{h}Gal(\bar{\eta}/\eta^{mu})}\rightarrow p_{\sigma*}\Ql{,X_{\sigma}}(\beta)^{\textup{h}Gal(\bar{\eta}/\eta)})\rightarrow fiber(\lambda) \rightarrow p_{\sigma*}\Phi^{(2)}_p(\Ql{,X}(\beta))^{\textup{h}Gal(\bar{\eta}/\eta)}
\end{equation}
They fit together in the following diagram
\begin{equation}
    \begindc{\commdiag}[18]
    \obj(-80,0)[1]{$i_{0}^*\mathcal{R}^{\ell,\vee}_{S}(\sing{Y})$}
    \obj(0,0)[2]{$fiber(\chi)\simeq fiber(\lambda)$}
    \obj(80,0)[3]{$p_{\sigma*}i_0^*j_{0*}j_{10!}j_{10}^*j_{10*}\Ql{,X_{\eta}(\beta)}$}
    \obj(0,40)[4]{$fiber(p_{\sigma*}\Ql{,X_{\sigma}}(\beta)^{\textup{h}Gal(\bar{\eta}/\eta^{mu})}\rightarrow p_{\sigma*}\Ql{,X_{\sigma}}(\beta)^{\textup{h}Gal(\bar{\eta}/\eta)})$}
    \obj(0,-40)[5]{$(p_{\sigma*}\Phi^{(2)}_p(\Ql{,X}(\beta)[-1])^{\textup{h}Gal(\bar{\eta}/\eta)}$}
    \mor{1}{2}{$$}
    \mor{2}{3}{$$}
    \mor{4}{2}{$$}
    \mor{2}{5}{$$}
    \mor(5,35)(75,5){$\xi$}
    \mor(-75,-5)(-5,-35){$$}
    \enddc
\end{equation}

The theorem follows by the previous lemma.
\end{proof}
